\documentclass[11pt,reqno]{amsart}
\usepackage{mathrsfs}
\usepackage{amsmath}
\usepackage{amssymb, amsmath}
\usepackage{color}

\numberwithin{equation}{section}

\allowdisplaybreaks

\let\al=\alpha

\let\e=\varepsilon

\let\lam=\lambda
\let\r=\rho

\let\f=\frac

\let\om=\omega
\let\G= \Gamma

\let\Om=\Omega

\let\na=\nabla

\let\pa=\partial

\def\cA{{\mathcal A}}
\def\cB{{\mathcal B}}

\def\cD{{\mathcal D}}
\def\cE{{\mathcal E}}

\def\cM{{\mathcal M}}

\def\cP{{\mathcal P}}

\def\cS{{\mathcal S}}

\def\R{\mathbf R}
\def\Z{\mathbf Z}

\def\no{\noindent}

\def\dv{\mbox{div}}
\def\dive{\mathop{\rm div}\nolimits}

\def\Ds{\langle D\rangle^s}
\def\eqdef{\buildrel\hbox{\footnotesize def}\over =}
\def\ef{\hphantom{MM}\hfill\llap{$\square$}\goodbreak}

\newcommand{\beq}{\begin{equation}}
\newcommand{\eeq}{\end{equation}}
\newcommand{\ben}{\begin{eqnarray}}
\newcommand{\een}{\end{eqnarray}}
\newcommand{\beno}{\begin{eqnarray*}}
\newcommand{\eeno}{\end{eqnarray*}}

\newtheorem{theorem}{Theorem}[section]
\newtheorem{definition}[theorem]{Definition}
\newtheorem{lemma}[theorem]{Lemma}
\newtheorem{proposition}[theorem]{Proposition}
\newtheorem{corol}[theorem]{Corollary}
\newtheorem{remark}[theorem]{Remark}

\begin{document}

\title[The incompressible Euler equations with free boundary]
{\small Local well-posedness and break-down criterion of the incompressible Euler equations with free boundary}

\author{Chao Wang}

\address{Beijing International Center For Mathematical Research, Peking University, 100871, P. R. China}
\email{wangchao@math.pku.edu.cn}

\author{Zhifei Zhang}
\address{School of Mathematical Sciences, Peking University, 100871, P. R. China}
\email{zfzhang@math.pku.edu.cn}

\author{Weiren Zhao}
\address{Department of Mathematics, Zhejiang University, 310027, P. R. China}
\email{zjzjzwr@126.com}

\author{Yunrui Zheng}
\address{Beijing International Center For Mathematical Research,  Peking University, 100871, P. R. China}
\email{ruixue@mail.ustc.edu.cn}

\date{\today}

\begin{abstract}
In this paper, we prove the local well-posedness of the free boundary problem for the incompressible Euler equations
in low regularity Sobolev spaces, in which the velocity is a Lipschtiz function and the free surface belongs to $C^{\f32+\varepsilon}$.
Moreover, we also present a Beale-Kato-Majda type break-down criterion of smooth solution
in terms of the mean curvature of the free surface, the gradient of the velocity and Taylor sign condition.
\end{abstract}

\maketitle

\section{Introduction}

\subsection{Presentation of the problem}

In this paper, we consider the motion of an ideal incompressible gravity fluid in a domain with free boundary of finite depth
\beno
\big\{(t,x,y)\in [0,T]\times\mathbf{R}^d\times\mathbf{R}: (x,y)\in \Omega_t\big\},
\eeno
where $\Omega_t$ is the fluid domain at time $t$ located by
\beno
\Omega_t=\big\{(x,y)\in\mathbf{R}^d\times\mathbf{R}:b(x)<y<\eta(t,x)\big\}.
\eeno

The motion of the fluid is described by the incompressible Euler equation
\beq\label{eq:euler}
\pa_t v+v\cdot \na_{x,y}v=-ge_{d+1}-\na_{x,y} P\qquad \hbox{in}\quad \Omega_t,\,\, t\geq 0,
\eeq
where $-ge_{d+1}=-g(0,\cdots,0,1)$ denotes the acceleration of gravity, $v=(v^1,...,v^{d+1})$
denotes the velocity field, and $P$ denotes the pressure.
The incompressibility of the fluid is expressed by
\beq \label{eq:euler-d}
\textrm{div}\, v=0 \qquad \hbox{in}\quad \Omega_t,\,\, t\geq 0.
\eeq

Assume that no fluid particles are transported across the surface. At the bottom, this is given by
\beq \label{eq:euler-b1}
v_n|_{y=b(X)}:=\textbf{n}_-\cdot\ v|_{y=b(X)}=0 \qquad \hbox{for}\quad t\ge 0, x\in\R^d
\eeq
where $\textbf{n}_-:=\frac{1}{\sqrt{1+|\na_X b|^2}}(\na_X b,-1)^T$ denotes the outward normal vector to
the lower boundary of $\Omega_t $. At the free surface, the boundary condition is kinematic and is given by
\beq\label{eq:euler-b2}
\pa_t \eta-\sqrt{1+|\na\eta|^2}v_n|_{y=\eta(t,x)}=0\quad \hbox{for}\quad t\ge 0, x\in\R^d,
\eeq
where $v_n=\textbf{n}_+\cdot v|_{y=\eta(t,x)}$, with $\textbf{n}_+:=\frac{1}{\sqrt{1+|\na
\eta|^2}}(-\na\eta, 1)^T$ denoting the outward normal vector to the free surface $\Sigma_t$.
In general, the pressure at the free surface is proportional to the mean curvature of the free surface, i.e.,
\ben\label{eq:euler-p}
P|_{y=\eta(t,x)}=-\kappa\na\cdot\Big(\frac{\na\eta}{\sqrt{1+|\na\eta|^2}}\Big)\qquad \hbox{for}\quad t\geq 0,
x\in \R^d,
\een
where $\kappa\ge 0$ is the surface tension coefficient.

In this paper, we consider the case without surface tension. Furthermore, we assume that
the bottom is flat(i.e., $b(x)=-1$) in order to simplify our presentation.
We also take the gravity constant $g=1$.

\subsection{Some known results}

Let us first review some known results concerning the water-wave equations without vorticity.
In the case when the surface tension is neglected and the motion of free surface is a
small perturbation of still water, one could check Nalimov
\cite{Nal}, Yosihara \cite{Yos} and Craig \cite{Craig} for the well-posedness of 2-D water-wave equations. In general,
the local well-posedness of the water-wave equations of infinite depth without surface tension
was solved by Wu \cite{Wu1,Wu2}, where she showed that the Taylor sign condition
\ben
-\f {\partial p} {\partial \textbf{n}}\big|_{y=\eta(t,x)}\ge c_0>0
\een
always holds as long as the free surface is no self-intersection. In \cite{AM1,AM2}, Ambrose and Masmoudi present a different proof.
Lannes \cite{Lan} first solves the water-wave equations of finite depth without surface tension in the framework of the Eulerian coordinates.
Ming and Zhang \cite{MZ} generalize Lannes's result to the case with surface tension.
In a series of works \cite{ABZ-DM, ABZ-AS, ABZ-IM}, Alazard, Burq and Zuily  use the tools of paradifferential operators and Strichartz estimates
to prove the local well-posedness of the water-wave equations in low regularity Sobolev spaces.

For small initial data, Wu \cite{Wu3} first proved the almost global well-posedness of 2-D water-wave equations.
Wu \cite{Wu4} and Germain, Masmoudi and Shatah \cite{GMS} proved the global well-posedness of 3-D water-wave equations by using different method.
Recently, Alarzard and Delort \cite{AD} and Ionescu and Pusateri \cite{IP} independently proved the global well-posedness of 2-D water-wave equations,
see also \cite{Ta1, Ta2} for a new proof based on the holomorphic coordinates.
On the other hand, Castro, Cordoba, Ferferman, Gancedo and Lopez-Fernandez \cite{CCF-AM1} showed that
there exists smooth initial data for the water-waves equations such that the solution overturns in finite time.
See \cite{CCF-AM2, CS-CMP} for the formation of the splash singularity.
Wu \cite{Wu5} also construct a class of self-similar solution for the 2-D water-wave equations without the gravity.

Now, we review some well-posedness results for the rotational water-wave equations. Christodoulou and Lindblad \cite{CL}
presented the a priori estimates of the incompressible Euler equations in a free domain diffeomorphic to a ball.
Later, Lindblad \cite{Lin} proved the local existence of smooth solution by using Nash-Moser iteration. Coutand and Shkoller \cite{CS-JM}
proved the local well-posedness of the incompressible Euler equations in both cases with surface tension and without surface tension by using the
lagrangian coordinates and a subtle mollification procedure. Zhang and Zhang \cite{ZZ} solves the incompressible Euler equations without surface tension
by using the framework of Clifford analysis introduced by Wu \cite{Wu2}.
Shatah and Zeng \cite{SZ1, SZ2, SZ3} solve this problem by deriving the evolution equations of geometry quantities, especially the mean curvature.

In this paper, we will first prove the local well-posedness of the rotational water-wave problem in low regularity Sobolev spaces,
and then present a break-down criterion to the obtained smooth solution in terms of physical quantity and geometrical quantity.
This work was motivated by Craig and Wayne's question proposed in \cite{CW}: ``{\bf How do solutions break down?}"

{\it There are several versions of this question, including `` What is the lowest exponent of a Sobolev space $H^s$ in which
one can produce an existence theorem local in time?" Or one could ask ``For which $\al$ is it true that, if one knows a priori that
$\sup_{[-T,T]}\|(\eta,\psi)\|_{C^\al}<+\infty$ and that $(\eta_0,\psi_0)\in C^\infty$, then the solution is fact $C^\infty$ over the time interval $[-T,T]$?"
$\cdots\cdots$It would be more satisfying to say that the solution fails to exist because the curvature of the surface has diverged at some point, or a
related geometrical and(or) physical statement.}

In the case without the vorticity and surface tension,
this question was solved by Alazard, Burq and Zuily for the low regularity well-posedness \cite{ABZ-IM},
and by Wang and Zhang for the break-down criterion \cite{WZ}.

\subsection{Main results}

The first main result of this paper is the local well-posedness of the water-wave equations in Sobolev spaces with low regularity,
where the regularity of the initial velocity is consistent with the classical local well-posedness result in $\R^{d+1}$ proved by
Kato and Ponce \cite{Kato}.

\begin{theorem}\label{thm:local}
Let $d\geq 1$ and $s>\f d 2+1$.
Assume that the initial data $(\eta_0, v_0)$ satisfies
\begin{eqnarray*}
\eta_0\in H^{s+\f 12}(\mathbf{R}^d),\quad v_0\in H^{s+\f12}(\Om_0).
\end{eqnarray*}
Furthermore, assume that there exist two positive constants $c_0>0$ and $h_0>0$ such that
\ben\label{ass:initial-P}
&&-(\pa_yP)(0,x,\eta_0(x))\ge c_0\quad\textrm{ for }x\in \R^d,\\
&&1+\eta_0(x)\ge h_0\quad \textrm{ for }x\in \R^d.\label{ass:initial-eta}
\een
Then there exists $T>0$ such that the system (\ref{eq:euler})--(\ref{eq:euler-p}) with the initial data
$(\eta_0, v_0)$ has a unique solution $(\eta, v)$ satisfying
\begin{gather*}
\eta\in C\big([0,T];H^{s+\f12}(\mathbf{R}^d)\big),\quad v\in C\big([0,T];H^{s+\f12}(\Om_t)\big).
\end{gather*}
\end{theorem}

\begin{remark}
The regularity of the initial velocity should be optimal by the recent strong ill-posedness result in the whole space proved by Bourgain and Li \cite{Bour}.
The regularity of the initial free surface could be further lowered by using the Strichartz type estimates, see \cite{ABZ-pre}
for the irrotational case.
\end{remark}

In a seminal paper \cite{BKM}, Beale, Kato and Majda showed that if $v$ is a smooth solution of the incompressible Euler equations in $[0,T)\times \R^3$ and satisfies
\beno
\int_0^T\|\na\times v(t)\|_{L^\infty(\R^3)}dt<+\infty,
\eeno
then the solution can be extended after $t=T$. The second main result of this paper is a Beale-Kato-Majda type blow-up criterion
for the free boundary problem of the incompressible Euler equations.

\begin{theorem}\label{thm:blow-up}
Let $s>\f d 2+1$ so that $s-\f12$ is an integer, and $(\eta, v, P)$ be the solution of the system (\ref{eq:euler})--(\ref{eq:euler-p})
in $[0,T]$ obtained in Theorem \ref{thm:local}.
If the solution $(\eta, v, P)$ satisfies
\beno
&&M(T)\eqdef \displaystyle\sup_{t\in [0,T]}\big(\|H(t)\|_{L^p\cap L^2}+\|v(t)\|_{W^{1,\infty}(\Om_t)}\big)<+\infty,\\
&&\displaystyle\inf_{(t,x,y)\in [0,T]\times \Sigma_t}-\frac {\pa P} {\pa \textbf{n}}(t,x,y)\ge c_0,\\
&&1+\eta(t,x)\ge h_0 \quad \textrm{ for }x\in \R^d,
\eeno
for some  $p>2d$ and $c_0>0, h_0>0$, then it holds that
\beno
\sup_{t\in [0,T]}E_s(t)\le C\big(E_s(0), M(T),T, c_0, h_0\big),
\eeno
Especially, the solution $(\eta, v)$ can be extended after $t=T$.
Here $H(t,x)$ is the mean curvature of the free surface and
\beno
E_s(t)\eqdef\|\eta(t)\|_{H^{s+\f12}}+\|v(t)\|_{H^{s+\f12}(\Om(t))}.
\eeno
\end{theorem}

\subsection{Main ideas}
We denote by $(V,B)$ the horizontal and vertical traces of the velocity on the free surface, i.e.,
\beno
V\triangleq (v^1,\cdots, v^d)|_{y=\eta},\quad B\triangleq v^{d+1}|_{y=\eta}.
\eeno
Introduce a good unknown $U=V+T_\zeta B$, where $\zeta=\na \eta$ and $T_\zeta$ is Bony's paraproduct.
We can derive the following evolution equation for $U$:
\ben\label{eq:U-intro}
D_t^2 U+ T_{a\lambda}U=f+f_\om.
\een
Here $D_t=\pa_t+T_V\cdot\na$, $T_{a\lambda}$ is an elliptic paradifferential operator of order one, and $f_\om$ is the nonlinear term induced by the vorticity.

Compared with the irrotational case, a main difficulty is that $f_\om$ lose one half derivative. More precisely, $f\in H^{s-\f12}$ but $f_\om \in H^{s-1}$.
Our key observation is that $D_tf_\om$ has the same regularity as $f_\om$, and $\|f_\om\|_{H^{s-1}}$ can be controlled by the lower order energy.
By using the following trick
\begin{align*}
\big\langle \langle D\rangle^{s-\f12} f_\om, \langle D\rangle^{s-\f12}D_tU\big\rangle
=&\f d {dt}\big\langle \langle D\rangle^{s-1}f_\om, \langle D\rangle^{s}U\big\rangle
-\big\langle \langle D\rangle^{s-1}D_tf_\om, \langle D\rangle^{s}U\big\rangle\\
&+\textrm{Lower order terms},
\end{align*}
we can obtain an energy inequality of the form
\ben\label{eq:E-intro}
E(t)\le \big\langle \langle D\rangle^{s-1}f_\om, \langle D\rangle^{s}U\big\rangle+E(0)+\int_0^t\cP(E(t'))dt',
\een
where $\cP$ is an increasing function. The term $\big\langle \langle D\rangle^{s-1}f_\om, \langle D\rangle^{s}U\big\rangle$ can be controlled by
\beno
\|f_\om\|_{H^{s-1}}\|U\|_{H^s}\le \cP(E_l(t))E(t)^\f12\le \cP(E_l(t))+\f12E(t).
\eeno
Here $E_l(t)$ is a lower order energy, which satisfies
\beno
E_l(t)\le E_l(0)+\int_0^t\cP(E(t'))dt'.
\eeno
This together with (\ref{eq:E-intro}) gives a close estimate for $E(t)$.

Compared with the work \cite{ABZ-IM}, a new technical ingredient is that we introduce Chemin-Lerner type Besov spaces in the elliptic estimates such that
we can obtain the maximal H\"{o}lder regularity estimates, which play an important role in the proof of break-down criterion.
Otherwise, if we just follow the framework of \cite{ABZ-IM},
it is possible to establish a similar break-down criterion, where $\|v(t)\|_{W^{1,\infty}(\Om_t)}$ is replaced by $\|v(t)\|_{C^{1,\al}(\Om_t)}$ for some $\al>0$.

For the free boundary problem, it is highly non trivial to obtain the existence of the solution from a priori estimates.
The main reason is that many special structures of the
system are used in the process of a priori estimates, however it is usually difficult to keep these structures for the approximate system.

To construct the approximate system, an immediate idea is that we use the equation (\ref{eq:U-intro}) of $U$
to construct the approximate system for the unknowns defined on the free surface. However, we find that it is
difficult to show that the limit system is equivalent to the original Euler system.
Instead, we still use the first order system to construct the iteration scheme.
In order to keep as more structures of the system as possible,
a key idea is that more unknowns and new equations are introduced such that
the structures are integrated into the new equations.
The construction of the iteration scheme are very tricky, where we also used another important observation:
the maximal regularity of the free surface is not needed for the estimate of the vorticity and the velocity,
and it is just used in the estimate of the pressure.

\section{Tools of paradifferential operators}

In this section, we introduce some basic results about the paradifferential operators from \cite{Me}(see also \cite{ABZ-IM}).

\subsection{Paradifferential operators}

Let us first introduce the definition of the symbol with limited spatial smoothness.
We denote by $W^{k,\infty}(\R^d)$ the usual Sobolev spaces for $k\in \mathbf{N}$, and
the H\"{o}lder space with exponent $k$ for $k\in (0,1)$.

\begin{definition}
Given $\mu\in [0,1]$ and $m\in \mathbf{R}$, we denote by $\Gamma^m_\mu(\mathbf{R}^d)$ the space of locally bounded functions $a(x,\xi)$ on
$\mathbf{R}^d\times \mathbf{R}^d \backslash\{0\}$, which are $C^\infty$ with respect to $\xi$ for $\xi\neq 0$ and such that, for all $\alpha\in \mathbf{N}^d$ and all
$\xi\neq 0$, the function $x\rightarrow \pa_\xi^\al a(x, \xi)$ belongs to $W^{\mu,\infty}$ and there exists a constant $C_\al$ such that
\beno
\|\pa_\xi^\al a(\cdot,\xi)\|_{W^{\mu,\infty}}\leq C_\al(1+|\xi|)^{m-|\al|}\quad \textrm{for any} \quad |\xi|\geq \f12.
\eeno
The semi-norm of the symbol is defined by
\beno
M^m_\mu(a)\eqdef\sup_{|\alpha|\leq 3d/2+1+\mu}\sup_{|\xi|\geq1/2} \|(1+|\xi|)^{|\alpha|-m}\pa^\alpha_\xi a(\cdot,\xi)\|_{W^{\mu,\infty}}.
\eeno
Especially, if $a$ is a function independent of $\xi$, then
\beno
M^m_\mu(a)=\|a\|_{W^{\mu,\infty}}.
\eeno
\end{definition}

Given a symbol $a$, the paradifferential operator $T_a$ is defined by
\ben\label{paradiff}
\widehat{T_au}(\xi)\eqdef (2\pi)^{-d}\int \chi(\xi-\eta,\eta)\widehat{a}(\xi-\eta,\eta)\psi(\eta)\widehat{u}(\eta) d \eta,
\een
where $\widehat a(\theta,\xi)$ is the Fourier transform of $a$ with respect to the first variable;
the $\chi(\theta,\xi)\in C^\infty(\R^d\times \R^d)$ is an admissible cut-off function: there exists $\e_1,\e_2$ such that $0<\e_1<\e_2$ and
\beno
\chi(\theta,\eta)=1 \quad\textrm{for}\quad |\theta|\leq \e_1 |\eta|, \quad \chi(\theta,\eta)=0 \quad\textrm{for} \quad|\theta|\geq \e_2 |\eta|,
\eeno
and such that for any $(\theta,\eta)\in \R^d\times \R^d$,
\beno
|\pa_\theta^\al \pa_\eta^\beta\chi(\theta,\eta)|\leq C_{\al,\beta}(1+|\eta|)^{-|\al|-|\beta|}.
\eeno
The cut-off function $\psi(\eta)\in C^\infty(\R^d)$ satisfies
\beno
\psi(\eta)=0 \quad\textrm{for}\quad |\eta|\leq 1, \quad \psi(\eta)=1 \quad\textrm{for}\quad |\eta|\geq 2.
\eeno
Here we will take the admissible cut-off function $\chi(\theta,\eta)$ as follows
\beno
\chi(\theta,\eta)=\sum_{k=0}^\infty\zeta_{k-3}(\theta)\varphi_k(\eta),
\eeno
where $\zeta(\theta)=1$ for $|\theta|\le 1.1$ and $\zeta(\theta)=0$ for $|\theta|\ge 1.9$; and
\beno
&&\zeta_k(\theta)=\zeta(2^{-k}\theta)\quad \textrm{for}\quad k\in\Z,\\
&&\varphi_0=\zeta,\quad \varphi_k=\zeta_k-\zeta_{k-1}\quad \textrm{for}\quad k\ge 1.
\eeno

We also introduce the Littlewood-Paley operators $\Delta_k, S_k$ defined by
\beno
&&\Delta_k u={\mathcal F}^{-1}\big(\varphi_k(\xi)\widehat u(\xi)\big)\quad\textrm{ for }k\ge 0, \quad \Delta_ku=0\quad\textrm{ for }k<0,\\
&&S_ku=\sum_{\ell\le k}\Delta_\ell u\quad\textrm{ for }k\in\Z.
\eeno
In the case when the function $a$ depends only on the first variable $x$ in $T_au$, we  take $\psi=1$. Then $T_au$ is just the usual
Bony's paraproduct  defined by
\ben\label{paraproduct}
T_au=\sum_{k}S_{k-3}a\Delta_ku.
\een
Furthermore, we have Bony's decomposition:
\ben\label{Bony}
au=T_au+T_ua+R(u,a),
\een
where the remainder term  $R(u,a)$ is defined by
\ben
&&R(u,a)=\sum_{|k-\ell|\le 2}\Delta_{k}a\Delta_\ell u.
\een

The following Berstein's inequality will be repeatedly used.

\begin{lemma}\label{lem:Berstein}
Let $1\le p\le q\le \infty, \al\in \mathbf{N}^d$. Then it holds that
\beno
&&\|\pa^\al S_k u\|_{L^q}\le C2^{kd(\f 1p-\f1q+|\al|)}\|S_ku\|_{L^p}\quad\textrm{ for }k\in\mathbf{N},\\
&&\|\Delta_ku\|_{L^q}\le C2^{kd(\f 1p-\f1q-|\al|)}\sup_{|\beta|=|\al|}\|\pa^\beta\Delta_ku\|_{L^p}\quad\textrm{ for }k\ge 1.
\eeno
\end{lemma}

\subsection{Functional spaces}

We introduce some functional spaces, which will be used throughout this paper.

\begin{definition}
Let $s\in \R$ and $p,q\in [1,\infty]$. The inhomogeneous Besov space $B^s_{p,q}(\R^d)$ is defined by
\beno
B^s_{p,q}(\R^d)\eqdef \Big\{f\in {\mathcal S}'(\R^d): \|f\|_{B^s_{p,q}}\triangleq\big(\sum_j2^{jsq}\|\Delta_jf\|_{L^p}^q\big)^\f1q<\infty\Big\}.
\eeno
\end{definition}

In the case of $p=q=2$, $B^{s}_{p,q}(\R^d)$ is just the usual Sobolev space $H^s(\R^d)$; In the case of $p=q=\infty$,
$B^{s}_{p,q}(\R^d)$ is the Zygmund space $C^s(\R^d)$.

Let $\cS=\R^d\times I$ with $I\subset \R$ an interval. We introduce the Sobolev space $H^s(\cS)$ on $\cS$.
When $s$ is an integer, $H^s(\cS)$ is just the usual Sobolev space. In general case, let $k=[s]$ and $\sigma=s-k\in (0,1)$.
The norm of $H^s(\cS)$ is defined by
\beno
\|u\|_{H^s(\cS)}\eqdef \sum_{\ell\le k}\|\na_{x,z}^\ell u\|_{L^2_z(I;H^\sigma)}
+\Big(\int_I\int_I\int_{\R^d}\f {|\na_{x,z}^k u(x,z)-\na_{x,z}^k u(x,z')|^2} {|z-z'|^{1+2\sigma}}dxdzdz'\Big)^\f12.
\eeno

In order to obtain the optimal elliptic regularity, let us introduce Chemin-Lerner type Besov space
$\widetilde{L}_z^q(I; B^s_{p,r}(\R^d))$, whose norm is defined by
\beno
\|f\|_{\widetilde{L}_z^q(I; B^s_{p,r})}\eqdef\Big(\sum_{k}2^{ksr}\|\Delta_k f\|_{L_z^q(I; L^p)}^r\Big)^\f1r.
\eeno
In the case of $p=r=\infty$, we denote it by $\widetilde{L}_z^q(I; C^s(\R^d))$; In the case of $p=q=r=2$,
we have $\widetilde{L}_z^q(I; B^s_{p,r}(\R^d))\equiv L^2_z(I;H^s(\R^d))$;
When $q=\infty, p=r=2$, we denote it by $\widetilde{L}_z^\infty(I; H^s(\R^d))$.
In this case, there holds
\beno
\|f\|_{{L}_z^\infty(I; H^s)}\le \|f\|_{\widetilde{L}_z^\infty(I; H^s)}.
\eeno
This kind of space was firstly introduced by Chemin and Lerner \cite{Che} to study the incompressible Navier-Stokes equations.

The following characterization of Sobolev space is very useful.

\begin{lemma}\label{lem:Sobolev}
Let $s\in \R$ and $c>0$. Suppose that $\{u_k\}_{k\in \mathbf{N}}$ is a sequence of functions in $L^2(\R^d)$ such that
(1)\,$u_0$ is spectrally supported in a ball $\{|\xi|\le c^{-1}\}$ and $u_k$ for $k>0$ is spectrally supported in an annulus $\{c2^k\le |\xi|\le  c^{-1}2^k\}$;
(2)\,$\{2^{ks}\|u_k\|_{L^2}\}_{k\in \mathbf{N}}\in\ell^2$. Then $u=\sum_k u_k\in H^s(\R^d)$ and
\beno
\|u\|_{H^s}\le C\Big(\sum_k2^{2ks}\|u_k\|_{L^2}^2\Big)^\f12.
\eeno
In addition, for $s>0$, it is sufficient to assume that
$u_k$ is spectrally supported in a ball $\{|\xi|\le  c^{-1}2^k\}$.
\end{lemma}

We also introduce the anisotropic Sobolev space $H^{s,\sigma}(\R^{d+1})$, whose norm is given by
\beno
\|u\|_{H^{s,\sigma}}\eqdef \|\langle D_x\rangle^\sigma u\|_{H^s(\R^{d+1})}, \quad x\in \R^d.
\eeno
We also have a similar characterization.

\begin{lemma}\label{lem:Sobolev-an}
For any $s, \sigma\in \R$, we have
\beno
\|u\|_{H^{s,\sigma}}^2\sim \sum_{\ell,j}2^{2\ell s}2^{2j\sigma}\|\Delta_\ell\Delta_j^h u\|_{L^2}^2.
\eeno
Here $\Delta_j^h$ is the Littlewood-Paley operator in the $x$ direction.
If $s,\sigma>0$, and the sequence $\big\{u_{\ell,j}\big\}_{\ell,j\in \mathbf{N}}$
is spectrally supported in $\{|\xi|\le  c^{-1}2^\ell\}\cap \{|\xi_h|\le  c^{-1}2^j\}$ for some $c>0$ and satisfies
$\|u_{\ell,j}\|_{L^2}\le c_{\ell,j}2^{-\ell s}2^{-j\sigma}$ with $\big\{c_{\ell,j}\big\}_{\ell,j\in \mathbf{N}}\in \ell^2$, then
$u=\sum_{\ell,j}u_{\ell,j}\in H^{s,\sigma}(\R^{d+1})$ satisfies
\beno
\|u\|_{H^{s,\sigma}}\le C\Big(\sum_{\ell,j}2^{2\ell s}2^{2j\sigma}\|u_{\ell,j}\|_{L^2}^2\Big)^\f12.
\eeno
\end{lemma}

\subsection{Symbolic calculus}

Let us recall the symbolic calculus and the boundedness in Sobolev space and Besov space of the paradifferential operators.

\begin{proposition}\label{prop:symbolic calculus}
Let $m, m'\in \R$.

\begin{itemize}

\item[1.]If $a\in\Gamma^m_0(\mathbf{R}^d)$, then for any $s\in\R$,
\beno
\|T_a\|_{H^s\rightarrow H^{s-m}}\leq CM_0^m(a).
\eeno

\item[2.] If $a\in\Gamma^{m}_{\rho}(\mathbf{R}^d), b \in\Gamma^{m'}_{\rho}(\mathbf{R}^d)$ for $\rho>0$, then for any $s\in\R$,
\beno
\|T_aT_b-T_{a\#b}\|_{H^s\rightarrow H^{s-m-m'+\rho}}\leq CM_\rho^{m_1}(a)M_0^{m'}(b)+CM_0^{m_1}(a)M_\rho^{m'}(b),
\eeno
where $a\#b=\sum_{|\al|<\rho}\pa_\xi^\al a(x,\xi)D_x^\al b(x,\xi), D_x=\f {\pa_x} i$.

\item[3.] If $a\in\Gamma^m_\rho(\mathbf{R}^d)$ for $\rho\in (0,1]$, then for any $s\in\R$,
\beno
\|T_{a^*}-(T_a)^*\|_{H^s\rightarrow H^{s-m+\rho}}\leq CM_\rho^m(a).
\eeno
\end{itemize}
Here $(T_a)^*$ is the adjoint operator of $T_a$, and $C$ is a constant independent of $a,b$.
\end{proposition}

\begin{proposition}\cite{WZ}\label{prop:Sym-Besov}
Let $m, m', s\in \R, q\in [1,\infty]$ and $\rho\in[0,1]$.

\begin{itemize}

\item[1.] If $a\in\Gamma^m_0(\mathbf{R}^d)$, then
\beno
\|T_a\|_{{B^s_{\infty,q}}\rightarrow B^{s-m}_{\infty,q}}\leq C M_0^m(a);
\eeno

\item[2.] If $a\in\Gamma^m_\rho(\mathbf{R}^d), b \in\Gamma^{m'}_\rho(\mathbf{R}^d)$, then
\beno
\|T_aT_b-T_{ab}\|_{B^s_{\infty,q}\rightarrow B^{s-m-m'+\rho}_{\infty,q}}\leq CM_\rho^m(a)M_0^{m'}(b)+ CM_0^m(a)M_\rho^{m'}(b).
\eeno

\end{itemize}
Here $C$ is a constant independent of $a,b$.
\end{proposition}

\begin{remark}\label{rem:symb}
If the symbol $a(x,\xi)$ satisfies
\beno
M^m_{-\mu}(a)\triangleq\sup_{|\alpha|\leq 3d/2+1}\sup_{|\xi|\geq1/2} \|(1+|\xi|)^{|\alpha|-m}\pa^\alpha_\xi a(\cdot,\xi)\|_{C^{-\mu}}<\infty
\eeno
for some $\mu>0$, then $T_a$ is bounded from $H^s(\R^d)$ to $H^{s-m-\mu}(\R^d)$ and $B^s_{\infty,q}(\R^d)$ to $B^{s-m-\mu}_{\infty,q}(\R^d)$
with the bound $M^m_{-\mu}(a)$.
\end{remark}

\begin{corol}\label{cor:symb}
Let $s, m_1,m_2,m_3\in \R$.
Suppose that $a\in \Gamma^{m_1}_{1},  b\in \Gamma^{m_2}_{2}, c\in \Gamma^{m_3}_{2}.$ Then we have
\beno
\big\|\big[T_a, [T_b, T_c]\big]\big\|_{H^s\rightarrow H^{s-m_1-m_2-m_3+2}}\le CM_{1}^{m_1}(a)M_{2}^{m_2}(b)M_{2}^{m_3}(c).
\eeno
\end{corol}

\no{\bf Proof.}\,It follows from Proposition \ref{prop:symbolic calculus} that
\beno
\|[T_b,T_c]-T_{p}+T_{p_1}\|_{H^s\rightarrow H^{s-m_2-m_3+2}}\le M_{2}^{m_2}(b)M_{2}^{m_3}(c).
\eeno
where $p(x,\xi)=\sum_{|\al|=1}\pa_\xi^\al b(x,\xi)D_x^\al c(x,\xi)$ and $p_1(x,\xi)=\sum_{|\al|=1}\pa_\xi^\al c(x,\xi)D_x^\al b(x,\xi)$.
Hence, it is sufficient to consider $[T_a, T_p]$ and $[T_a, T_{p_1}]$. Because of $p, p_1\in \Gamma_1^{m_2+m_3-1}$,
the corollary follows from Proposition \ref{prop:symbolic calculus}.\ef

\subsection{Tame estimates in Sobolev space}

Let us first recall some classical tame estimates. One can refer to \cite{BCD} for more general results.

\begin{lemma}\label{lem:remaider}
Let $s\in \R$, and $p,q\in [1,\infty]$. Then for any $\sigma>0$, we have
\beno
\|T_uv\|_{B^{s}_{p,q}}\le C\min\big(\|u\|_{C^{-\sigma}}\|v\|_{B^{s+\sigma}_{p,q}},\|u\|_{L^\infty}\|v\|_{B^{s}_{p,q}}\big).
\eeno
If $s_1+s_2>0$, then we have
\beno
\|R(u,v)\|_{H^{s_1+s_2-\f d2}}\le C\|u\|_{H^{s_1}}\|v\|_{H^{s_2}}.
\eeno
If $s>0$, then for any $\sigma\in \R$, we have
\beno
\|R(u,v)\|_{B^{s}_{p,q}}\le C\|u\|_{C^\sigma}\|v\|_{B^{s-\sigma}_{p,q}}.
\eeno
\end{lemma}

\no{\bf Proof.}\,The first two inequalities are classical, see \cite{BCD} for example.
We prove the third inequality. Recall that
\beno
R(u,v)=\sum_{|k-\ell|\le 2}\Delta_{k}u\Delta_\ell v.
\eeno
Hence, there exits some $N_0\in \textbf{N}$ so that
\beno
\Delta_jR(u,v)=\sum_{|k-\ell|\le 2;  k,\ell\ge j-N_0}\Delta_j\big(\Delta_{k}u\Delta_\ell v\big).
\eeno
It follows from Lemma \ref{lem:Berstein} that
\begin{align*}
\|\Delta_jR(u,v)\|_{L^p}\le& \sum_{|k-\ell|\le 2; k,\ell\ge j-N_0}\|\Delta_{k}u\|_{L^\infty}\|\Delta_\ell v\|_{L^p}\\
\le& C\|u\|_{C^{\sigma}}\sum_{\ell\ge j-N_0}2^{-\sigma\ell}\|\Delta_\ell v\|_{L^p}.
\end{align*}
This gives
\beno
2^{js}\|\Delta_jR(u,v)\|_{L^p}\le C\|u\|_{C^{\sigma}}\sum_{\ell\ge j-N_0}2^{-s(\ell-j)}2^{\ell(s-\sigma)}\|\Delta_\ell v\|_{L^p},
\eeno
which implies the second inequality by Young's inequality.\ef\medskip

A direct consequence of Lemma \ref{lem:remaider} is the following tame product estimate.

\begin{lemma}\label{lem:product}
Let $s\ge 0$. Then we have
\beno
&&\|fg\|_{H^s}\le C\big(\|f\|_{L^\infty}\|g\|_{H^s}+\|g\|_{L^\infty}\|f\|_{H^s}\big),\\
&&\|f\na g\|_{H^s}\le C\big(\|f\|_{L^\infty}\|\na g\|_{H^{s}}+\|g\|_{L^\infty}\|f\|_{H^{s+1}}\big).
\eeno
\end{lemma}

\begin{lemma}\cite{BCD}\label{lem:nonlinear}
Let $s>0, p,q\in [1,\infty]$ and $F$ be a smooth function with $F(0)=0$. Then we have
\beno
\|F(u)\|_{B^s_{p,q}}\le C(\|u\|_{L^\infty})\|u\|_{B^s_{p,q}}.
\eeno
Especially, for $p=q=2$, we have
\beno
\|F(u)\|_{H^s}\le C(\|u\|_{L^\infty})\|u\|_{H^s}.
\eeno
\end{lemma}

Using an extension argument, we deduce from Lemma \ref{lem:product} and Lemma \ref{lem:nonlinear} that

\begin{lemma}\label{lem:product-full}
Let $\cS=\R^d\times I$ with $I\subset \R$ an interval and $s\ge 0$. Then we have
\beno
&&\|uv\|_{H^s(\cS)}\le C\big(\|u\|_{L^\infty(\cS)}\|v\|_{H^s(\cS)}+\|v\|_{L^\infty(\cS)}\|u\|_{H^s(\cS)}\big),\\
&&\|u\na v\|_{H^s(\cS)}\le C\big(\|u\|_{L^\infty(\cS)}\|\na v\|_{H^s(\cS)}+\|v\|_{L^\infty(\cS)}\|u\|_{H^{s+1}(\cS)}\big).
\eeno
Let $F$ be a smooth function with $F(0)=0$. Then we have
\beno
\|F(u)\|_{H^s(\cS)}\le C(\|u\|_{L^\infty(\cS)})\|u\|_{H^{s}(\cS)}.
\eeno
\end{lemma}

Next we present a tame estimate in the anisotropic Sobolev space.

\begin{lemma}\label{lem:product-ans}
Let $s>0$ and $\sigma\in (0,1)$. Then it holds that for any $\epsilon>0$,
\beno
&&\|T_uv\|_{H^{s,\sigma}}+\|R(u,v)\|_{H^{s,\sigma}}\\
&&\le C\big(\|u\|_{L^\infty}+\|\langle D_x\rangle^{\sigma+\epsilon}u\|_{L^\infty_y(\R;L^2)}+\|\langle D_x\rangle^{\f d 2} u\|_{L^\infty_y(\R;L^2)}\big)
\|v\|_{H^{s,\sigma}}.
\eeno
\end{lemma}

\no{\bf Proof.}\,Using Bony's decomposition (\ref{Bony}), we write
\beno
T_uv=T(T^h+\overline{T}^h+R^h)(u,v).
\eeno
where $T(u,v)=T_uv, \overline{T}^h(u,v)=T^h(v,u)$, and $T^h_uv$ denote the paraproduct in the $x$ direction.

By using Lemma \ref{lem:Berstein} and Lemma \ref{lem:Sobolev-an}, it is easy to show that
\ben\label{eq:product-ans-1}
\|\big(TT^h+TR^h\big)(u,v)\|_{H^{s,\sigma}}\le C\|u\|_{L^\infty}\|v\|_{H^{s,\sigma}}.
\een
According to the definition of the paraproduct, we have
\beno
(T\overline{T}^h)(u,v)=\sum_{\ell,j}S_{\ell-3}\Delta_j^hu\Delta_\ell S^h_{j-3}v.
\eeno
For $\sigma<\f d 2$, we get by Lemma \ref{lem:Berstein} that
\begin{align*}
\|S_{\ell-3}\Delta_j^hu\Delta_\ell S^h_{j-3}v\|_{L^2(\R^{d+1})}
&\le \|S_{\ell-3}\Delta_j^hu\|_{L^\infty_y(\R;L^2)}\|\Delta_\ell S^h_{j-3}v\|_{L^2_y(\R;L^\infty)}\\
&\le C2^{j(\f d2-\sigma)}\|\Delta_j^hu\|_{L^\infty_y(\R;L^2)}\|\langle D_x\rangle^\sigma\Delta_\ell v\|_{L^2}\\
&\le  C2^{-j\sigma}\|\langle D_x\rangle^{\f d 2}\Delta_j^h u\|_{L^\infty_y(\R;L^2)}\|\langle D_x\rangle^\sigma\Delta_\ell v\|_{L^2}\\
&\le Cc_{\ell,j}2^{-\ell s}2^{-j\sigma}\|\langle D_x\rangle^{\f d 2} u\|_{L^\infty_y(\R;L^2)}\|v\|_{H^{s,\sigma}},
\end{align*}
and for $\sigma\ge \f d 2$ and any $\epsilon>0$, we have
\begin{align*}
\|S_{\ell-3}\Delta_j^h u\Delta_\ell S^h_{j-3} v\|_{L^2(\R^{d+1})}
&\le \|S_{\ell-3}\Delta_j^hu\|_{L^\infty_y(\R;L^2)}\|\Delta_\ell S^h_{j-3}v\|_{L^2_z(\R;L^\infty)}\\
&\le C2^{j\epsilon}\|\Delta_j^h u\|_{L^\infty_y(\R;L^2)}\|\langle D_x\rangle^{\f d 2-\epsilon}\Delta_\ell v\|_{L^2}\\
&\le Cc_{\ell,j}2^{-\ell s}2^{-j\sigma}\|\langle D_x\rangle^{\sigma+\epsilon}u\|_{L^\infty_y(\R;L^2)}
\|v\|_{H^{s,\sigma}},
\end{align*}
where $\|\{c_{\ell,j}\}\|_{\ell^2}\le 1$. Then Lemma \ref{lem:Sobolev-an} ensures that
\beno
&&\|(T\overline{T}^h)(u,v)\|_{H^{s,\sigma}}\le C\big(\|\langle D_x\rangle^{\sigma+\epsilon}u\|_{L^\infty_y(\R;L^2)}+\|\langle D_x\rangle^{\f d 2} u\|_{L^\infty_y(\R;L^2)}\big)
\|v\|_{H^{s,\sigma}},
\eeno
which together with (\ref{eq:product-ans-1}) gives the first part of the lemma. The proof of another part is similar.
\ef

Next, we give the material derivative estimates for $R(u,v)$.
\begin{lemma}\label{lem:remainder-Dt}
Let $\overline{\pa}_t \triangleq \pa_t+V\cdot\na$. Let $s>0$ and $\sigma_1\in (0,1),\sigma_2\in \R$ and $\sigma>1-\sigma_1$.
Then it holds that
\begin{align*}
\|\overline{\pa}_tR&(u(t), v(t))\|_{H^s}\le C\big(\|\overline{\pa}_t u\|_{C^{\sigma_1}}\|v\|_{H^{s-\sigma_1}}+\|\overline{\pa}_t v\|_{C^{\sigma_2}}\|u\|_{H^{s-\sigma_2}}\\
&+\|\na V\|_{L^\infty}\big(\|u\|_{C^{\sigma_1}}\|v\|_{H^{s-\sigma_1}}+\|v\|_{L^\infty}\|u\|_{H^s}\big)
+\|u\|_{C^{\sigma_1}}\|v\|_{L^\infty}\|V\|_{H^{s+\sigma}}\big),
\end{align*}
where $\|\overline{\pa}_t u\|_{C^{\sigma_1}}\|v\|_{H^{s-\sigma_1}}$ and $\|u\|_{C^{\sigma_1}}\|v\|_{H^{s-\sigma_1}}$  can also be replaced by
$\|\overline{\pa}_t u\|_{L^\infty}\|v\|_{H^{s}}$ and $\|u\|_{L^\infty}\|v\|_{H^{s}}$ respectively.
\end{lemma}

\no{\bf Proof.}\,A direct calculation gives
\begin{align*}
&\overline{\pa}_tR(u,v)=\sum_{|k-\ell|\le 2}\overline{\pa}_t(\Delta_k u\Delta_\ell V)\\
&=\sum_{|k-\ell|\le 2}\Big(\Delta_k(\overline{\pa}_tu)\Delta_\ell v+\Delta_k u\Delta_\ell(\overline{\pa}_tv)\\
&\qquad\quad-[\Delta_k,V]\cdot\na u\Delta_\ell v-\Delta_k u[\Delta_\ell, V]\cdot\na v\Big)\\
&=R\big(\overline{\pa}_t u, v\big)+R\big(u, \overline{\pa}_t v\big)
-\sum_{|k-\ell|\le 2}\big([\Delta_k,V]\cdot\na u\Delta_\ell v+\Delta_k u[\Delta_\ell, V]\cdot\na v\big).
\end{align*}
It follows from Lemma \ref{lem:remaider} that
\beno
&&\|R(\overline{\pa}_tu, v)\|_{H^s}\le C\|\overline{\pa}_t u\|_{C^{\sigma_1}}\|v\|_{H^{s-\sigma_1}},\\
&&\|R(u, \overline{\pa}_t v)\|_{H^s}\le C\|\overline{\pa}_t v\|_{C^{\sigma_2}}\|u\|_{H^{s-\sigma_2}}.
\eeno
Given $N_0\in \textbf{N}$ sufficiently large, we further decompose
\begin{align*}
\sum_{|k-\ell|\le 2}[\Delta_k,V]\cdot\na u\Delta_\ell v=&\sum_{|k-\ell|\le 2}[\Delta_k,S_{k-N_0}V]\cdot\na u\Delta_\ell v\\
&+\sum_{|k-\ell|\le 2}[\Delta_k,S^{k-N_0}V]\cdot\na u\Delta_\ell v,
\end{align*}
where $S^k=1-S_k$. By noting that $u$ has to be spectrally supported in  $\{|\xi|\sim 2^k\}$ and using the commutator estimate
\ben\label{eq:commutator}
\|[\Delta_k, u]\na v\|_{L^p}\le C\|\na u\|_{L^\infty}\|v\|_{L^p}\quad \textrm{for any } p\in [1,\infty],
\een
which follows from the identity
\begin{align*}
[\Delta_k,u]\na v(x)=&\int_{\R^d}\breve{\varphi}_k(x-x')(u(x')-u(x))\na v(x')dx'\\
=&\int_{\R^d}\na\breve{\varphi}_k(x-x')(u(x')-u(x))v(x')dx'\\
&-\int_{\R^d}\breve{\varphi}_k(x-x')\na u(x')v(x')dx',
\end{align*}
and $\|\breve{\varphi}_k\|_{L^1}+\|x\na\breve{\varphi}_k\|_{L^1}\le C$, then we can deduce that
\begin{align*}
\big\|[\Delta_k,S_{k-N_0}V]\cdot\na u\Delta_\ell v\big\|_{L^2}\le& C2^{-k\sigma_1}\|\na V\|_{L^\infty}\|u\|_{C^{\sigma_1}}\|\Delta_\ell v\|_{L^2}\\
\le& Cc_\ell2^{-\ell s}\|\na V\|_{L^\infty}\|u\|_{C^{\sigma_1}}\|v\|_{H^{s-\sigma_1}}
\end{align*}
with $\|\{c_\ell\}\|_{\ell^2}\le 1$.
Obviously, the term $[\Delta_k,S_{k-N_0}V]\cdot\na u\Delta_\ell v$ is spectrally supported in a ball $\{|\xi|\lesssim 2^k\}$ for $|k-\ell|\le 2$.
Then Lemma \ref{lem:Sobolev} ensures that
\beno
\Big\|\sum_{|k-\ell|\le 2}[\Delta_k,S_{k-N_0}V]\cdot\na u\Delta_\ell v\Big\|_{H^s}\le C\|\na V\|_{L^\infty}\|u\|_{C^{\sigma_1}}\|v\|_{H^{s-\sigma_1}}.
\eeno

Noticing that
\begin{align*}
&\|\Delta_k(S^{k-N_0}V\cdot\na u)\Delta_\ell v\|_{L^2}\le \|\Delta_k(S^{k-N_0}V\cdot\na S_ku)\Delta_\ell v\|_{L^2}\\
&\qquad\quad+\|\Delta_k(S^{k-N_0}\na\cdot V S^ku)\Delta_\ell v\|_{L^2}+\|\Delta_k\na\cdot(S^{k-N_0}VS^ku)\Delta_\ell v\|_{L^2}\\
&\le C2^{-\sigma_1k}\|\na V\|_{L^\infty}\|u\|_{C^{\sigma_1}}\|\Delta_\ell v\|_{L^2}.
\end{align*}
we can deduce from Lemma \ref{lem:Sobolev} that
\beno
\Big\|\sum_{|k-\ell|\le 2}\Delta_k(S^{k-N_0}V\cdot\na u)\Delta_\ell v\Big\|_{H^s}\le C\|\na V\|_{L^\infty}\|u\|_{C^{\sigma_1}}\|v\|_{H^{s-\sigma_1}}.
\eeno

On the other hand, we have
\begin{align*}
\Big\|\sum_{|k-\ell|\le 2}S^{k-N_0}V\cdot\na \Delta_k u\Delta_\ell v\Big\|_{H^s}
&\le C\sum_{|k-\ell|\le 2}\|S^{k-N_0}V\cdot\na \Delta_k u\Delta_\ell v\|_{H^{s}}\\
&\le C\sum_{|k-\ell|\le 2}2^{(1-\sigma_1)k}\|S^{k-N_0}V\|_{H^{s}}\|u\|_{C^{\sigma_1}}\|v\|_{L^\infty}\\
&\le C\sum_{|k-\ell|\le 2}2^{-k(\sigma-1+\sigma_1)}\|V\|_{H^{s+\sigma}}\|u\|_{C^{\sigma_1}}\|v\|_{L^\infty}\\
&\le C\|V\|_{H^{s+\sigma}}\|u\|_{C^{\sigma_1}}\|v\|_{L^\infty}.
\end{align*}
This shows that
\begin{align*}
\Big\|\sum_{|k-\ell|\le 2}[\Delta_k,S^{k-N_0}V]\cdot\na u\Delta_\ell v\Big\|_{H^s}
\le& C\big(\|\na V\|_{L^\infty}\|u\|_{C^{\sigma_1}}\|v\|_{H^{s-\sigma_1}}\\
&\quad+\|V\|_{H^{s+\sigma}}\|u\|_{C^{\sigma_1}}\|v\|_{L^\infty}\big).
\end{align*}
Obviously, $\|u\|_{C^{\sigma_1}}\|v\|_{H^{s-\sigma_1}}$ can be replaced by $\|u\|_{L^\infty}\|v\|_{H^{s}}$ in the above proof.

In a similar way, we can deduce that
\begin{align*}
\Big\|\sum_{|k-\ell|\le 2}\Delta_k u[\Delta_\ell, V]\cdot\na v\Big\|_{H^s}
\le& C\big(\|\na V\|_{L^\infty}\|v\|_{L^\infty}\|u\|_{H^s}\\
&\qquad+\|V\|_{H^{s+\sigma}}\|v\|_{L^\infty}\|u\|_{C^{\sigma_1}}\big).
\end{align*}
This completes the proof of the lemma.\ef

\subsection{Tame estimates in Chemin-Lerner spaces}

Let us recall the following lemmas from \cite{WZ}.

\begin{lemma}\label{lem:PP-Hs}
Let $s\in \R$ and $q,q_1,q_2\in [1,\infty]$ with $\f 1 q=\f1 {q_1}+\f1 {q_2}$. Then for any $s_1>0$, we have
\beno
\|T_gf\|_{\widetilde{L}^q_z(I;H^s)}\le C\min\big(\|g\|_{\widetilde{L}^{q_1}_z(I;L^\infty)}\|f\|_{\widetilde{L}^{q_2}_z(I;H^s)},\|g\|_{\widetilde{L}^{q_1}_z(I;C^{-s_1})}\|f\|_{\widetilde{L}^{q_2}_z(I;H^{s+s_1})}\big).
\eeno
\end{lemma}

\begin{lemma}\label{lem:PP-Holder}
Let $s\in \R$ and $q,q_1,q_2,r\in [1,\infty]$ with $\f 1 q=\f1 {q_1}+\f1 {q_2}$. Then for any $s_1>0$, we have
\beno
\|T_gf\|_{\widetilde{L}^q_z(I;B^s_{\infty,r})}\le C\min\big(\|g\|_{\widetilde{L}^{q_1}_z(I;L^\infty)}\|f\|_{\widetilde{L}^{q_2}_z(I;B^s_{\infty,r})},\|g\|_{\widetilde{L}^{q_1}_z(I;C^{-s_1})}\|f\|_{\widetilde{L}^{q_2}_z(I;B^{s+s_1}_{\infty,r})}\big).
\eeno
\end{lemma}

\begin{lemma}\label{lem:PR}
Let $q,q_1,q_2,r\in [1,\infty]$ with $\f 1 q=\f1 {q_1}+\f1 {q_2}$. Then for any $s>0$ and $s_1\in R$, we have
\beno
&&\|R(f,g)\|_{\widetilde{L}^q_z(I;H^s)}\le C\|g\|_{\widetilde{L}^{q_1}_z(I;C^{s_1})}\|f\|_{\widetilde{L}^{q_2}_z(I;H^{s-s_1})},\\
&&\|R(f,g)\|_{\widetilde{L}^q_z(I;B^s_{\infty,r})}\le C\|g\|_{\widetilde{L}^{q_1}_z(I;C^{s_1})}\|f\|_{\widetilde{L}^{q_2}_z(I;B^{s-s_1}_{\infty,r})}.
\eeno
If $s\le 0$ and $s_1+s_2>0$, then we have
\beno
&&\|R(f,g)\|_{\widetilde{L}^q_z(I;H^s)}\le C\|g\|_{\widetilde{L}^{q_1}_z(I;C^{s_1})}\|f\|_{\widetilde{L}^{q_2}_z(I;H^{s_2})},\\
&&\|R(f,g)\|_{\widetilde{L}^q_z(I;B^s_{\infty,r})}\le C\|g\|_{\widetilde{L}^{q_1}_z(I;C^{s_1})}\|f\|_{\widetilde{L}^{q_2}_z(I;C^{s_2})}.
\eeno
\end{lemma}


\subsection{Commutator estimates}

\begin{lemma}\label{lem:commu-Ds}
Let $m, \mu\in \R, s>0$ and $a\in \Gamma_\rho^m(\R^d)$ with $\rho\in (0,1]$. Then there holds
\beno
\big\|[\Ds, T_a]u\big\|_{H^\mu}\le CM_\rho^m(a)\|u\|_{H^{s+\mu+m-\rho}}.
\eeno
\end{lemma}

\no{\bf Proof.}\,We write
\beno
[\Ds, T_a]u=T_{\langle\xi\rangle^s}T_au-T_{\langle\xi\rangle^sa}u+(\Ds-T_{\langle\xi\rangle^s})T_au,
\eeno
then the proposition follows from Proposition \ref{prop:symbolic calculus} and the fact that
\beno
\|(\Ds-T_{\langle\xi\rangle^s})T_au\|_{H^\mu}\le C\|\Ds(1-\psi(D))T_au\|_{H^{\mu}}\le CM_0^m(a)\|u\|_{H^{-\mu'}}
\eeno
for any $\mu'>0$.\ef

%

\begin{lemma}\label{lem:comm-full}
Let $\cS=\R^d\times I$ with $I\subset \R$ an interval. Then it holds that for any integer $s\ge 1$,
\beno
\|[\na^s_{x,y}, v]\na_{x,y}u\|_{L^2(\cS)}\le C\big(\|\na_{x,y}v\|_{L^\infty(\cS)}\|u\|_{H^s(\cS)}+\|u\|_{L^\infty(\cS)}\|\na_{x,y}v\|_{H^s(\cS)}\big).
\eeno
\end{lemma}

\no{\bf Proof.}\,When $\cS=\R^{d+1}$, this lemma is classical. General case can be deduced by using an extension argument.\ef

%
%

\begin{proposition}\label{prop:commutator-tame}
Let $V\in C([0,T]; B^1_{\infty,1}(\mathbf{R}^d))$ and $p=p(t,x,\xi)$ be homogenous in $\xi$ of order $m$. Then for any $s\ge 0$,
\begin{align*}
&\big\|[T_p,\pa_t+T_V\cdot \nabla ]u(t)\big\|_{H^s}\\
&\leq C\Big(M^m_0(p)\|V(t)\|_{B^1_{\infty,1}}+M^m_{0}\big((\pa_t+T_V\cdot\nabla) p\big)\Big)\|u(t)\|_{H^{s+m}}.
\end{align*}
If $p(t)\in \Gamma^m_{\mu}$ for some $\mu>0$, $M^m_0(p)\|V(t)\|_{B^1_{\infty,1}}$ can be replaced by $M^m_{\mu}(p)\|V(t)\|_{W^{1,\infty}}$
in the case of $s>0$.
\end{proposition}

\no{\bf Proof.}\,The proof is motivated by Lemma 2.16 in \cite{ABZ-IM}. As in \cite{ABZ-IM}, it suffices to consider the case when $p=p(t, x)$ by decomposing $p$ into a sum of spherical harmonic.
In this case, $M_0^0(p)=\|p\|_{L^\infty}$.
A direct calculation gives
\begin{align*}
[\pa_t +T_V\cdot\nabla,T_p]u
=&T_{\pa_t p}u+T_V\cdot T_{\nabla p}u+T_V\cdot T_p\nabla u-T_pT_V\cdot\nabla u\\
=&T_{(\pa_t+T_V\cdot\nabla)p}u-\big(T_{T_V\cdot\nabla p}-T_{V}\cdot T_{\nabla p}\big)u\\
&+\big(T_V\cdot T_p\nabla u-T_pT_V\cdot\nabla u\big).
\end{align*}

We infer from Lemma \ref{lem:remaider} that
\ben\label{eq:comm-est1}
\|T_{(\pa_t+T_V\cdot\nabla)p}u\|_{H^s}\le C\|(\pa_t+T_V\cdot\nabla)p\|_{L^\infty}\|u\|_{H^{s}}.
\een

Recalling the definition $S^ju=u-S_ju$, we decompose $T_V\cdot T_{\nabla p}u$ as
\begin{align*}
&T_V\cdot T_{\nabla p}u=\sum_{j,k} S_{j-3}V\cdot\Delta_j\big(S_{k-3}(\nabla p)\Delta_k u\big)\\
&=\sum_k S_{k-3}V\cdot S_{k-3}(\nabla p)\Delta_k u+\sum_{j,k}\big(S_{j-3}V-S_{k-3}V\big)\cdot\Delta_j\big(S_{k-3}(\nabla p)\Delta_k u\big)\\
&=\sum_k S_{k-3}(V\cdot\nabla p)\Delta_k u-\sum_kS_{k-3}(S^{k-3}V\cdot\na p)\Delta_k u\\
&\quad-\sum_k[S_{k-3},S_{k-3}V]\cdot\na p\Delta_k u+\sum_{j,k}\big(S_{j-3}V-S_{k-3}V\big)\cdot\Delta_j\big(S_{k-3}(\nabla p)\Delta_k u\big)\\
&\triangleq T_{V\cdot\na p}u+I_1+I_2+I_3.
\end{align*}

We get by Lemma \ref{lem:Berstein} that
\beno
&&\|\big(S_{j-3}V-S_{k-3}V\big)\cdot\Delta_j\big(S_{k-3}(\nabla p)\Delta_k u\big)\|_{L^2}\\
&&\quad\le C2^{-k(s-1)}\|S_{j-3}V-S_{k-3}V\|_{L^\infty}\|p\|_{L^\infty}\|u\|_{H^s}.
\eeno
Note that the summation index $(j,k)$ in $I_3$ satisfies $|k-j|\le N_0$ for some $N_0\in \mathbf{N}$,
and $\big(S_{j-3}V-S_{k-3}V\big)\cdot\Delta_j\big(S_{k-3}(\nabla p)\Delta_k u\big)$ is spectrally supported in a ball $\{|\xi|\lesssim 2^j\}$.
Then Lemma \ref{lem:Sobolev} ensures that for $s>0$,
\begin{align*}
\|I_3\|_{H^s}\le C\|V\|_{W^{1,\infty}}\|p\|_{L^\infty}\|u\|_{H^s}
\end{align*}
by using the fact that
\beno
\sum_{|j-k|\le N_0}2^k\|S_{j-3}V-S_{k-3}V\|_{L^\infty}\le C\|V\|_{W^{1,\infty}}.
\eeno
For $s=0$, we have
\begin{align*}
\|I_3\|_{L^2}\le C\|V\|_{B^1_{\infty,1}}\|p\|_{L^\infty}\|u\|_{L^2}.
\end{align*}

By Lemma \ref{lem:Berstein} again, we have
\begin{align*}
\|S_{k-3}(S^{k-3}V\cdot\na p)\Delta_k u\|_{L^2}\le& C2^{-ks}c_k\big(\|S^{k-3}(\na V)\|_{L^\infty}+2^k\|S^{k-3}V\|_{L^\infty}\big)\|p\|_{L^\infty}\|u\|_{H^s}\\
\le& C2^{-ks}c_k\|V\|_{W^{1,\infty}}\|p\|_{L^\infty}\|u\|_{H^s},
\end{align*}
with $\|\{c_k\}\|_{\ell^2}\le 1$, and $S_{k-3}(S^{k-3}V\cdot\na p)\Delta_k u$
is spectrally supported in an annulus $\{|\xi|\sim 2^k\}$. Then Lemma \ref{lem:Sobolev} ensures that
\begin{align*}
\|I_1\|_{H^s}\le C\|p\|_{L^\infty}\|V\|_{W^{1,\infty}}\|u\|_{H^s}.
\end{align*}

Noticing that  $[S_{k-3},S_{k-3}V]\cdot\na p\Delta_k u$ is spectrally supported in a ball $\{|\xi|\lesssim 2^k\}$
and by (\ref{eq:commutator}),
\beno
\|[S_{k-3},S_{k-3}V]\cdot\na p\Delta_k u\|_{L^2}\le Cc_k2^{-ks}\|V\|_{W^{1,\infty}}\|p\|_{L^\infty}\|u\|_{H^s},
\eeno
with $\|\{c_k\}\|_{\ell^2}\le 1$, we infer from Lemma \ref{lem:Sobolev} that for $s>0$,
\begin{align*}
\|I_2\|_{H^s}\le C\|p\|_{L^\infty}\|V\|_{W^{1,\infty}}\|u\|_{H^s}.
\end{align*}

In the case of $s=0$, we need to decompose $I_2$ as
\beno
I_2=\sum_k[S_{k-3},S_{k-N_0}V]\cdot\na p\Delta_k u+\sum_k[S_{k-3},(S_{k-3}-S_{k-N_0})V]\cdot\na p\Delta_k u,
\eeno
where we take $N_0$ big enough so that $[S_{k-3},S_{k-N_0}V]\cdot\na p\Delta_k u$ is spectrally supported in an annulus $\{|\xi|\sim 2^j\}$.
Then we have
\beno
\|I_2\|_{L^2}\le C\|p\|_{L^\infty}\|V\|_{B^1_{\infty,1}}\|u\|_{L^2}.
\eeno

Putting the estimates of $I_1, I_2$ and $I_3$ together, we deduce that
\beno
&&\|\big(T_V\cdot T_{\nabla p}-T_{V\cdot\nabla p}\big)u\|_{H^s}\le C\|p\|_{L^\infty}\|V\|_{W^{1,\infty}}\|u\|_{H^s}(s>0),\\
&&\|\big(T_V\cdot T_{\nabla p}-T_{V\cdot\nabla p}\big)u\|_{L^2}\le C\|p\|_{L^\infty}\|V\|_{B^1_{\infty,1}}\|u\|_{L^2}.
\eeno

By Bony's decomposition (\ref{Bony}), we have
\begin{align*}
\|T_{V\cdot\nabla p-T_V\cdot\na p}u\|_{H^s}&\le \|V\cdot\nabla p-T_V\cdot\na p\|_{L^\infty}\|u\|_{H^s}\\
&\le C\|p\|_{L^\infty}\|V\|_{B^{1}_{\infty,1}}\|u\|_{H^s},
\end{align*}
and if $p\in C^\mu$, we have
\begin{align*}
\|T_{V\cdot\nabla p-T_V\cdot\na p}u\|_{H^s}\le C\|p\|_{C^\mu}\|V\|_{W^{1,\infty}}\|u\|_{H^s}.
\end{align*}
This shows that for $s>0$
\ben\label{eq:comm-est2}
\|\big(T_{T_V\cdot\nabla p}-T_{V}\cdot T_{\nabla p}\big)u\|_{H^s}\le C\min\big(\|p\|_{L^\infty}\|V\|_{B^{1}_{\infty,1}},\|p\|_{C^\mu}\|V\|_{W^{1,\infty}}\big)\|u\|_{H^s},
\een
and for $s=0$,
\ben\label{eq:comm-est3}
\|\big(T_{T_V\cdot\nabla p}-T_{V}\cdot T_{\nabla p}\big)u\|_{H^s}\le \|p\|_{L^\infty}\|V\|_{B^{1}_{\infty,1}}\|u\|_{H^s}.
\een

Next we decompose $T_V\cdot T_p \nabla u$ as
\begin{align*}
T_V\cdot T_p \nabla u=&\sum_{j,k}S_{j-3}V\Delta_j\big(S_{k-3}p\cdot\Delta_k\nabla u\big)\\
=&\sum_{k}S_{k-3}VS_{k-3}p\cdot\Delta_k\nabla u\\
&+\sum_{j,k}\big(S_{j-3}V-S_{k-3}V\big)\Delta_j\big(S_{k-3}p\cdot\Delta_k\nabla u\big)\\
\triangleq&\sum_{k}S_{k-3}VS_{k-3}p\cdot\Delta_k\nabla u+II_1.
\end{align*}
On the other hand, we have
\beno
&&T_pT_V\cdot \nabla u=\sum_j S_{j-3}p\Delta_j\big(S_{j-3}V\cdot \nabla u\big)\\
&&\qquad+\sum_{j,k} S_{j-3}p\Delta_j\big((S_{k-3}-S_{j-3})V\cdot\nabla \Delta_k u\big)\\
&&= \sum_j S_{j-3}VS_{j-3}p\cdot\Delta_j \nabla u +\sum_jS_{j-3}p[\Delta_j,S_{j-3}V]\cdot\nabla u\\
&&\quad+\sum_{j,k} S_{j-3}p\cdot\Delta_j\big((S_{k-3}-S_{j-3})V\cdot\nabla \Delta_k u\big)\\
&&\triangleq \sum_j S_{j-3}VS_{j-3}p\cdot \na \Delta_j u+II_2+II_3.
\eeno
This gives
\beno
T_V\cdot T_p \nabla u-T_pT_V\cdot \nabla u=II_1-II_2-II_3.
\eeno

Notice that the summation index $(j,k)$ in $II_3$ satisfies $|k-j|\le N_0$ for some $N_0\in \mathbf{N}$.
Similar to $I_3$, we can deduce that
\beno
\|II_3\|_{H^s}\le C \|p\|_{L^\infty}\|V\|_{W^{1,\infty}}\|u\|_{H^s}.
\eeno
For $II_1$, we have
\beno
&&\|II_1\|_{H^s}\le C \|p\|_{L^\infty}\|V\|_{W^{1,\infty}}\|u\|_{H^s}(s>0),\\
&&\|II_1\|_{L^2}\le C\|p\|_{L^\infty}\|V\|_{B^1_{\infty,1}}\|u\|_{L^2}.
\eeno

We decompose $II_2$ as
\begin{align*}
II_2=&\sum_jS_{j-3}p\big(\Delta_j(T_{S_{j-3}V}\cdot\na u)-T_{S_{j-3}V}\cdot\na\Delta_j u\big)\\
&+\sum_jS_{j-3}p\Delta_j\big(S_{j-3}V\cdot\na u-T_{S_{j-3}V}\cdot\na u\big)\\
&-\sum_jS_{j-3}p\big(S_{j-3}V\cdot\na\Delta_j u-T_{S_{j-3}V}\cdot\na\Delta_j u\big)\\
\triangleq& II_2^1+II_2^2+II_2^3.
\end{align*}
Similar to $I_2$, we have
\beno
&&\|II_2^1\|_{H^s}\le C \|p\|_{L^\infty}\|V\|_{W^{1,\infty}}\|u\|_{H^s}(s>0),\\
&&\|II_2^1\|_{L^2}\le C \|p\|_{L^\infty}\|V\|_{B^1_{\infty,1}}\|u\|_{L^2}.
\eeno
Using Bony's decomposition (\ref{Bony}), we decompose $II_2^2$ as
\beno
II_2^2=\sum_jS_{j-3}p\Delta_j\big(R(S_{j-3}V,\na u)+T_{S^{j-N_0}\na u}\cdot {S_{j-3}V}\big)
\eeno
for some $N_0\in \mathbf{N}$. It is easy to prove that for $s>0$,
\beno
\|S_{j-3}p\Delta_j\big(R(S_{j-3}V,\na u)+T_{S^{j-N_0}\na u}\cdot {S_{j-3}V}\big)\|_{L^2}
\le Cc_j2^{-js}\|p\|_{L^\infty}\|V\|_{W^{1,\infty}}\|u\|_{H^s}
\eeno
with $\|\{c_j\}\|_{\ell^2}\le 1$. Then Lemma \ref{lem:Sobolev} ensures that for $s>0$,
\beno
\|II_2^2\|_{H^s}\le C \|p\|_{L^\infty}\|V\|_{W^{1,\infty}}\|u\|_{H^s}.
\eeno
For $s=0$, it is obvious that
\beno
\|II_2^2\|_{L^2}\le C \|p\|_{L^\infty}\|V\|_{B^1_{\infty,1}}\|u\|_{L^2}.
\eeno
Similarly, we have
\beno
&&\|II_2^3\|_{H^s}\le C \|p\|_{L^\infty}\|V\|_{W^{1,\infty}}\|u\|_{H^s}(s>0),\\
&&\|II_2^3\|_{L^2}\le C \|p\|_{L^\infty}\|V\|_{B^1_{\infty,1}}\|u\|_{L^2}.
\eeno
This proves that
\beno
&&\|T_V\cdot T_p\nabla u-T_p T_V\cdot\nabla u\|_{H^s}\le C\|p\|_{L^\infty}\|V\|_{W^{1,\infty}}\|u\|_{H^s}(s>0),\\
&&\|T_V\cdot T_p\nabla u-T_p T_V\cdot\nabla u\|_{L^2}\le C\|p\|_{L^\infty}\|V\|_{B^1_{\infty,1}}\|u\|_{L^2},
\eeno
which together with (\ref{eq:comm-est1})-(\ref{eq:comm-est3}) give the proposition.\ef

\begin{remark}\label{rem:comm}
If the symbol $p(t,x,\xi)$ satisfies
\beno
M^m_{-\mu}(\pa_t+T_V\cdot\na p)\triangleq\sup_{|\alpha|\leq 3d/2+1}\sup_{|\xi|\geq1/2} \|(1+|\xi|)^{|\alpha|-m}\pa^\alpha_\xi(\pa_t+T_V\cdot\na p)\|_{C^{-\mu}}<\infty
\eeno
for some $\mu>0$, then there holds
\begin{align*}
&\big\|[T_p,\pa_t+T_V\cdot \nabla ]u(t)\big\|_{H^s}\\
&\leq CM^m_0(p)\|V(t)\|_{B^1_{\infty,1}}\|u(t)\|_{H^{s+m}}+CM^m_{-\mu}\big((\pa_t+T_V\cdot\nabla) p\big)\|u(t)\|_{H^{s+m+\mu}},
\end{align*}
which can be seen from (\ref{eq:comm-est1}).
\end{remark}

The following proposition will be used in the proof of well-posedness. Although the estimate is very rough, it is sufficient for our application.

\begin{proposition}\label{prop:commutator-Dt}
Assume that $V\in C([0,T]; W^{2,\infty}), \pa_t V\in C([0,T]; W^{1,\infty})$, and the symbol $p=p(t,x,\xi)$ is homogenous in $\xi$ of order $m\in \R$.
Then it holds that
\begin{align*}
\big\|(\pa_t+T_V\cdot \nabla)&[T_p,\pa_t+T_V\cdot \nabla ]u(t)\big\|_{H^s}\le C\big(1+\|V(t)\|_{W^{2,\infty}}+\|\pa_tV(t)\|_{W^{1,\infty}}\big)^2\\
&\times\big(M^m_2(p)+M^m_1(\pa_tp)+M^m_0(\pa_t^2p)\big)\|\big((\pa_t+T_V\cdot \nabla)  u(t), u(t)\big)\|_{H^{m+s}}.
\end{align*}
\end{proposition}

\no{\bf Proof.}\,Note that
\beno
[T_p,\pa_t+T_V\cdot \nabla ]u(t)=-\big(T_{\pa_tp}+T_VT_{\nabla p}\big) u-\big(T_VT_p-T_pT_V\big)\nabla u.
\eeno
A direct calculation gives
\beno
&&(\pa_t+T_V\cdot \nabla)[T_p,\pa_t+T_V\cdot \nabla]u(t)\\
&&=-(\pa_t+T_V\cdot \nabla)\big((T_{\pa_tp}+T_VT_{\nabla p}) u+(T_VT_p-T_pT_V)\nabla u\big)\\
&&=-\big[\pa_t+T_V\cdot \nabla,T_{\pa_tp}+T_VT_{\nabla p}\big] u+ (T_{\pa_tp}+T_VT_{\nabla p})(\pa_t+T_V\cdot \nabla)u\\
&&\quad-\big[\pa_t+T_V\cdot \nabla,(T_VT_p-T_pT_V)\big](\nabla u)-[T_p,T_V](\pa_t+T_V\cdot \nabla)(\nabla u).
\eeno
It follows from Lemma \ref{lem:remaider} that
\beno
&&\big\|(T_{\pa_tp}+T_VT_{\nabla p})(\pa_t+T_V\cdot \nabla)u\big\|_{H^s}\\
&&\le C\big(M^m_0(\partial_tp)+\|V(t)\|_{L^\infty}M^m_1(p)\big)\|(\pa_t+T_V\cdot \nabla)u\|_{H^{m+s}}.
\eeno
We decompose
\beno
&&\big[\pa_t+T_V\cdot \nabla,T_{\pa_tp}+T_VT_{\nabla p}\big]u\\
&&=\big[\pa_t+T_V\cdot \nabla,T_{\pa_tp}\big]u+\big[\pa_t+T_V\cdot \nabla,T_V\big]T_{\na p}u+T_V\big[\pa_t+T_V\cdot \nabla,T_{\na p}\big]u.
\eeno
Then it follows from Proposition \ref{prop:commutator-tame} and Proposition \ref{prop:symbolic calculus} that
\beno
&&\big\|\big[\pa_t+T_V\cdot \nabla,T_{\pa_tp}+T_VT_{\nabla p}\big]u\|_{H^s}\\
&&\le C\big(1+\|V(t)\|_{C^{1+\e}}+\|\pa_t V\|_{L^\infty}\big)^2\big(M^m_0(\na_{t,x}p)+M^m_0(\nabla_{t,x}^2p)\big)\|u(t)\|_{H^{m+s}}.
\eeno

We decompose
\beno
[T_p,T_V](\pa_t+T_V\cdot \nabla)(\nabla u)=[T_p,T_V]\na (\pa_t+T_V\cdot \nabla)u-[T_p,T_V]T_{\na V}\cdot \na u,
\eeno
which along with Proposition \ref{prop:symbolic calculus} gives
\beno
&&\big\|[T_p,T_V](\pa_t+T_V\cdot \nabla)(\nabla u)\big\|_{H^s}\\
&&\le C\|V(t)\|_{W^{1,\infty}}M_1^m(p)\big(\|(\pa_t+T_V\cdot \nabla)u\|_{H^{m+s}}+\|\na V(t)\|_{L^{\infty}}\|u\|_{H^{m+s}}\big).
\eeno

We write
\beno
\big[\pa_t,(T_VT_p-T_pT_V)\big](\nabla u)=[T_{\pa_tV},T_p](\na u)+[T_{V},T_{\pa_t p}](\na u),
\eeno
from which and Proposition \ref{prop:symbolic calculus}, we deduce that
\beno
\big\|[\pa_t,(T_VT_p-T_pT_V)](\nabla u)\big\|_{H^s}
\le C\|(V,\pa_t V)\|_{W^{1,\infty}}\big(M_1^m(p)+M_1^m(\pa_t p)\big)\|u\|_{H^{m+s}}.
\eeno

Finally, we deduce from Corollary \ref{cor:symb} that
\beno
\big\|\big[T_V\cdot \nabla,[T_V,T_p]\big](\nabla u)\big\|_{H^s}\le C\|V\|_{W^{2,\infty}}M_2^m(p)\|u\|_{H^{m+s}}.
\eeno
Putting the above estimates together, the proposition is proved.\ef

%

\section{Parabolic evolution equation}

Let $I=[z_0, z_1]$ be a finite interval. We denote by $\G_\r^m(I\times \R^d)$ the space of symbols $a(z;x,\xi)$ satisfying
\beno
{\cM}^m_\rho(a)\eqdef\sup_{z\in I}\sup_{|\alpha|\leq 3d/2+1+\rho}\sup_{|\xi|\geq1/2} \|(1+|\xi|)^{|\alpha|-m}\pa^\alpha_\xi a(z;\cdot,\xi)\|_{W^{\rho,\infty}}
<+\infty.
\eeno

We consider the parabolic evolution equation
\ben\label{eq:parabolic evolution}
\left\{
\begin{array}{l}
\pa_z w+T_a w=f,\\
{w|_{z=z_0}=w_0,}
\end{array}
\right.
\een
where the symbol $a\in \G^1_\r(I\times \R^d)$ is elliptic in the sense that there exists $c_1>0$ such that for any $z\in I, (x,\xi)\in \R^d\times \R^d$,
it holds that
\ben
\textrm{Re}\,a(z;x,\xi)\geq c_1|\xi|.
\een

For the elliptic estimates, we need to introduce two kinds of spaces. The first kind of spaces are intended for Sobolev elliptic estimates:
\beno
&&X^\sigma(I)\eqdef  \widetilde{L}^\infty_{z}(I;H^\sigma(\mathbf{R}^d))\cap L^2_{z}(I;H^{\sigma+\f12}(\mathbf{R}^d)),\\
&&Y^\sigma(I)\eqdef \widetilde{L}^1_{z}(I;H^\sigma(\mathbf{R}^d))+L^2_{z}(I;H^{\sigma-\f12}(\mathbf{R}^d)).
\eeno
The second kind of spaces are intended for H\"{o}lder elliptic estimates:
\beno
&&{X}^{\sigma}_{p,q}(I)\eqdef \widetilde{L}_z^p(I;B^{\sigma+\f1 p}_{\infty,q}),\\
&&{Y}^{\sigma}_q(I)\eqdef\widetilde{L}_z^1(I;B^{\sigma}_{\infty,q})+\widetilde{L}_z^2(I;B^{\sigma-\f12}_{\infty,q})+\widetilde{L}_z^\infty(I;B^{\sigma-1}_{\infty,q}).
\eeno

Let us recall the following parabolic estimates in Chemin-Lerner type spaces, which have been essentially proved in \cite{WZ}.

\begin{proposition}\label{prop:parabolic estimate}
Let $a\in \Gamma^1_\rho(I\times \R^d)$ for some $\rho>0$. Assume that $w_0\in H^\sigma$ and $f\in Y^\sigma(I)$ for $\sigma\in \R$.
Then there exists a unique solution $w\in X^\sigma(I)$ of (\ref{eq:parabolic evolution}) satisfying
\beno
\|w\|_{X^\sigma(I)} \leq C\big({\cM}^1_{\rho}(a)\big)\big( \|w_0\|_{H^\sigma}+\|f \|_{Y^\sigma(I)}\big),
\eeno
where $C$ is an increasing function depending on $c_1$ and $|I|$.
\end{proposition}

\begin{proposition}\label{prop:parabolic estimate-maximal}
Let $\sigma\in \R, p, q\in [1,\infty]$ and $a\in \Gamma^1_\rho(I\times \R^d)$ for some $\rho>0$. Assume that $w\in X^\sigma_{p,q}(I)$ is a solution of (\ref{eq:parabolic evolution}).
Then for any $\delta>0$,
\beno
\|w\|_{X^\sigma_{p,q}(I)}\leq C\big({\cM}^1_{\rho}(a)\big)
\big( \|w_0\|_{B^\sigma_{\infty,q}}+\|f\|_{Y^\sigma_{p,q}(I)}+\|w\|_{\widetilde{L}_z^p(I;C^{-\delta})}\big),
\eeno
where $Y^\sigma_{p,q}(I)=\widetilde{L}_z^p(I;B^{\sigma-1+\f1p}_{\infty,q})$ and $C$ is an increasing function depending on $c_1$ and $|I|$.
\end{proposition}

\section{Elliptic estimates in a strip}

The goal of this section is to establish the elliptic estimates in Sobolev spaces and Besov spaces.
These estimates will be used to estimate the Dirichlet-Neumann operator, the velocity and the pressure.

Throughout this section, we assume that $\eta\in C^{\f32+\e}(\R^d)$ for some $\e>0$ and there exists some $h_0>0$ such that
\ben\label{ass:height}
1+\eta(x)\ge h_0\quad \textrm{ for }x\in \R^d.
\een
Let $\cS\triangleq\big\{(x,y): x\in \R^d,-1< y<\eta(x)\big\}$ be a strip.

In the sequel, we denote by $K_\eta=C(\|\eta\|_{C^{\f 32+\e}})$ an increasing function depending on $h_0$,
which may change from line to line; $\na=(\pa_{x_1},\cdots,\pa_{x_d}), D=\f \na i$ and $\Delta=\Delta_x$.

\subsection{Elliptic boundary problem}

We consider the elliptic equations in $\cS$:
\begin{equation}\label{eq:elliptic}
\Delta_{x,y} \phi=g\quad\textrm{ in }\quad \cS,
\end{equation}
with the Dirichlet boundary condition
\ben\label{eq:elliptic-D}
\phi|_{y=\eta}=f,\quad \phi|_{y=-1}=f_b,
\een
or the mixed boundary condition
\ben\label{eq:elliptic-N}
\phi|_{y=\eta}=f,\quad \pa_y\phi|_{y=-1}=0.
\een

Given $f,f_b\in H^\f12(\R^d)$ and $g\in H^{-1}(\cS)$, the existence and uniqueness of weak solution
for the elliptic equation (\ref{eq:elliptic})-(\ref{eq:elliptic-D}) or (\ref{eq:elliptic})-(\ref{eq:elliptic-N})
can be deduced by using Riesz theorem(see \cite{ABZ-IM} for example). Moreover, the solution $\phi\in H^1(\cS)$ satisfies
\ben\label{eq:elliptic-H1}
\|\phi\|_{H^1(\cS)}\le C(\|\eta\|_{W^{1,\infty}},h_0)\big(\|(f,f_b)\|_{H^\f12}+\|g\|_{H^{-1}(\cS)}\big).
\een

As a warm-up, we first establish the elliptic estimates in the flat strip $\cS=[-1,0]\times \R^d$ by using the parabolic estimates.
The same ideas will be used to deal with the general case in the next subsections. In this subsection, we denote that $I=[-1,0]$.

\begin{proposition}\label{prop:elliptic-flat-Hs}
Let $\phi$ be a solution of (\ref{eq:elliptic})-(\ref{eq:elliptic-D}). Then for any $\sigma\ge -\f12$, there holds
\begin{eqnarray*}
&&\|\nabla_{x,y}\phi\|_{X^{\sigma}(I)}\leq C\big(\|\na_{x,y}\phi\|_{L^2(\cS)}+\|(f,f_b)\|_{H^{\sigma+1}}+\|g\|_{Y^{\sigma}(I)}\big),\\
&&\|\na_{x,y}\phi\|_{H^{\sigma+\f12}({\cS})}
\leq C\big(\|\na_{x,y}\phi\|_{L^2(\cS)}+\|(f,f_b)\|_{H^{\sigma+1}}+\|g\|_{H^{\sigma-\f12}({\cS})}\big).
\end{eqnarray*}
\end{proposition}

\no{\bf Proof.} We split the proof into three steps.

{\bf Step 1.} Estimates on  [-1/2, 0]

Let $w=\chi(y)(\pa_y-|D|)\phi$, where $\chi(y)$ is a smooth function satisfying $\chi(-1)=0$ and $\chi(y)=1$ for $y\in [-\f12,0]$.
Then $w$ satisfies
\beno
(\pa_y+|D|)w=\chi(y)g+\chi'(y)(\pa_y+|D|)\phi\triangleq g_1,\quad w(-1)=0.
\eeno
Then for any $y\in [-1,0]$, we have
\beno
w(x,y)=\int_{-1}^ye^{(y'-y)|D|}g_1(x,y')dy'.
\eeno
Using the fact(see \cite{Zhang} for example) that there exists $c>0$ such that for any $y\le 0$ and $p\in [1,\infty]$, there holds
\ben\label{eq:heat}
\|e^{y|D|}\Delta_jf\|_{L^p}\le Ce^{cy2^j}\|\Delta_j f\|_{L^p},
\een
we can easily deduce that for $\sigma_1=\min(\sigma, \f12)$,
\beno
\|w\|_{X^{\sigma_1}(I)}\le C\big(\|g\|_{Y^\sigma(I)}+\|\na_{x,y}\phi\|_{L^2(\cS)}\big).
\eeno
On the other hand, we have
\beno
(\pa_y-|D|)\phi=w(y)\quad \textrm{for }y\in [-\f12,0],\quad \phi(0)=f.
\eeno
Hence, for $y\in [-\f12,0]$,
\beno
\phi(x,y)=e^{y|D|}f+\int_0^ye^{-(y'-y)|D|}w(y')dy',
\eeno
from which and (\ref{eq:heat}), we infer that
\begin{align}
\|\na_{x,y}\phi\|_{X^{\sigma_1}([-\f12,0])}\le& C\big(\|f\|_{H^{\sigma_1+1}}+\|w\|_{X^{\sigma_1}(I)}\big)\nonumber\\
\le& C\big(\|f\|_{H^{\sigma_1+1}}+\|g\|_{Y^\sigma(I)}+\|\na_{x,y}\phi\|_{L^2(\cS)}\big).\label{eq:elliptic-est1-1}
\end{align}

 { \bf Step 2.} Estimates on  $[-1, -1/2]$

If we let $w=\chi(y)(\pa_y+|D|)\phi$, where $\chi(y)$ is a smooth function satisfying $\chi(0)=1$ and $\chi(y)=1$ for $y\in [-1,-\f12]$.
Then we have
\beno
&&(\pa_y-|D|)w=\chi(y)g+\chi'(y)(\pa_y-|D|)\phi,\quad w(0)=0,\\
&&(\pa_y+|D|)\phi=w\quad \textrm{for }y\in [-1,-\f12],\quad \phi(-1)=f_b.
\eeno
Then a similar argument leading to (\ref{eq:elliptic-est1-1}) gives
\begin{align*}
\|\na_{x,y}\phi\|_{X^{\sigma_1}([-1,-\f12])}
\le C\big(\|f_b\|_{H^{\sigma_1+1}}+\|g\|_{Y^\sigma(I)}+\|\na_{x,y}\phi\|_{L^2(\cS)}\big),
\end{align*}
which along with (\ref{eq:elliptic-est1-1}) gives
\begin{align}
\|\na_{x,y}\phi\|_{X^{\sigma_1}(I)}
\le C\big(\|(f, f_b)\|_{H^{\sigma_1+1}}+\|g\|_{Y^\sigma(I)}+\|\na_{x,y}\phi\|_{L^2(\cS)}\big).\label{eq:elliptic-est1-2}
\end{align}
A bootstrap argument will ensure that the inequality (\ref{eq:elliptic-est1-2}) holds for $\sigma_1=\sigma$.

{\bf Step 3}. Proof of the second result

To show the second inequality, it suffices to consider the estimate of the normal derivative.
Let $\sigma+\f12=k+\sigma_1$ for $k\in \textbf{N}$ and $\sigma_1\in [0,1)$. Using the interpolation inequality and the elliptic equation,
we deduce that
\begin{align*}
\|\na^k\pa_y\phi\|_{H^{\sigma_1}(\cS)}\le& C\big(\|\pa_y\phi\|_{L^2(I;H^{\sigma+\f12})}+\|\pa_y^2\phi\|_{L^2(I;H^{\sigma-\f12})}\big)\\
\le& C\big(\|\na_{x,y}\phi\|_{L^2(I;H^{\sigma+\f12})}+\|g\|_{H^{\sigma-\f12}({\cS})}\big)\\
\le& C\big(\|(f, f_b)\|_{H^{\sigma_1+1}}+\|g\|_{H^{\sigma-\f12}({\cS})}+\|\na_{x,y}\phi\|_{L^2(\cS)}\big).
\end{align*}
Let us assume that for $\ell\in [1,k]$, there holds
\begin{align*}
\|\na^{k-\ell+1}\pa_y^{\ell}\phi\|_{H^{\sigma_1}(\cS)}
\le C\big(\|(f, f_b)\|_{H^{\sigma_1+1}}+\|g\|_{H^{\sigma-\f12}({\cS})}+\|\na_{x,y}\phi\|_{L^2(\cS)}\big).
\end{align*}
Then the estimate $\|\pa_y^{k+1}\phi\|_{H^{\sigma_1}(\cS)}$ can be deduced from
\begin{align*}
\|\pa_y^{k+1}\phi\|_{H^{\sigma_1}(\cS)}\le& C\big(\|\pa_y^{k-1}\Delta\phi\|_{L^2(I;H^{\sigma_1})}+\|g\|_{H^{\sigma-\f12}({\cS})}\big)\\
\le& C\big(\|\pa_y\phi\|_{L^2(I;H^{\sigma+\f12})}+\|\pa_y^k\phi\|_{L^2(I;H^{\sigma_1+1})}+\|g\|_{H^{\sigma-\f12}({\cS})}\big),
\end{align*}
and an induction assumption. The proof is finished.\ef

\medskip

Next, we give the estimates of $\phi$ in the Besov space.

\begin{proposition}\label{prop:elliptic-flat-Holder}
Let $\phi$ be a solution of (\ref{eq:elliptic})-(\ref{eq:elliptic-D}) with $g=g_1+\pa_yg_2$.
Let $\sigma\in \R, q, r\in [1,\infty]$. Then it holds that for any $\delta>0$,
\begin{align*}
\|\na_{x,y}\phi\|_{{X}^{\sigma}_{r,q}(I)}
\le C\big(\|(f,f_b)\|_{B^{\sigma+1}_{\infty,q}}+\|g_1\|_{{Y}^{\sigma}_q(I)}+\|g_2\|_{{X}^{\sigma}_{r,q}(I)}+\|\na_{x,y}\phi\|_{L^2_z(I;C^{-\delta})}\big).
\end{align*}
\end{proposition}

\no{\bf Proof.}\,Let $w=\chi(y)(\pa_y-|D|)\phi$, where $\chi(y)$ is a smooth function satisfying $\chi(-1)=0$ and $\chi(y)=1$ for $y\in [-\f12,0]$.
Then $w$ satisfies
\beno
(\pa_y+|D|)(w-\chi(y)g_2)=\chi(y)g_1+\chi'(y)(\pa_y+|D|)\phi-\chi(y)|D|g_2-\pa_y\chi g_2 \triangleq g
\eeno
with $w(-1)=0$.
Then for any $y\in [-1,0]$, we have
\beno
(w-\chi(y)g_2 )(x,y)=\int_{-1}^ye^{(y'-y)|D|}g(y')dy'.
\eeno
from which and (\ref{eq:heat}), we infer that for $\sigma_1=\min(-\delta+\f13,\sigma)$,
\begin{align*}
\|w\|_{{X}^{\sigma_1}_{r,q}(I)}
\le C\big(\|g_1\|_{{Y}^{\sigma}_q(I)}+\|g_2\|_{{X}^{\sigma}_{r,q}(I)}+\|\na_{x,y}\phi\|_{L^2_z(I;C^{-\delta})}\big).
\end{align*}
On the other hand, for $y\in [-\f12,0]$,
\beno
\phi(x,y)=e^{y|D|}f+\int_0^ye^{-(y'-y)|D|}w(y')dy',
\eeno
from which and (\ref{eq:heat}), we infer that
\begin{align*}
\|\na_{x,y}\phi\|_{X^{\sigma_1}_{r,q}([-\f12,0])}\le& C\big(\|f\|_{B^{\sigma+1}_{\infty,q}}+\|w\|_{X^{\sigma_1}_{r,q}(I)}\big)\\
\le& C\big(\|f\|_{B^{\sigma+1}_{\infty,q}}+\|g_1\|_{{Y}^{\sigma}_q(I)}+\|g_2\|_{{X}^{\sigma}_{r,q}(I)}+\|\na_{x,y}\phi\|_{L^2_z(I;C^{-\delta})}\big).
\end{align*}
Similarly, we can deduce that
\begin{align*}
\|\na_{x,y}\phi\|_{X^{\sigma_1}_{r,q}([-1,-\f12])}
\le C\big(\|f\|_{B^{\sigma+1}_{\infty,q}}+\|g_1\|_{{Y}^{\sigma}_q(I)}+\|g_2\|_{{X}^{\sigma}_{r,q}(I)}+\|\na_{x,y}\phi\|_{L^2_z(I;C^{-\delta})}\big).
\end{align*}
Then the proposition follows by using a bootstrap argument.\ef

\subsection{Flatten the boundary and paralinearization}

Motivated by \cite{Lan, ABZ-IM}, we flatten the boundary of $\cS$ by a regularized mapping:
\beno
(x,z)\in \R^d\times I\longmapsto (x,\rho_\delta(x,z))\in\cS,
\eeno
where $I=[-1,0]$ and $\rho_\delta$ with $\delta>0$ is given by
\begin{eqnarray}\label{mapping}
\rho_\delta(x,z)=z+(1+z)e^{\delta z|D|}\eta(x).
\end{eqnarray}

In the sequel, we denote $\overline{\cS}=\R^d\times I$ with $I=[-1,0]$.
For a function $f(x,y)$ defined in $\cS$, we denote $\widetilde{f}(x,z)=f(x,\rho_\delta(x,z))$.

For any $z\in [-1,0]$, we have
\begin{align*}
&\pa_z\rho_\delta=1+\eta(x)+\big(e^{\delta z|D|}-1\big)\eta(x)+(1+z)\delta e^{\delta z|D|}|D|\eta(x).
\end{align*}
It is easy to show that
\beno
&&\|\big(e^{\delta z|D|}-1\big)\eta+\delta(1+z) e^{\delta z|D|}|D|\eta\|_{L^\infty}\\
&&\le \delta\int_z^0\|e^{\delta z'|D|}|D|\eta\|_{L^\infty}dz'+\delta\|e^{\delta z|D|}|D|\eta(x)\|_{L^\infty}\\
&&\le C\delta\||D|\eta\|_{L^\infty}\le C\delta\|\eta\|_{C^{1+\e}}.
\eeno
Hence, we can take $\delta$ small enough depending only on $\|\eta\|_{C^{1+\e}}$ and $h_0$
so that
\ben
\pa_z\rho_\delta(x,z)\ge \f {h_0} 2\quad \textrm{for }(x,z)\in \overline{\cS}.
\een

We have the following regularity information for the regularized map $\rho_\delta$,
which can be easily verified by the definition of $\rho_\delta$.

\begin{lemma}\label{lem:map-regularity}
For any $\sigma\ge 0$, there holds
\beno
&&\|(\nabla\rho_\delta,\pa_z\rho_\delta-1)\|_{X^{\sigma-\f12}(I)}+\|(\nabla\rho_\delta,\pa_z\rho_\delta-1)\|_{H^{\sigma}(\overline{\cS})}\leq K_\eta\|\eta\|_{H^{\sigma+\f12}}.
\eeno
Let $p, q\in [1,\infty]$ and $\sigma\in \R$. Then we have
\beno
\|(\nabla\rho_\delta,\pa_z\rho_\delta-1)\|_{\widetilde{L}^p_z(I;B^{\sigma}_{\infty,q})}\le K_\eta\|\eta\|_{B^{\sigma+1-\f1p}_{\infty,q}}.
\eeno
\end{lemma}

We set $v(x,z)=\phi(x,\rho_\delta(x,z))$. It is easy to find that $v$ satisfies
\begin{equation}\label{eq:ellitic-flat}
\pa^2_z v+\al \Delta v+\beta\cdot\nabla \pa_z v-\gamma \pa_z v= F_{0},
\end{equation}
where $F_{0}= \alpha g(x,\rho_\delta(x,z))\triangleq F_{00}+\na_{z}G_0$ and the coefficients $\al,\beta,\gamma$ are defined by
\ben\label{eq:elliptic coefficients}
\alpha=\f{(\pa_z \rho_\delta)^2}{1+|\nabla \rho_\delta|^2},\quad \beta=-2\f{\pa_z\rho_\delta \nabla \rho_\delta}{1+|\nabla \rho_\delta|^2},
\quad\gamma=\f{1}{\pa_z\rho_\delta}(\pa_z^2\rho_\delta+\alpha\Delta\rho_\delta+\beta\cdot\nabla\pa_z\rho_\delta).
\een

We collect the following Sobolev and H\"{o}lder regularity information for the elliptic coefficients,
which can be proved by using Bony's decomposition (\ref{Bony}),
Lemma \ref{lem:PP-Hs}-Lemma \ref{lem:PR} and Lemma \ref{lem:map-regularity}(see Lemma 4.4 in \cite{WZ}).

\begin{lemma}\label{lem:coeff}
If $\sigma>\f32$, then it holds that
\beno
&&\|\alpha -1\|_{X^{\sigma-\f12}(I)}+\|\beta\|_{X^{\sigma-\f12}(I)}+\|\gamma\|_{X^{\sigma-\f32}(I)}
\leq K_\eta\|\eta\|_{H^{\sigma+\f12}},\\
&&\|\alpha-1\|_{\widetilde{L}^1_z(I; H^{\sigma+\f12})}+\|\beta\|_{\widetilde{L}^1_z(I; H^{\sigma+\f12})}
+\|\gamma\|_{\widetilde{L}^1_z(I; H^{\sigma-\f12})}\le K_\eta\|\eta\|_{H^{\sigma+\f12}},\\
&&\|\alpha\|_{\widetilde{L}_z^\infty(I; C^{\f12+\e})}+\|\beta\|_{\widetilde{L}_z^\infty(I; C^{\f12+\e})}+\|\gamma\|_{\widetilde{L}_z^\infty(I; C^{-\f12+\e})}\le K_\eta,\\
&&\|\alpha\|_{\widetilde{L}_z^2(I; C^{1+\e})}+\|\beta\|_{\widetilde{L}_z^2(I; C^{1+\e})}
+\|\gamma\|_{\widetilde{L}_z^2(I; C^{\e})}\leq K_\eta.
\eeno
\end{lemma}

\begin{lemma}\label{lem:coeff2}
It holds that
\beno
&&\|({\al}-1,{\beta})\|_{H^\sigma(\overline{\cS})} \le K_\eta\|\eta\|_{H^{\sigma+\f12}}\quad \textrm{for }\sigma\ge 0,\\
&&\|{\gamma}\|_{H^{\sigma-1}(\overline{\cS})} \le K_\eta\|\eta\|_{H^{\sigma+\f12}}\quad \textrm{for }\sigma\ge 1.
\eeno

\end{lemma}

\no{\bf Proof.}\,The lemma follows from Lemma \ref{lem:product-full} and Lemma \ref{lem:map-regularity}.\ef
\vspace{0.2cm}

In order to obtain the tame elliptic estimates, we paralinearize the elliptic equation (\ref{eq:ellitic-flat}) as
\ben\label{eq:elliptic-pl}
\pa^2_z v+T_\al\Delta v+T_\beta\cdot\nabla \pa_z v=F_{0}+F_1+F_2,
\een
with $F_1, F_2$ given by
\beno
F_1=\gamma \pa_zv,\quad F_2=(T_\al-\al)\Delta v+(T_\beta-\beta)\cdot\nabla \pa_z v.
\eeno
Following \cite{ABZ-IM}, we decouple the equation (\ref{eq:elliptic-pl}) into a forward and a backward parabolic
evolution equations:
\begin{eqnarray}\label{eq:elliptic-decouple}
(\pa_z-T_a)(\pa_z-T_A)v= F_{0 }+F_1+F_2+F_3\triangleq F,
\end{eqnarray}
where
\beno
&&a=\f12\big(-i\beta\cdot\xi-\sqrt{4\al|\xi|^2-(\beta\cdot\xi)^2}\big),\\
&&A=\f12\big(-i\beta\cdot\xi+\sqrt{4\al|\xi|^2-(\beta\cdot\xi)^2}\big),\\
&&F_3=(T_aT_A-T_\al\Delta)v-(T_a+T_A+T_\beta\cdot\nabla )\pa_zv-T_{\pa_z A}v.
\eeno

The symbols $a, A$ satisfy

\begin{lemma}\label{lem:symbol}
It holds that
\beno
&&a(z;x,\xi)\cdot A(z;x,\xi)=-\al(x,z)|\xi|^2,\\
&&a(z;x,\xi)+ A(z;x,\xi)=-i\beta(x,z)\cdot \xi.
\eeno
Moreover, $a, A\in {\cM}^1_{\f12+\e}$ with the bounds
\beno
{\cM}^1_{\f12+\e}(a)\le K_\eta,\quad {\cM}^1_{\f12+\e}(A)\le K_\eta.
\eeno
\end{lemma}

\no{\bf Proof.} The two equalities are obvious. Note that
\beno
4\al|\xi|^2-(\beta\cdot\xi)^2\ge c_2|\xi|^2
\eeno
for some $c_2>0$ depending only on $\|\eta\|_{C^{\f32+\e}}$.
Then the second statement follows from Lemma \ref{lem:nonlinear} and Lemma \ref{lem:coeff}.\ef\vspace{0.2cm}

\subsection{Elliptic estimates in Sobolev space}

We first present the trace estimate in terms of the $H^1$ norm of the solution.

\begin{lemma}\label{lem:elliptic-H1}
Assume that $v\in H^1(\overline{\cS})$ is a solution of (\ref{eq:ellitic-flat}) in $\overline{\cS}$. Then we have
\begin{align*}
&\|\na v\|_{X^{-\f12}(I)}\leq C\| v\|_{H^1(\overline{\cS})},\\
&\|\pa_z v\|_{X^{-\f12}(I)}\leq K_\eta\big(\|F_0\|_{Y^{-\f12}(I)}+\| v\|_{H^1(\overline{\cS})}\big).
\end{align*}
\end{lemma}

\no{\bf Proof.}\, Let $v_j=\Delta_jv$. Let $\chi(z)$ be a smooth function supported in $[-\f34,0]$ and $\chi(z)=1$ on $[-\f12,0]$. For any $z\in [-\f12,0]$, we have
\begin{align*}
\big\langle\nabla v_j(z),\nabla v_j(z)\big\rangle_{H^{-\f12}}
\leq& 2\int^z_{-1}\chi(z')\big\langle\pa_{z'}\nabla v_j(z'),\nabla v_j(z')\big\rangle_{H^{-\f12}}dz'\\
&+\int^z_{-1}\chi'(z')\big\langle\nabla v_j(z'),\nabla v_j(z')\big\rangle_{H^{-\f12}}dz'\\
\le& C\|\nabla_{x,z}v_j\|^2_{L^2(\overline{\cS})}+C\|v\|^2_{L^2(\overline{\cS})}
\end{align*}
which implies that
\begin{align*}
\|\na v\|_{\widetilde{L}^\infty(-\f12,0;H^{-\f12})}\leq C\| v\|_{H^1(\overline{\cS})}.
\end{align*}
By the equation (\ref{eq:ellitic-flat}) and Lemma 4.8 in \cite{WZ}, we get
\begin{align*}
\big\langle\pa_z {v_j(z)},\pa_z v_j(z)\big\rangle_{H^{-\f12}}
=&2\int^z_{-1}\chi(z')\big\langle\pa_{z'}^2 v_j(z'), \pa_z v_j(z')\big\rangle_{H^{-\f12}}dz'\\
&+\int^z_{-1}\chi'(z')\big\langle\pa_{z'}v_j(z'), \pa_z v_j(z')\big\rangle_{H^{-\f12}}dz'\\
=&2\int^z_{-1}\chi(z')\big\langle \Delta_j(F_0-\al\Delta v+\beta\nabla \pa_zv-\gamma\pa_z v),\pa_z v_j\big\rangle_{H^{-\f12}}dz'\\
&+\int^z_{-1}\chi'(z')\big\langle\pa_{z'}v_j(z'), \pa_z v_j(z')\big\rangle_{H^{-\f12}}dz'\\
\leq& Cc_j\Big(\|F_0\|_{Y^{-\f12}(I)}+\|\dv_x(\al \nabla v+\beta\pa_z v)\|_{L^2_z(I;H^{-1})}\\
&+\|\nabla\al\nabla v+\nabla\beta \pa_z v+\gamma\pa_z v\|_{\widetilde{L}^1_z(I; H^{-\f12})}+\|\nabla_{x,z}v\|_{L^2(I\times\R^d)}\Big)\|\pa_zv\|_{X^{-\f12}(I)}\\
\le& K_\eta c_j\big(\|F_0\|_{Y^{-\f12}(I)}+\|v\|_{H^1(\overline{\cS})}\big)\|\pa_zv\|_{X^{-\f12}(I)},
\end{align*}
with $\|\{c_j\}\|_{\ell^1}\le 1$, which implies that
\begin{align*}
\|\pa_z v\|_{\widetilde{L}^\infty(-\f12,0;H^{-\f12})}^2\leq K_\eta\big(\|F_0\|_{Y^{-\f12}(I)}+\|\nabla_{x,z}v\|_{L^2(\overline{\cS})}\big)\|\na_{x,z}v\|_{X^{-\f12}(I)}.
\end{align*}
The same estimates hold  for $z\in [-1,-\f12]$ by  a similar cut-off argument. This gives
\beno
\|\na_{x,z}v\|_{X^{-\f12}(I)}^2\le K_\eta\big(\|F_0\|_{Y^{-\f12}(I)}+\|\nabla_{x,z}v\|_{L^2(\overline{\cS})}\big)\|\na_{x,z}v\|_{X^{-\f12}(I)},
\eeno
which implies the desired result.\ef

\begin{remark}\label{rem:ellip-H1}
If $\eta\in H^s(\R^d)$ for $s>\f d 2+1$, then we have
\beno
\|\na_{x,z}v\|_{X^{-\f12}(I)}\leq C(\|\eta\|_{H^s}, h_0)\big(\|F_0\|_{Y^{-\f12}(I)}+\| v\|_{H^1(\overline{\cS})}\big).
\eeno
Indeed, we have $\na\al, \na \beta, \gamma\in H^{s-\f32}(\overline{\cS})$, which implies by Lemma \ref{lem:remaider} that
\beno
\|\nabla\al\nabla v+\nabla\beta \pa_z v+\gamma\pa_z v\|_{{L}^1_z(I; H^{-\f12})}\le C(\|\eta\|_{H^s}, h_0)\|\nabla_{x,z}v\|_{L^2(\overline{\cS})}.
\eeno
\end{remark}

The following elliptic estimates will be used to estimate the velocity in the proof of well-posedness. Here and in what follows, we denote
by $P_{\sigma,\eta}=C\big(\|\eta\|_{H^{\sigma}}\big)$ an increasing function depending on $h_0$, which may be different from line to line.

Let us first recall the following elliptic estimate for tangential derivatives, which has been essentially proved in \cite{ABZ-IM}.
\begin{proposition}\label{prop:elliptic-tan-nontame}
Let $v\in H^1(\overline{\cS})$ be a solution of (\ref{eq:ellitic-flat}) in $\overline{\cS}$ with $v(x,0)=f(x)$ and $v(x,-1)=f_b(x)$.
Assume that $\eta\in H^s(\R^d)$ for $s>\f d 2+\f 32$.
Then for any $\sigma\in [-\f12, s-\f12]$, it holds that
\begin{align*}
&\|\nabla_{x,z} v\|_{X^{\sigma-\f12}(I)}\leq P_{s,\eta}\big(\|\na_{x,z}v\|_{L^2(\overline{\cS})}+\|(f, f_b)\|_{H^{\sigma+\f12}}+\|F_{0}\|_{L^2_z(I:H^{\sigma-1})}\big),\\
&\|\nabla_{x,z} v\|_{X^{\sigma-\f12}([a,0])}\leq P_{s,\eta}\big(\|\na_{x,z}v\|_{L^2(\overline{\cS})}+\|f\|_{H^{\sigma+\f12}}+\|F_{0}\|_{L^2_z(I:H^{\sigma-1})}\big),
\end{align*}
for any $a\in (-1,0)$.
\end{proposition}

Now we present the elliptic estimate for full derivatives.

\begin{proposition}\label{prop:elliptic-boun}
Let $v\in H^1(\overline{\cS})$ be a solution of (\ref{eq:ellitic-flat}) in $\overline{\cS}$ with $v(x,0)=f(x)$ and $v(x,-1)=f_b(x)$.
Assume that $\eta\in H^s(\R^d)$ for $s>\f d 2+\f 32$. Then it holds that
\begin{align*}
\|\nabla_{x,z}v\|_{H^{s-\f12}{(\overline{\cS})}}\leq& P_{s,\eta}\big(\|\na_{x,z}v\|_{L^2(\overline{\cS})}+
\|(f,f_b)\|_{H^{s}}+\|F_{0}\|_{H^{s-\f32}(\overline{\cS})}\big).
\end{align*}

\end{proposition}

We need the following lemma.

\begin{lemma}\label{lem:product-app1}
Let $\sigma=s-\f12$ and $\sigma_1=\sigma-[\sigma]$. Then for any positive integer $k\le [\sigma]$, there holds
\beno
&&\|\pa_z^{k-1}(\al \Delta v+\beta\cdot\nabla \pa_z v-\gamma \pa_z v)\|_{L^2_z(I;H^{\sigma_1})}\nonumber\\
&&\le P_{s,\eta}\Big(\|\langle D\rangle^{\sigma_1} \na{v}\|_{H^{k}(\overline{\cS})}+\|\langle D\rangle^{\sigma_1+\f12}\pa_z{v}\|_{H^{k-1}(\overline{\cS})}
+\|\na_{x,z}{v}\|_{L^\infty(\overline{\cS})}\\
&&\qquad\qquad+\|\langle D\rangle^{\sigma_1+\f12}\na_{x,z}{v}\|_{L^\infty_z(I;L^2)}
+\|\langle D\rangle^{\f d 2}\na_{x,z}{v}\|_{L^\infty_z(I;L^2)}\Big).
\eeno
\end{lemma}

\no{\bf Proof.}\,For $k=1$, we have by Lemma \ref{lem:product} and Lemma \ref{lem:coeff}  that
\beno
&&\|\al \Delta v+\beta\cdot\nabla \pa_z v-\gamma \pa_z v\|_{L^2_z(I;H^{\sigma_1})}\\
&&\le C\big(\|(\al-1,\beta)\|_{L^\infty}+1\big)\|\langle D\rangle^{\sigma_1}\na {v}\|_{H^{1}(\overline{\cS})}
+\|\gamma\|_{L^\infty_z(I;C^{-\f12})}\|\langle D\rangle^{\sigma_1+\f12}\pa_z {v}\|_{L^2(\overline{\cS})}\\
&&\qquad+C\|\na_{x,z}{v}\|_{L^\infty(\overline{\cS})}\big(\|(\al-1,\beta)\|_{L^2(I;H^{1+\sigma_1})}+\|\gamma\|_{L^2(I;H^{\sigma_1})}\big)\\
&&\le P_{s,\eta}\big(\|\na_{x,z}{v}\|_{L^\infty(\overline{\cS})}+\|\langle D\rangle^{\sigma_1}\na {v}\|_{H^{1}(\overline{\cS})}
+\|\langle D\rangle^{\sigma_1+\f12}\pa_z {v}\|_{L^2(\overline{\cS})}\big).
\eeno
Here we used $\sigma_1+\f32\le s$ and $s>\f d 2+\f32$.

Next we consider the case of $k\ge 2$.
Let $(\overline{v}, \overline{\al}, \overline{\beta},\overline{\gamma})$ be an extension of $(v, \al,\beta,\gamma)$ to $\R^{d+1}$ so that
\beno
&&\|\na_{x,z}\overline{v}\|_{H^{\sigma}(\R^{d+1})}\leq C\|\na_{x,z}{v}\|_{H^{\sigma}(\overline{\cS})},\\
&&\|(\overline{\al}-1, \overline{\beta})\|_{H^{\sigma}(\R^{d+1})}\le C\|({\al}-1, {\beta})\|_{H^{\sigma}(\overline{\cS})},\\
&&\|\overline{\gamma}\|_{H^{\sigma-1}(\R^{d+1})}\le C\|{\gamma}\|_{H^{\sigma-1}(\overline{\cS})}.
\eeno
Using Bony's decomposition (\ref{Bony}), we write
\begin{align*}
\overline{\beta}\na\pa_z\overline{v}=&T_{\overline{\beta}}\na\pa_z\overline{v}+T_{\na\pa_z\overline{v}}\overline{\beta}
+R(\overline{\beta},\na\pa_z\overline{v})\\
=&T_{\overline{\beta}}\na\pa_z\overline{v}+\pa_z T_{\na\overline{v}}\overline{\beta}-T_{\na\overline{v}}\pa_z\overline{\beta}
+R(\overline{\beta},\na\pa_z\overline{v}),
\end{align*}
from which and Lemma \ref{lem:product-ans}, we infer that for any $\epsilon>0$
\begin{align*}
\|\overline{\beta}\na\pa_z\overline{v}\|_{H^{k-1,\sigma_1}}
&\le C\big(\|\overline{\beta}\|_{L^\infty}+\|\langle D\rangle^{\sigma_1+\epsilon}\overline{\beta}\|_{L^\infty_z(\R;L^2)}
+\|\langle D\rangle^{\f d 2}\overline{\beta}\|_{L^\infty_z(\R;L^2)}\big)
\|\na\overline{v}\|_{H^{k,\sigma_1}}\\
&+C\big(\|\na_{x,z}\overline{v}\|_{L^\infty}+\|\langle D\rangle^{\sigma_1+\epsilon}\na\overline{v}\|_{L^\infty_z(\R;L^2)}
+\|\langle D\rangle^{\f d 2}\na\overline{v}\|_{L^\infty_z(\R;L^2)}\big)
\|\overline{\beta}\|_{H^{k,\sigma_1}}.
\end{align*}
Due to $\sigma_1+\epsilon+1\le s$ for some $\epsilon\in (0,1]$ and $s>\f d2+\f32$, it follows from Lemma \ref{lem:coeff2} that
\beno
\|\overline{\beta}\|_{L^\infty}+\|\langle D\rangle^{\sigma_1+\epsilon}\overline{\beta}\|_{L^\infty_z(\R;L^2)}
+\|\langle D\rangle^{\f d 2}\overline{\beta}\|_{L^\infty_z(\R;L^2)}+\|\overline{\beta}\|_{H^{k,\sigma_1}}
\le P_{s,\eta}.
\eeno
This shows that
\begin{align*}
\|\langle D\rangle^{\sigma_1}({\beta}\na\pa_z{v})\|_{H^{k-1}(\overline{\cS})}\le
P_{s,\eta}\big(&\|\langle D\rangle^{\sigma_1} \na{v}\|_{H^{k}(\overline{\cS})}+\|\na_{x,z}{v}\|_{L^\infty(\overline{\cS})}\\
&+\|\langle D\rangle^{\sigma_1+\epsilon}\na{v}\|_{L^\infty_z(I;L^2)}
+\|\langle D\rangle^{\f d 2}\na{v}\|_{L^\infty_z(I;L^2)}\big).
\end{align*}
Similarly, we have
\begin{align*}
\|\langle D\rangle^{\sigma_1}(\al \Delta v)\|_{H^{k-1}(\overline{\cS})}\le P_{s,\eta}\big(&\|\langle D\rangle^{\sigma_1} \na{v}\|_{H^{k}(\overline{\cS})}
+\|\na_{x,z}{v}\|_{L^\infty(\overline{\cS})}\\
&+\|\langle D\rangle^{\sigma_1+\epsilon}\na{v}\|_{L^\infty_z(I;L^2)}
+\|\langle D\rangle^{\f d 2}\na{v}\|_{L^\infty_z(I;L^2)}\big).
\end{align*}

While, we can see from the proof of Lemma \ref{lem:product-ans} that
\begin{align*}
\|\langle D\rangle^{\sigma_1}\overline{\gamma}\pa_z\overline{v}&\|_{H^{k-1,\sigma_1}}\le C\big(\|\overline{\gamma}\|_{L^\infty_z(\R;C^{-\f12})}
+\|\langle D\rangle^{\sigma_1+\epsilon}\overline{\gamma}\|_{L^2}
+\|\langle D\rangle^{\f d 2}\overline{\gamma}\|_{L^2}\big)
\|\langle D\rangle^{\sigma_1+\f12}\pa_z\overline{v}\|_{H^{k-1}}\\
&+C\big(\|\na_{x,z}\overline{v}\|_{L^\infty}+\|\langle D\rangle^{\sigma_1+\epsilon}\pa_z\overline{v}\|_{L^\infty_y(\R;L^2)}
+\|\langle D\rangle^{\f d 2}\pa_z\overline{v}\|_{L^\infty_y(\R;L^2)}\big)\|\overline{\gamma}\|_{H^{k-1,\sigma_1}}.
\end{align*}
We know from Lemma \ref{lem:coeff} and Lemma \ref{lem:coeff2} that
\beno
\|\overline{\gamma}\|_{L^\infty_z(\R; C^{-\f12})}+\|\langle D\rangle^{\sigma_1+\epsilon}\overline{\gamma}\|_{L^2}
+\|\langle D\rangle^{\f d 2}\overline{\gamma}\|_{L^2}+\|\overline{\gamma}\|_{H^{k-1,\sigma_1}}\le P_{s,\eta}.
\eeno
Then we infer that
\begin{align*}
\|{\gamma}\pa_z{v}\|_{H^{k-1}(\overline{\cS})}
\le P_{s,\eta}\big(&\|\langle D\rangle^{\sigma_1}\pa_z{v}\|_{H^{k-\f12}(\overline{\cS})}+\|\na_{x,z}{v}\|_{L^\infty(\overline{\cS})}\\
&+\|\langle D\rangle^{\sigma_1+\epsilon}\pa_z{v}\|_{L^\infty_z(I;L^2)}
+\|\langle D\rangle^{\f d 2}\pa_z{v}\|_{L^\infty_z(I;L^2)}\big).
\end{align*}

Summing up the above estimates, we conclude the lemma.\ef\medskip

\no{\bf Proof of Proposition \ref{prop:elliptic-boun}.}\,Let $\sigma=s-\f12$. Proposition \ref{prop:elliptic-tan-nontame} ensures that
\begin{align}
&\|\nabla_{x,z} v\|_{L^\infty_z(I;H^{\sigma-\f12})}+\|\nabla_{x,z} v\|_{L^2_z(I;H^{\sigma})}\nonumber\\
&\leq P_{s,\eta}\big(\|\na_{x,z}v\|_{L^2(\overline{\cS})}+\|(f, f_b)\|_{H^{\sigma+\f12}}+\|F_{0}\|_{L^2_z(I:H^{\sigma-1})}\big).\label{eq:ell-Hs-tan2}
\end{align}
Assume that for $1\leq \ell\le [\sigma]-1$, there holds
\begin{align}
\|\nabla_{x,z}^{\ell+1}v\|_{L^2_z(I;H^{\sigma-\ell})}\leq& P_{s,\eta}\big(\|\na_{x,z}v\|_{L^2(\overline{\cS})}+
\|(f,f_b)\|_{H^{\sigma+\f12}}+\|F_{0}\|_{H^{\sigma-1}(\overline{\cS})}\big).
\label{eq:ell-Hs-induct2}
\end{align}
We will prove that the inequality holds for $\ell=[\sigma]\triangleq k$. Let $\sigma_1=\sigma-k$.
Using the equation (\ref{eq:ellitic-flat}) and Lemma \ref{lem:product-app1}, we deduce that
\begin{align}
&\|\pa_z^{k+1} v\|_{L^2_z(I;H^{\sigma_1})}\leq \|\pa_z^{k-1}(\al \Delta v+\beta\cdot\nabla \pa_z v-\gamma \pa_z v)\|_{L^2_z(I;H^{\sigma_1})}
+\|\pa_z^{k-1} F_{0}\|_{L^2_z(I;H^{\sigma_1})}\nonumber\\
&\le P_{s,\eta}\Big(\|\langle D\rangle^{\sigma_1} \na{v}\|_{H^{k}(\overline{\cS})}+\|\langle D\rangle^{\sigma_1+\f12}\pa_z{v}\|_{H^{k-1}(\overline{\cS})}
+\|\na_{x,z}{v}\|_{L^\infty_z(I;H^{\sigma-\f12})}\nonumber\\
&\qquad+\|\langle D\rangle^{\sigma_1+\f12}\na_{x,z}{v}\|_{L^\infty_z(I;L^2)}
+\|\langle D\rangle^{\f d 2}\na_{x,z}{v}\|_{L^\infty_z(I;L^2)}\Big)+\|F_{0}\|_{H^{\sigma-1}(\overline{\cS})}.\label{eq:ell-Hs-f1}
\end{align}
Note that $\sigma_1+\f12\le \sigma-\f12, \f d2<\sigma-\f12$. This together with (\ref{eq:ell-Hs-tan2}) and (\ref{eq:ell-Hs-induct2})
ensures that the inequality (\ref{eq:ell-Hs-induct2}) holds for $\ell=[\sigma]$.

Using Lemma \ref{lem:product-full}, Lemma \ref{lem:coeff2} and the interpolation, we obtain
\begin{align*}
\|\pa_z^{k+1} v\|_{H^{\sigma_1}(\overline{\cS})}\leq& \|\pa_z^{k-1}(\al \Delta v+\beta\cdot\nabla \pa_z v-\gamma \pa_z v)\|_{H^{\sigma_1}(\overline{\cS})}+\|\pa_z^{k-1} F_{0}\|_{H^{\sigma_1}(\overline{\cS})}\nonumber\\
\le& P_{s,\eta}\big(\|\na\na_{x,z}v\|_{H^{k-1+\sigma_1}(\overline{\cS})}+\|\langle D\rangle^{\f12}\pa_{z}v\|_{H^{k-1+\sigma_1}(\overline{\cS})}\big)
+\|F_{0}\|_{H^{\sigma-1}(\overline{\cS})}\\
\le&P_{s,\eta}\big(\|\na_{x,z}v\|_{L^2_z(I;H^{k+\sigma_1})}+\|\pa_{z}^{k+1} v\|_{L^2_z(I;H^{\sigma_1})}\big)
+\|F_{0}\|_{H^{\sigma-1}(\overline{\cS})},
\end{align*}
which together with (\ref{eq:ell-Hs-tan2}) and (\ref{eq:ell-Hs-f1}) implies the proposition.\ef

\subsection{Tame elliptic estimates}

To prove the break-down criterion, we need to establish the tame elliptic estimates.
In this subsection, we assume that $\eta\in H^{s+\f12}(\R^d)$ for $s>\f d 2+1$.

The following elliptic estimate will be used to estimate the Drichlet-Neumann operator and the pressure.

\begin{proposition}\label{prop:elliptic-upboun}
Let $I_0=[a,0]$ for $a\in (-1,0)$.
Assume that $v\in H^{s+1}(\overline{\cS})$ is a solution of (\ref{eq:ellitic-flat}) in $\overline{\cS}$ with $v(x,0)=f(x)$.
Then there exists an interval $I_1\subset (-1,0)$ so that for all $\sigma\in [-\f12, s-\f12]$,
\begin{eqnarray*}
\|\nabla_{x,z} v\|_{X^{\sigma}(I_0)}\leq K_\eta\big(\|\na_{x,z}v\|_{L^2(\overline{\cS})}+
\|f\|_{H^{\sigma+1}}+\|F_{00}\|_{Y^{\sigma}(I)}+\|G_0\|_{X^{\sigma}(I )}+\|\eta\|_{H^{s+\f12}}\|\nabla_{x,z}v\|_{L^{\infty}(\widetilde{\cS})}\big).
\end{eqnarray*}
Here $\widetilde{S}=\R^d\times I_1$. Moreover, if $\sigma<s-\f12$, $\|\nabla_{x,z} v\|_{L^{\infty}(\widetilde{\cS})}$ can be replaced by $\|\nabla_{x,z} v\|_{L^{\infty}_z(I_1;C^{0})}$.
\end{proposition}

The proof of the proposition need the following lemma, which can be proved by using Bony's decomposition, Lemma \ref{lem:PP-Hs} and Lemma \ref{lem:PR}.
Since the proof is almost the same as those in \cite{WZ}, we omit the details.

\begin{lemma}\label{lem:F1}
Let $I_1\subseteq I$ be an interval and $\widetilde{S}=\R^d\times I_1$. For any $\sigma\leq s-\f12$, it holds that
\beno
&&\|F_1\|_{Y^{\sigma}(I_1)}\leq K_\eta\big(\|\pa_zv\|_{L_z^2(I_1; H^{\sigma})}+\|\pa_z v\|_{L^{\infty}(\widetilde{\cS})}\|\eta\|_{H^{s+\f12}}\big),\\
&&\|F_2\|_{Y^{\sigma}(I_1)}\leq K_\eta\|\nabla_{x,z}v\|_{L^\infty(\widetilde{\cS})}\|\eta\|_{H^{s+\f12}},\\
&&\|F_3\|_{Y^{\sigma}(I_1)}\leq K_\eta\|\nabla v\|_{L_z^2(I_1; H^{\sigma})}\quad \textrm{for any}\quad\sigma\in \R.
\eeno
If $\sigma<s-\f12$, $\|\nabla_{x,z} v\|_{L^{\infty}(\widetilde{\cS})}$ can be replaced by $\|\nabla_{x,z} v\|_{L^{\infty}_z(I_1;C^{0})}$.
\end{lemma}

\medskip

\no{\bf Proof of Proposition \ref{prop:elliptic-upboun}}.
As in \cite{ABZ-IM}, the proof uses the induction argument.
Suppose that there exits $I_2=[\xi_0,0]\subseteq I_1$ with $\xi_0\in (-1,a)$ such that for some $\sigma\in [-\f12,s-\f12)$, there holds
\begin{eqnarray*}
\|\nabla_{x,z} v\|_{X^{\sigma}(I_2)}\leq K_\eta\big(\|\na_{x,z}v\|_{L^2(\overline{\cS})}+
\|f\|_{H^{\sigma+1}}+\|F_{00}\|_{Y^{\sigma}(I)}+\|G_0\|_{X^{\sigma}(I )}+\|\eta\|_{H^{s+\f12}}\|\nabla_{x,z} v\|_{L^{\infty}(\widetilde{\cS})}\big).
\end{eqnarray*}
This is indeed true for $\sigma=-\f12$ by Lemma \ref{lem:elliptic-H1}.
We will show that the above inequality still holds when $\sigma$ and $I_2$ are replaced by $\sigma+\delta_1$ and $I_3$,
where $\sigma+\delta_1\le s-\f12$ with $\delta_1\in (0,\f12]$, and $I_3=[\xi_1,0]$ with $\xi_1=\xi_0+\f12(a-\xi_0)\in (\xi_0,a)$.

Let $\chi$ be a smooth function such that $\chi(\xi_0)=0$ and $\chi(z)=1$ for $z\ge \xi_1$.
Set $w=\chi(z)(\pa_z-T_A)v-\chi(z)G_0$. Then $w$ satisfies
\beno
\pa_zw-T_aw=F',\quad w(\xi_0)=0,
\eeno
where $F'=\chi(z)\big(F_{00}+F_1+F_2+F_3\big)+\chi'(z)((\pa_z-T_A)v-G_0)+\chi(z)T_aG_0$.

We deduce from Proposition \ref{prop:symbolic calculus} and Lemma \ref{lem:symbol} that for $\delta_1\le \f12$,
\beno
\|(\pa_z-T_A)v\|_{L^2{(I_2;H^{\sigma+\delta_1})}}\le K_\eta\|\na_{x,z}v\|_{L^2{(I_2;H^{\sigma+\delta_1})}},
\eeno
which along with Lemma \ref{lem:F1} gives
\begin{align}
\|F'\|_{Y^{\sigma+\delta_1}(I_2)}\le& \|F_0\|_{Y^{\sigma+\delta_1}(I_2)}+K_\eta\|G_0\|_{X^{\sigma}(I_2)}\nonumber\\
&+K_\eta\big(\|\nabla_{x,z} v\|_{L^\infty(I_2\times \R^d)}\|\eta\|_{H^{s+\f12}}+\|\nabla_{x,z}v\|_{X^\sigma(I_2)}\big).\label{eq:F'}
\end{align}
Then it follows from Proposition \ref{prop:parabolic estimate} that
\begin{align}\label{eq:elliptic-Hs-2}
&\|w\|_{X^{\sigma+\delta_1}(I_2)}\leq K_\eta\|F\|_{Y^{\sigma+\delta_1}(I_2)}\nonumber\\
&\le K_\eta \big(\|F_0\|_{Y^{r+\delta_1}(I)}+\|G_0\|_{X^{\sigma}(I )}+\|\nabla_{x,z}v\|_{X^\sigma(I_2)}
+\|\nabla_{x,z} v\|_{L^\infty(I_2\times \R^d)}\|\eta\|_{H^{s+\f12}}\big).
\end{align}

To obtain the estimate of $v$, we consider the backward parabolic equation
\beno
\pa_zv-T_Av=w+G_0\quad \textrm{in }I_3,\quad v|_{z=0}=f.
\eeno
Take $\na$ to the equation of $v$ to get
\beno
(\pa_z-T_A)\na v=\na (w+G_0)+T_{\na A}v,\quad \na v|_{z=0}=\na f.
\eeno
By Remark \ref{rem:symb} and Lemma \ref{lem:symbol}, we have
\beno
\|T_{\na A}v\|_{L^2_z(I_3; H^{\sigma+\delta_1-\f12})}\le K_\eta\|\na v\|_{L^2_z(I_3;H^{\sigma+\delta_1})}.
\eeno
Then Proposition \ref{prop:parabolic estimate}, (\ref{eq:elliptic-Hs-2}) and the induction assumption ensure that
\begin{align*}
&\|\na v\|_{X^{\sigma+\delta_1}(I_3)}\leq K_\eta\big(\|w\|_{X^{\sigma+\delta_1}(I_3)}+\|G_0\|_{X^{\sigma+\delta_1}(I_3)}+\|f\|_{H^{\sigma+1+\delta_1}}+\|\na v\|_{L^2(I_3;H^{\sigma+\delta_1})}\big)\\
&\leq K_\eta\big(\|f\|_{H^{\sigma+1+\delta_1}}+\|F_0\|_{Y^{\sigma+\delta_1}(I)}+\|G_0\|_{X^{\sigma+\delta_1}(I_3)}+\|\nabla_{x,z}v\|_{X^\sigma(I_2)}
+\|\nabla_{x,z} v\|_{L^\infty(I_2\times \R^d)}\|\eta\|_{H^{s+\f12}}\big)\\
&\le K_\eta\big(\|\na_{x,z}v\|_{L^2(\overline{\cS})}+\|f\|_{H^{\sigma+\delta_1+1}}+\|F_0\|_{Y^{\sigma+\delta_1}(I)} +\|G_0\|_{X^{\sigma+\delta_1}(I_3)}
 +\|\nabla_{x,z}v\|_{L^\infty(I_2\times \R^d)}\|\eta\|_{H^{s+\f12}}\big).
\end{align*}
The same estimate holds for $\pa_z v$ by using $\pa_zv=T_Av+w+G_0$.

In the case when $\sigma+\delta_1<s-\f12$, $\|\nabla_{x,z} v\|_{L^{\infty}(I_2\times \R^d)}$ in (\ref{eq:F'})
can be replaced by $\|\nabla_{x,z} v\|_{L^{\infty}_z(I_2;B^{0}_{\infty,\infty})}$ by Lemma \ref{lem:F1}, and so does the final result.
\ef

\medskip

The following elliptic estimate will be used to estimate the velocity  in the proof of break-down criterion.

\begin{proposition}\label{prop:elliptic-boun-b}
Assume that $v\in H^{s+\f12}(\overline{\cS})$ is a solution of (\ref{eq:ellitic-flat}) in $\overline{\cS}$ with $v(x,0)=f(x)$ and $v(x,-1)=f_b(x)$.
Let $k=s-\f12$ be an integer. Then it holds that
\begin{align*}
\|\nabla_{x,z}v\|_{H^k(\overline{\cS})}\leq K_\eta\big(&\|\na_{x,z}v\|_{L^2(\overline{\cS})}+
\|(f,f_b)\|_{H^{s}}\\
&+\|F_{0}\|_{H^{k-1}(\overline{\cS})}+\|\eta\|_{H^{s+\f12}}\|\nabla_{x,z} v\|_{L^{\infty}(\overline{\cS})}\big).
\end{align*}
\end{proposition}

Let us first present the tame estimate for the tangential derivatives.

\begin{lemma}\label{lem:elliptic-tan}
Assume that $v\in H^{s+1}(\overline{\cS})$ is a solution of (\ref{eq:ellitic-flat}) in $\overline{\cS}$ with $v(x,0)=f(x)$ and $v(x,-1)=f_b(x)$.
Then for all $\sigma\in [-\f12, s-\f12]$ , it holds that
\begin{eqnarray*}
\|\nabla_{x,z} v\|_{X^{\sigma}(I)}\leq K_\eta\big(\|\na_{x,z}v\|_{L^2(\overline{\cS})}+
\|(f, f_b)\|_{H^{\sigma+1}}+\|F_{0}\|_{Y^{\sigma}(I)}+\|\eta\|_{H^{s+\f12}}\|\nabla_{x,z} v\|_{L^{\infty}(\overline{\cS})}\big).
\end{eqnarray*}
\end{lemma}

\no{\bf Proof}. Because the proof is similar to Proposition  \ref{prop:elliptic-upboun},  we just present a sketch.
Assume that  for some $\sigma\in [-\f12,s-\f12)$, there holds
\begin{eqnarray*}
\|\nabla_{x,z} v\|_{X^{\sigma}(I)}\leq K_\eta\big(\|\na_{x,z}v\|_{L^2(\overline{\cS})}+
\|f\|_{H^{\sigma+1}}+\|F_{0}\|_{Y^{\sigma}(I)}+\|\eta\|_{H^{s+\f12}}\|\nabla_{x,z} v\|_{L^{\infty}(\overline{\cS})}\big).
\end{eqnarray*}
We show that the above inequality still holds when $\sigma$  is replaced by $\sigma+\delta_1$,
where $\sigma+\delta_1\le s-\f12$ with $\delta_1\in (0,\f12]$.

Let $\chi$ be a smooth function such that $\chi(-1)=0$ and $\chi(z)=1$ for $z\in [-\f12,0]$.
Set $w=\chi(z)(\pa_z-T_A)v$. Then $(w,v)$ satisfies
\beno
&&\pa_zw-T_aw=\widetilde{F},\quad w(-1)=0,\\
&&\pa_zv-T_Av=w,\quad v(0)=f,
\eeno
where $\widetilde{F}=\chi(z)\big(F_0+F_1+F_2+F_3\big)+\chi'(z)(\pa_z-T_A)v$.
Using the same argument of Proposition \ref{prop:elliptic-upboun}, we can deduce that
\begin{align*}
&\|\na_{x,z}v\|_{X^{\sigma+\delta_1}([-\f12,0])}\\
&\le K_\eta\big(\|\na_{x,z}v\|_{L^2(\overline{\cS})}+\|f\|_{H^{\sigma+\delta_1+1}}+\|F_0\|_{Y^{\sigma+\delta_1}(I)}
 +\|\nabla_{x,z}v\|_{L^\infty(\overline{\cS})}\|\eta\|_{H^{s+\f12}}\big).
\end{align*}

On the other hand, the equation (\ref{eq:elliptic-pl}) can also be decoupled into a  backward and forward parabolic
evolution equations, i.e.,
\beno
(\pa_z-T_A)(\pa_z-T_a)v= F_{0 }+F_1+F_2+\widetilde{F}_3,
\eeno
where
\beno
\widetilde{F}_3=(T_AT_a-T_\al\Delta)v-(T_a+T_A+T_\beta\cdot\nabla )\pa_zv-T_{\pa_z a}v.
\eeno
Now we let $\chi$ be a smooth function such that $\chi(0)=0$ and $\chi(z)=1$ for $z\in [-1,-\f12]$.
Set $w=\chi(z)(\pa_z-T_a)v$. Then $(w,v)$ satisfies
\beno
&&\pa_zw-T_Aw=\widetilde{F},\quad w(0)=0,\\
&&\pa_zv-T_av=w,\quad v(-1)=f_b,
\eeno
A similar argument ensures that
\begin{align*}
&\|\na_{x,z}v\|_{X^{\sigma+\delta_1}([-1-\f12])}\\
&\le K_\eta\big(\|\na_{x,z}v\|_{L^2(\overline{\cS})}+\|f_b\|_{H^{\sigma+\delta_1+1}}+\|F_0\|_{Y^{\sigma+\delta_1}(I)}
 +\|\nabla_{x,z}v\|_{L^\infty(\overline{\cS})}\|\eta\|_{H^{s+\f12}}\big).
\end{align*}
This together with the estimates in the interval $[-\f12,0]$ gives the lemma.\ef
\medskip

\no{\bf Proof of Proposition \ref{prop:elliptic-boun-b}.}\,Lemma \ref{lem:elliptic-tan} ensures that
\begin{align}
\|\nabla_{x,z} v\|_{L^\infty(I;H^{k-\f12})}+&\|\nabla_{x,z} v\|_{L^2(I;H^{k})}\leq K_\eta\big(\|\na_{x,z}v\|_{L^2(\overline{\cS})}+\|(f, f_b)\|_{H^{s}}\nonumber\\
&+\|F_{0}\|_{L^2_z(I:H^{k-1})}+\|\eta\|_{H^{s+\f12}}\|\nabla_{x,z}v\|_{L^{\infty}(\overline{\cS})}\big).\label{eq:ell-Hs-tan}
\end{align}
Let us assume that for $0\le \ell \le k-1$, there holds
\begin{align}
\|\nabla_{x,z}^{\ell+1} v\|_{L^2(I;H^{k-\ell})}\leq& K_\eta\Big(\|\na_{x,z}v\|_{L^2(\overline{\cS})}+
\|(f,f_b)\|_{H^{s}}\nonumber\\
&\quad+\sum_{\ell'\le \ell-1}\|\nabla_{x,z}^{\ell'}F_{0}\|_{L^2(\overline{\cS})}+\|\eta\|_{H^{s+\f12}}\|\nabla_{x,z} v\|_{L^{\infty}(\overline{\cS})}\Big).
\label{eq:ell-Hs-induct}
\end{align}

Next we show that the inequality holds for $\ell=k$.
Using the equation (\ref{eq:ellitic-flat}), we get
\ben\label{eq:ell-v-normal}
 \|\pa_z^{k+1} v\|_{L^2(\overline{\cS})}\leq \|\pa_z^{k-1}(\al \Delta v+\beta\cdot\nabla \pa_z v-\gamma \pa_z v)\|_{L^2(\overline{\cS})}
 +\|\pa_z^{k-1} F_{0}\|_{L^2(\overline{\cS})}.
\een

As in the proof of Proposition \ref{prop:elliptic-boun},
let $(\overline{v}, \overline{\al}, \overline{\beta},\overline{\gamma})$ be an extension of $(v, \al,\beta,\gamma)$ to $\R^{d+1}$
keeping the corresponding Sobolev norm.
We infer from Lemma \ref{lem:product} that
\begin{align*}
\|(\overline{\al}-1)\Delta \overline{v}\|_{H^{k-1}}\le& C\|\overline{\al}-1\|_{L^\infty}\|\na^2\overline{v}\|_{H^{k-1}}
+C\|\na\overline{v}\|_{L^\infty}\|\overline{\al}-1\|_{H^{k}}\\
\le&C\|{\al}-1\|_{L^\infty(\overline{\cS})}\|\na^2{v}\|_{H^{k-1}(\overline{\cS})}
+C\|\na{v}\|_{L^\infty(\overline{\cS})}\|{\al}-1\|_{H^{k}(\overline{\cS})}
\end{align*}
Similarly, we have
\begin{align*}
\|\overline{\beta}\na\pa_z\overline{v}\|_{H^{k-1}}
\le C\|{\beta}\|_{L^\infty(\overline{\cS})}\|\na\na_{x,z}{v}\|_{H^{k-1}(\overline{\cS})}
+C\|\na_x{v}\|_{L^\infty(\overline{\cS})}\|{\beta}\|_{H^{k}(\overline{\cS})}.
\end{align*}

Using Bony's decomposition (\ref{Bony}), we write
\beno
\overline{\gamma}\pa_z\overline{v}=T_{\overline{\gamma}}\pa_z\overline{v}+T_{\pa_z\overline{v}}\overline{\gamma}+R(\pa_z\overline{v},\overline{\gamma}).
\eeno
It follows from Lemma \ref{lem:Berstein} that
\begin{align*}
\|S_{j-3}\overline{\gamma}\Delta_j(\pa_z\overline{v})\|_{L^2}\le& C\|\overline{\gamma}\|_{L^\infty_z(\R;C^{-\f12})}
\|\langle D\rangle^\f12\Delta_j(\pa_z\overline{v})\|_{L^2}\\
\le& Cc_j2^{-j(k-1)}\|\overline{\gamma}\|_{L^\infty_z(\R;C^{-\f12})}\|\langle D\rangle^\f12\pa_z\overline{v}\|_{H^{k-1}}
\end{align*}
with $\|\{c_j\}\|_{\ell^2}\le 1$. Then Lemma \ref{lem:Sobolev} ensures that
\begin{align*}
\|T_{\overline{\gamma}}\pa_z\overline{v}\|_{H^{k-1}}\le& C\|\overline{\gamma}\|_{L^\infty_z(\R;C^{-\f12})}\|\langle D\rangle^\f12\pa_z\overline{v}\|_{H^{k-1}}.
\end{align*}
By Lemma \ref{lem:remaider}, we have
\beno
\|T_{\pa_z\overline{v}}\overline{\gamma}\|_{H^{k-1}}\le C\|\pa_z{v}\|_{L^\infty(\overline{\cS})}\|{\gamma}\|_{H^{k-1}(\overline{\cS})},
\eeno
and for $k>1$, we have
\beno
\|R(\pa_z\overline{v},\overline{\gamma})\|_{H^{k-1}}\le C\|\pa_z{v}\|_{L^\infty(\overline{\cS})}\|{\gamma}\|_{H^{k-1}(\overline{\cS})}.
\eeno
This proves that for $k>1$,
\beno
&&\|\overline{\gamma}\pa_z\overline{v}\|_{H^{k-1}}\le C\|\pa_z{v}\|_{L^\infty(\overline{\cS})}\|{\gamma}\|_{H^{k-1}(\overline{\cS})}
+C\|{\gamma}\|_{L^\infty_z(I;C^{-\f12})}\|\langle D\rangle^\f12\pa_z{v}\|_{H^{k-1}(\overline{\cS})}.
\eeno
For $k=1$, it is obvious that
\beno
\|\overline{\gamma}\pa_z\overline{v}\|_{L^2}\le
C\|{\gamma}\|_{L^2_z(I;L^\infty)}\|\pa_z{v}\|_{L^\infty_z(I;L^2)}.
\eeno

The above estimates together with Lemma \ref{lem:coeff} and Lemma \ref{lem:coeff2} ensure that
\beno
&&\|\pa_z^{k-1}(\al \Delta v+\beta\cdot\nabla \pa_z v-\gamma \pa_z v)\|_{L^2(\overline{\cS})}\\
&&\le C\|\na_{x,z}{v}\|_{L^\infty(\overline{\cS})}\big(\|(\al-1, {\beta})\|_{H^{k}(\overline{\cS})}+\|{\gamma}\|_{H^{k-1}(\overline{\cS})}\big)
+C\|({\al},\beta)\|_{L^\infty(\overline{\cS})}\|\na\na_{x,z}{v}\|_{H^{k-1}(\overline{\cS})}\\
&&\qquad+C\|{\gamma}\|_{L^\infty_z(I;C^{-\f12})}\|\langle D\rangle^\f12\pa_z{v}\|_{H^{k-1}(\overline{\cS})}+C\|{\gamma}\|_{L^2_z(I;L^\infty)}^2\|\pa_z{v}\|_{L^\infty_z(I;L^2)}\\
&&\le C\|\na_{x,z}{v}\|_{L^\infty(\overline{\cS})}\|\eta\|_{H^{s}}
+K_\eta\big(\|\na\na_{x,z}{v}\|_{H^{k-1}(\overline{\cS})}+\|\langle D\rangle^\f12\pa_z{v}\|_{H^{k-1}(\overline{\cS})}+\|\pa_z{v}\|_{L^\infty_z(I;L^2)}\big),
\eeno
from which and (\ref{eq:ell-v-normal}), we infer that
\begin{align*}
\|\pa_z^{k+1}& v\|_{L^2(\overline{\cS})}\le C\|\na_{x,z}{v}\|_{L^\infty(\overline{\cS})}\|\eta\|_{H^{s}}\\
&+K_\eta\big(\|\na\na_{x,z}{v}\|_{H^{k-1}(\overline{\cS})}+\|\langle D\rangle^\f12\pa_z{v}\|_{H^{k-1}(\overline{\cS})}+\|\pa_z{v}\|_{L^\infty_z(I;L^2)}\big)+\|\pa_z^{k-1}F_{0}\|_{L^2(\overline{\cS})}.
\end{align*}
which along with (\ref{eq:ell-Hs-tan}) and (\ref{eq:ell-Hs-induct}) implies the proposition by the interpolation.
\ef

\subsection{Elliptic estimates in Besov space}

In this subsection, we present the elliptic estimates in Besov space which will be used in the proof of break-down criterion.

\begin{proposition}\label{prop:elliptic Holder est}
Let $I_0=[a,0]$ and $I_1=[b,0]$ with $b<a$. Let $q\in [1,\infty], r\in [0,\f12]$.
Assume that $v\in X^r_q(I_1)$ is a solution of the elliptic equation
\begin{equation}
\pa^2_z v+\al \Delta v+\beta\cdot\nabla \pa_z v-\gamma \pa_z v= F_{00}+\pa_zG_0\quad \textrm{in  }\overline{\cS}\nonumber
\end{equation}
with $v(x,0)=f(x)$. Then it holds that for any $\delta>0$,
\begin{align*}
\|\na_{x,z} v\|_{{X}^{r}_{q}(I_0)}
\le K_\eta\Big(&\|f\|_{B^{r+1}_{\infty,q}}+\|F_{00}\|_{{Y}^{r}_q(I)}+\|G_0\|_{{X}^{r}_{q}(I)}\\
&+\|\na_{x,z}v\|_{{L}^\infty_z(I_0; C^{-\delta})}+\|\na_{x,z}v\|_{\widetilde{L}^2_z(I_1\setminus I_0; C^{-\f \e 2})}
\Big).
\end{align*}
Here ${X}^{r}_{q}(I_0)=\widetilde{L}_z^\infty(I_0;B^{r}_{\infty,q})$.
\end{proposition}

Let us recall the following lemma from \cite{WZ}, which can be proved by using Bony's decomposition (\ref{Bony}),
Lemma \ref{lem:PP-Holder} and Lemma \ref{lem:PR}.

\begin{lemma}\label{lem:F1-2}
Let $I_0=[a,0]$ and $I_1=[b,0]$ with $b<a$.
Then for any $r\leq \f12$ and $q\in [1,\infty]$, it holds that
\begin{align*}
\|F_1\|_{\widetilde{L}_z^2(I_1; B^{r-\f12}_{\infty,q})}
\leq& K_\eta\big(\|\nabla_{x,z} v\|_{\widetilde{L}_z^\infty(I_0; C^{-\f\e 2})}+\|\pa_z v\|_{\widetilde{L}_z^\infty(I_0; B^{r-\f12}_{\infty,q})}\\
&\qquad+\|\nabla_{x,z} v\|_{\widetilde{L}_z^2(I_1\setminus I_0; C^{-\f\e 2})}\big),\\
\|F_2\|_{\widetilde{L}_z^2(I_1; B^{r-\f12}_{\infty,q})}\leq& K_\eta\big(\|\nabla_{x,z} v\|_{\widetilde{L}_z^\infty(I_0; C^{-\f\e 2})}+\|\nabla_{x,z} v\|_{\widetilde{L}_z^2(I_1\setminus I_0; C^{-\f\e 2})}\big),\\
\|F_3\|_{\widetilde{L}_z^2(I_1; B^{r-\f12}_{\infty,q})}\leq& K_\eta\big(\|\nabla_{x,z} v\|_{\widetilde{L}_z^\infty(I_0; C^{-\f\e 2})}+\|\nabla_{x,z} v\|_{\widetilde{L}_z^2(I_1\setminus I_0; C^{-\f\e 2})}\big).
\end{align*}
\end{lemma}
\medskip

\no{\bf Proof of Proposition \ref{prop:elliptic Holder est}.}\,Let $\chi$ be a smooth function supported in $I_1$
and $\chi(z)=1$ for $z\in I_0$. Set $w=\chi(z)(\pa_z-T_A)v$. Then $(w,v)$ satisfies
\beno
&&\pa_zw-T_aw=\widetilde{F},\quad w(b)=0,\\
&&\pa_zv-T_Av=w\quad \textrm{in }[a,0],\quad v(0)=f,
\eeno
where $\widetilde{F}=\chi(z)\big(F_{00}+\pa_zG_0+F_1+F_2+F_3\big)+\chi'(z)(\pa_z-T_A)v$.

Let $w_1=w-\chi(z)G_0$. Then $w_1$ satisfies
\beno
\pa_zw_1-T_aw_1=\widetilde{F}-\pa_z(\chi(z)G_0)-\chi(z)T_aG_0,\quad w(b)=0.
\eeno
Then Proposition \ref{prop:parabolic estimate-maximal} together with Lemma \ref{lem:F1-2} and Proposition \ref{prop:Sym-Besov} ensures that
\begin{align*}
\|w-G_0\|_{{X}^{r}_{2,q}(I_0)}
\le& K_\eta\big(\|\widetilde{F}-\chi(z)\pa_z G_0\|_{{Y}^{r}_q(I_1)}+\|G_0\|_{{X}^{r}_q(I_1)}+\|w\|_{\widetilde{L}_z^2(I_1; C^{-\delta})}\big)\nonumber\\
\le& {K_\eta}\big(\|F_{00}\|_{{Y}^{r}_q(I)}+\|G_0\|_{{X}^{r}_q(I_1)}+\|\nabla_{x,z} v\|_{\widetilde{L}_z^\infty(I_0;B^{r-\f12}_{\infty,q})}\\
&\qquad+\|\nabla_{x,z} v\|_{\widetilde{L}_z^\infty(I_0; C^{-\f \e2})}+\|\nabla_{x,z} v\|_{\widetilde{L}_z^2(I_1\setminus I_0; C^{-\f \e2})}\big).
\end{align*}
Note that $(\pa_z-T_A)\na v=\na (w-G_0)+\na G_0+T_{\na A}v$ on $I_0$. It follows from Proposition \ref{prop:parabolic estimate-maximal},
Lemma \ref{lem:symbol} and Remark \ref{rem:symb} that
\begin{align*}
\|\na v\|_{{X}^{r}_q(I_0)}&\leq K_\eta\big(\|f\|_{B^{r+1}_{\infty,q}}+\|w-G_0\|_{\widetilde{L}_z^2(I_0;B^{r+\f12}_{\infty,q})}+\|G_0\|_{{X}^{r}_q(I_1)} +\|\na_{x,z}v\|_{\widetilde{L}_z^\infty(I_0;B^{r-\f12}_{\infty,q})}
+\|\na_{x,z}v\|_{\widetilde{L}^\infty_z(I_0; C^{-\delta})}\big)\\
&\le K_\eta\big(\|f\|_{B^{r+1}_{\infty,q}}+\|F_{00}\|_{{Y}^{r}_q(I)}+\|G_0\|_{{X}^{r}_q(I_1)}
+\|\nabla_{x,z} v\|_{\widetilde{L}_z^\infty(I_0;B^{r-\f12}_{\infty,q})}\\
&\qquad\qquad+\|\nabla_{x,z} v\|_{\widetilde{L}_z^\infty(I_0; C^{-\f \e2})}+\|\nabla_{x,z} v\|_{\widetilde{L}_z^2(I_1\setminus I_0; C^{-\f \e2})}\big).
\end{align*}
The same estimate for $\pa_z v$ can be deduced by using $\pa_zv=T_Av+w$ and Proposition \ref{prop:Sym-Besov}. Thus, we obtain
\begin{align*}
\|\na_{x,z}v\|_{{X}^{r}_q(I_0)}
\le&  K_\eta\big(\|f\|_{B^{r+1}_{\infty,q}}+\|F_{00}\|_{{Y}^{r}_q(I)}+\|G_0\|_{{X}^{r}_q(I_1)}
+\|\nabla_{x,z} v\|_{\widetilde{L}_z^\infty(I_0;B^{r-\f12}_{\infty,q})}\\
&\qquad+\|\nabla_{x,z} v\|_{\widetilde{L}_z^\infty(I_0; C^{-\f \e2})}+\|\nabla_{x,z} v\|_{\widetilde{L}_z^2(I_1\setminus I_0; C^{-\f \e2})}\big).
\end{align*}
This together with the interpolation inequality
\beno
\|\nabla_{x,z} v\|_{\widetilde{L}_z^\infty(I_0; B^{r-\f12}_{\infty,q}\cap C^{-\f \e2})}\le K_\eta\|\nabla_{x,z} v\|_{{L}_z^\infty(I_0; C^{-\delta})}
+\f1 {4K_\eta}\|\na_{x,z}v\|_{\widetilde{L}_z^\infty(I_0;B^{r}_{\infty,q})}
\eeno
implies the proposition.\ef


\subsection{Interior $W^{1,p}$ estimate}

We consider the elliptic equation
\ben\label{eq:elliptic-inter}
\Delta_{x,y}\phi=\textrm{div}_{x,y}g\quad \textrm{in}\quad\cS.
\een
Given a point $X_0=(x_0,y_0)\in \R^d\times \R$, let $B_{r}(X_0)$ be a ball of radius $r$ centered at $X_0$.
We have the following interior $W^{1,p}$ estimate.

\begin{proposition}\label{prop:holder-inter}
Suppose that $\phi\in H^1(\cS)$ is a solution of (\ref{eq:elliptic-inter}) with $g\in L^p(\cS)$ for $p\ge 2$.
Let $I_x=[-1+c_1{h_0}, \eta(x)-c_2{h_0}]$, where $c_1,c_2>0$ and $c_1+c_2<1$.
Then there exits $\delta_1>0$ depending only on $\|\eta\|_{C^{1+\varepsilon}}, h_0$ and $c_1, c_2$ so that $B_{\delta_1}(X_0)\subseteq \cS$ for any $X_0\in \{(x,y): x\in \R^d, y\in I_x\}$
and
\beno
\|\phi\|_{W^{1,p}(B_{{\delta_1}/2}(X_0))}\le K_\eta\big(\|g\|_{L^p(B_{\delta_1}(X_0))}+\|\na_{x,y}\phi\|_{L^2(B_{\delta_1}(X_0))}\big).
\eeno
\end{proposition}

\no{\bf Proof.}\,Given any point $X_0\in \R^d\times I_x$.
Let $\overline{\phi}$ be a solution of the elliptic equation in $ B_{\delta_1}(X_0)$:
\beno
\Delta_{x,y}\overline{\phi}=\textrm{div}_{x,y}{g}\quad \textrm{in  } B_{\delta_1}(X_0), \quad \overline{\phi}|_{\pa B_{\delta_1}(X_0)}=0.
\eeno
The classical $W^{1,p}$ elliptic estimate ensures that
\beno
\|\overline{\phi}\|_{W^{1,p}(B_{\delta_1}(X_0))}\le C(\delta_1)\|{g}\|_{L^p(B_{\delta_1}(X_0))}.
\eeno
Then by using the interior gradient estimate of the harmonic function $\phi-\overline{\phi}$, we deduce that
\begin{align*}
\|\phi\|_{W^{1,p}(B_{{\delta_1}/2}(X_0))}\le& \|\phi-\overline{\phi}\|_{W^{1,p}(B_{{\delta_1}/2}(X_0))}
+\|\overline{\phi}\|_{W^{1,p}(B_{{\delta_1}/2}(X_0))}\\
\le& C(\delta_1)\big(\|\na_{x,y}(\phi-\overline{\phi})\|_{L^2(B_{\delta_1}(X_0))}+\|\overline{\phi}\|_{W^{1,p}(B_{\delta_1}(X_0))}\big)\\
\le& C(\delta_1)\big(\|\na_{x,y}\phi\|_{L^2(B_{\delta_1}(X_0))}+\|g\|_{L^p(B_{\delta_1}(X_0))}\big).
\end{align*}
The proof is finished.\ef

\section{Dirichlet-Neumann operator}

In this section, we assume that $\eta\in C^{\f32+\e}(\R^d)$ for some $\e>0$ and satisfies (\ref{ass:height}).
We will use some notations from section 4.

\subsection{Definition and paralinearization}

We consider the boundary value problem
\ben\label{eq:elliptic-DN}
\left\{
\begin{array}{l}
\Delta_{x,y} \phi=0\qquad\textrm{ in}\quad \cS,\\
\phi|_{y=\eta(x)}=f,\quad \phi|_{y=-1}=0.
\end{array}
\right.
\een
where $\cS=\big\{(x,y): x\in \R^d,-1< y<\eta(x)\big\}$.
Given $f\in H^\f12(\R^d)$, the existence of the variation solution $\phi$ with $\na_{x,y}\phi\in L^2(\cS)$ can be deduced by using
Riesz theorem, see \cite{ABZ-IM} for example. Moreover, it holds that
\ben\label{eq:DN-H1}
\|\na_{x,y}\phi\|_{L^2(\cS)}\le C(\|\eta\|_{W^{1,\infty}},h_0)\|f\|_{H^\f12}.
\een
Let $\eta\in H^s(\R^d)$ for $s>\f d 2+\f3 2$. This together with  Proposition \ref{prop:elliptic-boun} yields that
for any $\sigma\in [-\f12,s-1]$,
\begin{align}\label{eq:DN-v-est}
\|\widetilde{\phi}\|_{H^{\sigma+\f32}(\overline{\cS})}\le P_{s,\eta}\|f\|_{H^{\sigma+1}}.
\end{align}

\begin{definition}
Given $\eta, f, \phi$ as above, the Dirichlet-Neumann(DN) operator $G(\eta)$ is defined by
\beno
G(\eta)f\eqdef\sqrt{1+|\nabla \eta|^2}\pa_n \phi\big|_{y=\eta}.
\eeno
\end{definition}

The DN operator is a positive self-adjoint operator. More precisely, for any $f,g\in H^\f12(\R^d)$, we have
\beno
&&\langle G(\eta)f,g\rangle=\langle f,G(\eta)g\rangle,\\
&&\langle G(\eta)f,f\rangle=\|\na_{x,y}\phi\|_{L^2(\cS)}\ge c\|f\|_{H^\f12}^2.
\eeno

In terms of $\widetilde{\phi}$, the Dirichlet-Neumann operator $G(\eta)$ can be written as
\beno
G(\eta)f=\Big(\f{1+|\nabla \rho_\delta|^2}{\pa_z \rho_\delta}\pa_z \widetilde{\phi}-\nabla \rho_\delta \cdot \nabla \widetilde{\phi}\Big)\Big|_{z=0}.
\eeno

Following \cite{ABZ-IM}, we first paralinearize $G(\eta)$. We set
\beno
\zeta_1=\f{1+|\nabla \rho_\delta|^2}{\pa_z \rho_\delta}\big|_{z=0},\quad \zeta_2=\na\rho_\delta\big|_{z=0}.
\eeno
It is easy to show that
\ben\label{eq:DN-zeta-Hs}
\|\zeta_1-1\|_{H^{s-1}}+\|\zeta_2\|_{H^{s-1}}\le K_\eta\|\eta\|_{H^s}.
\een

Using Bony's decomposition (\ref{Bony}), we decompose $G(\eta)$ as
\begin{align*}
G(\eta)f=&\pa_z\widetilde{\phi}+T_{\zeta_1-1}\pa_z\widetilde{\phi}+T_{\pa_z \widetilde{\phi}}(\zeta_1-1)+R(\zeta_1-1,\pa_z \widetilde{\phi})-T_{i\zeta_2\cdot \xi}\widetilde{\phi}\\
&-T_{\nabla \widetilde{\phi}}\cdot\zeta_2-R(\zeta_2,\nabla \widetilde{\phi})\big|_{z=0}.
\end{align*}
Replacing $\pa_z \widetilde{\phi}$ by $T_A\widetilde{\phi}$, we get
\ben\label{DN-para}
G(\eta)f=T_{\lambda}f+R(\eta)f,
\een
where $\lambda=\zeta_1A-i\zeta_2\cdot \xi\big|_{z=0}$ with
\ben\label{eq:DN-A}
A=\f{1}{2}({-i\beta\cdot\xi}+\sqrt{4\al|\xi|^2-(\beta\cdot\xi)^2}),
\een
and $R(\eta)$ is the remainder of DN operator given by
\begin{align}
R(\eta)f
=&\Big[\big(T_{\zeta_1}T_A-T_{\zeta_1 A}\big)\widetilde{\phi}-T_{\zeta_1}(\pa_z-T_A)\widetilde{\phi}\nonumber\\
 &+\big(S_2(\pa_z \widetilde{\phi})+T_{\pa_z \widetilde{\phi}}(\zeta_1-1)+R(\zeta_1-1,\pa_z \widetilde{\phi})-T_{\nabla \widetilde{\phi}}\cdot\zeta_2-R(\nabla \widetilde{\phi},\zeta_2)\big)\Big]\bigg|_{z=0}\nonumber\\
\triangleq & R_1(\eta)f+R_2(\eta)f+R_3(\eta)f.\label{eq:DN-R}
\end{align}

\subsection{Sobolev estimate of the remainder}

In this subsection, we present Sobolev estimates for the remainder in the case when the boundary is smooth,
which will used in the proof of the uniform estimates of the approximate solutions.

\begin{proposition}\label{prop:DN-remainder}
Assume that $\eta\in H^s(\R^d)$ for $s>\f d 2+2$. Then it holds that
\beno
\|R(\eta)f\|_{H^{s-1}}\le P_{s,\eta}\|f\|_{H^{s-1}}.
\eeno
\end{proposition}

\no{\bf Proof.}\,Thanks to the fact that for $s>\f d 2+2$,
\beno
M_1^0(\zeta_1)+\cM_1^1(A)\le P_{s,\eta},
\eeno
we deduce from Proposition \ref{prop:symbolic calculus} and (\ref{eq:DN-v-est}) that
\begin{align*}
\|R_1(\eta)f\|_{H^{s-1}}\le& P_{s,\eta}\|\na \widetilde{\phi}\|_{L^\infty_z(I;H^{s-2})}\le P_{s,\eta}\|f\|_{H^{s-1}}.
\end{align*}
By Lemma \ref{lem:remaider} and (\ref{eq:DN-zeta-Hs}), we have
\begin{align*}
\|R_3(\eta)f\|_{H^{s-1}}\le& C\|\na_{x,z}\widetilde{\phi}\|_{L^\infty_z(I;L^2)}+P_{s,\eta}\|\nabla_{x,z}\widetilde{\phi}\|_{{L}^\infty(\overline{\cS})}\\
\le& P_{s,\eta}\|\na_{x,z}\widetilde{\phi}\|_{L^\infty_z(I;H^{s-2})} \le P_{s,\eta}\|f\|_{H^{s-1}}.
\end{align*}

For $R_2(\eta)$, we get by Lemma \ref{lem:remaider} that
\begin{align}
\|R_2(\eta)f\|_{H^{s-1}}&\le P_{s,\eta}\|(\pa_z-T_A)\widetilde{\phi}\|_{L^\infty_z([-\f12,0];H^{s-1})}.\label{eq:DN-R2-Hs}
\end{align}
While by the proof of Proposition \ref{prop:elliptic-upboun}, we know that $w=\chi(z)(\pa_z-T_A)\widetilde{\phi}$ satisfies
\beno
\pa_zw-T_aw=F',\quad w(-1)=0.
\eeno
where $F'=\chi(z)\big(F_1+F_2+F_3\big)+\chi'(z)(\pa_z-T_A)\widetilde{\phi}$ and $\chi$ is a smooth function satisfying $\chi(-1)=0$ and $\chi(z)=1$ for $z\in [-\f12,0]$.
Then it follows from Proposition \ref{prop:parabolic estimate}, Lemma \ref{lem:DN-F-est} and (\ref{eq:DN-v-est}) that
\begin{align*}
\|(\pa_z-T_A)\widetilde{\phi}\|_{X^{s-1}([-\f12,0])}\le& \|w\|_{X^{s-1}(I)}\\
\le& P_{s,\eta}\Big(\sum_{i=1}^3\|F_i\|_{L^2_z(I;H^{s-\f32})}+\|\na_{x,z}\widetilde{\phi}\|_{X^{s-2}(I)}\Big)\\
\le& P_{s,\eta}\|f\|_{H^{s-1}},
\end{align*}
which along with (\ref{eq:DN-R2-Hs}) gives
\begin{align*}
\|R_2(\eta)f\|_{H^{s-1}}\le P_{s,\eta}\|f\|_{H^{s-1}}.
\end{align*}

Putting the estimates of $R_i(\eta)f$ together concludes  the proposition.\ef

\begin{lemma}\label{lem:DN-F-est}
Let $s>\f d 2+2$. It holds that for $i=1,2,3$,
\beno
\|F_i\|_{L^2_z(I;H^{s-\f32})}\le P_{s,\eta}\|\na_{x,z}\widetilde{\phi}\|_{X^{s-2}(I)}.
\eeno
\end{lemma}

\no{\bf Proof.}\,Recall that $F_1=\gamma\pa_z\widetilde{v}$. It follows from Lemma \ref{lem:product-full} that
\begin{align*}
\|F_1\|_{L^2_z(I;H^{s-\f32})}\le& C\|\gamma\|_{L^\infty}\|\pa_z\widetilde{\phi}\|_{L^2_z(I;H^{s-\f32})}
+C\|\pa_z\widetilde{\phi}\|_{L^\infty}\|\gamma\|_{L^2_z(I;H^{s-\f32})}\\
\le& C\|\gamma\|_{L^\infty_z(I;H^{s-2})}\|\pa_z\widetilde{\phi}\|_{L^2_z(I;H^{s-\f32})}
+C\|\pa_z\widetilde{\phi}\|_{L^\infty_z(I;H^{s-2})}\|\gamma\|_{L^2_z(I;H^{s-\f32})}\\
\le& P_{s,\eta}\|\na_{x,z}\widetilde{\phi}\|_{X^{s-2}(I)}.
\end{align*}

Recall that
\beno
F_2=(T_\al-\al)\Delta\widetilde{\phi}+(T_\beta-\beta)\cdot\nabla \pa_z \widetilde{\phi}.
\eeno
Then we get by Lemma \ref{lem:remaider} that
\begin{align*}
\|F_2\|_{L^2_z(I;H^{s-\f32})}\le& C\|\na_{x,z}\widetilde{\phi}\|_{L^\infty}\big(1+\|\al-1\|_{L^2_z(I;H^{s-\f12})}+\|\beta\|_{L^2_z(I;H^{s-\f12})}\big)\\
\le& P_{s,\eta}\|\na_{x,z}\widetilde{\phi}\|_{X^{s-2}(I)}.
\end{align*}

Recall that
\beno
F_3=(T_aT_A-T_\al\Delta)\widetilde{\phi}-(T_a+T_A+T_\beta\cdot\nabla )\pa_z\widetilde{\phi}-T_{\pa_z A}\widetilde{\phi}.
\eeno
For $s>\f d 2+2$, we have
\beno
\cM_1^0(a)+\cM_1^1(A)\le P_{s,\eta}.
\eeno
Then we deduce from Lemma \ref{lem:symbol} and Proposition \ref{prop:symbolic calculus} that
\begin{align*}
\|F_3\|_{L^2_z(I;H^{s-\f32})}\le P_{s,\eta}\|\na_{x,z}\widetilde{\phi}\|_{X^{s-2}(I)}.
\end{align*}
The proof is finished.\ef

\subsection{Tame estimate of the remainder}

In this subsection, we present tame estimates for the remainder in the case when the boundary  has more one half derivative than $f$.
The result will be used in the proof of the break-down criterion.

\begin{proposition}\label{prop:DN-Hs}
Assume that $\eta\in H^{s+\f12}(\R^d)$ for $s>\f d 2+1$.
Let $I_1$ be as in Proposition \ref{prop:elliptic-upboun} and $\widetilde{S}=\R^d\times I_1$. It holds that
\beno
&&\|R(\eta)f\|_{{H^{s-\f12}}}
\leq K_\eta\big(\|f\|_{H^s}+\|\nabla_{x,z}\widetilde{\phi}\|_{{L}^\infty(\widetilde{\cS})}\|\eta\|_{H^{s+\f12}}\big),\\
&&\|R(\eta)f\|_{{H^{s-1}}}
\leq  K_\eta\big(\|f\|_{H^{s-\f12}}+\|\nabla_{x,z}\widetilde{\phi}\|_{L^\infty_z(I_1;C^0)}\|\eta\|_{H^{s+\f12}}\big).
\eeno
\end{proposition}

\no{\bf Proof.}\,Note that $A\in \Gamma^1_{\f12+\e}(I\times \R^d)$ and $\zeta_1\in \Gamma^0_{\f12+\e}(\R^d)$ with the bound
\ben\label{eq:DN-Azeta-est}
{\cM}^1_{\f12+\e}(A)+{M}^0_{\f12+\e}(\zeta_1)\le K_\eta.
\een
Then we get by Proposition \ref{prop:symbolic calculus}, (\ref{eq:DN-zeta-Hs}), Proposition \ref{prop:elliptic-upboun} and (\ref{eq:DN-H1}) that
\begin{align*}
\|R_1(\eta)f\|_{H^{s-\f12}}\le& K_\eta\|\na \widetilde{\phi}\|_{L^\infty_z([z_0,0];H^{s-1})}\\
\le& K_\eta\big(\|\na_{x,z}\widetilde{\phi}\|_{L^2(\overline{\cS})}+\|f\|_{H^s}+\|\nabla_{x,z} \widetilde{\phi}\|_{L^{\infty}_z(I_1;C^{0})}\|\eta\|_{H^{s+\f12}}\big)\\
\le& K_\eta\big(\|f\|_{H^s}+\|\nabla_{x,z}\widetilde{\phi}\|_{L^{\infty}_z(I_1;C^{0})}\|\eta\|_{H^{s+\f12}}\big).
\end{align*}
Similarly, we deduce from Lemma \ref{lem:remaider} that
\begin{align*}
\|R_3(\eta)f\|_{H^{s-\f12}}\le& K_\eta\|\nabla_{x,z}\widetilde{\phi}\|_{{L}^\infty(\widetilde{\cS})}\|\eta\|_{H^{s+\f12}}+\|\pa_z\widetilde{\phi}\|_{L^\infty_z([z_0,0];H^{s-1})}\\
\le& K_\eta\big(\|f\|_{H^s}+\|\nabla_{x,z}\widetilde{\phi}\|_{{L}^\infty(\widetilde{\cS})}\|\eta\|_{H^{s+\f12}}\big).
\end{align*}

While by the proof of (\ref{eq:elliptic-Hs-2}), we see that
\begin{align*}
\|(\pa_z-T_A)\widetilde{\phi}\|_{L^\infty_z([z_0,0];H^{s-\f12})}
\le& K_\eta\big(\|\nabla_{x,z}\widetilde{\phi}\|_{X^{s-1}([z_0,0])}+\|\nabla_{x,z}\widetilde{\phi}\|_{{L}^\infty(\widetilde{\cS})}\|\eta\|_{H^{s+\f12}}\big)\\
\le&  K_\eta\big(\|f\|_{H^s}+\|\nabla_{x,z}\widetilde{\phi}\|_{{L}^\infty(\widetilde{\cS})}\|\eta\|_{H^{s+\f12}}\big).
\end{align*}
Then we get by Lemma \ref{lem:remaider} that
\begin{align*}
\|R_2(\eta)f\|_{H^{s-\f12}}&\le K_\eta\|(\pa_z-T_A)\widetilde{\phi}\|_{L^\infty_z([z_0,0];H^{s-\f12})}\\
&\le K_\eta\big(\|f\|_{H^s}+\|\nabla_{x,z}\widetilde{\phi}\|_{{L}^\infty(\widetilde{\cS})}\|\eta\|_{H^{s+\f12}}\big).
\end{align*}
This completes the proof of the first inequality.
The second inequality can be proved similarly.\ef

\subsection{H\"{o}lder estimate of the remainder}

\begin{proposition}\label{prop:DN-Holder}
Let $I_0=[-\f12,0]$ and $I_1=[-\f34,0]$.
It holds that for any $\delta>0$,
\beno
\|R(\eta)f\|_{{C^{\f12}}}\leq K_\eta\big(\|f\|_{C^1}+\|\na_{x,z}\widetilde{\phi}\|_{\widetilde{L}^\infty_z(I_0;C^{-\delta})}+
\|\na_{x,z}\widetilde{\phi}\|_{\widetilde{L}^2_z(I_1\backslash I_0;C^{-\f \e 2})}\big).
\eeno
\end{proposition}
\no{\bf Proof.}\,Thanks to (\ref{eq:DN-Azeta-est}),
we get by Proposition \ref{prop:Sym-Besov} that
\begin{align*}
\|R_1(\eta)f\|_{C^\f12}\le& K_\eta\|\na\widetilde{\phi}\|_{L^\infty_z(I_0; C^0)}.
\end{align*}
We can see from the proof of Proposition \ref{prop:elliptic Holder est} that
\beno
\|(\pa_z-T_A)\widetilde{\phi}\|_{\widetilde{L}^\infty_z([-\f14,0];C^\f12)}\le K_\eta(\|\nabla_{x,z}\widetilde{\phi}\|_{\widetilde{L}_z^\infty(I_0; C^0)}+
\|\na_{x,z}\widetilde{\phi}\|_{\widetilde{L}^2_z(I_1\backslash I_0;C^{-\f \e 2})}),
\eeno
which together with Lemma \ref{lem:remaider} gives
\begin{align*}
\|R_2(\eta)f\|_{C^\f12}&\le K_\eta\|(\pa_z-T_A)\widetilde{\phi}\|_{L^\infty_z([-\f14,0];C^\f12)}
\le K_\eta(\|\nabla_{x,z}\widetilde{\phi}\|_{\widetilde{L}_z^\infty(I_0; C^0)}+
\|\na_{x,z}\widetilde{\phi}\|_{\widetilde{L}^2_z(I_1\backslash I_0;C^{-\f \e 2})}).
\end{align*}
By Lemma \ref{lem:remaider} again, we get
 \begin{align*}
\|R_3(\eta)f\|_{C^\f12}
&\le K_\eta\|{\na_{x,z} \widetilde{\phi}}\|_{L^\infty_z(I; C^0)}.
\end{align*}
This together with Proposition \ref{prop:elliptic Holder est} shows that
\begin{align*}
\|R(\eta)f\|_{C^\f12}&\le K_\eta(\|\nabla_{x,z}\widetilde{\phi}\|_{\widetilde{L}_z^\infty(I_0; C^0)}+
\|\na_{x,z}\widetilde{\phi}\|_{\widetilde{L}^2_z(I_1\backslash I_0;C^{-\f \e 2})})\\
&\le K_\eta\big(\|f\|_{B^1_{\infty,\infty}}+\|\na_{x,z}\widetilde{\phi}\|_{\widetilde{L}^\infty_z(I_0;C^{-\delta})}+
\|\na_{x,z}\widetilde{\phi}\|_{\widetilde{L}^2_z(I_1\backslash I_0;C^{-\f \e 2})}\big).
\end{align*}
This finishes the proof.\ef

\section{New formulation and paralinearization}

\subsection{New formulation}

First of all, we derive the evolution equations for the free surface and the trace of the velocity on the boundary.
We denote
\beno
&&V\triangleq \big(v^1,\cdots, v^d\big)\big|_{y=\eta},\quad V_b\triangleq \big(v^1,\cdots, v^d\big)\big|_{y=-1},\quad B\triangleq v^{d+1}\big|_{y=\eta},\\
&&a\triangleq-\pa_y P|_{y=\eta},\quad \zeta\triangleq\na\eta.
\eeno
Using the fact that for any function $f=f(t,x,y)$,
\begin{align}
(\pa_t+V\cdot\nabla)(f|_{y=\eta})=\big(\pa_tf+v\cdot\nabla_{x,y}f\big)\big|_{y=\eta},\nonumber
\end{align}
we deduce from (\ref{eq:euler}) that
\ben
&&(\pa_t+V\cdot\nabla)B=a-1,\label{eq:B}\\
&&(\pa_t+V\cdot\nabla)V+a\zeta=0,\label{eq:V}\\
&&(\pa_t+V_b\cdot\nabla)V_b+\na P|_{y=-1}=0.\label{eq:Vb}
\een

Let $\om$ be the vorticity of the fluid, which is defined by
\beno
\om=\big(\om_{i,j}\big)_{1\le i,j\le d+1},\quad \om_{i,j}=\pa_{x_i}v^j-\pa_{x_j}v^i.
\eeno
The motion of the fluid is determined by the vorticity equation
\ben\label{eq:vorticity}
\om_t+v\cdot\na_{x,y}\om=\om\otimes\na_{x,y} v\quad \textrm{in} \quad \Om_t.
\een
Here $(\om\otimes\na_{x,y} v)_{i,j}=\om_{k,i}\pa_j v^k+\om_{k,j}\pa_i v^k$.

The velocity $v$ can be recovered from the vorticity $\om$ by solving the elliptic equation
\beq\label{eq:v-vorticity}
\left\{
\begin{array}{ll}
\Delta_{x, y} v=-\na_{x,y}\times\om\quad \textrm{in} \quad \Omega_t,\\
v|_{y=\eta}=(V, B),\quad v|_{y=-1}=(V_b,0).
\end{array}\right.
\eeq
where $\om=\big(\om_{i,j}\big)_{1\le i,j\le d+1}$ with $\om_{i,j}=\pa_{x_i}v^j-\pa_{x_j}v^i$.
The pressure $P$ of the fluid is determined by solving the following elliptic equation:
\begin{equation}\label{eq:pressure}
\left\{
\begin{aligned}
&-\Delta_{x,y} P=\pa_i{v}^j\pa_j{v}^i\quad \textrm{in} \quad \Om_t,\\
&P|_{y=\eta}=0,\quad \pa_yP|_{y=-1}=-1.
\end{aligned}
\right.
\end{equation}

We decompose the velocity $v$ into the irrotational part $v_{ir}$ and the rotational part $v_{\om}$, i.e.,
\beq\label{eq:v-ir}
\left\{
\begin{array}{ll}
\Delta_{x, y} v_{ir}=0\quad \textrm{in} \quad \Omega_t,\\
v_{ir}|_{y=\eta}=(V, B),\quad v_{ir}|_{y=-1}=0,
\end{array}\right.
\eeq
and
\beq\label{eq:v-vor}
\left\{
\begin{array}{ll}
\Delta_{x, y} v_{\om}=-\na_{x,y}\times\om\quad \textrm{in} \quad \Om_t,\\
v_{\om}|_{y=\eta}=0,\quad v_\om|_{y=-1}=(V_b, 0).
\end{array}\right.
\eeq

It follows from (\ref{eq:euler-b2}) that
\ben\label{eq:zeta}
(\pa_t+V\cdot\nabla)\zeta=\na B-\sum_{j}\na V_j\zeta^j.
\een
Then a direct calculation yields
\begin{align*}
&\pa_{x_i}B-\pa_{x_i}V_j\pa_{x_j}\eta\\
&=\pa_{x_i}v^{d+1}+\pa_{x_i}\eta\pa_yv^{d+1}-\pa_{x_j}\eta\big(\pa_{x_i}v^j+\pa_{x_i}\eta\pa_yv^j\big)\big|_{y=\eta}\\
&=\big(\pa_y v^i-\pa_{x_j}v^i\cdot \pa_{x_j}\eta\big)+\pa_{x_i}\eta\big(\pa_y v^{d+1}-\pa_{x_j}\eta\pa_{x_j}v^{d+1}\big)\\
&\quad+\big(\om_{i,d+1}-\pa_{x_j}\eta\om_{ij}+\pa_{x_i}\eta\pa_{x_j}\eta\om_{j,d+1}\big)\big|_{y=\eta}\\
&=\big(\pa_y v_{ir}^i-\pa_{x_j}v_{ir}^i\cdot \pa_{x_j}\eta\big)+\pa_{x_i}\eta\big(\pa_y v_{ir}^{d+1}-\pa_{x_j}\eta\pa_{x_j}v_{ir}^{d+1}\big)\\
&\quad+\big(\pa_y v_{\om}^i-\pa_{x_j}v_{\om}^i\cdot \pa_{x_j}\eta\big)+\pa_{x_i}\eta\big(\pa_y v_{\om}^{d+1}-\pa_{x_j}\eta\pa_{x_j}v_{\om}^{d+1}\big)\\
&\quad+\big(\om_{i,d+1}-\pa_{x_j}\eta\om_{ij}+\pa_{x_i}\eta\pa_{x_j}\eta\om_{j,d+1}\big)\big|_{y=\eta}\\
&=G(\eta)V_i+\pa_{x_i}\eta G(\eta)B+R_\om^i,
\end{align*}
where
\begin{align*}
R_\om^i\triangleq&\big(\pa_y v_{\om}^i-\pa_{x_j}v_{\om}^i\cdot \pa_{x_j}\eta\big)+\pa_{x_i}\eta\big(\pa_y v_{\om}^{d+1}-\pa_{x_j}\eta\pa_{x_j}v_{\om}^{d+1}\big)\\
&+\big(\om_{i,d+1}-\pa_{x_j}\eta\om_{ij}+\pa_{x_i}\eta\pa_{x_j}\eta\om_{j,d+1}\big)\big|_{y=\eta}.
\end{align*}
Thus, $\zeta=\na \eta$ satisfies
\ben\label{eq:zeta-new}
(\pa_t+V\cdot\nabla)\zeta=G(\eta)V+\zeta G(\eta)B+R_\om.
\een
The term $R_\om$ induced by the vorticity will lead the system to lose one half derivative.

\subsection{Paralinearization}

We paralinearize the system (\ref{eq:B}), (\ref{eq:V}) and (\ref{eq:zeta-new}). For this end, we introduce so called good unknown $U=V+T_\zeta B$.
Applying Bony's decomposition and (\ref{DN-para}) to (\ref{eq:B}), (\ref{eq:V}) and (\ref{eq:zeta-new}), we obtain
\ben\label{eq:ZK-para}
\left\{
\begin{array}{l}
(\pa_t+T_V\cdot\na)V+T_a\zeta+T_\zeta(\pa_t+T_V\cdot\na)B=h_1,\\
(\pa_t+T_V\cdot\na)\zeta=T_\lambda U+h_2+R_\om,
\end{array}\right.
\een
where
\beno
&& h_1\triangleq(T_V-V)\cdot\nabla V-R(a,\zeta)+T_\zeta(T_V-V)\cdot\nabla B,\\
&&h_2\triangleq(T_V-V)\cdot\na \zeta+[T_\zeta, T_\lambda]B+(\zeta-T_\zeta)T_\lambda B+R(\eta)V+\zeta R(\eta)B.
\eeno
Let $D_t\triangleq \pa_t+T_V\cdot\na$. In terms of good unknown, the first equation of (\ref{eq:ZK-para}) can be rewritten as
\ben\label{eq:U}
D_t U+T_a\zeta=h_1+[D_t,T_\zeta]B.
\een
Taking $D_t$ on both sides of (\ref{eq:U}), we get
\beno
D_t^2 U+ T_a D_t\zeta=D_th_1+[T_a, D_t]\zeta+D_t[D_t, T_\zeta]B,
\eeno
which along with the second equation of (\ref{eq:ZK-para}) gives
\ben\label{eq:ZK-para-2}
D_t^2 U+ T_{a\lambda}U=f+f_\om,
\een
where $(f,f_\om)$ is given by
\beno
f\triangleq D_th_1+(T_{a\lambda}-T_\lambda T_a)U+[T_a,D_t]\zeta+D_t[D_t, T_\zeta]B-T_ah_2,\quad f_\om\triangleq-T_aR_\om.
\eeno

Similarly, we have
\ben\label{eq:ZK-para-3}
D_t^2\zeta+ T_{a\lambda}\zeta=g+g_\om,
\een
where $(g,g_\om)$ is given by
\beno
g\triangleq D_th_2+[D_t, T_\lambda]U+T_\lambda(h_1+D_t[D_t, T_\zeta]B)+(T_{a\lambda}-T_\lambda T_a)\zeta,\quad g_\om\triangleq D_tR_\om.
\eeno

\section{Estimate of the pressure}

The pressure $P$ satisfies
\ben\label{eq:pressure equ}
\left\{
\begin{aligned}
&-\Delta_{x,y} P=\pa_i\big({v}^j\pa_j{v}^i\big)=\pa_i v^j\pa_j v^i\quad \textrm{in} \quad \Om_t,\\
&P|_{y=\eta}=0,\quad \pa_yP|_{y=-1}=-1.
\end{aligned}
\right.
\een
Here $v$ is the velocity. In this section, we denote $P_1\triangleq P+y$.

\subsection{$H^2$ estimate of the pressure}

\begin{lemma}\label{lem:pressure-L2}
Let $P$ be a solution of (\ref{eq:pressure equ}). It holds that
\beno
\|\na_{x,y}P_1\|_{L^2(\Om_t)}\le K_\eta\big(\|\na_{x,y}v\|_{L^\infty(\Om_t)}\|v\|_{L^2(\Om_t)}+\|\eta\|_{H^\f12}\big).
\eeno
\end{lemma}

\no{\bf Proof.}\,Let $\overline{\eta}(t,x,y)=\chi(y)e^{(y-\eta(t,x))|D|}\eta$ for $y\le \eta(t,x)$,
where $\chi$ is a smooth function satisfying $\chi(0)=1$ and $\chi(y)=0$ for $y\le -1+h_0/2$.
Then $P_2\triangleq P+y-\overline{\eta}$ satisfies
\beno
\left\{
\begin{aligned}
&-\Delta_{x,y}P_2=\Delta_{x,y}\overline{\eta}+\pa_i\big({v}^j\pa_j{v}^i\big)\quad \textrm{in} \quad \Om_t,\\
&P_2|_{y=\eta}=0,\quad \pa_yP_2|_{y=-1}=0.
\end{aligned}
\right.
\eeno
Thanks to $v^{d+1}|_{y=-1}=0$, we get by integration by parts that
\begin{align*}
\int_{\Om_t}|\na_{x,y}P_2|^2dxdy&=\int_{\Om_t}\big(\Delta_{x,y}\overline{\eta}+\pa_i\pa_j({v}^j{v}^i)\big)P_2dxdy\\
&=-\int_{\Om_t}\big(\pa_j({v}^j{v}^h)+\pa_y({v}^h{v}^{d+1})\big)\cdot\na P_2dxdy\\
&\quad-\int_{\Om_t}\pa_y(v^{d+1}v^{d+1})\pa_yP_1dxdy-\int_{\Om_t}\na_{x,y}\overline{\eta}\cdot\na_{x,y}P_2dxdy\\
&\le C\big(\|\na_{x,y}v\|_{L^\infty(\Om_t)}\|v\|_{L^2(\Om_t)}+\|\na_{x,y}\overline{\eta}\|_{L^2(\Om_t)}\big)
\|\na_{x,y}P_2\|_{L^2(\Om_t)},
\end{align*}
from which, we deduce
\beno
\|\na_{x,y}P_2\|_{L^2(\Om_t)}\le C\big(\|\na_{x,y}v\|_{L^\infty(\Om_t)}\|v\|_{L^2(\Om_t)}+\|\na_{x,y}\overline{\eta}\|_{L^2(\Om_t)}\big).
\eeno
Then the lemma follows by uisng the fact that
\beno
\|\na_{x,y}\overline{\eta}\|_{L^2(\Om_t)}\le C(\|\eta\|_{W^{1,\infty}},h_0)\|\eta\|_{H^\f12}.
\eeno
This completes the proof.
\ef

\medskip

Next, we give the higher estimates of pressure:
\begin{lemma}\label{lem:pressure-H2}
Let $P$ be a solution of (\ref{eq:pressure equ}). It holds that
\beno
\|\na_{x,y}P_1\|_{H^1(\Om_t)}\le K_\eta\big(\|\na_{x,y}v\|_{L^\infty(\Om_t)}\|v\|_{H^1(\Om_t)}+\|\eta\|_{H^\f32}\big).
\eeno
\end{lemma}

\no{\bf Proof.}\,Let $P_2$ be as in the proof of Lemma \ref{lem:pressure-L2}. Then $\widetilde{P}_2$ satisfies
\ben\label{eq:P1-flat}
\left\{
\begin{array}{l}
\pa^2_z \widetilde{P}_2+\al \Delta \widetilde{P}_2+\beta\cdot\nabla \pa_z \widetilde{P}_2-\gamma \pa_z \widetilde{P}_2= F_{0},\\
\widetilde{P}_2|_{z=0}=0,\quad \pa_z\widetilde{P}_2|_{z=-1}=0,
\end{array}\right.
\een
where $F_{0}= \alpha\widetilde{\pa_iv^j}\widetilde{\pa_jv^i}+\al \widetilde{\Delta_{x,y}\overline{\eta}}$. It follows from (\ref{eq:P1-flat}) that
\begin{align*}
&\int_{\overline{\cS}}\pa_z^2\widetilde{P}_2\Delta \widetilde{P}_2dxdz+\int_{\overline{\cS}}\al|\Delta \widetilde{P}_2|^2dxdz
+\int_{\overline{\cS}}\beta\cdot\nabla \pa_z \widetilde{P}_2\Delta \widetilde{P}_2dxdz\\
&=\int_{\overline{\cS}}\gamma\pa_z \widetilde{P}_2\Delta \widetilde{P}_2dxdz+\int_{\overline{\cS}}F_0\Delta \widetilde{P}_2dxdz.
\end{align*}
By integration by parts, we get
\begin{align*}
\int_{\overline{\cS}}\pa_z^2\widetilde{P}_2\Delta \widetilde{P}_2dxdz=&-\int_{\overline{\cS}}\pa_z^2\na\widetilde{P}_2\cdot\na \widetilde{P}_2dxdz
=\int_{\overline{\cS}}|\na\pa_z\widetilde{P}_2|^2dxdz.
\end{align*}
Thanks to the definition of $\al,\beta$, it is easy to see that there exists $c>0$ depending on $\|\eta\|_{W^{1,\infty}}, h_0$ so that
\begin{align*}
&\int_{\overline{\cS}}\pa_z^2\widetilde{P}_2\Delta \widetilde{P}_2dxdz+\int_{\overline{\cS}}\al|\Delta \widetilde{P}_2|^2dxdz
+\int_{\overline{\cS}}\beta\cdot\nabla \pa_z \widetilde{P}_2\Delta \widetilde{P}_2dxdz\\
&\ge c\int_{\overline{\cS}}\big(|\Delta \widetilde{P}_2|^2+|\na\pa_z \widetilde{P}_2|^2\big)dxdz.
\end{align*}
Hence, we obtain
\begin{align*}
\int_{\overline{\cS}}\big(|\Delta \widetilde{P}_1|^2+|\na\pa_z \widetilde{P}_1|^2\big)dxdz\le C\big(\|\gamma\pa_z \widetilde{P}_1\|_{L^2(\overline{\cS})}+\|F_0\|_{L^2(\overline{\cS})}\big)\|\Delta \widetilde{P}_1\|_{L^2(\overline{\cS})}.
\end{align*}
It follows from Lemma \ref{lem:coeff} that
\begin{align*}
&\|\gamma\pa_z \widetilde{P}_2\|_{L^2(\overline{\cS})}\le C\|\gamma\|_{L^2_z(I:L^\infty)}\|\pa_z \widetilde{P}_2\|_{L^\infty_z(I:L^2)}
\le K_\eta\|\pa_z\widetilde{P}_2\|_{L^2(\overline{\cS})}^\f12\|\pa_z\na \widetilde{P}_1\|_{L^2(\overline{\cS})}^\f12.
\end{align*}
This shows that for any $\epsilon>0$,
\begin{align}
\|\Delta \widetilde{P}_2\|_{L^2(\overline{\cS})}+\|\na\pa_z\widetilde{P}_2\|_{L^2(\overline{\cS})}
\le& \|F_0\|_{L^2(\overline{\cS})}+K_\eta\|\pa_z\widetilde{P}_1\|_{L^2(\overline{\cS})}+\epsilon\|\pa_z\na \widetilde{P}_2\|_{L^2(\overline{\cS})}.\label{eq:P1-est1}
\end{align}
Using the equation (\ref{eq:P1-flat}), we infer that
\begin{align*}
\|\pa_z^2\widetilde{P}_2\|_{L^2(\overline{\cS})}\le K_\eta\big(&\|\Delta \widetilde{P}_2\|_{L^2(\overline{\cS})}+\|\na\pa_z\widetilde{P}_2\|_{L^2(\overline{\cS})}\\
&+\|\pa_z\widetilde{P}_2\|_{L^2(\overline{\cS})}^\f12\|\pa_z^2\widetilde{P}_2\|_{L^2(\overline{\cS})}^\f12
+\|F_0\|_{L^2(\overline{\cS})}\big),
\end{align*}
which along with (\ref{eq:P1-est1}) gives by taking $\epsilon$ small that
\begin{align*}
\|\na_{x,z}^2\widetilde{P}_2\|_{L^2(\Om_t)}\le& K_\eta\big(\|F_0\|_{L^2(\overline{\cS})}+\|\pa_z\widetilde{P}_1\|_{L^2(\overline{\cS})}\big)\\
\le&  K_\eta\big(\|\eta\|_{H^\f32}+\|\na_{x,y}v\|_{L^\infty(\Om_t)}\|\na_{x,y}v\|_{L^2(\Om_t)}+\|\pa_z\widetilde{P}_1\|_{L^2(\overline{\cS})}\big).
\end{align*}
Then the lemma follows by using Lemma \ref{lem:pressure-L2}.\ef\medskip

Next let us turn to $H^1$ estimate of $(\pa_t+v\cdot\na_{x,y})P\triangleq \cD_t P$. Using the equation (\ref{eq:euler}), a direct calculation gives
\beq\label{eq:pressure-Dt}
\left\{
\begin{array}{ll}
&\Delta_{x, y} \cD_t P=\pa_k P\Delta_{x, y} v^k+G,\\
&\cD_t P|_{y=\eta}=0,\quad \pa_y\cD_t P|_{y=-1}=\na\cdot v^{h}+\pa_y v^h\cdot\na P,
\end{array}\right.
\eeq
where
\begin{align*}
G=&4\delta^{ij}\pa_iv^k\pa_{j}\pa_kP-2(\pa_iv^j)(\pa_jv^k)\pa_kv^i\\
=&4\pa_k(\pa_iv^k\pa_iP)-2(\pa_iv^j)(\pa_jv^k)\pa_kv^i.
\end{align*}
Recall $\om_{i,k}=\pa_iv^k-\pa_kv^i$. We have
\beno
\pa_k P\cdot \Delta_{x, y} v^k=\pa_i(\pa_k P\om_{i,k})-\pa_i\pa_kP\om_{i,k}=\pa_i(\pa_k P\om_{i,k}).
\eeno
\begin{lemma}\label{lem:pressure-Dt-L2}
Let $\cD_t P$ be a solution of (\ref{eq:pressure-Dt}). It holds that
\beno
\|\na_{x,y}\cD_t P\|_{L^2(\Om_t)}\le K_\eta\big(1+\|\na_{x,y}v\|_{L^\infty(\Om_t)}^2\big)\big(\|v\|_{H^1(\Om_t)}+\|\eta\|_{H^\f32}\big).
\eeno
\end{lemma}

\no{\bf Proof.}\,  We get by integration by parts that
\begin{align*}
&\int_{\Om_t}|\na_{x,y}\cD_t P|^2dxdy+\int_{\R^d}\pa_y\cD_t P \cD_t P|_{y=-1}dx\\
&=-\int_{\Om_t}\big(4(\pa_iv^k\pa_iP)\pa_k\cD_t P+\pa_k P\om_{i,k}\pa_i\cD_t P+2(\pa_iv^j)(\pa_jv^k)\pa_kv^i\cD_t P\big)dxdy\\
&\quad+\int_{\R^d}(\pa_k P\om_{d+1,k})\cD_t P+(\pa_iv^{d+1}\pa_iP)\cD_t Pdx\big|_{y=-1}\\
&\le C\big(\|\na_{x,y}v\|_{L^2(\Om_t)}+\|\na_{x,y}v\|_{L^\infty(\Om_t)}\|\na_{x,y}P_1\|_{L^2(\Om_t)}\big)\|\na_{x,y}\cD_t P\|_{L^2(\Om_t)}\\
&\quad+C\|\na_{x,y}v\|_{L^\infty(\Om_t)}^2\|\na_{x,y}v\|_{L^2(\Om_t)}\|\cD_t P\|_{L^2(\Om_t)}\\
&\quad+C\|\na_{x,y}v\|_{L^\infty(\Om_t)}\|\na_{x,y}P_1(\cdot,-1)\|_{L^2(\R^d)}\|\cD_t P(\cdot,-1)\|_{L^2}+\int_{\R^d}\na\cdot v^{h}\cD_t P|_{y=-1}dx.
\end{align*}
Here we used $\dv v=0$. By the boundary condition of (\ref{eq:pressure-Dt}), we have
\begin{align*}
\int_{\R^d}\pa_y\cD_t P \cD_t P-\na\cdot v^{h}\cD_t P|_{y=-1}dx=&\int_{\R^d}(\pa_y v^h\cdot\na P)\cD_t P|_{y=-1}dx\\
\le& \|\na_{x,y}v\|_{L^\infty(\Om_t)}\|\na_{x,y}P_1(\cdot,-1)\|_{L^2}\|\cD_t P(\cdot,-1)\|_{L^2}\\
\le& K_\eta\|\na_{x,y}v\|_{L^\infty(\Om_t)}\|\na_{x,y}P_1\|_{H^1(\Om_t)}\|\cD_t P\|_{H^1(\Om_t)}
\end{align*}
Thanks to $\cD_t P|_{y=\eta}=0$ so that
$$
\|\cD_t P\|_{L^2(\Om_t)}\le K_\eta\|\cD_t P\|_{H^1(\Om_t)}.
$$
we deduce that
\begin{align*}
\|\na_{x,y}\cD_t P\|_{L^2(\Om_t)}\le K_\eta\big(1+\|\na_{x,y}v\|_{L^\infty(\Om_t)}^2\big)\big(\|\na_{x,y}P_1\|_{H^1(\Om_t)}+\|\na_{x,y}v\|_{L^2(\Om_t)}\big),
\end{align*}
which together with Lemma \ref{lem:pressure-H2} concludes the lemma.\ef

\subsection{H\"{o}lder estimate of the pressure}

In the sequel, we denote
\ben\label{def:A(t)}
A(t)\triangleq 1+\|v(t)\|_{W^{1,\infty}(\Om_t)}+\|v(t)\|_{H^1(\Om_t)}+\|\eta(t)\|_{H^\f32}.
\een

First, we give the H\"older estimate of the pressure.
\begin{lemma}\label{lem:pressure-holder}
Let $I_0=[-\f34,0]$ and $P$ be a solution of (\ref{eq:pressure equ}). Then it holds that
\beno
\|\nabla_{x, z} \widetilde{P}\|_{\widetilde{L}^\infty_z(I_0;C^\f12)}\leq K_\eta A(t)^2.
\eeno
\end{lemma}

\no{\bf Proof.}\,Apply Proposition \ref{prop:elliptic Holder est} with $F_0=\alpha\widetilde{\pa_iv^j}\widetilde{\pa_jv^i}$ and $G_0=0$ to obtain
\begin{align*}
\|\nabla_{x, z} \widetilde{P}\|_{\widetilde{L}^\infty_z(I_0;C^\f12)}\le& K_\eta\big(\|F_0\|_{\widetilde{L}^2_z(I;C^0)}+\|\nabla_{x, z} \widetilde{P}\|_{\widetilde{L}^\infty_z(I_0;C^{-\delta})}
+\|\nabla_{x, z} \widetilde{P}\|_{\widetilde{L}^2_z(I_1\setminus I_0;C^{-\f{\varepsilon}2)}}\big),
\end{align*}
where $I_1=[-\f78, 0]$ and $\delta>0$ is taken so that $-\delta+\f d 2\le -\f12$.
It is obvious that
\beno
\|F_0\|_{L^2_z(I;C^0)}\le K_\eta\|\na_{x,y}v\|_{L^\infty(\Om_t)}^2.
\eeno
It follows from Lemma \ref{lem:elliptic-H1} and Lemma \ref{lem:pressure-L2} that
\begin{align*}
\|\nabla_{x,z}\widetilde{P}\|_{\widetilde{L}^\infty_z(I_0;C^{-\delta})}
\le& \|\nabla_{x,z}\widetilde{P_1}\|_{\widetilde{L}^\infty_z(I_0;C^{-\delta})}+\|\nabla_{x,z}\rho_\delta \|_{\widetilde{L}^\infty_z(I_0;C^{-\delta})}\\
\le& C\|\nabla_{x,z}\widetilde{P_1}\|_{L^\infty_z(I_0;H^{-\f12})}+C\|\eta\|_{L^\infty}\\
\le& K_\eta\big(\|\na_{x,y}P_1\|_{L^2(\Om_t)}+\|F_0\|_{L^2_z(I;L^2)}\big)+C\|\eta\|_{L^\infty}\\
\le& K_\eta\big(\|\na_{x,y}v\|_{L^\infty(\Om_t)}\|v\|_{H^1(\Om_t)}+\|\eta\|_{H^\f12}+1\big).
\end{align*}
Take $c_1, c_2\in (0,1)$ depending on $K_\eta, a$ such that $\rho_\delta(x, I_1\setminus I_0)\in \cS_1=\big\{(x,y): y\in [-1+c_1h_0, \eta(x)-c_2h_0]\big\}$.
Let $\delta_1$ be as in Proposition \ref{prop:holder-inter}. Then for $p$ big enough depending on $\varepsilon$ and $d$,  we have
\begin{align*}
\|\nabla_{x, z} \widetilde{P}\|_{\widetilde{L}^2_z(I_1\setminus I_0;C^{-\f{\varepsilon}2)}}\le&
C(\delta_1)\|\nabla_{x, z} \widetilde{P}_1\|_{{L}^2_z(I;L^2)}+C\|\eta\|_{W^{1,\infty}}\\
&\quad+\sup_{X_0\in \cS_1}\|\nabla_{x,y}P\|_{L^{p}(B_{\delta_1}(X_0))},
\end{align*}
which together with  Proposition \ref{prop:holder-inter} and Lemma \ref{lem:pressure-L2} implies that
\begin{align}
\|\nabla_{x, z} \widetilde{P}\|_{\widetilde{L}^2_z(I_1\setminus I_0;C^{-\f{\varepsilon}2)}}
&\le K_\eta\big(\|\na_{x,y}v\|_{L^\infty(\Om_t)}^2+\|\nabla_{x,y}P_1\|_{L^2(\Om_t)}\big)+K_\eta\nonumber\\
&\le K_\eta A(t)^2.\label{eq:pressure-est4}
\end{align}

Putting the above estimates together concludes the lemma.\ef

\medskip

Next, we give the estimate of $\nabla_{x, z} \widetilde{\cD_t P}$ in Besov space:
\begin{lemma}\label{lem:pressure-Dt-holder}
Let $I_0=[-\f34,0]$ and $\cD_t P$ be a solution of (\ref{eq:pressure-Dt}). Then it holds that
\beno
\|\nabla_{x, z} \widetilde{\cD_t P}\|_{\widetilde{L}^\infty_z(I_0;B^0_{\infty,1})}\leq K_\eta\big(1+\|{\na_{x,z}\widetilde{v}}\|_{L^\infty_z(I;B^0_{\infty,1})}\big)A(t)^2.
\eeno
\end{lemma}

\no{\bf Proof.}\,We denote
\begin{align}
F_0=&\sum_{i=1}^d\pa_i\big(\al\big(\widetilde{\pa_k P}\widetilde{\om_{i,k}}+4\widetilde{\pa_kv^i}\widetilde{\pa_kP}\big)\big)-2\al(\widetilde{\pa_iv^j})(\widetilde{\pa_jv^k})\widetilde{\pa_kv^i}\nonumber\\
&+\Big(\sum_{i=1}^d(\pa_z\al_i-\pa_i\al)\big(\widetilde{\pa_k P}\widetilde{\om_{i,k}}+4\widetilde{\pa_kv^i}\widetilde{\pa_kP}\big)\nonumber\\
&-\pa_z\al_{d+1}\big(\widetilde{\pa_k P}\widetilde{\om_{d+1,k}}+4\widetilde{\pa_kv^{d+1}}\widetilde{\pa_kP}\big)\Big)\triangleq F_0^1+F_0^2+F_0^3,\label{eq:F0}\\
G_0=&-\sum_{i=1}^d\al_i\big(\widetilde{\pa_k P}\widetilde{\om_{i,k}}+4\widetilde{\pa_kv^i}\widetilde{\pa_kP}\big)+\al_{d+1}\big(\widetilde{\pa_k P}\widetilde{\om_{d+1,k}}
+4\widetilde{\pa_kv^{d+1}}\widetilde{\pa_kP}\big),\label{eq:G0}
\end{align}
where $\al_1=\al\f {\pa_i\rho_\delta} {\pa_z \rho_\delta}$ for $i=1,\cdots d$ and $\al_{d+1}=\f {\al} {\pa_z \rho_\delta}$.
Then $\cD_t P$ satisfies
\beno
\pa^2_z \cD_t P+\al \Delta \cD_t P+\beta\cdot\nabla \pa_z \cD_t P-\gamma \pa_z \cD_t P= F_{0}+\pa_zG_0.
\eeno
We apply Proposition \ref{prop:elliptic Holder est} to obtain
\begin{align*}
\|\nabla_{x, z} \widetilde{\cD_t P}\|_{\widetilde{L}^\infty_z(I_0;B^0_{\infty,1})}\le K_\eta\big(&\|F_0\|_{Y^0_1(I)}+\|G_0\|_{\widetilde{L}^\infty_z(I_0;B^0_{\infty,1})}+\|\nabla_{x, z} \widetilde{\cD_t P}\|_{\widetilde{L}^\infty_z(I_0;C^{-\delta})}\\
&+\|\nabla_{x, z} \widetilde{\cD_t P}\|_{\widetilde{L}^2_z(I_1\setminus I_0;C^{-\f{\varepsilon}2)}}\big),
\end{align*}
where $I_1=[-\f78, 0]$. By Lemma \ref{lem:pressure-F0-G0} and Lemma \ref{lem:pressure-holder}, we have
\begin{align}
&\|F_0\|_{Y^0_1(I)}+\|G_0\|_{\widetilde{L}^\infty_z(I_0;B^0_{\infty,1})}\nonumber\\
&\le K_\eta\|{\na_{x,z}\widetilde{v}}\|_{\widetilde{L}^\infty_z(I;B^0_{\infty,1})}
\big(\|\na_{x,y}v\|_{L^\infty(\Om_t)}^2+\|{\na_{x,z}\widetilde{P}}\|_{\widetilde{L}^\infty_z(I;C^\f12)}\big)\nonumber\\
&\le K_\eta\|{\na_{x,z}\widetilde{v}}\|_{\widetilde{L}^\infty_z(I;B^0_{\infty,1})}A(t)^2.\label{eq:pressure-F0-G0}
\end{align}
We get by Lemma \ref{lem:pressure-Dt-L2} and Lemmma \ref{lem:elliptic-H1} that
\begin{align*}
\|\na \widetilde{\cD_t P}\|_{\widetilde{L}^\infty_z(I_0;C^{-\delta})}\le& C\|\na\widetilde{\cD_t P}\|_{L^\infty_z(I_0;H^{-\f12})}\le C\|\na_{x,y}\cD_t P\|_{L^2(\Om_t)}\\
\le& K_\eta\big(1+\|\na_{x,y}v\|_{L^\infty(\Om_t)}\big)A(t)^2.
\end{align*}
On the other hand, for $z\in [-1,0]$,
\begin{align*}
\pa_z \widetilde{\cD_t P}(x,z)&=\int_{-1}^z\pa_z^2\widetilde{\cD_t P}(x,z')dz'\\
&=G(z)-G(-1)+\int_{-1}^z\big(F_0-\al\Delta \widetilde{\cD_t P}+\beta\nabla \pa_z\widetilde{\cD_t P}-\gamma\pa_z \widetilde{\cD_t P}\big)dz'.
\end{align*}
From Lemma 4.8 in \cite{WZ} and Lemma \ref{lem:pressure-Dt-L2}, we know that
\begin{align*}
\|\al\Delta \widetilde{\cD_t P}-\beta\nabla \pa_z\widetilde{\cD_t P}+\gamma\pa_z \widetilde{\cD_t P}\|_{L^1_z(-1,0;H^{-1})}\le& K_\eta\|\na_{x,z}\widetilde{\cD_t P}\|_{L^2(\overline{\cS})}\\
\le& K_\eta\big(1+\|\na_{x,y}v\|_{L^\infty(\Om_t)}\big)A(t),
\end{align*}
which along with (\ref{eq:pressure-F0-G0}) implies
\beno
\|\pa_z\widetilde{\cD_t P}\|_{\widetilde{L}^\infty_z(I_0;C^{-\delta})}\le K_\eta\big(1+\|{\na_{x,z}\widetilde{v}}\|_{\widetilde{L}^\infty_z(I;B^0_{\infty,1})}\big)A(t)^2.
\eeno
A similar argument leading to (\ref{eq:pressure-est4}) yields
\begin{align*}
&\|\nabla_{x, z} \widetilde{\cD_t P}\|_{\widetilde{L}^2_z(I_1\setminus I_0;C^{-\f{\varepsilon}2)}}\\
&\le K_\eta\big(\|\na_{x,y}v\|_{L^\infty(\Om_t)}\|\na_{x,y}P\|_{L^\infty(\Om_t)}+\|\na_{x,y}v\|_{L^\infty(\Om_t)}^3+\|\na_{x,y}\cD_t P\|_{L^2(\Om_t)}\big)\\
&\le K_\eta\big(1+\|\na_{x,y}v\|_{L^\infty(\Om_t)}\big)A(t)^2.
\end{align*}

Putting the above estimates together concludes the lemma.
\ef

Next, we give the estimates of $F_0$ and $G_0$.

\begin{lemma}\label{lem:pressure-F0-G0}
Let $F_0$ and $G_0$ be given by (\ref{eq:F0}) and (\ref{eq:G0}) respectively. Then we have
\beno
&&\|F_0\|_{Y^0_1(I)}+\|G_0\|_{\widetilde{L}^\infty_z(I;B^0_{\infty,1})}\\
&&\le K_\eta\|{\na_{x,z}\widetilde{v}}\|_{\widetilde{L}^\infty_z(I;B^0_{\infty,1})}
\big(\|\na_{x,y}v\|_{L^\infty(\Om_t)}^2+\|{\na_{x,z}\widetilde{P}}\|_{L^\infty_z(I;C^\f12)}\big).
\eeno
\end{lemma}

\no{\bf Proof.}\,We can deduce from Lemma \ref{lem:PP-Holder}, Lemma \ref{lem:PR}  that for $\delta<0$,
\beno
&&\|fg\|_{\widetilde{L}^p_z(I;B^\delta_{\infty,1})}\le C\|f\|_{\widetilde{L}^p_z(I;B^0_{\infty,1})}\|g\|_{L^\infty(\overline{\cS})},\\
&&\|fg\|_{\widetilde{L}^p_z(I;B^0_{\infty,1})}\le C\|f\|_{\widetilde{L}^p_z(I;B^0_{\infty,1})}\|g\|_{\widetilde{L}^\infty_z(I;B^0_{\infty,1})},
\eeno
which together with Lemma \ref{lem:coeff} imply that
\beno
&&\|G_0\|_{\widetilde{L}^\infty_z(I;B^0_{\infty,1})}\le K_\eta\|\na_{x,y}v\|_{L^\infty(\Om_t)}\|{\na_{x,z}\widetilde{P}}\|_{\widetilde{L}^\infty_z(I;C^\f12)},\\
&&\|F_0^1\|_{\widetilde{L}^\infty_z(I;B^{-1}_{\infty,1})}\le K_\eta\|{\na_{x,z}\widetilde{v}}\|_{\widetilde{L}^\infty_z(I;B^0_{\infty,1})}\|{\na_{x,z}\widetilde{P}}\|_{\widetilde{L}^\infty_z(I;C^\f12)},\\
&&\|F_0^2\|_{\widetilde{L}^2_z(I_0;B^{-\f12}_{\infty,1})}\le K_\eta\|\na_{x,y}v\|_{L^\infty(\Om_t)}^2\|{\na_{x,z}\widetilde{v}}\|_{\widetilde{L}^\infty_z(I;B^0_{\infty,1})},\\
&&\|F_0^3\|_{\widetilde{L}^2_z(I;B^{-\f12}_{\infty,1})}\le K_\eta\|\na_{x,y}v\|_{L^\infty(\Om_t)}\|{\na_{x,z}\widetilde{P}}\|_{\widetilde{L}^\infty_z(I;C^\f12)}.
\eeno
Here we also used the fact that
\beno
\|\widetilde{\na_{x,y}v}\|_{\widetilde{L}^\infty_z(I;B^0_{\infty,1})}\le K_\eta\|{\na_{x,z}\widetilde{v}}\|_{\widetilde{L}^\infty_z(I;B^0_{\infty,1})}.
\eeno

This gives the lemma by the definition of $Y^0_1(I)$.\ef

\subsection{Sobolev estimate of the pressure}

\begin{lemma}\label{lem:pressure-Hs}
Let $I_0=[-\f12,0]$ and $P$ be a solution of (\ref{eq:pressure equ}). Then it holds that
\beno
&&\|\nabla_{x, z}\widetilde{P_1}\|_{X^{s-\f12}(I_0)}\le K_\eta A(t)^2\big(\|\na_{x,z}\widetilde{v}\|_{L^2_z(I;H^{s-1})}+\|\eta\|_{H^{s+\f12}}\big),\\
&&\|\na_{x,y}P_1(\cdot,-1)\|_{H^s}\le K_\eta A(t)^2\big(\|\na_{x,z}\widetilde{v}\|_{H^{s-\f12}(\overline{\cS})}+\|\eta\|_{H^{s+\f12}}\big).
\eeno
\end{lemma}

\no{\bf Proof.}\,Apply Proposition \ref{prop:elliptic-upboun} with $F_0=\al(\widetilde{\pa_iv^j}\widetilde{\pa_jv^i})$ and $f=0$ to obtain
\begin{align*}
\|\nabla_{x, z}\widetilde{P_1}\|_{X^{s-\f12}(I_0)}\le K_\eta\big(&\|\na_{x,z}\widetilde{P_1}\|_{L^2(\overline{\cS})}+\|F_0\|_{Y^{s-\f12}(I)}\\
&\quad+\|\eta\|_{H^{s+\f12}}\|\na_{x,z}\widetilde{P_1}\|_{L^\infty(\R^d\times I_1)}\big)
\end{align*}
with $I_1=[-\f34,0]$. By Lemma \ref{lem:product} and Lemma \ref{lem:coeff}, we have
\begin{align*}
\|F_0\|_{Y^{s-\f12}(I)}\le& \|\widetilde{\pa_iv^j}\widetilde{\pa_jv^i}\|_{L^2_z(I;H^{s-1})}+\|(\al-1)(\widetilde{\pa_iv^j}\widetilde{\pa_jv^i})\|_{L^2_z(I;H^{s-1})}\\
\le& K_\eta\|\widetilde{\pa_iv^j}\|_{L^\infty(\overline{\cS})}\|\widetilde{\pa_jv^i}\|_{L^2_z(I;H^{s-1})}
+C\|\widetilde{\pa_iv^j}\|_{L^\infty(\overline{\cS})}^2\|\al-1\|_{L_z^2(I;H^{s-1})}\\
\le& K_\eta\big(\|\na_{x,y}v\|_{L^\infty(\Om_t)}\|\na_{x,z}\widetilde{v}\|_{L^2_z(I;H^{s-1})}+\|\na_{x,y}v\|_{L^\infty(\Om_t)}^2\|\eta\|_{H^{s-\f12}}\big),
\end{align*}
which along with Lemma \ref{lem:pressure-L2} and Lemma \ref{lem:pressure-holder} gives the first inequality.

In fact, the first inequality also holds with $[a,0]$ for $a>-1$ instead of $I_0$. This implies there exists $y_0\in (-1,-1+h_0)$  so that
\beno
\|P_1(\cdot,y_0)\|_{H^s}\le K_\eta A(t)^2\big(\|\na_{x,z}\widetilde{v}\|_{L^2_z(I;H^{s-1})}+\|\eta\|_{H^{s+\f12}}\big).
\eeno
Let $\cS_1=\R^d\times [-1,y_1]$ where $y_1<y_0$. Then the standard elliptic estimate ensures that
\begin{align*}
\|\na_{x,y}P_1\|_{H^{s+\f12}(\cS_1)}\le& C\big(\|P_1(\cdot,y_0)\|_{H^s}+\|\na_{x,y}v\na_{x,y}v\|_{H^{s-\f12}(\cS_1)}\big)\\
\le&  K_\eta A(t)^2\big(\|\na_{x,z}\widetilde{v}\|_{L^2_z(I;H^{s-1})}+\|\na_{x,y}{v}\|_{H^{s-\f12}({\cS_1})}+\|\eta\|_{H^{s+\f12}}\big),
\end{align*}
which together with the fact that
\beno
\|\na_{x,y}v\|_{H^{s-\f12}(\cS_1)}\le K_\eta\big(\|\na_{x,z}\widetilde{v}\|_{H^{s-\f12}(\overline{\cS})}+\|\eta\|_{H^{s+\f12}}\big),
\eeno
implies the second inequality by the trace theorem.\ef

\subsection{Estimate of $a$}

Recall that $a(t,x)=-(\pa_y P)(t,x,\eta(t,x))$.

\begin{lemma}\label{lem:a-holder}
It holds that
\begin{align*}
&\|a\|_{C^{\f12}}\leq K_\eta A(t)^2,\\
&\|D_ta\|_{B^0_{\infty,1}}\le K_\eta\big(1+\|{\widetilde{v}}\|_{L_z^\infty(I;B^1_{\infty,1})}\big)A(t)^2.
\end{align*}
\end{lemma}

\no{\bf Proof.}\,The first inequality is a direct consequence of Lemma \ref{lem:pressure-holder} and Lemma \ref{lem:coeff2} by using
$\pa_yP=\f {\pa_z \widetilde{P}} {\pa_z\rho_\delta}.$ To show the second inequality, we write
\begin{align*}
\pa_t a+T_V\cdot\na a=&\pa_y(\pa_t+v\cdot\na_{x,y}P)-\pa_yv\cdot\na_{x,y}P|_{y=\eta}+(T_V-V)\cdot\na a\\
=&\f {\pa_z \widetilde{\cD_t P}} {\pa_z\rho_\delta}-\f {\pa_z \widetilde{v}} {\pa_z\rho}\cdot
\big(\na \widetilde{P}-\f {\na \rho_\delta} {\pa_z\rho_\delta}\pa_z \widetilde{P}, \f {\pa_z \widetilde{P}} {\pa_z\rho}\big)\big|_{z=0}+(T_V-V)\cdot\na a\\
\triangleq&I_1+I_2+I_3.
\end{align*}
Using Lemma \ref{lem:remaider}, Lemma \ref{lem:nonlinear} and Lemma \ref{lem:coeff2}, it is easy to show that
\beno
&&\|I_1\|_{B^{0}_{\infty,1}}\le K_\eta\|\pa_z\widetilde{\cD_t P}\|_{L_z^\infty(I;B^0_{\infty,1})}, \\
&&\|I_2\|_{B^{0}_{\infty,1}}\le K_\eta\|\na_{x,z}\widetilde{P}\|_{L_z^\infty(I;C^\f12)}\|{\na_{x,z}\widetilde{v}}\|_{L_z^\infty(I;B^0_{\infty,1})},\\
&&\|I_3\|_{B^0_{\infty,1}}\le C\|V\|_{B^1_{\infty,1}}\|a\|_{C^\f12}\le C\|{\widetilde{v}}\|_{L_z^\infty(I;B^1_{\infty,1})}\|a\|_{C^\f12},
\eeno
which along with Lemma \ref{lem:pressure-holder} and Lemma \ref{lem:pressure-Dt-holder} lead to the second inequality.
\ef

\begin{lemma}\label{lem:a-Hs}
It holds that
\begin{align*}
\|a-1\|_{H^{s-\f12}}\le K_\eta A(t)^2\big(\|\na_{x,z}\widetilde{v}\|_{L^2_z(I;H^{s-1})}+\|\eta\|_{H^{s+\f12}}\big).
\end{align*}
\end{lemma}

\no{\bf Proof.}\,Note that $a-1=\f {1} {\pa_z\rho_\delta}\pa_z\widetilde{P_1}\big|_{z=0}$. Then we deduce from Lemma \ref{lem:product},
Lemma \ref{lem:nonlinear} and Lemma \ref{lem:coeff2} that
\beno
\|a-1\|_{H^{s-\f12}}\le K_\eta\big(\|\pa_z\widetilde{P_1}\|_{L^\infty_z([-\f12,0];H^{s-\f12})}+\|\eta\|_{H^{s+\f12}}\big),
\eeno
which together with Lemma \ref{lem:pressure-Hs} yields the result.\ef

\section{Estimate of the velocity}

The velocity $v$ satisfies
\beno
\left\{
\begin{array}{ll}
\Delta_{x, y} v=-\na_{x,y}\times\om\quad \textrm{in} \quad \Omega_t,\\
v|_{y=\eta}=(V, B),\quad v|_{y=-1}=(V_b,0).
\end{array}\right.
\eeno
Here $\om$ is the vorticity.

\subsection{Sobolev estimate of the velocity}
The following $H^1$ estimate is classical:
\ben\label{eq:v-H1}
\|v\|_{H^1(\Om_t)}\le K_\eta\big(\|\om\|_{L^2(\Om_t)}+\|(V,B,V_b)\|_{H^\f12}\big).
\een
This together with Proposition \ref{prop:elliptic-upboun} and Proposition \ref{prop:elliptic-boun-b} ensures that

\begin{lemma}\label{lem:v-Hs}
Let $I_0=[a,0]$ for $a\in (-1,0)$. It holds that
\begin{align*}
\|\na_{x,z}\widetilde{v}\|_{X^{s-1}(I_0)}\le K_\eta\big(\|(V,B,V_b)\|_{H^s}+\|\widetilde{\om}\|_{H^{s-\f12}(\overline{\cS})}+\|\na_{x,y}v\|_{L^\infty(\Om_t)}\|\eta\|_{H^{s+\f12}}\big).
\end{align*}
If $s-\f12$ is an integer, we have
\begin{align*}
\|\na_{x,z}\widetilde{v}\|_{H^{s-\f12}(\overline{\cS})}\le K_\eta\big(\|(V,B,V_b)\|_{H^s}+\|\widetilde{\om}\|_{H^{s-\f12}(\overline{\cS})}+\|\na_{x,y}v\|_{L^\infty(\Om_t)}\|\eta\|_{H^{s+\f12}}\big).
\end{align*}
\end{lemma}

\subsection{The estimate of the irrotational part}

The irrotational part $v_{ir}$ of the velocity is defined by
\beno
\left\{
\begin{array}{ll}
\Delta_{x, y} v_{ir}=0\quad \textrm{in} \quad \Omega_t,\\
v_{ir}|_{y=\eta}=(V, B),\quad v_{ir}|_{y=-1}=0.
\end{array}\right.
\eeno

\begin{lemma}\label{lem:vir}
Let $I_0=[-\f34,0]$. It holds that
\begin{align*}
&\|\na_{x,z}\widetilde{v_{ir}}\|_{\widetilde{L}^\infty_z(I_0; C^0)}\le K_\eta A(t),\\
&\|\na_{x,z}\widetilde{v_{ir}}\|_{\widetilde{L}^\infty_z(I_0;B^0_{\infty,1})}\le K_\eta\big(\|\widetilde{v}\|_{L^\infty(I;B^{1}_{\infty,1})}+\|v\|_{H^1(\Om_t)}\big),
\end{align*}
where $A(t)$ is given by (\ref{def:A(t)}).
\end{lemma}

\no{\bf Proof.}\,Let $I_1=[-\f78,0]$. It follows from Proposition \ref{prop:elliptic Holder est} that
\begin{align*}
\|\na_{x,z}\widetilde{v_{ir}}\|_{\widetilde{L}^\infty_z(I_0; C^0)}\le K_\eta\big(\|(V,B)\|_{C^1}+\|\na_{x,z}\widetilde{v}_{ir}\|_{\widetilde{L}^\infty_z(I_0; C^{-\delta})}+\|\na_{x,z}\widetilde{v}_{ir}\|_{\widetilde{L}^2_z(I_1\setminus I_0; C^{-\f \varepsilon2})}\big)
\end{align*}
for any $\delta>0$. Take $\delta>0$ so that $-\delta+\f d 2\le -\f 1 2$. Then we get by Lemma \ref{lem:elliptic-H1} that
\begin{align*}
\|\na_{x,z}\widetilde{v_{ir}}\|_{\widetilde{L}^\infty_z(I_0; C^{-\delta})}\le& C\|\na_{x,z}\widetilde{v}_{ir}\|_{L^\infty(I_0;H^{-\f12})}\\
\le& K_\eta\|\na_{x,y}v_{ir}\|_{L^2(\Om_t)}\le K_\eta\|(V,B)\|_{H^\f12}.
\end{align*}
While, Proposition \ref{prop:holder-inter} implies that
\beno
\|\na_{x,z}\widetilde{v_{ir}}\|_{\widetilde{L}^2_z(I_1\setminus I_0; C^{-\f \varepsilon2})}\le K_\eta\|\na_{x,y}v_{ir}\|_{L^2(\Om_t)}\le K_\eta\|(V,B)\|_{H^\f12}.
\eeno
This shows that
\begin{align*}
\|\na_{x,z}\widetilde{v_{ir}}\|_{\widetilde{L}^\infty_z(I_0; C^0)}\le& K_\eta\big(\|(V,B)\|_{C^1}+\|(V,B)\|_{H^\f12}\big)\le K_\eta A(t).
\end{align*}
Similar argument leads to
\begin{align*}
\|\na_{x,z}\widetilde{v_{ir}}\|_{\widetilde{L}^\infty_z(I_0;B^0_{\infty,1})}\le& K_\eta\big(\|(V,B)\|_{B^1_{\infty,1}}+\|(V,B)\|_{H^\f12}\big)\nonumber\\
\le& K_\eta\big(\|\widetilde{v}\|_{L^\infty(I;B^{1}_{\infty,1})}+\|v\|_{H^1(\Om_t)}\big).
\end{align*}
This finishes the proof.
\ef

With Lemma \ref{lem:vir}, we deduce from Proposition \ref{prop:DN-Hs} and Proposition \ref{prop:DN-Holder} that
\begin{proposition}\label{prop:remainder}
It holds that
\begin{align*}
&\|R(\eta)(V,B)\|_{C^\f12}\le K_\eta A(t),\\
&\|R(\eta)(V,B)\|_{H^{s-1}}\le K_\eta\big(\|(V,B)\|_{H^{s-\f12}}+A(t)\|\eta\|_{H^{s+\f12}}\big),\\
&\|R(\eta)(V,B)\|_{H^{s-\f12}}\le K_\eta\big(\|(V,B)\|_{H^s}+(\|\widetilde{v}\|_{L^\infty(I;B^{1}_{\infty,1})}+A(t))\|\eta\|_{H^{s+\f12}}\big).
\end{align*}
\end{proposition}

\medskip

\subsection{The estimate of the rotational part}
The rotational part $v_{\om}$ of the velocity is defined by
\beno
\left\{
\begin{array}{ll}
\Delta_{x, y} v_{\om}=-\na_{x,y}\times\om\quad \textrm{in} \quad \Omega_t,\\
v_{\om}|_{y=\eta}=0,\quad v_{\om}|_{y=-1}=(V_b,0).
\end{array}\right.
\eeno

\begin{lemma}\label{lem:vr}
Let $I_0=[-\f34,0]$. It holds that
\begin{align*}
&\|\na_{x,z}\widetilde{v_{\om}}\|_{\widetilde{L}^\infty_z(I_0; C^0)}\le K_\eta A(t),\\
&\|\na_{x,z}\widetilde{v_{\om}}\|_{\widetilde{L}^\infty_z(I_0;B^0_{\infty,1})}\le K_\eta\big(\|\widetilde{v}\|_{L^\infty(I;B^{1}_{\infty,1})}+\|v\|_{H^1(\Om_t)}\big),\\
&\|\na_{x,z}\widetilde{v_{\om}}\|_{X^{s-1}([-\f12,0])}\le K_\eta\big(\|V_b\|_{H^\f12}+\|\widetilde{\om}\|_{H^{s-\f12}(\overline{\cS})}+A(t)\|\eta\|_{H^{s}}\big).
\end{align*}
\end{lemma}

\no{\bf Proof.}\, The first two inequalities follows from Lemma \ref{lem:vir} and the fact $v_\om=v-v_{ir}$.
Then we get by Proposition \ref{prop:elliptic-upboun} that
\begin{align*}
\|\na_{x,z}\widetilde{v_{\om}}\|_{X^{s-1}([-\f12,0])}\le& K_\eta\big(\|\na_{x,y}v_\om\|_{L^2(\Om_t)}
+\|\widetilde{\om}\|_{H^{s-\f12}(\overline{\cS})}+\|\na_{x,y}\widetilde{v_{\om}}\|_{\widetilde{L}^\infty_z(I_0;C^0)}\|\eta\|_{H^{s}}\big),
\end{align*}
which along with the first inequality yields the third inequality.\ef\medskip

Next we show that $(\pa_t+v\cdot\na_{x,y})v_\om\triangleq \cD_t v_\om$ has the similar estimates.
It is easy to verify that $\cD_t v_\om$ satisfies
\beno
\left\{
\begin{array}{ll}
\Delta_{x, y} \cD_tv_{\om}=G_\om\quad \textrm{in} \quad \Omega_t,\\
\cD_t v_{\om}|_{y=\eta}=0,\quad \cD_tv_{\om}|_{y=-1}=(\dot V_b,0).
\end{array}\right.
\eeno
where $\dot V_b=(\pa_t+V_b\cdot\na)V_b$ and
\beno
G_\om=-\na_{x,y}\times\cD_t \om+\na_{x,y}\cdot(\om\cdot \na_{x,y}\times v)+ \na\times(\om \cdot \na v_\om)+2\pa_i(\na_{k}v^i\cdot\na_{k}v_\om).
\eeno
It is easy to find that $G_\om$ can be rewritten $\widetilde{G}_\om= \widetilde{G}_\om^0+ \na_{ z} \widetilde{G}_\om^1 $ where $G_\om^1$ satisfies that (by Lemma \ref{lem:product-full} and Lemma \ref{lem:map-regularity}),
\begin{align*}
&\|\widetilde{G}_\om^0\|_{L^2(I_0;H^{s-\f32})}+\|\widetilde{G}_\om^1\|_{L^2(I_0;H^{s-\f12})}\\
& \le K_\eta\big(\|\na_{x,y}v\|_{L^\infty(\Om_t)}^2+
\|\na_{x,y}v\|_{L^\infty(\Om_t)}\|\na_{x,y}v_\om\|_{L^\infty(\Om_t)}\big)\|\eta\|_{H^{s+\f12}}\\
&+K_\eta\big(\|\na_{x,y}v\|_{L^\infty(\Om_t)}
+\|\na_{x,y}v_\om\|_{L^\infty(\Om_t)}\big)\|\na_{x,z}\widetilde{v}\|_{L^2(I_0;H^{s-\f12})}\\
&+K_\eta\|\na_{x,y}v\|_{L^\infty(\Om_t)}\|\na_{x,z}\widetilde{v_\om}\|_{L^2(I_0;H^{s-\f12})}.
\end{align*}

\begin{lemma}\label{lem:vr-Dt}
It holds that
\begin{align*}
&\|\na_{x,z}\widetilde{\cD_t v_\om}\|_{\widetilde{L}^\infty([-\f12,0];C^0)}\le K_\eta\big(1+\|\widetilde{v}\|_{L^\infty(I;B^{1}_{\infty,1})}\big)A(t)^2,\\
&\|\na_{x,z}\widetilde{\cD_t v_{\om}}\|_{X^{s-1}([-\f12,0])}\le K_\eta\big(1+\|\widetilde{v}\|_{L^\infty(I;B^{1}_{\infty,1})}\big)A(t)^2\\
&\qquad\qquad\qquad\qquad\qquad\times\big(\|\eta\|_{H^{s+\f12}}+\|(V,B,V_b)\|_{H^{s}}+\|\widetilde{\om}\|_{H^{s-\f12}}\big).
\end{align*}
\end{lemma}

\no{\bf Proof.}\,Thanks to $\cD_t\om=\om\cdot\na v$, we get
\beno
\|\cD_t\om\|_{L^2(\Om_t)}\le \|\na_{x,y}v\|_{L^\infty(\Om_t)}\|\na_{x,y}v\|_{L^2(\Om_t)},
\eeno
and by (\ref{eq:Vb}) and Lemma \ref{lem:pressure-H2}, we have
\beno
\|\dot V_b\|_{H^\f12}\le K_\eta\|\na_{x,y}P\|_{H^1(\Om_t)}\le K_\eta A(t)^2.
\eeno
Then it is easy to show that
\begin{align*}
\|\na_{x,y}\cD_t v_\om\|_{L^2(\Om_t)}\le& K_\eta\big(\|\dot V_b\|_{H^\f12}+\|\na_{x,y}v\|_{L^\infty(\Om_t)}(\|\na_{x,y}v\|_{L^2(\Om_t)}+\|\na_{x,y}v_\om\|_{L^2(\Om_t)})\big)\\
\le& K_\eta A(t)^2.
\end{align*}
Then a similar argument leading to Lemma \ref{lem:pressure-Dt-holder} ensures that
\begin{align*}
\|\na_{x,z}\widetilde{\cD_t v_\om}\|_{L^\infty([-\f34,0];C^0)}\le& K_\eta\big(A(t)^2+\|\na_{x,y}v_\om\|_{L^\infty(\Om_t)}A(t)\big)\nonumber\\
\le& K_\eta\big(1+\|\widetilde{v}\|_{L^\infty(I;B^{1}_{\infty,1})}\big)A(t)^2.
\end{align*}
Let $I_0=[-\f12,0]$ and $I_1=[-\f34,0]$.
Then Proposition \ref{prop:elliptic-upboun} together with Lemma \ref{lem:v-Hs} and Lemma \ref{lem:vr} implies that
\begin{align*}
\|\na_{x,z}\widetilde{\cD_tv_{\om}}\|_{X^{s-1}(I_0)}\le& K_\eta\big(\|\na_{x,y}\cD_t v_\om\|_{L^2(\Om_t)}+ \|G_\om^0\|_{L^2(I_0;H^{s-\f32})}+ \|G_\om^1\|_{L^2(I_0;H^{s-\f12})}\\
&\qquad\qquad\quad+\|\na_{x,z}\widetilde{\cD_tv_\om}\|_{L^\infty(I_1;C^0)}\|\eta\|_{H^{s+\f12}}\big)\\
\le& K_\eta\big(1+\|\widetilde{v}\|_{L^\infty(I;B^{1}_{\infty,1})}\big)A(t)^2\\
&\quad\times\big(\|\eta\|_{H^{s+\f12}}+\|(V,B,V_b)\|_{H^{s}}+\|\widetilde{\om}\|_{H^{s-\f12}}\big).
\end{align*}
The proof is finished.\ef

\section{Proof of break-down criterion}
In this section, we assume that $(\eta, v)$ is a solution of the system (\ref{eq:euler})--(\ref{eq:euler-p}) obtained in Theorem \ref{thm:local}
in $\Om_t=\big\{(x,y)\in\mathbf{R}^d\times\mathbf{R}:-1<y<\eta(t,x)\big\}$ for $t\in [0,T]$.
We denote
\beno
&&A(t)\triangleq 1+\|v(t)\|_{W^{1,\infty}(\Om_t)}+\|v(t)\|_{H^1(\Om_t)}+\|\eta(t)\|_{H^\f32}+\f{1}{c_0},\\
&&\cB(t)\triangleq 1+\|\widetilde{v}\|_{L^\infty_z(I;B^1_{\infty,1})}.
\eeno

\subsection{The $H^1$ energy estimate}

We have the following basic energy law for the system (\ref{eq:euler})--(\ref{eq:euler-p}).

\begin{lemma}\label{lem:basic energy}
For any $t\in [0,T]$, there holds
\beno
E(t)=E(0),\quad E(t)\triangleq\|v(t)\|^2_{L^2(\Om_t)}+\|\eta(t)\|_{L^2}^2.
\eeno
\end{lemma}

\no{\bf Proof.}\,By (\ref{eq:euler}), (\ref{eq:euler-b2}) and integration by parts, we get
\begin{align*}
\f d {dt}\int_{\Om_t}|v(t,x,y)|^2dxdy=&\int_{\R^d}\pa_t\eta|v|^2dx+2\int_{\Om_t}\pa_tv\cdot vdxdy\\
=&\int_{\R^d}\pa_t\eta|v|^2dx-2\int_{\Om_t}\big(v\cdot\na_{x,y}v+\na(P+y)\big)\cdot vdxdy\\
=&\int_{\R^d}\pa_t\eta|v|^2dx-\int_{\R^d}v\cdot \textbf{n}_+(|v|^2+2\eta)\sqrt{1+|\na\eta|^2}dx\\
=&-2\int_{\R^d}\pa_t\eta(t,x)\eta(t,x)dx=-\f d {dt}\int_{\R^d}|\eta(t,x)|^2dx.
\end{align*}
This shows that
\beno
\f d {dt}\Big(\int_{\Om_t}|v(t,x,y)|^2dxdy+\int_{\R^d}|\eta(t,x)|^2dx\Big)=0.
\eeno
Hence, $E(t)=E(0)$ for $t\in [0,T]$.\ef

\begin{lemma}\label{lem:velocity-H1}
It holds that for any $t\in [0,T]$,
\beno
\f {d} {dt}\|\na_{x,y}v(t)\|_{L^2(\Om_t)}\le K_\eta\big(\|v(t)\|_{W^{1,\infty}(\Om_t)}\|v(t)\|_{H^1(\Om_t)}+\|\eta(t)\|_{H^\f32}\big).
\eeno
\end{lemma}

\no{\bf Proof.}\,Similar to the proof of Lemma \ref{lem:basic energy}, we have
\begin{align}
&\f d {dt}\int_{\Om_t}|\na_{x,y}v(t,x,y)|^2dxdy\nonumber\\
&=\int_{\R^d}\pa_t\eta|\na_{x,y}v|^2dx+2\int_{\Om_t}\pa_t\na_{x,y}v\cdot\na_{x,y}vdxdy\nonumber\\
&=\int_{\R^d}\pa_t\eta|\na_{x,y}v|^2dx-\int_{\Om_t}v\cdot\na_{x,y}|\na_{x,y}v|^2dxdy\nonumber\\
&\quad-2\int_{\Om_t}\big(\na_{x,y}v\cdot\na_{x,y}v+\na_{x,y}^2(P+y)\big)\cdot\na_{x,y}vdxdy\nonumber\\
&=-2\int_{\Om_t}\big(\na_{x,y}v\cdot\na_{x,y}v+\na_{x,y}^2(P+y)\big)\cdot\na_{x,y}vdxdy\nonumber\\
&\le 2\|\na_{x,y}v\|_{L^\infty(\Om_t)}\|\na_{x,y}v\|_{L^2(\Om_t)}^2+\|\na_{x,y}^2(P+y)\|_{L^2(\Om_t)}\|\na_{x,y}v\|_{L^2(\Om_t)},\nonumber
\end{align}
which along with Lemma \ref{lem:pressure-H2} yields the result.\ef

\subsection{Energy estimate of the trace of the velocity and the free surface}

Let us first present the lower order energy estimate.

\begin{proposition}\label{prop:energy-lower}
It holds that
\begin{align*}
&\f d {dt}\big(\|(V,B)\|_{H^{s-\f12}}^2+\|V_b\|_{H^s}^2+\|\eta\|_{H^s}^2\big)\\
&\le K_\eta A(t)^2\Big(\|(V,B,V_b)\|_{H^s}^2+\|\na_{x,z}\widetilde{v}\|_{H^{s-\f12}(\overline{\cS})}^2+
\|\eta\|_{H^{s+\f12}}^2\Big).
\end{align*}
\end{proposition}

\no{\bf Proof.}\,Recall that $\eta(t,x)$ satisfies
\beno
\pa_t\eta+V\cdot\na \eta=B.
\eeno
Make $H^s$ energy estimate to obtain
\begin{align*}
\f12\f {d} {dt}\|\eta(t)\|_{H^s}^2\le&-\big\langle\langle D\rangle^s(V\cdot\na \eta), \langle D\rangle^s\eta\big\rangle+\|B\|_{H^s}\|\eta\|_{H^s}.
\end{align*}
We write
\begin{align*}
\big\langle\langle D\rangle^s(V\cdot\na \eta), \langle D\rangle^s\eta\big\rangle
=&\big\langle\langle D\rangle^s(T_V\cdot\na \eta), \langle D\rangle^s\eta\big\rangle+\big\langle\langle D\rangle^s((V-T_V)\cdot\na \eta), \langle D\rangle^s\eta\big\rangle\\
=&\big\langle[\langle D\rangle^s,T_V]\cdot\na\eta, \langle D\rangle^s\eta\big\rangle-\big\langle T_{\na\cdot V}\langle D\rangle^s\eta, \langle D\rangle^s\eta\big\rangle\\
&+\big\langle\langle D\rangle^s((V-T_V)\cdot\na \eta), \langle D\rangle^s\eta\big\rangle.
\end{align*}
Then we deduce from Lemma \ref{lem:commu-Ds} and Lemma \ref{lem:remaider} that
\beno
\big\langle\langle D\rangle^s(V\cdot\na \eta), \langle D\rangle^s\eta\big\rangle
\le C\|V\|_{W^{1,\infty}}\|\eta\|_{H^s}^2+K_\eta\|V\|_{H^s}\|\eta\|_{H^s}.
\eeno
This shows that
\begin{align}\label{eq:eta-Hs}
\f {d} {dt}\|\eta(t)\|_{H^s}^2\le& K_\eta A(t)\big(\|(V,B)\|_{H^s}^2+\|\eta\|_{H^s}^2\big).
\end{align}

Recall that $(V,B, V_b)$ satisfies
\beno
&&(\pa_t+V\cdot\nabla)B=a-1,\\
&&(\pa_t+V\cdot\nabla)V+a\zeta=0,\\
&&(\pa_t+V_b\cdot\nabla)V_b+\na P|_{y=-1}=0.
\eeno
In a similar way leading to (\ref{eq:eta-Hs}), we deduce
\begin{align*}
&\f12\f {d} {dt}\big(\|(V,B)\|_{H^{s-\f12}}^2+\|V_b\|_{H^s}^2\big)\\
&\le C\|(\na V,\na B)\|_{L^\infty}\|(V,B)\|_{H^{s-\f12}}^2+\|a-1\|_{H^{s-\f12}}\|B\|_{H^{s-\f12}}\\
&\qquad+\|a\zeta\|_{H^{s-\f12}}\|V\|_{H^{s-\f12}}+C\|\na V_b\|_{L^\infty}\|V_b\|_{H^{s}}^2+\|\na P|_{y=-1}\|_{H^{s}}\|V_b\|_{H^{s}},
\end{align*}
which together with Lemma \ref{lem:a-Hs} and Lemma \ref{lem:pressure-Hs} yields
\begin{align*}
&\f {d} {dt}\big(\|(V,B)\|_{H^{s-\f12}}^2+\|V_b\|_{H^s}^2\big)\\
&\le K_\eta A(t)^2\big(\|(V,B,V_b)\|_{H^s}^2+\|\na_{x,z}\widetilde{v}\|_{H^{s-\f12}(\overline{\cS})}^2+
\|\eta\|_{H^{s+\f12}}^2\big),
\end{align*}
from which and (\ref{eq:eta-Hs}), we deduce the proposition.\ef\medskip

Next we present the high order energy estimate.

\begin{proposition}\label{prop:energy-boundary}
Let $U$ be a solution of (\ref{eq:ZK-para-2}). Then there holds
\begin{align*}
&\f12\f d {dt}\big(\|D_t U\|_{H^{s-\f12}}^2+\|T_{\sqrt{a\lambda}}U\|_{H^{s-\f12}}^2\big)\\
&\le \big\langle (f+f_\om)_{s-1/2}, (D_t U)_{s-1/2}\big\rangle+K_\eta\cB(t)A(t)^2\big(\|D_t U\|_{H^{s-\f12}}^2+\|U\|_{H^s}^2\big).
\end{align*}
Here we denote $f_s\triangleq \langle D\rangle^s f$.
\end{proposition}

\no{\bf Proof.}\,It follows from (\ref{eq:ZK-para-2}) that $(D_t U)_{s-1/2}$ satisfies
\begin{align*}
D_t(D_t U)_{s-1/2}+T_{a\lambda}U_{s-1/2}=&[D_t, \langle D \rangle^{s-1/2}]D_tU+[T_{a\lambda}, \langle D \rangle^{s-1/2}]U\nonumber\\
&+(f+f_\om)_{s-1/2}\triangleq F.
\end{align*}
Taking $L^2$ inner product with $(D_t U)_{s-1/2}$, we get
\begin{align}
\big\langle D_t(D_t U)_{s-1/2}, (D_t U)_{s-1/2}\big\rangle+\big\langle T_{a\lambda}U_{s-1/2}, (D_t U)_{s-1/2}\big\rangle
=\big\langle F, (D_t U)_{s-1/2}\big\rangle.\label{eq:U-E2}
\end{align}
By integration by parts, we get
\begin{align}
&\big\langle D_t(D_t U)_{s-1/2}, (D_t U)_{s-1/2}\big\rangle\nonumber\\
&=\f12\f d {dt}\big\langle (D_t U)_{s-1/2}, (D_t U)_{s-1/2}\big\rangle-\f12\big\langle \na\cdot V(D_t U)_{s-1/2}, (D_t U)_{s-1/2}\big\rangle\nonumber\\
&\qquad+\big\langle(T_V\cdot\na-V\cdot\na)(D_t U)_{s-1/2}, (D_t U)_{s-1/2}\big\rangle\nonumber\\
&\ge \f12\f d {dt}\big\langle (D_t U)_{s-1/2}, (D_t U)_{s-1/2}\big\rangle-C\|V\|_{B^1_{\infty,1}}\|D_t U\|_{H^{s-\f12}}^2,\label{eq:U-E3}
\end{align}
where we used the inequality
\beno
\|(T_V\cdot\na-V\cdot\na)(D_t U)_{s-1/2}\|_{L^2}\le C\|V\|_{B^1_{\infty,1}}\|D_t U\|_{H^{s-\f12}}.
\eeno
Similarly, we have
\begin{align*}
&\big\langle T_{a\lambda}U_{s-1/2}, (D_t U)_{s-1/2}\big\rangle\\
&=\big\langle (T_{\sqrt{a\lambda}})^*T_{\sqrt{a\lambda}}U_{s-1/2}, (D_t U)_{s-1/2}\big\rangle
+\big\langle\big(T_{a\lambda}-(T_{\sqrt{a\lambda}})^*T_{\sqrt{a\lambda}}\big)U_{s-1/2}, (D_t U)_{s-1/2}\big\rangle\\
&=\big\langle (T_{\sqrt{a\lambda}}U_{s-1/2}, T_{\sqrt{a\lambda}}(D_t U)_{s-1/2}\big\rangle
+\big\langle\big(T_{a\lambda}-(T_{\sqrt{a\lambda}})^*T_{\sqrt{a\lambda}}\big)U_{s-1/2}, (D_t U)_{s-1/2}\big\rangle\\
&=\f12\f d {dt}\big\langle T_{\sqrt{a\lambda}}U_{s-1/2}, T_{\sqrt{a\lambda}}U_{s-1/2}\big\rangle
-\f12\big\langle \na\cdot VT_{\sqrt{a\lambda}}U_{s-1/2}, T_{\sqrt{a\lambda}}U_{s-1/2}\big\rangle\\
&\quad+\big\langle T_{\sqrt{a\lambda}}U_{s-1/2}, [T_{\sqrt{a\lambda}}\langle D\rangle^{s-\f12}, D_t]U\big\rangle
+\big\langle T_{\sqrt{a\lambda}}U_{s-1/2}, (T_V\cdot\na-V\cdot\na)T_{\sqrt{a\lambda}}U_{s-1/2}\big\rangle\\
&\quad+\big\langle\big(T_{a\lambda}-(T_{\sqrt{a\lambda}})^*T_{\sqrt{a\lambda}}\big)U_{s-1/2}, (D_t U)_{s-1/2}\big\rangle.
\end{align*}
It follows from Proposition \ref{prop:symbolic calculus} that
\beno
&&\big\langle \na\cdot VT_{\sqrt{a\lambda}}U_{s-1/2}, T_{\sqrt{a\lambda}}U_{s-1/2}\big\rangle
\le C\|\na V\|_{L^\infty}M_0^\f12(\sqrt{a\lambda})^2\|U\|_{H^s}^2,\\
&&\big\langle\big(T_{a\lambda}-(T_{\sqrt{a\lambda}})^*T_{\sqrt{a\lambda}}\big)U_{s-1/2}, (D_t U)_{s-1/2}\big\rangle
\le CM_0^\f12(\sqrt{a\lambda})^2\|U\|_{H^s}\|D_tU\|_{H^{s-\f12}},\\
&&\big\langle T_{\sqrt{a\lambda}}U_{s-1/2}, (T_V\cdot\na-V\cdot\na)T_{\sqrt{a\lambda}}U_{s-1/2}\big\rangle
\le M_0^\f12(\sqrt{a\lambda})^2\|V\|_{B^{1}_{\infty,\infty}}\|D_tU\|_{H^{s-\f12}}^2.
\eeno
We get by Proposition \ref{prop:commutator-tame} that
\begin{align*}
&\big\langle T_{\sqrt{a\lambda}}U_{s-1/2}, [T_{\sqrt{a\lambda}}\langle D\rangle^{s-\f12}, D_t]U\big\rangle\\
&\le CM_0^\f12(\sqrt{a\lambda})\big(M_0^\f12(\sqrt{a\lambda})\|V\|_{B^{1}_{\infty,\infty}}+M_0^\f12(D_t\sqrt{a\lambda})\big)\|U\|_{H^s}^2.
\end{align*}
This proves that
\begin{align}
\big\langle& T_{a\lambda}U_{s-1/2}, (D_t U)_{s-1/2}\big\rangle\ge \f12\f d {dt}\big\langle T_{\sqrt{a\lambda}}U_{s-1/2}, T_{\sqrt{a\lambda}}U_{s-1/2}\big\rangle\nonumber\\
&-C\|V\|_{B^{1}_{\infty,\infty}}M_0^\f12(\sqrt{a\lambda})^2\big(\|U\|_{H^s}^2+\|D_tU\|_{H^{s-\f12}}^2\big)-M_0^\f12(D_t\sqrt{a\lambda})\|U\|_{H^s}^2\nonumber\\
&\quad-CM_0^\f12(\sqrt{a\lambda})^2\|U\|_{H^s}\|D_tU\|_{H^{s-\f12}}.\label{eq:U-E4}
\end{align}

By Lemma \ref{lem:commu-Ds} and Proposition \ref{prop:symbolic calculus}, we have
\begin{align}
&\|[D_t, \langle D \rangle^{s-1/2}]D_tU\|_{L^2}\le C\|V\|_{W^{1,\infty}}\|D_tU\|_{H^{s-\f12}},\label{eq:U-E5}\\
&\|[T_{a\lambda}, \langle D \rangle^{s-1/2}]U\|_{L^2}\le CM_{\f12}^1(a\lambda)\|U\|_{H^s}.\label{eq:U-E6}
\end{align}
Then it follows from (\ref{eq:U-E2})-(\ref{eq:U-E6}) that
\begin{align*}
&\f12\f d {dt}\big(\|D_t U\|_{H^{s-\f12}}^2+\|T_{\sqrt{a\lambda}}U\|_{H^{s-\f12}}^2\big)\le \langle (f+f_\om)_{s-1/2}, (D_t U)_{s-1/2}\big\rangle\\
&\quad+C\big(1+\|V\|_{B^1_{\infty,\infty}}\big)\big(1+M_0^\f12(\sqrt{a\lambda})^2+M_{\f12}^1(a\lambda)\big)\big(\|D_t U\|_{H^{s-\f12}}^2+\|U\|_{H^s}^2\big)\\
&\qquad+CM_0^\f12(D_t\sqrt{a\lambda})\|U\|_{H^s}^2,
\end{align*}
while, by Lemma \ref{lem:a-holder} and Lemma \ref{lem:symbol}, we have
\begin{align*}
&M_0^\f12(\sqrt{a\lambda})^2+M_{\f12}^1(a\lambda)\le K_\eta A(t)^2,\\
&M_0^\f12(D_t\sqrt{a\lambda})\le K_\eta\cB(t)A(t)^2.
\end{align*}
Then the proposition follows easily.\ef

\subsection{Energy estimate of the vorticity}

Using the equation (\ref{eq:vorticity}), it is easy to see that the vorticity $\widetilde{\om}(t,x,z)$ satisfies
\beno
\pa_t \widetilde{\om}+\overline{v}\cdot \nabla_{x,z}\widetilde{\om}=\widetilde{\om}^h\cdot\big(\na\widetilde{v}-\f {\na \rho_\delta} {\pa_z\rho_\delta}\pa_z\widetilde{v}\big)
+\widetilde{\om}^{d+1}\f {\pa_z \widetilde{v}} {\pa_z\rho_\delta}\triangleq F,
\eeno
where $\widetilde{v}^h=\big(\widetilde{v}^1, \cdots, \widetilde{v}^d\big), \widetilde{\om}^h=\big(\widetilde{\om}^1, \cdots, \widetilde{\om}^d\big)$ and
\beno
\overline{v}=\big(\widetilde{v}^h,\f{1}{\pa_z \rho_\delta}(\widetilde{v}^{d+1}-\pa_t\rho_\delta-\widetilde{v}^h\cdot\na\rho_\delta)\big).
\eeno

\begin{proposition}\label{prop:vorticity}
Let $s>\f d 2+1$ and $s-\f12$ be an integer. Then we have
\begin{align*}
\f d {dt}\|\widetilde{\om}(t)\|_{H^{s-\f12}(\overline{\cS})}^2\le&K_\eta\cB(t)A(t)\big(\|\widetilde{\om}\|_{H^{s-\f12}(\overline{\cS})}^2
+\|\eta\|_{H^{s}}^2+\|\widetilde{v}\|_{H^{s+\f12}(\overline{\cS})}^2\big).
\end{align*}
\end{proposition}

We need the following lemma.

\begin{lemma}\label{lem:v-trans}
Let $s>1+\f d 2$. Then we have
\beno
&&\|\na_{x,z}\overline{v}\|_{H^{s-\f12}(\overline{\cS})}\le K_\eta\big(\|\widetilde{v}\|_{H^{s+\f12}(\overline{\cS})}
+\cB(t)\|\eta\|_{H^{s}}\big),\\
&&\|\na_{x,z}\overline{v}\|_{L^{\infty}(\overline{\cS})}\le K_\eta\cB(t).
\eeno
\end{lemma}

\no{\bf Proof.}\,Thanks to the definition of $\overline{v}$, it suffices to consider $\overline{v}^{d+1}$.
Let $\phi=\widetilde{v}^{d+1}-\pa_t\rho_\delta-\widetilde{v}^h\cdot\na\rho_\delta$.
Let $\Delta^\delta_{x,z}=\Delta+\delta^{-2}\pa_z^2$ and $\na_{x,z}^\delta=(\na, \delta^{-1}\pa_z)$. Using the fact that
\beno
\Delta^\delta_{x,z}\rho_\delta=2\delta^{-1}e^{\delta z|D|}|D|\eta\triangleq f_\eta,
\eeno
we find that
\begin{align*}
\Delta^\delta_{x,z}\phi=
&\Delta^\delta_{x,z}\widetilde{v}^{d+1}-\pa_tf_\eta-\Delta^\delta_{x,z}\widetilde{v}^h\cdot\na\rho_\delta
-2\na^\delta_{x,z}\widetilde{v}^h\cdot\na\na_{x,z}^\delta\rho_\delta-\widetilde{v}\cdot\na f_\eta\triangleq F_{\eta,v}
\end{align*}
together with the boundary condition $\phi=0$ on $z=0$ and $z=-1$.

By (\ref{eq:euler-b2}), we have
\beno
\|\pa_tf_\eta\|_{H^{s-\f32}(\overline{\cS})}\le C\|\pa_t\eta\|_{H^{s-1}}\le C\big(\|\na \eta\|_{L^\infty}\|V\|_{H^{s-1}}+\|V\|_{L^\infty}\|\eta\|_{H^s}+\|B\|_{H^{s-1}}\big),
\eeno
which together with Lemma \ref{lem:product-full}, Lemma \ref{lem:map-regularity} implies that
\beno
\|F_{\eta,v}\|_{H^{s-\f32}(\cS)}\le K_\eta\big(\|\widetilde{v}\|_{H^{s+\f12}(\overline{\cS})}
+\|\widetilde{v}\|_{W^{1,\infty}(\overline{\cS})}\|\eta\|_{H^s}\big).
\eeno
Then Proposition \ref{prop:elliptic-flat-Hs} ensures that
\ben\label{eq:phi-hs}
\|\na_{x,z}\phi\|_{H^{s-\f12}(\overline{\cS})}\le K_\eta\big(\|\widetilde{v}\|_{H^{s+\f12}(\overline{\cS})}
+\|\widetilde{v}\|_{W^{1,\infty}(\overline{\cS})}\|\eta\|_{H^s}\big).
\een

We write
\begin{align*}
F_{\eta,v}=&\na_{x,z}^\delta\cdot\big(\na_{x,z}^\delta \widetilde{v}_{d+1}-\na_{x,z}^\delta \widetilde{v}^h\cdot\na \rho_\delta+(\widetilde{v}^h,0)f_\eta\big)\\
&+\big(-\pa_t f_\eta-\na_{x,z}^\delta \widetilde{v}^h\cdot \na_{x,z}^\delta\na\rho_\delta+\na\cdot \widetilde{v}^hf_\eta\big)\\
=&\na_{x,z}^\delta\cdot F_{\eta,v}^1+F_{\eta,v}^2.
\end{align*}
By Lemma \ref{lem:PP-Holder} and Lemma \ref{lem:PR}, we have
\beno
&&\|T_{\na_{x,z}^\delta \widetilde{v}^h}\cdot \na\rho_\delta\|_{L_z^\infty(I;B^{0}_{\infty,1})}\le C\|\na_{x,z}\widetilde{v}\|_{L^\infty(\overline{\cS})}\|\na\rho_\delta\|_{L_z^\infty(I;B^{0}_{\infty,1})}
,\\
&&\|T_{\na\rho_\delta}{\na_{x,z}^\delta\widetilde{v}^h}\|_{L_z^\infty(I;B^{0}_{\infty,1})}\le C\|\na\rho_\delta\|_{L^\infty(\overline{\cS})}\|\na_{x,z}\widetilde{v}\|_{L_z^\infty(I;B^{0}_{\infty,1})}
,\\
&&\|R(\na_{x,z}^\delta \widetilde{v}^h, \na\rho_\delta)\|_{L_z^\infty(I;B^{0}_{\infty,1})}
\le C\|\na_{x,z}\widetilde{v}\|_{L^\infty(\overline{\cS})}\|\na\rho_\delta\|_{L_z^\infty(I;B^{\epsilon}_{\infty,1})},
\eeno
for any $\epsilon>0$, which gives
\beno
\|\na_{x,z}^\delta \widetilde{v}^h\cdot \rho_\delta\|_{L_z^\infty(I;B^{0}_{\infty,1})}
\le K_\eta\|\widetilde{v}\|_{L^\infty_z(I;B^1_{\infty,1})}.
\eeno
The same estimate holds for $(\widetilde{v}^h,0)f_\eta$. Hence,
\beno
\|F^1_{\eta,v}\|_{L_z^\infty(I;B^{0}_{\infty,1})}\le K_\eta\|\widetilde{v}\|_{L^\infty_z(I;B^1_{\infty,1})}.
\eeno
By Lemma \ref{lem:PP-Holder} and Lemma \ref{lem:PR} again, we have
\beno
&&\|T_{\na_{x,z}^\delta \widetilde{v}^h}\na_{x,z}^\delta\na\rho_\delta\|_{\widetilde{L}_z^2(I;B^{-\f12}_{\infty,1})}
\le C\|\na_{x,z}\widetilde{v}\|_{L^\infty(\overline{\cS})}\|\na_{x,z}^\delta\rho_\delta\|_{\widetilde{L}_z^2(I;B^{\f12}_{\infty,1})},\\
&&\|T_{\na_{x,z}^\delta\na\rho_\delta}{\na_{x,z}^\delta \widetilde{v}^h}\|_{\widetilde{L}_z^2(I;B^{-\f12}_{\infty,1})}
\le C\|\na_{x,z}\widetilde{v}\|_{L^\infty_z(I;B^{0}_{\infty,1})}\|\na_{x,z}^\delta\na\rho_\delta\|_{\widetilde{L}_z^2(I;B^{-\f12}_{\infty,1})},\\
&&\|R(\na_{x,z}^\delta\na\rho_\delta,{\na_{x,z}^\delta \widetilde{v}^h})\|_{\widetilde{L}_z^2(I;B^{-\f12}_{\infty,1})}
\le C\|\na_{x,z}\widetilde{v}\|_{L^\infty(\overline{\cS})}\|\na_{x,z}^\delta\rho_\delta\|_{\widetilde{L}_z^2(I;C^{1+\epsilon})},
\eeno
for any $\epsilon>0$, which gives
\beno
\|{\na_{x,z}^\delta \widetilde{v}^h}\na_{x,z}^\delta\na\rho_\delta\|_{\widetilde{L}_z^2(I;B^{-\f12}_{\infty,1})}
\le K_\eta\|\na_{x,z}\widetilde{v}\|_{L^\infty(\overline{\cS})}.
\eeno
The same estimate holds for $\na\cdot \widetilde{v}^hf_\eta$ and $\pa_tf_\eta$. Hence,
\beno
\|F_{\eta,v}^2\|_{\widetilde{L}_z^2(I;B^{-\f12}_{\infty,1})}\le K_\eta\|\widetilde{v}\|_{W^{1,\infty}(\overline{\cS})}.
\eeno
Then we apply Proposition \ref{prop:elliptic-flat-Holder} to conclude that
\ben\label{eq:phi-hol-1}
&&\|\phi\|_{L^\infty_z(I;B^1_{\infty,1})}\le K_\eta\|\widetilde{v}\|_{L^\infty_z(I;B^1_{\infty,1})}.
\een

Next we turn to the estimate of $\overline{v}^{d+1}$, which satisfies
\begin{align*}
\Delta_{x,z}^\delta \overline{v}^{d+1}=&\Big(\f {2|\na_{x,z}^\delta\pa_z\rho_\delta|^2} {(\pa_z\rho_\delta)^3}-\f {\pa_z f_\eta} {(\pa_z\rho_\delta)^2}\Big)\phi
+2\na_{x,z}^\delta\big(\f 1 {\pa_z\rho_\delta}\big)\cdot\na_{x,z}\phi+\f 1 {\pa_z\rho_\delta}\Delta_{x,z}^\delta\phi\triangleq G_{\eta,v}
\end{align*}
with $\overline{v}_{d+1}|_{z=0}=0$ and  $\overline{v}_{d+1}|_{z=-1}=0$.
By Lemma \ref{lem:product-full} and Lemma \ref{lem:map-regularity}, we have
\begin{align*}
&\|G_{\eta,v}\|_{H^{s-\f32}(\overline{\cS})}\le  K_\eta\big(\|\phi\|_{H^{s+\f12}(\overline{\cS})}
+\|\na_{x,z}\phi\|_{L^\infty(\overline{\cS})}\|\eta\|_{H^{s}}\big),
\end{align*}
from which, (\ref{eq:phi-hs})-(\ref{eq:phi-hol-1}) and Proposition \ref{prop:elliptic-flat-Hs}, we deduce the first inequality.

We write
\begin{align*}
G_{\eta,v}=&\Big(\f {2|\na_{x,z}^\delta\pa_z \rho_\delta|^2} {(\pa_z\rho_\delta)^3}-\f {\pa_zf_\eta} {(\pa_z\rho_\delta)^2}\Big)\phi
+\na_{x,z}^\delta\big(\f 1 {\pa_z\rho_\delta}\big)\cdot\na_{x,z}\phi+\na_{x,z}^\delta\cdot\big(\f 1 {\pa_z\rho_\delta}\na_{x,z}^\delta\phi\big)\\
\triangleq& G_{\eta,v}^1+\na_{x,z}^\delta\cdot G_{\eta,v}^2.
\end{align*}
Then a similar argument leading to (\ref{eq:phi-hol-1}) yields
\begin{align*}
\|G_{\eta,v}^1\|_{\widetilde{L}_z^2(I;B^{-\f12}_{\infty,1})}+\|G^2_{\eta,v}\|_{L_z^\infty(I;B^{0}_{\infty,1})}\le K_\eta
\|\na_{x,z}\phi\|_{L^\infty_z(I;B^0_{\infty,1})}.
\end{align*}
Then Proposition \ref{prop:elliptic-flat-Holder} together with (\ref{eq:phi-hol-1}) gives the second inequality.\ef

\medskip

Now we are in position to prove Proposition \ref{prop:vorticity}.
\medskip

{\bf Proof of Proposition \ref{prop:vorticity}:}
Let $k\in [0,s-\f12]$ be an interger. Then we have
\begin{align*}
\f d {dt}\|\na_{x,z}^k\widetilde{\om}\|_{L^2(\overline{\cS})}^2=2\big\langle-\na_{x,z}^k(\overline{v}\cdot \nabla_{x,z}\widetilde{\om})+\na_{x,z}^kF, \na_{x,z}^k\widetilde{\om}\big\rangle.
\end{align*}

First of all, we have by Lemma \ref{lem:product-full} that
\begin{align*}
\|\na_{x,z}^kF\|_{L^2(\overline{\cS})}\le& K_\eta\|\na_{x,z}\widetilde{v}\|_{L^\infty(\overline{\cS})}\big(\|\widetilde{\om}\|_{H^{s-\f12}(\overline{\cS})}
+\|\na_{x,z}\widetilde{v}\|_{H^{s-\f12}(\overline{\cS})}\big)\\
&\quad+K_\eta\|\na_{x,z}\widetilde{v}\|_{L^\infty(\overline{\cS})}^2\|\eta\|_{H^s}.
\end{align*}
Thanks to $\overline{v}^3=0$ on $z=0$ and $z=-1$, we have
\begin{align*}
\big\langle\na_{x,z}^k(\overline{v}\cdot \nabla_{x,z}\widetilde{\om}), \na_{x,z}^k\widetilde{\om}\big\rangle
=\big\langle[\na_{x,z}^k, \overline{v}]\cdot \nabla_{x,z}\widetilde{\om}, \na_{x,z}^k\widetilde{\om}\big\rangle
-\big\langle\textrm{div}_{x,z}\overline{v}\nabla_{x,z}^k\widetilde{\om}, \na_{x,z}^k\widetilde{\om}\big\rangle.
\end{align*}
Then we deduece from Lemma \ref{lem:comm-full} and Lemma \ref{lem:v-trans} that
\begin{align*}
\big\langle\na_{x,z}^k(\overline{v}\cdot \nabla_{x,z}\widetilde{\om}), \na_{x,z}^k\widetilde{\om}\big\rangle
\le& C\|\na_{x,z}\widetilde{v}\|_{L^\infty(\overline{\cS})}\|\widetilde{\om}\|_{H^{s-\f12}(\overline{\cS})}\|\na_{x,z}\overline{v}\|_{H^{s-\f12}(\overline{\cS})}\\
&+C\|\na_{x,z}\overline{v}\|_{L^\infty(\overline{\cS})}\|\widetilde{\om}\|_{H^{s-\f12}(\overline{\cS})}^2\\
\le& K_\eta A(t)\cB(t)\big(\|\widetilde{\om}\|_{H^{s-\f12}(\overline{\cS})}^2+\|\eta\|_{H^{s}}^2+\|\widetilde{v}\|_{H^{s+\f12}(\overline{\cS})}^2\big).
\end{align*}

Summing up, we conclude the proposition.
\ef

\subsection{Nonlinear estimates}
Recall that the nonlinear term $f$ is given by
\beno
f=D_th_1+(T_{a\lambda}-T_aT_\lambda)U+[T_a, D_t]\zeta+D_t[D_t, T_\zeta]B-T_ah_2,
\eeno
where
\beno
&& h_1=(T_V-V)\cdot\nabla V-R(a,\zeta)+T_\zeta(T_V-V)\cdot\nabla B,\\
&&h_2=(T_V-V)\cdot\na \zeta+[T_\zeta, T_\lambda]B+(\zeta-T_\zeta)T_\lambda B+R(\eta)V+\zeta R(\eta)B.
\eeno

\begin{lemma}\label{lem:h1}
It holds that
\beno
\|h_1\|_{H^{s-\f12}}\le K_\eta A(t)^2\big(\|V\|_{H^{s-\f12}}+\|\eta\|_{H^{s}}\big).
\eeno
\end{lemma}

\no{\bf Proof.}\,It follows from Lemma \ref{lem:remaider} that
\begin{align*}
\|h_1\|_{H^{s-\f12}}\le C\|\na V\|_{L^\infty}\|V\|_{H^{s-\f12}}+C\|a\|_{C^\f12}\|\zeta\|_{H^{s-1}}+K_\eta\|\na B\|_{L^\infty}\|V\|_{H^{s-\f12}},
\end{align*}
which along with Lemma \ref{lem:a-holder} gives the lemma.\ef

\begin{lemma}\label{lem:h1-Dt}
It holds that
\beno
\|D_th_1\|_{H^{s-\f12}}\le K_\eta \cB(t)A(t)^2\big(\|(V,B)\|_{H^s}+\|\na_{x,z}\widetilde{v}\|_{L^2_z(I;H^{s-\f12})}+\|\eta\|_{H^{s+\f12}}\big).
\eeno
\end{lemma}

\no{\bf Proof.}\,
First, we consider the second term of the $h_1$.
We denote $\overline{\pa}_t=\pa_t+V\cdot\na$. By Lemma \ref{lem:remainder-Dt}, we get
\begin{align*}
\|&\overline{\pa}_tR(a-1, \zeta)\|_{H^{s-\f12}}\le C\big(\|\overline{\pa}_ta\|_{L^\infty}\|\eta\|_{H^{s+\f12}}+\|\overline{\pa}_t\zeta\|_{L^\infty}\|a-1\|_{H^{s-\f12}}\\
&+\|\na V\|_{L^\infty}\big(\|a\|_{L^\infty}\|\eta\|_{H^{s+\f12}}+\|\zeta\|_{L^\infty}\|a-1\|_{H^{s-\f12}}\big)
+\|V\|_{H^s}\|a\|_{L^\infty}\|\zeta\|_{C^{\f 12+\e}}\big).
\end{align*}
Using the fact that
\beno
\overline{\pa}_ta=D_ta+(V-T_V)\cdot\na a,
\eeno
we get by Lemma \ref{lem:remaider} and Lemma \ref{lem:a-holder} that
\begin{align*}
\|\overline{\pa}_ta\|_{L^\infty} \le \|D_ta\|_{L^\infty}+\|V\|_{W^{1,\infty}}\|a\|_{L^\infty}\le K_\eta \cB(t)A(t)^2.
\end{align*}
Note that $\overline{\pa}_t\zeta=\na B-\na V\cdot\na\eta$, hence,
\beno
\|\overline{\pa}_t\zeta\|_{L^\infty}\le K_\eta A(t).
\eeno
Then by Lemma \ref{lem:a-holder}, we obtain
\begin{align*}
\|\overline{\pa}_tR(a-1, \zeta)\|_{H^{s-\f12}}\le K_\eta \cB(t)A(t)^2\big(\|\na_{x,z}\widetilde{v}\|_{L^2_z(I;H^{s-\f12})}+\|V\|_{H^s}+\|\eta\|_{H^{s+\f12}}\big).
\end{align*}
which implies that
\begin{align*}
\|{D}_tR(a-1, \zeta)\|_{H^{s-\f12}}\le& \|\overline{\pa}_tR(a-1, \zeta)\|_{H^{s-\f12}}+C\|R(a-1,\zeta)\|_{C^\f12}\|V\|_{H^s}\\
\le &K_\eta \cB(t)A(t)^2\big(\|\na_{x,z}\widetilde{v}\|_{L^2_z(I;H^{s-\f12})}+\|V\|_{H^s}+\|\eta\|_{H^{s+\f12}}\big).
\end{align*}

\medskip
First, we consider the first and third terms of the $h_1$. We write
\beno
D_t(T_V-V)\cdot \na V=-[D_t,T_{\na V}]V-T_{\na V}D_tV-D_tR(\na V, V).
\eeno
By Remark \ref{rem:comm}, we have
\begin{align*}
\|[D_t,T_{\na V}]V\|_{H^{s-\f12}}\le C\|V\|_{B^1_{\infty,1}}\|V\|_{H^{s-\f12}}+C\|D_t\na V\|_{C^{-\f12}}\|V\|_{H^s}.
\end{align*}
While, using the equation
\beno
D_t\na V=\na D_tV-T_{\na V}\cdot\na V=\na(a\zeta)+\na(T_V-V)\cdot\na V-T_{\na V}\cdot\na V,
\eeno
we deduce from Lemma \ref{lem:remaider} and Lemma \ref{lem:a-holder} that
\beno
\|D_t\na V\|_{C^{-\f12}}\le K_\eta\|a\|_{C^\f12}+C\|V\|_{W^{1,\infty}}^2\le K_\eta A(t)^2.
\eeno
This shows that
\beno
\|[D_t,T_{\na V}]V\|_{H^{s-\f12}}\le C\|V\|_{B^1_{\infty,1}}\|V\|_{H^{s-\f12}}+K_\eta A(t)^2\|V\|_{H^s}.
\eeno
We get by Lemma \ref{lem:remaider} that
\beno
\|T_{\na V}D_tV\|_{H^{s-\f12}}\le C\|V\|_{W^{1,\infty}}\|D_tV\|_{H^{s-\f12}}.
\eeno
It follows from Lemma \ref{lem:remainder-Dt} with $u=V$ and $v=\na V$ that
\begin{align*}
\|\overline{\pa}_tR(\na V, V)\|_{H^{s-\f12}}\le& C\big(\|\overline{\pa}_t\na V\|_{C^{-\f12}}+\|\overline{\pa}_tV\|_{C^\f12}+\|V\|_{W^{1,\infty}}^2\big)\|V\|_{H^s}\\
\le & K_\eta A(t)^2\|V\|_{H^s},
\end{align*}
which ensures that
\begin{align*}
\|{D}_tR(\na V, V)\|_{H^{s-\f12}}&\le \|\overline{\pa}_tR(\na V, V)\|_{H^{s-\f12}}+\|R(\na V, V)\|_{C^\f12}\|V\|_{H^s}\\
&\le K_\eta A(t)^2\|V\|_{H^s}.
\end{align*}
Thus, we obtain
\beno
\|D_t(T_V-V)\cdot\na V\|_{H^{s-\f12}}\le K_\eta\cB(t)A(t)^2\|V\|_{H^s}.
\eeno

We write
\beno
D_tT_\zeta(T_V-V)\cdot\nabla B=[D_t,T_\zeta](T_V-V)\cdot\nabla B+T_\zeta D_t(T_V-V)\cdot\nabla B
\eeno
In a similar way as the above, we can deduce that
\beno
\|D_tT_\zeta(T_V-V)\cdot\nabla B\|_{H^{s-\f12}}\le K_\eta\cB(t)A(t)^2\|(V,B)\|_{H^s}.
\eeno

Putting the above estimates together gives the lemma.\ef

\begin{lemma}\label{lem:h2}
It holds that
\beno
\|h_2\|_{H^{s-\f12}}\le K_\eta\cB(t)A(t)\big(\|(V,B)\|_{H^s}+\|\eta\|_{H^{s+\f12}}\big).
\eeno
\end{lemma}

\no{\bf Proof.}\,By Lemma \ref{lem:remaider}, we get
\beno
\|(T_V-V)\cdot\nabla \zeta\|_{H^{s-\f12}}\le C\|\zeta\|_{C^\f 12}\|V\|_{H^{s}}\le K_\eta\|V\|_{H^{s}}.
\eeno
By Proposition \ref{prop:symbolic calculus} and Lemma \ref{lem:symbol}, we have
\beno
\|[T_\zeta, T_\lambda]B\|_{H^{s-\f12}}\le K_\eta\|B\|_{H^s}.
\eeno
It follows from  Proposition \ref{prop:remainder} that
\beno
\|R(\eta)V\|_{H^{s-\f12}}\le  K_\eta\big(\|V\|_{H^s}+(\|\widetilde{v}\|_{L^\infty(I;B^{1}_{\infty,1})}+A(t))\|\eta\|_{H^{s+\f12}}\big).
\eeno
By Lemma \ref{lem:product} and Proposition \ref{prop:remainder}, we get
\begin{align*}
\|\zeta R(\eta)B\|_{H^{s-\f12}}\le& C\|R(\eta)B\|_{L^\infty}\|\eta\|_{H^{s+\f12}}+C\|\zeta\|_{L^\infty}\|R(\eta)B\|_{H^{s-\f12}}\\
\le &K_\eta \big(1+\|\widetilde{v}\|_{L^\infty(I;B^{1}_{\infty,1})}+A(t)\big)\big(\|B\|_{H^s}+\|\eta\|_{H^{s+\f12}}\big).
\end{align*}
By Lemma \ref{lem:remaider} and Proposition \ref{prop:Sym-Besov}, we have
\begin{align*}
\|(\zeta-T_\zeta)T_\lambda B\|_{H^{s-\f12}}\le C\|T_\lambda B\|_{L^\infty}\|\eta\|_{H^{s+\f12}}\le K_\eta\|B\|_{B^1_{\infty,1}}\|\eta\|_{H^{s+\f12}}.
\end{align*}
This completes the proof of the lemma.\ef

\begin{lemma}\label{lem:f1-est}
Let $f_1=f-D_t[D_t,T_\zeta]B$. It holds that
\beno
\|f_1\|_{H^{s-\f12}}\le K_\eta\cB(t)A(t)^3\big(\|(V,B)\|_{H^s}+\|\na_{x,z}\widetilde{v}\|_{L^2_z(I;H^{s-\f12})}+\|\eta\|_{H^{s+\f12}} \big).
\eeno
\end{lemma}

\no{\bf Proof.}\,By Proposition \ref{prop:symbolic calculus}, Lemma \ref{lem:a-holder} and Lemma \ref{lem:symbol}, we get
\begin{align*}
\|(T_{a\lambda}-T_aT_\lambda)U\|_{H^{s-\f12}}\le& CM^1_\f12(\lambda)M^0_\f12(a)\|U\|_{H^s}\\
\le& K_\eta A(t)^2\|U\|_{H^s}\leq K_\eta A(t)^2\|(V,B)\|_{H^s} .
\end{align*}
It follows from Proposition \ref{prop:commutator-tame} and Lemma \ref{lem:a-holder} that
\begin{align*}
\|[T_a, D_t]\zeta\|_{H^{s-\f12}}\le& C\big(M_0^0(a)\|V\|_{B^1_{\infty,1}}+M_0^0(D_ta)\big)\|\eta\|_{H^{s+\f12}}\\
\le& K_\eta\cB(t)A(t)^2\|\eta\|_{H^{s+\f12}}.
\end{align*}
By Lemma \ref{lem:remaider}, Lemma \ref{lem:a-holder} and Lemma \ref{lem:h2}, we get
\begin{align*}
\|T_ah_2\|_{H^{s-\f12}}\le& C\|a\|_{L^\infty}\|h_2\|_{H^{s-\f12}}\\
\le& K_\eta\cB(t)A(t)^3\big(\|(V,B)\|_{H^s}+\|\eta\|_{H^{s+\f12}}\big).
\end{align*}
Finally, the estimate for $D_th_1$ follows from Lemma \ref{lem:h2}.\ef\vspace{0.2cm}

Next we estimate $f_\om=-T_aR_\om$, where $R_\om$ is defined by
\begin{align*}
R_\om^i=&\big(\pa_y v_{\om}^i-\pa_{x_j}v_{\om}^i\cdot \pa_{x_j}\eta\big)+\pa_{x_i}\eta\big(\pa_y v_{\om}^{d+1}-\pa_{x_j}\eta\pa_{x_j}v_{\om}^{d+1}\big)\\
&+\big(\om_{i,d+1}-\pa_{x_j}\eta\om_{ij}+\pa_{x_i}\eta\pa_{x_j}\eta\om_{j,d+1}\big)\big|_{y=\eta}.
\end{align*}

\begin{lemma}\label{lem:fom}
It holds that
\beno
\|f_\om\|_{H^{s-1}}\le K_\eta A(t)^3\big(\|V_b\|_{H^\f12}+\|\widetilde{\om}\|_{H^{s-\f12}(\overline{\cS})}+\|\eta\|_{H^{s}}\big).
\eeno
\end{lemma}

\no{\bf Proof.}\,Let $I_0=[-\f12,0]$. By Lemma \ref{lem:product} and Lemma \ref{lem:remaider}, it is easy to see that
\begin{align*}
\|R_\om\|_{H^{s-1}}\le& K_\eta\big(\|\na_{x,z}\widetilde{v_\om}\|_{L^\infty(I_0;H^{s-1})}+\|\widetilde{\om}\|_{L^\infty(I_0;H^{s-1})}\big)\\
&+C\big(\|\na_{x,z}\widetilde{v_\om}\|_{L^\infty_z(I_0; C^0)}+\|\om\|_{L^\infty(\Om_t)}\big)\|\eta\|_{H^{s}},
\end{align*}
from which and Lemma \ref{lem:vr}, we deduce that
\begin{align}\label{eq:Rom-Hs}
\|R_\om\|_{H^{s-1}}\le& K_\eta A(t)\big(\|V_b\|_{H^\f12}+\|\widetilde{\om}\|_{H^{s-\f12}(\overline{\cS})}+\|\eta\|_{H^{s}}\big).
\end{align}
Then the lemma follows from Lemma \ref{lem:remaider} and Lemma \ref{lem:a-holder}.
\ef

\begin{lemma}\label{lem:fom-Dt}
It holds that
\beno
\|D_t f_\om\|_{H^{s-1}}\le K_\eta\cB(t)A(t)^4\big(\|(V,B,V_b)\|_{H^{s}}+\|\eta\|_{H^{s+\f12}}+\|\widetilde{\om}\|_{H^{s-\f12}(\overline{\cS})}\big).
\eeno
\end{lemma}

\no{\bf Proof.}\,Let $\overline{\pa}_t=\pa_t+V\cdot\na$. Then $\overline{\pa}_t R_\om$ could be written as
\begin{align*}
\overline{\pa}_tR_\om=&g_1(\na\eta)(\pa_t+v\cdot\na_{x,y})\na_{x,y}v_\om+g_2(\na\eta)\overline{\pa}_t\na\eta\na_{x,z}v_\om\\
&+g_3(\na\eta)(\pa_t+v\cdot\na_{x,y})\om+g_4(\na\eta)\overline{\pa}_t\na\eta\om\big|_{y=\eta},
\end{align*}
for some smooth functions $g_i(i=1,2,3,4)$. Then it follows from Lemma \ref{lem:product} that
\begin{align*}
\|\overline{\pa}_tR_\om&\|_{H^{s-1}}\le K_\eta\big(\|(\pa_t+v\cdot\na_{x,y})\na_{x,y}v_\om|_{y=\eta}\|_{H^{s-1}}+\|(\pa_t+v\cdot\na_{x,y})\om|_{y=\eta}\|_{H^{s-1}}\big)\\
&+K_\eta\big(\|(\pa_t+v\cdot\na_{x,y})\na_{x,y}v_\om|_{y=\eta}\|_{C^0}+\|(\pa_t+v\cdot\na_{x,y})\om\|_{L^\infty(\Om_t)}\big)\|\eta\|_{H^{s+\f12}}\\
&+K_\eta\|\overline{\pa}_t\na \eta\|_{L^\infty}\big(\|\na_{x,y}v_\om|_{y=\eta}\|_{H^{s-1}}+\|\om|_{y=\eta}\|_{H^{s-1}}\big)\\
&+K_\eta\big(\|\na_{x,y}v_\om|_{y=\eta}\|_{L^\infty}+\|\om\|_{L^\infty(\Om_t)}\big)\|\overline{\pa}_t\na\eta\|_{H^{s-1}}\\
&+K_\eta\big(\|\na_{x,y}v_\om|_{y=\eta}\|_{L^\infty}+\|\om\|_{L^\infty(\Om_t)}\big)\|\overline{\pa}_t\na\eta\|_{L^\infty}\|\eta\|_{H^{s+\f12}}.
\end{align*}

Let $I_0=[-\f12,0]$. By Lemma \ref{lem:product} again, we get
\begin{align*}
\|(\pa_t+v&\cdot\na_{x,y})\na_{x,y}v_\om|_{y=\eta}\|_{H^{s-1}}
\le K_\eta\Big(\|\na_{x,z}\widetilde{\dot v_\om}\|_{X^{s-1}(I_0)}+\|\na_{x,z}\widetilde{v}_\om\|_{L^\infty_z(I_0;L^\infty)}\|\na_{x,z}\widetilde{v}\|_{X^{s-1}(I_0)}
\\&+\|\na_{x,z}\widetilde{v}_\om\|_{_{X^{s-1}(I_0)}}\|\na_{x,y}v\|_{L^\infty(\Om_t)}
+\|\na_{x,z}\widetilde{v}_\om\|_{L^\infty_z(I_0;L^\infty)}\|\na_{x,y}v\|_{L^\infty(\Om_t)}\|\eta\|_{H^s}\Big).
\end{align*}
Here $\dot v_\om=(\pa_t+v\cdot\na_{x,y})v_\om$.
Using the equation $(\pa_t+v\cdot\na_{x,y})\om=\om\cdot\na_{x,y}v$ and Lemma \ref{lem:product}, we infer that
\begin{align*}
\|(\pa_t+v\cdot\na_{x,y})\om|_{y=\eta}\|_{H^{s-\f12}}\le K_\eta\big(&\|\na_{x,y}v\|_{L^\infty(\Om_t)}\|\na_{x,z}\widetilde{v}\|_{X^{s-1}(I_0)}\\
&\quad+\|\na_{x,y}v\|_{L^\infty(\Om_t)}^2\|\eta\|_{H^s}\big).
\end{align*}
Using the equation $\overline{\pa}_t\na\eta=\na B-\na V\cdot\na \eta$ and Lemma \ref{lem:product}, we get
\begin{align*}
\|\overline{\pa}_t\na\eta\|_{H^{s-1}}\le K_\eta\big(\|(V,B)\|_{H^s}+\|\na V\|_{L^\infty}\|\eta\|_{H^s}\big).
\end{align*}
On the other hand, we have
\begin{align*}
&\|(\pa_t+v\cdot\na_{x,y})\na_{x,y}v_\om|_{y=\eta}\|_{C^0}\\
&\qquad\le K_\eta\big(\|\na_{x,z}\widetilde{\dot v_\om}\|_{L^\infty_z(I_0;C^0)}+\|\na_{x,z}\widetilde{v_\om}\|_{L^\infty_z(I_0;L^\infty)}\|\na_{x,y}v\|_{L^\infty(\Om_t)}\big),\\
&\|(\pa_t+v\cdot\na_{x,y})\om|_{y=\eta}\|_{L^\infty}\le K_\eta\|\na_{x,y}v\|_{L^\infty(\Om_t)}^2,\\
&\|\overline{\pa}_t\na\eta|_{y=\eta}\|_{L^\infty}\le K_\eta\|(\na B,\na V)\|_{L^\infty}.
\end{align*}

Summing up the above estimates, we apply Lemma \ref{lem:vr} and Lemma \ref{lem:vr-Dt} to obtain
\begin{align*}
\|\overline{\pa}_tR_\om\|_{H^{s-1}}\le K_\eta\cB(t)A(t)^2\big(\|(V,B,V_b)\|_{H^{s}}+\|\eta\|_{H^{s+\f12}}+\|\widetilde{\om}\|_{H^{s-\f12}(\overline{\cS})}\big),
\end{align*}
which implies that
\begin{align*}
\|{D}_tR_\om\|_{H^{s-1}}\le& \|\overline{\pa}_tR_\om\|_{H^{s-1}}+C\|R_\om\|_{L^\infty}\|V\|_{H^s}\\
\le& K_\eta\cB(t)A(t)^2\big(\|(V,B,V_b)\|_{H^{s}}+\|\eta\|_{H^{s+\f12}}+\|\widetilde{\om}\|_{H^{s-\f12}(\overline{\cS})}\big).
\end{align*}
which along with Proposition \ref{prop:commutator-tame}, Lemma \ref{lem:a-holder} and (\ref{eq:Rom-Hs}) gives
\begin{align*}
\|{D}_tf_\om\|_{H^{s-1}}\le& \|[D_t,T_a]R_\om\|_{H^{s-1}}+\|T_aD_tR_\om\|_{H^{s-1}}\\
\le& K_\eta\cB(t)A(t)^4\big(\|(V,B,V_b)\|_{H^{s}}+\|\eta\|_{H^{s+\f12}}+\|\widetilde{\om}\|_{H^{s-\f12}(\overline{\cS})}\big).
\end{align*}
The proof is finished.\ef

\begin{lemma}\label{lem:fom-inn}
It holds that
\begin{align*}
\big\langle(f_\om)_{s-1/2}, (D_t U)_{s-1/2}&\big\rangle\le \f d {dt}\big\langle (f_\om)_{s-1}, U_{s}\big\rangle\\
 +&K_\eta\cB(t)A(t)^4\big(\|(V,B,V_b)\|_{H^{s}}+\|\eta\|_{H^{s+\f12}}+\|\widetilde{\om}\|_{H^{s-\f12}(\overline{\cS})}\big)\|U\|_{H^s}.
\end{align*}
\end{lemma}

\no{\bf Proof.}\,A direct calculation yields
\begin{align*}
\big\langle(f_\om)_{s-1/2}, (D_t U)_{s-1/2}\big\rangle
=&-\big\langle(f_\om)_{s-1/2}, [T_V\cdot\na, \langle D\rangle^{s-1/2}]U\big\rangle+\big\langle(f_\om)_{s-1/2}, D_t U_{s-1/2}\big\rangle\\
=&-\big\langle(f_\om)_{s-1/2}, [T_V\cdot\na, \langle D\rangle^{s-1/2}]U\big\rangle+\f d {dt}\big\langle (f_\om)_{s-1}, U_{s}\big\rangle\\
&-\big\langle [T_V\cdot\na, \langle D\rangle^{s-1/2}]f_\om, U_{s-1/2}\big\rangle-\big\langle (D_tf_\om)_{s-1}, U_{s}\big\rangle\\
&+\big\langle((T_V\cdot\na)^*+T_V\cdot\na)(f_\om)_{s-1/2}, U_{s-1/2}\big\rangle.
\end{align*}
By Lemma \ref{lem:commu-Ds} and Proposition \ref{prop:symbolic calculus}, we have
\beno
&&\|[T_V\cdot\na, \langle D\rangle^{s-\f12}]U\|_{H^{\f12}}\le C\|V\|_{W^{1,\infty}}\|U\|_{H^s},\\
&&\|[T_V\cdot\na, \langle D\rangle^{s-\f12}]f_\om\|_{H^{-\f12}}\le C\|V\|_{W^{1,\infty}}\|f_\om\|_{H^{s-1}},\\
&&\|(T_V\cdot\na)^*+T_V\cdot\na)(f_\om)_{s-1/2}\|_{H^{-\f12}}\le C\|V\|_{W^{1,\infty}}\|f_\om\|_{H^{s-1}}.
\eeno
Then the lemma follows from Lemma \ref{lem:fom} and Lemma \ref{lem:fom-Dt}.\ef

\begin{lemma}\label{lem:f1-inn}
It holds that
\begin{align*}
\big\langle &(D_t[D_t,T_\zeta]B)_{s-\f12},(D_t U)_{s-1/2}\big\rangle
=\f d {dt}\Big(\big\langle g_{s-1/2}, h_{s-1/2}\big\rangle+\f12\big\langle g_{s-1/2}, g_{s-1/2}\big\rangle\Big)\\
&\quad+K_\eta \cB(t)A(t)^3\big(\|D_tU\|_{H^{s-\f12}}+\|(V,B)\|_{H^s}+\|\na_{x,z}\widetilde{v}\|_{L^2_z(I;H^{s-\f12})}+\|\eta\|_{H^{s+\f12}}\big)^2.
\end{align*}
where $g=[D_t,T_\zeta]B$ and $h=D_tU-[D_t,T_\zeta]B$.
\end{lemma}

\no{\bf Proof.}\,A direct calculation yields
\begin{align*}
&\big\langle (D_tg)_{s-1/2}, (D_t U)_{s-1/2}\big\rangle=\big\langle (D_tg)_{s-1/2}, h_{s-1/2}\big\rangle+\big\langle (D_tg)_{s-1/2}, g_{s-1/2}\big\rangle\\
&=\f d {dt}\Big(\big\langle g_{s-1/2}, h_{s-1/2}\big\rangle+\f12\big\langle g_{s-1/2}, g_{s-1/2}\big\rangle\Big)\\
&\quad-\big\langle [T_V\cdot\na, \langle D\rangle^{s-1/2}]g, h_{s-1/2}\big\rangle+\big\langle g_{s-1/2}, (D_t)^*h_{s-1/2}\big\rangle\\
&\quad-\big\langle [T_V\cdot\na, \langle D\rangle^{s-1/2}]g, g_{s-1/2}\big\rangle+\f12\big\langle(T_V\cdot\na)^*+T_V\cdot\na)g_{s-1/2}, g_{s-1/2}\big\rangle\\
&\ge \f d {dt}\Big(\big\langle g_{s-1/2}, h_{s-1/2}\big\rangle+\f12\big\langle g_{s-1/2}, g_{s-1/2}\big\rangle\Big)\\
&\quad-\|[T_V\cdot\na, \langle D\rangle^{s-1/2}]g\|_{L^2}\big(\|h\|_{H^{s-\f12}}+\|g\|_{H^{s-\f12}}\big)\\
&\quad-\|(T_V\cdot\na)^*+T_V\cdot\na)g_{s-1/2}\|_{L^2}\|g\|_{H^{s-\f12}}-\|g\|_{H^s}\|(D_t)^*h_{s-1/2}\|_{H^{-\f12}}.
\end{align*}

It follows from Proposition \ref{prop:commutator-tame} that
\begin{align*}
&\|g\|_{H^{s-\f12}}\le C\big(\|V\|_{W^{1,\infty}}+\|D_t\zeta\|_{L^\infty}\big)\|B\|_{H^{s-\f12}}\le K_\eta A(t)\|B\|_{H^{s-\f12}},\\
&\|g\|_{H^{s}}\le C\big(\|V\|_{W^{1,\infty}}+\|D_t\zeta\|_{L^\infty}\big)\|B\|_{H^{s}}\le K_\eta A(t)\|B\|_{H^{s}},
\end{align*}
from which, Proposition \ref{prop:symbolic calculus} and Lemma \ref{lem:commu-Ds}, we infer
\beno
&&\|[T_V\cdot\na, \langle D\rangle^{s-1/2}]g\|_{L^2}\le K_\eta A(t)^2\|B\|_{H^{s-\f12}},\\
&&\|((T_V\cdot\na)^*+T_V\cdot\na)g_{s-1/2}\|_{L^2}\le K_\eta A(t)^2\|B\|_{H^{s-\f12}},\\
&&\|h\|_{H^{s-\f12}}\le \|D_tU\|_{H^{s-\f12}}+K_\eta A(t)\|B\|_{H^{s-\f12}}.
\eeno
We have by Lemma \ref{lem:commu-Ds} and Lemma \ref{lem:remaider} that
\begin{align*}
\|(D_t)^*h_{s-1/2}\|_{H^{-\f12}}\le &\|D_th\|_{H^{s-1}}+\|[T_V\cdot\na,\langle D\rangle^{s-\f12}]h\|_{L^2}+\|T_{\na\cdot V}h_{s-1/2}\|_{L^2}\\
\le& \|D_th\|_{H^{s-1}}+C\|V\|_{W^{1,\infty}}\|h\|_{H^{s-\f12}}\\
\le& \|D_th\|_{H^{s-1}}+K_\eta A(t)^2\big(\|D_tU\|_{H^{s-\f12}}+\|B\|_{H^{s-\f12}}\big).
\end{align*}
While, by the equation (\ref{eq:U}), we find
\beno
D_th=D_th_1-[D_t,T_a]\zeta-T_aD_t\zeta.
\eeno
Then we get by Lemma \ref{lem:h1-Dt}, Proposition \ref{prop:commutator-tame} and Lemma \ref{lem:a-holder} that
\beno
\|D_th\|_{H^{s-1}}\le  K_\eta \cB(t)A(t)^2\big(\|(V,B)\|_{H^s}+\|\na_{x,z}\widetilde{v}\|_{L^2_z(I;H^{s-\f12})}+\|\eta\|_{H^{s+\f12}}\big).
\eeno
This shows that
\begin{align*}
\|(D_t)^*h_{s-1/2}\|_{H^{-\f12}}\le  K_\eta \cB(t)A(t)^2\big(&\|D_tU\|_{H^{s-\f12}}+\|(V,B)\|_{H^s}\\
&+\|\na_{x,z}\widetilde{v}\|_{L^2_z(I;H^{s-\f12})}+\|\eta\|_{H^{s+\f12}}\big).
\end{align*}

Summing up the above estimates, we conclude the lemma.\ef

\subsection{Energy functional}

We introduce an energy functional $\cE_s(t)$ defined by
\begin{align}
\cE_{s}(t)\eqdef& \|v(t)\|_{H^1(\Om(t))}+\|\eta(t)\|_{H^s}+\|(V,B)(t)\|_{H^{s-\f12}}+\|V_b(t)\|_{H^s}\nonumber\\
&+\|\widetilde{\om}(t)\|_{H^{s-\f12}(\overline{\cS})}+\|T_{\sqrt{a\lambda}}U(t)\|_{H^{s-\f12}}+\|D_tU(t)\|_{H^{s-\f12}}.\label{def:energy-t}
\end{align}

\begin{proposition}\label{prop:energy-rela}
It holds that
\beno
\|U\|_{H^s}+\|(V,B)\|_{H^s}+\|\eta\|_{H^{s+\f12}}+\|\widetilde{v}\|_{H^{s+\f12}(\overline{\cS})}\le K_\eta A(t)^3\cE_s(t).
\eeno
\end{proposition}

\no{\bf Proof.}\,It follows from Lemma \ref{lem:v-Hs} that
\beno
\|\widetilde{v}\|_{H^{s+\f12}(\overline{\cS})}\le K_\eta A(t)\big(\|(V,B,V_b)\|_{H^s}+\|\widetilde{\om}\|_{H^{s-\f12}(\overline{\cS})}+\|\eta\|_{H^{s+\f12}}\big).
\eeno
This together with Lemma \ref{lem:eta-E} and Lemma \ref{lem:VB-E} yields the result.\ef

\begin{lemma}\label{lem:eta-E}
It holds that
\beno
\|\eta\|_{H^{s+\f12}}\le K_\eta A(t)^2\cE_s(t).
\eeno
\end{lemma}

\no{\bf Proof.}\,Recall that $D_tU=-T_a\zeta+h_1+[D_t,T_\zeta]B$. Due to $a\ge c_0$, we have
\begin{align*}
\zeta=&T_{a^{-1}}T_a\zeta+(T_{a^{-1}}T_a-1)\zeta\\
=&T_{a^{-1}}\big(-D_tU+h_1+[D_t,T_\zeta]B\big)+(T_{a^{-1}}T_a-1)\zeta
\end{align*}
which along with Proposition \ref{prop:symbolic calculus}, Lemma \ref{lem:h1}, Proposition \ref{prop:commutator-tame} and Lemma \ref{lem:a-holder} yields
\begin{align*}
\|\zeta\|_{H^{s-\f12}}\le K_\eta\|D_tU\|_{H^{s-\f12}}+K_\eta A(t)^2\big(\|\eta\|_{H^{s}}+\|(V,B)\|_{H^{s-\f12}}\big).
\end{align*}
This gives the lemma by recalling $\zeta=\na\eta$.\ef

\begin{lemma}\label{lem:VB-E}
It holds that
\beno
\|U\|_{H^s}+\|(V,B)\|_{H^s}\le K_\eta A(t)^3\cE_s(t).
\eeno
\end{lemma}

\no{\bf Proof.}\,Recall that $U=V+T_\zeta B$. Hence,
\beno
\|V\|_{H^s}\leq \|U\|_{H^s}+\|T_\zeta B\|_{H^{s}}\leq K_\eta\big(\|U\|_{H^s}+\|B\|_{H^{s}}\big).
\eeno
So, it is sufficient to consider $B$. We have
\beno
\dv U=\dv V+T_\zeta\cdot\nabla B+T_{\dv \zeta}B.
\eeno
On the other hand, we have
\begin{align*}
\dv V=&\sum_{i=1}^d\pa_iv^i+\pa_i\eta\pa_yv^i\big|_{y=\eta}\\
=&-\pa_yv^{d+1}+\na\eta\cdot\na v^{d+1}\big|_{y=\eta}+\pa_i\eta\om_{d+1,i}\big|_{y=\eta}\\
=&-G(\eta)B-\pa_yv^{d+1}_\om+\na \eta\cdot\na v^{d+1}_\om\big|_{y=\eta}+\pa_i\eta\om_{d+1,i}\big|_{y=\eta}\\
\triangleq &-G(\eta)B+V_\om.
\end{align*}
Then we deduce that
\begin{align*}
\dv U=&\dv V+T_\zeta\cdot \nabla B+T_{\dv \zeta}B\\
=&-G(\eta)B+V_\om+T_\zeta\cdot \nabla B+T_{\dv \zeta}B\\
=&-T_\lambda B- R(\eta)B+V_\om+T_\zeta\cdot \nabla B+T_{\dv \zeta}B\\
=&-T_q B- R(\eta)B+V_\om+T_{\dv \zeta}B,
\end{align*}
where the symbol $q= \lambda-i\zeta \cdot \xi$. Thus, it follows from Proposition \ref{prop:remainder}, Lemma \ref{lem:vr} and Lemma \ref{lem:eta-E} that
\begin{align*}
\|T_qB\|_{H^{s-1}}\le& \|U\|_{H^s}+\|R(\eta)B\|_{H^{s-1}}+\|V_\om\|_{H^{s-1}}+K_\eta\|B\|_{H^{s-1}}\\
\le& \|U\|_{H^s}+K_\eta A(t)\big(\|(B,V_b)\|_{H^{s-1}}+\|\om\|_{H^{s-\f12}(\overline{\cS})}+\|\eta\|_{H^{s+\f12}}\big)\\
\le& \|U\|_{H^s}+K_\eta A(t)^3\cE(t).
\end{align*}
On the other hand, we get by Proposition \ref{prop:symbolic calculus} and Lemma \ref{lem:a-holder} that
\begin{align*}
\|U\|_{H^s}\le& \|T_{(\sqrt{a\lambda})^{-1}}T_{\sqrt{a\lambda}}U\|_{H^s}+\|(T_{(\sqrt{a\lambda})^{-1}}T_{\sqrt{a\lambda}}-1)U\|_{H^{s}}\\
\le& K_\eta A(t)^2\|T_{\sqrt{a\lambda}}U\|_{H^{s-\f12}}\|+K_\eta A(t)^2\|U\|_{H^{s-\f12}}\\
\le& K_\eta A(t)^2\|T_{\sqrt{a\lambda}}U\|_{H^{s-\f12}}+K_\eta A(t)^2\|(V,B)\|_{H^{s-\f12}},
\end{align*}
from which and Proposition \ref{prop:symbolic calculus}, we infer that
\begin{align*}
\|B\|_{H^s}\le& K_\eta\|T_qB\|_{H^{s-1}}+\|(T_{q^{-1}}T_q-1)B\|_{H^{s}}\\
\le& K_\eta A(t)^3\cE_s(t)+K_\eta\|B\|_{H^{s-\f12}}\le K_\eta A(t)^3\cE_s(t).
\end{align*}
The proof is finished.\ef

\subsection{Proof of Theorem \ref{thm:blow-up}}

We first recover the regularity of the free surface from the mean curvature.

\begin{lemma}\label{lem:eta-Holder}
Assume that the mean curvature $\kappa\in L^2\cap L^p(\R^d)$ for some $p>d$. Then we have
\beno
\|\eta\|_{H^2}+\|\eta\|_{C^{2-\f d p}}\le C\big(\|\eta\|_{L^2}, \|\na\eta\|_{L^\infty}, \|H\|_{L^2\cap L^p}\big).
\eeno
\end{lemma}

\no{\bf Proof.}\,The estimate of $\|\eta\|_{C^{2-\f d p}}$ has been proved in \cite{WZ}.
Let $\eta_\ell=\pa_\ell\eta$ and $a_{ij}=(1+|\na\eta|^2)^{-\f32}\big((1+|\na\eta|^2)\delta_{ij}-\pa_i\eta\pa_j\eta\big)$. A direct calculation gives
\beno
\pa_j\big(a_{ij}\pa_i\eta_\ell\big)=\pa_\ell H.
\eeno
It is easy to verify that the matrix $\big(a_{ij}\big)$ is uniformly elliptic with the elliptic constants depending on $\|\na\eta\|_{L^\infty}$, which implies
that $\|\na \eta_\ell\|_{L^2}\le C(\|\na\eta\|_{L^\infty})\|H\|_{L^2}$.\ef

\begin{lemma}\label{lem:A-eta}
It holds that
\beno
\sup_{t\in [0,T]}\big(A(t)+\|\eta(t)\|_{C^{\f32+\e}}\big)\le C\big(T, M(T),\|v_0\|_{H^1(\Om_0)},\|\eta_0\|_{H^s}\big),
\eeno
where $C$ is an increasing function depending on $h_0$.
\end{lemma}

\no{\bf Proof.}\,Recall that $\zeta=\na\eta$ satisfies
\beno
\pa_t\zeta+V\cdot\na \zeta=\na B+\na V\cdot\zeta,
\eeno
which implies
\beno
\sup_{t\in [0,T]}\|\na\eta(t)\|_{L^\infty}\le C\big(T, M(T), \|\na\eta_0\|_{L^\infty}\big),
\eeno
which along with Lemma \ref{lem:basic energy} and Lemma \ref{lem:eta-Holder} implies that
\ben\label{eq:break-est1}
\sup_{t\in [0,T]}\big(\|\eta(t)\|_{H^2}+\|\eta(t)\|_{C^{\f 32+\e}}\big)\le C\big(T, M(T),\|v_0\|_{L^2(\Om_0)},\|\eta_0\|_{H^s}\big).
\een
By Lemma \ref{lem:basic energy} and Lemma \ref{lem:velocity-H1}, we have
\begin{align*}
\f d {dt}\big(\|v\|_{H^1(\Om_t)}^2+\|\eta\|_{L^2}^2\big)\le K_\eta\big(\|v(t)\|_{W^{1,\infty}(\Om_t)}\|v(t)\|_{H^1(\Om_t)}+\|\eta(t)\|_{H^\f32}\big),
\end{align*}
from which and (\ref{eq:break-est1}), we deduce the lemma. \ef

\begin{lemma}\label{lem:log-sobolev}
Let $s>\f d2+1$ and $f\in H^s(\R^d)$. Then there holds
\beno
\|f\|_{B^1_{\infty,1}}\le C\big(1+\|f\|_{W^{1,\infty}}\big)\ln(e+\|f\|_{H^s}).
\eeno
\end{lemma}

\no{\bf Proof.}\,Given an integer $N$, we get by Lemma \ref{lem:Berstein} that
\begin{align*}
\|f\|_{B^{1}_{\infty,1}}=&\sum_{j\ge -1}^N2^j\|\Delta_j f\|_{L^\infty}+\sum_{j>N}2^j\|\Delta_j f\|_{L^\infty}\\
\le& C(N+1)\|f\|_{W^{1,\infty}}+C\sum_{j>N}2^{(1+\f d2)j}\|\Delta_j f\|_{L^2}\\
\le& C(N+1)\|f\|_{W^{1,\infty}}+C2^{-N(s-1-\f d2)}\|f\|_{H^s}.
\end{align*}
Take $N$ so that $2^{-N(s-1-\f d2)}\|f\|_{H^s}\sim 1$(i.e., $N\sim \ln(e+\|f\|_{H^s})$). Then the lemma follows easily.
\ef

\medskip

Now we are in position to prove Theorem \ref{thm:blow-up}.
We denote
\beno
\cP_T\triangleq\cP\big(T, M(T), \|\na\eta_0\|_{L^\infty}\big)
\eeno
for some increasing function $\cP$ depending on $c_0, h_0$, which may change from line to line. By Lemma \ref{lem:A-eta}, $K_\eta$ and $A(t)$
is bounded by $\cP_T$. By Proposition \ref{prop:energy-rela}, we have
\beno
\|U\|_{H^s}+\|(V,B)\|_{H^s}+\|\eta\|_{H^{s+\f12}}+\|\widetilde{v}\|_{H^{s+\f12}(\overline{\cS})}\le P_T\cE_s(t).
\eeno

We first deduce from Proposition \ref{prop:energy-boundary} that
\begin{align*}
&\f d {dt}\big(\|D_tU\|_{H^{s-\f12}}^2+\|T_{\sqrt{a\lambda}}U\|_{H^{s-\f12}}^2\big)
\le \cP_T\cB(t)\cE_s(t)^2+2\langle (f+f_\om)_{s-1/2}, (D_t U)_{s-1/2}\big\rangle,
\end{align*}
which along with Lemma \ref{lem:f1-est}, Lemma \ref{lem:fom-inn} and Lemma \ref{lem:f1-inn} gives
\begin{align*}
&\f d {dt}\Big(\|D_tU\|_{H^{s-\f12}}^2+\|T_{\sqrt{a\lambda}}U\|_{H^{s-\f12}}^2+\big\langle (f_\om)_{s-1}, U_{s}\big\rangle\\
&\qquad+\big\langle g_{s-1/2}, h_{s-1/2}\big\rangle+\f12\big\langle g_{s-1/2}, g_{s-1/2}\big\rangle\Big)\leq
\cP_T\cB(t)\cE_s(t)^2,
\end{align*}
where $g=[D_t,T_\zeta]B$ and $h=D_tU-[D_t,T_\zeta]B$. By the proof of Lemma \ref{lem:f1-inn} and Lemma \ref{lem:fom}, we know
\begin{align*}
&\big\langle (f_\om)_{s-1}, U_{s}\big\rangle\le K_\eta A(t)^3\big(\|V_b\|_{H^{\f12}}+\|\widetilde{\om}\|_{H^{s-\f12}}+\|\eta\|_{H^s}\big)\|U\|_{H^s},\\
&\big\langle g_{s-1/2}, h_{s-1/2}\big\rangle\le K_\eta A(t)\|B\|_{H^{s-\f12}}\|D_tU\|_{H^{s-\f12}}+K_\eta A(t)^3\|B\|_{H^{s-\f12}}^2,\\
&\big\langle g_{s-1/2}, g_{s-1/2}\big\rangle\le K_\eta A(t)^2\|B\|_{H^{s-\f12}}^2.
\end{align*}
This shows that
\begin{align*}
&\|D_tU(t)\|_{H^{s-\f12}}+\|T_{\sqrt{a\lambda}}U(t)\|_{H^{s-\f12}} \le \cP_T\cE_s(0)+\cP_0\int_0^t\cB(t')\cE_s(t')dt'\\
&+\cP_T\big(\|(V,B)(t)\|_{H^{s-\f12}}+\|V_b(t)\|_{H^s}+\|\widetilde{\om}(t)\|_{H^{s-\f12}(\overline{\cS})}+\|\eta(t)\|_{H^s}\big),
\end{align*}
which together with Proposition \ref{prop:vorticity}, Proposition \ref{prop:energy-lower} and Lemma \ref{lem:log-sobolev} gives
\begin{align*}
\cE_s(t)\le& \cP_0\cE_s(0)+\cP_T\int_0^t \cB(t')\cE_s(t')dt'\\
\le &\cP_T\cE_s(0)+\cP_T\int_0^t\ln(e+\cE_s(t'))\cE_s(t')dt'.
\end{align*}
Then Gronwall's inequality ensures that
\beno
\cE_s(t)\le C\big(E_s(0), \cP_T\big).
\eeno
This completes the proof of Theorem \ref{thm:blow-up}.\ef

\section{Iteration scheme and symmetrization}

We begin with the proof of local well-posedness from this section.
We first establish the local well-posedness result for sufficiently smooth data.
In the last section, we extend it to the low regularity data.
Although there has been a lot of work \cite{Lin, SZ3, ZZ} devoted to the local well-posedness of the Euler equations with free surface for
smooth data, they do not work in the Eulerian coordinates and do not consider the case of finite depth to our knowledge.

\begin{theorem}\label{thm:local-smooth}
Let $d\geq 1$ and $s>\f d 2+10$ be an integer.
Assume that the initial data $(\eta_0, v_0)$ satisfies
\begin{eqnarray*}
\eta_0\in H^{s+1/2}(\mathbf{R}^d),\quad v_0\in H^{s}(\Om_0).
\end{eqnarray*}
Furthermore, we assume that there exist two positive constants $c_0>0$ and $h_0>0$ such that
\beno
&&-(\pa_yP)(0,x,\eta_0(x))\ge c_0\quad\textrm{ for }x\in \R^d,\\
&&1+\eta_0(x)\ge h_0\quad \textrm{ for }x\in \R^d.\label{ass:initial-eta}
\eeno
Then there exists $T>0$ such that the system (\ref{eq:euler})--(\ref{eq:euler-p}) with the initial data
$(\eta_0, v_0)$ has a unique solution $(\eta, v)$ satisfying
\begin{gather*}
\eta\in C\big([0,T];H^{s+\f12}(\mathbf{R}^d)\big),\quad v\in C\big([0,T];H^{s}(\Om_t)\big).
\end{gather*}
\end{theorem}

The proof of Theorem \ref{thm:local-smooth} is conducted in the following four section.
In this section, we construct a sequence of approximate solutions by an iteration.
The next section is devoted to the uniform estimates of the approximate solutions.
The last two sections are devoted to show that the approximate sequence is a Cauchy
sequence and converges to the solution of the Euler system (\ref{eq:euler})--(\ref{eq:euler-p}).

\subsection{Iteration scheme}

We construct the approximate solution by an iteration. Assume that the initial data $(\eta_0, v_0)$ is smooth.
Assume that we construct a smooth solution $\big(V^n, B^n, V_b^n,\eta^n, P^n, \om^n\big)$ and
$\big(V^n_1, B^n_1, V_{b,1}^n,\eta^n_1\big)$ in $n$-th iteration. Here $\om^n$
is a function defined in $\Om_t^n=\big\{(x,y): -1<y<\eta_1^n(t,x)\big\}$ and $P^n$ is a function defined in $\big\{(x,y): -1<y<\eta^n(t,x)\big\}$.
We will construct the solution $\big(V^{n+1}, B^{n+1}, V_b^{n+1}, \eta^{n+1}, P^{n+1}, \om^{n+1})$
and $\big(V^{n+1}_1, B^{n+1}_1, V_{b,1}^{n+1}, \eta^{n+1}_1\big)$ in the $(n+1)$-th iteration by the following scheme.

We still denote by $D_t\triangleq\pa_t+T_{V^n}\cdot\na$ for the simplicity of notation. Let
\beno
a^n=-\pa_yP^n|_{y=\eta^n},\quad \lambda^n=\lambda(\eta^n_1),\quad \zeta^n=\na\eta^n,\quad \zeta^n_1=\na\eta^n_1.
\eeno
We first introduce the evolutional system on the trace of velocity and the free surface in the $(n+1)$-th iteration.
\begin{equation}\label{eq:iteration-velocity}
\left\{
\begin{aligned}
&D_tV^{n+1}=-T_{\zeta^{n}}a^{n}-T_{a^{n}}\zeta^{n+1}-R(\zeta^n, a^n)+(T_{V_1^{n}}-V_1^{n})\cdot\nabla V^{n},\\
&D_tB^{n+1}= a^{n}-1+(T_{V_1^{n}}-V_1^{n})\cdot\nabla B^{n},\\
&(\pa_t+V_b^n\cdot\na)V_b^{n+1}=-\na P^n|_{y=-1},\\
&D_t\zeta^{n+1}=T_{\lam^n}(V^{n+1}+T_{\zeta^n}B^{n+1})+(T_{V_1^{n}}-V_1^{n})\cdot\zeta^n+[T_{\zeta^n},T_{\lam^n}]B^{n}_1\\
&\qquad+(\zeta^{n} -T_{\zeta^{n} })T_{\lam^n}B^n+R(\eta^n_1)V^{n}_1+\zeta^n R(\eta^n_1)B^{n}_1+R^n_\om,\\
&\big(V^{n+1}, B^{n+1}, \zeta^{n+1}, V^{n+1}_b\big)\big|_{t=0}=\big(V_0,B_0, \na\eta_0, V_{b,0}\big),
\end{aligned}\right.
\end{equation}
where $R^n_\om$ is given by
\begin{align*}
(R_\om^n)^i=&\big(\pa_y (v_{\om}^n)^i-\pa_{x_j}(v_{\om}^n)^i\cdot \pa_{x_j}\eta^n_1\big)+\pa_{x_i}\eta^n_1\big(\pa_y (v_{\om}^n)^{d+1}
-\pa_{x_j}\eta^n_1\pa_{x_j}(v_{\om}^n)^{d+1}\big)\\
&+\big(\om_{i,d+1}^n-\pa_{x_j}\eta^n_1\om_{ij}^n+\pa_{x_i}\eta^n_1\pa_{x_j}\eta^n_1\om_{j,d+1}^n\big)\big|_{y=\eta^n_1}.
\end{align*}
Here $v_\om^n$ solves the following elliptic equation in $\Om^n_t$:
\beno
\left\{
\begin{array}{ll}
-\Delta_{x, y} v_{\om}^n=\na_{x,y}\times\om^n\quad \textrm{in} \quad \Om_t^n,\\
v_{\om}^n|_{y=\eta^n_1}=0,\quad v_\om^n|_{y=-1}=(V_{b,1}^n, 0).
\end{array}\right.
\eeno
Given $(V^{n+1}, B^{n+1}, V^{n+1}_b)$, we introduce a new boundary velocity $(V^{n+1}_1, B^{n+1}_1, V^{n+1}_{b,1})$ defined by
\begin{align}\label{eq:iteration-velocity-new}
\left\{
\begin{array}{l}
\big(\pa_t+T_{V^{n+1}}\cdot\na\big)\big(V_1^{n+1}, B^{n+1}_1\big)=D_t\big(V^{n+1}, B^{n+1}\big),\\
\big(\pa_t+V_{b}^{n+1}\cdot\na\big)V^{n+1}_{b,1}=-\na P^n|_{y=-1},\\
\big(V^{n+1}_1, B^{n+1}_1, V^{n+1}_{b,1}\big)\big|_{t=0}=\big(V_0,B_0, V_{b,0}\big).
\end{array}\right.
\end{align}
A key property is that $(V^{n+1}_1,B^{n+1}_1)$ has the same regularity as $\big(\pa_t+T_{V^{n+1}}\cdot\na\big)\big(V_1^{n+1}, B^{n+1}_1\big)$ .
While, $\big(\pa_t+T_{V^{n+1}}\cdot\na\big)\big(V^{n+1}, B^{n+1}\big)$ will lose one derivative.

Given $(V^{n+1}, B^{n+1}, \zeta^{n+1})$, we define $\eta^{n+1}$ by
\ben\label{eq:iteration-eta}
-\Delta\eta^{n+1}+\eta^{n+1}=-\dive\zeta^{n+1}+ {\eta}^{n+1}_1,
\een
where ${\eta}^{n+1}_1$ is determined by
\begin{equation}
\left\{
\begin{aligned}
&(\pa_t+V^{n+1}\cdot\nabla) {\eta}_1^{n+1}=B^{n+1}_1\\
& {\eta}^{n+1}_1|_{t=0}=\eta_0.
\end{aligned}\right.
\end{equation}

Given $\eta^{N+1}_1$, let $\Om_t^{n+1}=\big\{(x,y):-1\leq y\leq \eta^{n+1}_1(t,x)\big\}$.
The vorticity $\om^{n+1}$ in $(n+1)$-th iterative is given by solving the nonlinear vorticity equation in the known domain $\Om_t^{n+1}$:
\begin{equation}\label{eq:iterative-vorticity}
\left\{
\begin{aligned}
&\pa_t {\om}^{n+1}+ \big((v^{n+1})^h\cdot\na+(v_1^{n+1})^{d+1}\pa_y\big){\om}^{n+1}= {\om}^{n+1}\cdot\nabla_{x,y}{v}^{n+1}_1,\\
& {\om}^{n+1}\mid_{t=0}=\om_0,
\end{aligned}
\right.
\end{equation}
where the velocity $v^{n+1}$ is given by
\begin{equation}\label{eq:iterative-velocity}
\left\{
\begin{aligned}
&-\Delta { v}^{n+1}=\nabla_{x,y}\times{\om}^{n+1}\quad\textrm{in}\quad\Om^{n+1}_t,\\
& {v}^{n+1}\mid_{y= {\eta}^{n+1}_1}=(V^{n+1},B^{n+1}),\\
&  v^{n+1} |_{y=-1}= (V^{n+1}_{b},0),
\end{aligned}
\right.
\end{equation}
and the velocity $v^{n+1}_1$ is given by
\begin{equation}\label{eq:iterative-velocity-new}
\left\{
\begin{aligned}
&-\Delta { v}_1^{n+1}=\nabla_{x,y}\cdot{\om}^{n+1}\quad\textrm{in}\quad\Om^{n+1}_t,\\
& {v}^{n+1}_1\mid_{y= {\eta}^{n+1}_1}=(V^{n+1}_1,B^{n+1}_1),\\
&  v^{n+1}_1|_{\Gamma_b}= (V^{n+1}_{b,1},0).
\end{aligned}
\right.
\end{equation}

Finally, we need to construct the pressure $P^{n+1}$ in a smoother domain $\widetilde{\Om}^{n+1}_t=\big\{(x,y):-1<y<\eta^{n+1}(t,x)\big\}$:
\begin{equation}\label{eq:iteration-pressure}
\left\{
\begin{aligned}
&-\Delta P^{n+1}=\big(\pa_i( {v_1}^j)^{n+1}\pa_j ( {v_1}^i)^{n+1}\big)\circ \Phi^{n+1}_1\circ(\Phi^{n+1})^{-1},\\
&P^{n+1}\mid_{y=\eta^{n+1}}=0,\quad (\pa_yP^{n+1})\mid_{y=-1}=-1,
\end{aligned}
\right.
\end{equation}
where the map $\Phi^{n+1}$ and $\Phi^{n+1}_1$ are given by
\beno
&&\Phi^{n+1}:(x,z)\in \cS\longmapsto(x,\rho_\delta^{n+1}(x,z))\in \widetilde{\Om}^{n+1}_t,\\
&&\Phi^{n+1}_1:(x,z)\in \cS\longmapsto(x,\rho_{\delta,1}^{n+1}(x,z))\in {\Om}^{n+1}_t.
\eeno
Here $\rho_\delta^{n+1}=\rho_{\delta,\eta^{n+1}}$ and  $\rho_{\delta,1}^{n+1}=\rho_{\delta,\eta^{n+1}_1}$ with $\rho_{\delta,\eta}(x,z)=z+(1+z)e^{\delta z|D|}\eta$.

\subsection{Symmetrization}

We introduce a good unknown $U^{n+1}=V^{n+1}+T_{\zeta^n}B^{n+1}$.
It follows from (\ref{eq:iteration-velocity}) that
\begin{equation}\label{eq:iterative-U+zeta}
\left\{
\begin{aligned}
&D_tU^{n+1}+T_{a^n}\zeta^{n+1}=h_1^{n}+[D_t,T_{\zeta^n}]B^{n+1},\\
&D_t\zeta^{n+1}=T_{\lam^n}U^{n+1}+h_2^n+R^n_\om,
\end{aligned}\right.
\end{equation}
where $(h_1^n, h_2^n)$ is given by
\begin{align*}
h_1^n=&(T_{V^n_1}-V^n_1)\cdot\nabla V^n-R(a^n,\zeta^n)+T_{\zeta^n}(T_{V^n_1}-V^n_1)\cdot\nabla B^n,\\
h_2^n=&(T_{V^n_1}-V^n_1)\cdot\na \zeta^n+[T_{\zeta^n}, T_{\lambda^n}]B_1^n+(\zeta^n-T_{\zeta^n})T_{\lambda^n} B^n\\
&+R(\eta^n_1)V^n_1+\zeta^n R(\eta^n_1)B^n_1.
\end{align*}
A direct calculation gives
\ben\label{eq:iteration-zeta}
D_t^2\zeta^{n+1}+ T_{a^n\lambda^n}\zeta^{n+1}=f^{n}_1+f_2^n+D_tf_3^n+D_tf_4^n+D_tf_\om^n,
\een
where $(f^n_1, f_2^n, f_3^n, f_4^n, f_\om^n)$ are given by
\begin{align*}
f^n_1=& D_t\big((T_{V^n_1}-V^n_1)\cdot\na \zeta^n+[T_{\zeta^n}, T_{\lambda^n}]B_1^n\big)+(T_{a^n\lambda^n}-T_{\lambda^n}T_{a^n})\zeta^{n+1},\\
f^n_2=&[D_t, T_{\lambda^n}]U^{n+1}+T_{\lambda^n}\big(h_1^n+[D_t, T_{\zeta^n}]B^{n+1}\big),\\
f_3^n=&(\zeta^n-T_{\zeta^n})T_{\lambda^n} B^n+(\zeta^n-T_{\zeta^n})R(\eta^n_1)B^n_1,\\
f_4^n=&R(\eta^n_1)V^n_1+T_{\zeta^n}R(\eta^n_1)B^n_1,\\
f_\om^n=&R_\om^n.
\end{align*}

The local existence of smooth solution for the approximate system (\ref{eq:iteration-velocity})--(\ref{eq:iteration-pressure})
can be proved by using the theory of symmetric hyperbolic system and elliptic equations.
Here we ignore the proof.

\section{Uniform energy estimates}

\subsection{Set-up}

Throughout this section, we denote
\beno
D_t\triangleq\pa_t+T_{V^n}\cdot\na,\quad D_t^u\triangleq\pa_t+V^n\cdot\na,\quad D_t^b\triangleq\pa_t+V^n_{b}\cdot\na,\quad {\cD}_t\triangleq\pa_t+{v}^n_2\cdot\na_{x,y},
\eeno
where ${v}^n_2=\big((v^n)^h, (v_1^n)^{d+1}\big)$. We denote
\beno
\widetilde{D}_t\triangleq\pa_t+T_{V^{n+1}}\cdot\na,\quad\widetilde{D}_t^u\triangleq\pa_t+V^{n+1}\cdot\na, \quad \widetilde{D}_t^{b}\triangleq\pa_t+V^{n+1}_b\cdot\na,
\quad\widetilde{\cD}_t\triangleq\pa_t+v^{n+1}_2\cdot\na_{x,y},
\eeno
where $v^{n+1}_2=\big((v^{n+1})^h, (v_1^{n+1})^{d+1}\big)$. For a function $f(x,y)$ defined on $\big\{(x,y): -1<y<\eta(x)\big\}$,
we denote $\widetilde{f}(x,z)\triangleq f(x,\rho_\delta(x,z))$, where $\rho_{\delta}(x,z)=z+(1+z)e^{\delta z|D|}\eta$.

Let $E_i(i=1,\cdots,6)$ be some constants determined later. Assume that there exists $T$ independent of $n$ determined later such that the solution in the $n$-th iteration satisfies
\begin{itemize}
\item[H1.] For any $t\in [0,T]$, there holds
\begin{align*}
&\big\|(V^n, B^n, V_b^n,V^n_1, B^n_1, V^n_{b,1})(t)\big\|_{H^{s-\f12}}+\|(\eta^n,\eta^n_1)(t)\|_{H^{s-\f12}}+\|\widetilde{\om}^n(t)\|_{H^{s-1}(\overline{\cS})}\le E_1;
\end{align*}

\item[H2.] For any $t\in [0,T]$, there holds
\begin{align*}
&\|\na_{x,z}\widetilde{v_1^n}(t)\|_{H^{s-1}(\overline{\cS})}+\|\widetilde{v_\om^n}(t)\|_{H^s(\overline{\cS})}+\|a^n(t)\|_{H^{s-\f32}}+\|\na P^{n}(t,\cdot,-1)\|_{H^{s-\f32}}\le E_2;
\end{align*}

\item[H3.] For any $t\in [0,T]$, there holds
\begin{align*}
\|a^n(t)\|_{H^{s-\f12}}+\|\na P^{n}(t,\cdot,-1)\|_{H^{s-\f12}}\le E_3^1,\quad \|\eta^n(t)\|_{H^{s+\f12}}\le E_3^2;
\end{align*}

\item[H4.] For any $t\in [0,T]$, there holds
\begin{align*}
&\|D_t(V^n_1,B^n_1)(t)\|_{H^{s-\f32}}+\|\pa_t\eta^n(t)\|_{H^{s-\f32}}+\|D_t^u\eta_1^n(t)\|_{H^{s-\f12}}\\
&\quad+\|\widetilde{{\cD}_t{\om}^n}(t)\|_{H^{s-1}(\overline{\cS})}+\|\widetilde{{\cD}_tv_\om^n}(t)\|_{X^{s-\f12}([-\f12,0])}\le E_4;
\end{align*}

\item[H5.] For any $t\in [0,T]$, there holds
\begin{align*}
&\|\big(D_t(V^n_1, B^n_1), D_t^bV^n_{b,1}\big)(t)\|_{H^{s-\f12}}+\|\pa_t^2\eta^n(t)\|_{H^{s-\f52}}+\|\pa_ta^n(t)\|_{H^{s-\f32}}\\
&+\|(D_t^u)^2\eta_1^n(t)\|_{H^{s-\f12}}+\|\widetilde{\cD_t^2\om^n}(t)\|_{H^{s-1}(\overline{\cS})}+\|\na_{x,z}\widetilde{{\cD}_tv^{n}_1}\|_{H^{s-1}(\overline{\cS})}\\
&\quad+\|\na\pa_tP^{n}(t,\cdot,-1)\|_{H^{s-\f32}}\le E_5;
\end{align*}

\item[H6.]  For any $t\in [0,T]$, there holds
\begin{align*}
&\|\pa_t\big(D_t(V^n_1, B^n_1),D_t^bV^n_{b,1}\big)(t)\|_{H^{s-\f32}}+\|\widetilde{\cD_t^2v^n_\om}(t)\|_{X^{s-\f32}([-\f12,0])}\le E_6;
\end{align*}

\item[H7.]  For any $(t,x)\in [0,T]\times \R^d$, there holds
\begin{align*}
&a^n(t,x)=-\pa_y P^n\big|_{y=\eta^n(t,x)}\ge \f {c_0} 2;
\end{align*}

\item[H8.]  For any $(t,x)\in [0,T]\times \R^d$, there holds
\begin{align*}
\eta^n(t,x)+1\ge \f {h_0} 2,\quad \eta^n_1(t,x)+1\ge \f {h_0} 2.
\end{align*}
\end{itemize}

The purpose of this section is to show that (H1)--(H8) also hold for the solution in the $(n+1)$-th iteration.

In the sequel, we denote $\cA_k=\cA_k\big(E_1,\cdots, E_k\big)$ for $k=1,\cdots,6$,
and $P_{s,\eta_1,\cdots,\eta_k}=P_{s,\eta_1,\cdots, \eta_k}\big(\|\eta_1\|_{H^{s}},\cdots, \|\eta_k\|_{H^{s}}\big)$,
where $\cA_k$ and $P_{s,\eta_1,\cdots, \eta_k}$ are some nondecreasing functions depending on $c_0, h_0$ and may change from line to line.
We also denote by $\cP(\cdot,\cdots,\cdot)$ some increasing function depending on $c_0, h_0$, which may be different from line to line.

\subsection{Energy functional}

We introduce the energy functional $\cE^{n+1}(t)$ defined by
\begin{align*}
\cE^{n+1}(t)\eqdef& \cE^{n+1}_1(t)+\cE^{n+1}_2(t),
\end{align*}
where $\cE^{n+1}_1(t)$ and $\cE^{n+1}_2(t)$ are given by
\begin{align*}
\cE^{n+1}_1(t)=&\|D_t\zeta^{n+1}\|_{H^{s-1}}+\|T_{\sqrt{a^n\lambda^n}}\zeta^{n+1}\|_{H^{s-1}},\\
\cE^{n+1}_2(t)=&\|\widetilde{\om}^{n+1}\|_{H^{s-1}(\overline{\cS})}+\|\zeta^{n+1}\|_{H^{s-\f32}}
+\|(V^{n+1}, B^{n+1}, V^{n+1}_b)\|_{H^{s-\f12}}\\
&+\|(V^{n+1}_1, B^{n+1}_1, V^{n+1}_{b,1}, \eta_1^{n+1})\|_{H^{s-\f12}}.
\end{align*}

Using the equations, we can establish the following regularity information in terms of $\cE^{n+1}(t)$ for the solution in the $(n+1)$-th iteration.

\begin{lemma}\label{lem:energy-relation}
It holds that
\begin{align*}
\|\eta^{n+1}\|_{H^{s+\f12}}\le \cA_2\cE^{n+1}(t),\quad \|U^{n+1}\|_{H^{s-\f12}}\le \cA_1\cE^{n+1}_2(t).
\end{align*}
\end{lemma}

\no{\bf Proof.}\,We write
\beno
\zeta^{n+1}=T_{(\sqrt{a^n\lambda^n})^{-1}}T_{\sqrt{a^n\lambda^n}}\zeta^{n+1}+\big(T_{(\sqrt{a^n\lambda^n})^{-1}}T_{\sqrt{a^n\lambda^n}}-1\big)\zeta^{n+1},
\eeno
from which and Proposition \ref{prop:symbolic calculus}, we infer that
\begin{align*}
\|\zeta^{n+1}\|_{H^{s-\f12}}\le& \cA_2\|\zeta^{n+1}\|_{H^{s-\f32}}+\cA_2\|T_{\sqrt{a^n\lambda^n}}\zeta^{n+1}\|_{H^{s-1}}\\
\le& \cA_2\cE^{n+1}(t).
\end{align*}
Recall that $\eta^{n+1}$ satisfies
\beno
-\Delta\eta^{n+1}+\eta^{n+1}=-\dive\zeta^{n+1}+ {\eta}^{n+1}_1,
\eeno
which implies that
\beno
\|\eta^{n+1}\|_{H^{s+\f12}}\le \|\zeta^{n+1}\|_{H^{s-\f12}}+\|\eta^{n+1}_1\|_{H^{s-\f32}}\le \cA_2\cE^{n+1}(t).
\eeno

The second inequality follows from $U^{n+1}=V^{n+1}+T_{\zeta^n}B^{n+1}$.\ef

\begin{lemma}\label{lem:U-n+1-Dt}
It holds that
\beno
\|D_tU^{n+1}\|_{H^{s-\f12}}\le \cA_1+\cA_4\cE^{n+1}(t).
\eeno
\end{lemma}

\no{\bf Proof.}\,Recall that
\ben\label{eq:DtU-n+1}
D_tU^{n+1}+T_{a^n}\zeta^{n+1}=h_1^{n}+[D_t,T_{\zeta^n}]B^{n+1}.
\een
By Lemma \ref{lem:remaider}, we have
\beno
\|h_1^n\|_{H^{s-\f12}}\le \cA_1,
\eeno
which along with Proposition \ref{prop:commutator-tame} and  Lemma \ref{lem:energy-relation} gives
\begin{align*}
\|D_tU^{n+1}\|_{H^{s-\f12}}\le& \cA_1+\cA_2\|\zeta^{n+1}\|_{H^{s-\f12}}+\cA_4\|B^{n+1}\|_{H^{s-\f12}}\\
\le& \cA_1+\cA_4\cE^{n+1}(t).
\end{align*}
The proof is finished.\ef

\begin{lemma}\label{lem:eta-n+1-t2}
It holds that
\beno
&&\|\widetilde{D}_t^u\eta^{n+1}_1\|_{H^{s-\f12}}\le \cE^{n+1}_2(t),\\
&&\|(\widetilde{D}_t^u)^2\eta^{n+1}_1\|_{H^{s-\f12}}\le \cA_3\big(1+\cE^{n+1}_2(t)\big)^2,\\
&&\|D_t\zeta^{n+1}\|_{H^{s-\f32}}\le \cA_2\big(1+\cE^{n+1}_2(t)\big),\\
&&\|\pa_t\eta^{n+1}\|_{H^{s-\f32}}\le \cA_2\big(1+\cE^{n+1}_2(t)\big)^2,\\
&&\|\pa_t^2\eta^{n+1}\|_{H^{s-\f52}}\le \cA_4\big(1+\cE^{n+1}_2(t)\big)^3.
\eeno
\end{lemma}

\no{\bf Proof.}\,The first two inequalities are obvious. Recall that
\beno
&&D_t\zeta^{n+1}=T_{\lam^n}(V^{n+1}+T_{\zeta^n}B^{n+1})+(T_{V_1^{n}}-V_1^{n})\cdot\zeta^n+[T_{\zeta^n},T_{\lam^n}]B^{n}_1\\
&&\qquad\qquad+(\zeta^{n} -T_{\zeta^{n} })T_{\lam^n}B^n+R(\eta^n_1)V^{n}_1+\zeta^n R(\eta^n_1)B^{n}_1+R^n_\om.
\eeno
Then it follows from Proposition \ref{prop:symbolic calculus} and Proposition \ref{prop:DN-remainder} that
\begin{align*}
\|D_t\zeta^{n+1}\|_{H^{s-\f32}}\le \cA_1\cE_2^{n+1}(t)+\cA_2.
\end{align*}
For $\pa_t\eta^{n+1}$, we get by the elliptic estimate that
\begin{align*}
\|\pa_t\eta^{n+1}\|_{H^{s-\f32}}\le& C\big(\|\pa_t\zeta^{n+1}\|_{H^{s-\f52}}+\|\pa_t\eta^{n+1}_1\|_{H^{s-\f72}}\big)\\
\le& \cA_2+\cA_1\cE^{n+1}_2(t)+C\cE^{n+1}_2(t)^2.
\end{align*}
Using Proposition \ref{prop:symbolic calculus} and Proposition \ref{prop:DN-remainder-Dt}, we can also deduce that
\begin{align*}
\|\pa_tD_t\zeta^{n+1}\|_{H^{s-\f52}}\le \cA_4\cE_2^{n+1}(t)+\cA_4,
\end{align*}
which implies the estimate for $\pa_t^2\eta^{n+1}$.\ef

\begin{lemma}\label{lem:VB-n+1-Dt}
It holds that
\beno
&&\|(D_tV^{n+1}, \widetilde{D}_tV^{n+1}_1 )\|_{H^{s-\f12}}\le \cA_3\big(1+\cE^{n+1}(t)\big),\\
&&\|(D_tB^{n+1},\widetilde{D}_tB^{n+1}_1,\widetilde{D}_t^bV^{n+1}_{b,1})\|_{H^{s-\f12}}\le \cA_3,\\
&&\|(D_t^{b}V_b^{n+1}, \widetilde{D}_t^bV^{n+1}_{b,1})\|_{H^{s-\f32}}\le \cA_2,\\
&&\|(\widetilde{D}_tV^{n+1}_1, \widetilde{D}_tB^{n+1}_1)\|_{H^{s-\f32}}\le \cA_2\big(1+\cE_2^{n+1}(t)\big),\\
&&\|(\pa_t\widetilde{D}_t(V^{n+1}_1,B^{n+1}_1),\pa_t\widetilde{D}_t^{b}V_{b,1}^{n+1})\|_{H^{s-\f32}}\le \cA_5\big(1+\cE_2^{n+1}(t)\big).
\eeno
\end{lemma}

\no{\bf Proof.}\,Recall that $D_tV^{n+1}=\widetilde{D}_tV^{n+1}_1$ and
\beno
D_tV^{n+1}=-T_{\zeta^{n}}a^{n}-T_{a^{n}}\zeta^{n+1}-R(\zeta^n, a^n)+(T_{V_1^{n}}-V_1^{n})\cdot\nabla V^{n}.
\eeno
It follows from Lemma \ref{lem:remaider} and Lemma \ref{lem:energy-relation} that
\beno
\|(D_tV^{n+1}, \widetilde{D}_tV^{n+1}_1 )\|_{H^{s-\f12}}\le \cA_3+\cA_2\|\zeta^{n+1}\|_{H^{s-\f12}}
\le\cA_3+\cA_2\cE^{n+1}(t).
\eeno
The proof of the other inequalities is similar. \ef

\subsection{Estimate of the velocity}
In the sequel, we assume that (H8) holds for $(\eta^{n+1},\eta^{n+1}_1)$, i.e.,
\beno
\eta^{n+1}(t,x)+1\ge \f {h_0} 2,\quad \eta_1^{n+1}(t,x)+1\ge \f {h_0} 2.
\eeno

Recall that $v^{n+1}$ satisfies
\beno
\left\{
\begin{aligned}
&-\Delta_{x,y}{ v}^{n+1}=\nabla_{x,y}\times{\om}^{n+1}\quad\textrm{in}\quad\Om^{n+1}_t,\\
&{v}^{n+1}\mid_{y= {\eta}^{n+1}_1}=(V^{n+1},B^{n+1}),\quad v^{n+1} |_{y=-1}= (V^{n+1}_{b},0),
\end{aligned}
\right.
\eeno
and $v_1^{n+1}$ satisfies
\beno
\left\{
\begin{aligned}
&-\Delta_{x,y}{v_1}^{n+1}=\nabla_{x,y}\times{\om}^{n+1}\quad\textrm{in}\quad\Om^{n+1}_t,\\
&{v}^{n+1}_1\mid_{y= {\eta}^{n+1}_1}=(V^{n+1}_1,B^{n+1}_1),\quad v^{n+1}_1|_{y=-1}= (V^{n+1}_{b,1},0).
\end{aligned}
\right.
\eeno

First of all, we apply Proposition \ref{prop:elliptic-boun} to obtain
\begin{align}
\|\na_{x,z}\widetilde{v^{n+1}}\|_{H^{s-1}(\overline{\cS})}\le&
P_{s-\f12,\eta^{n+1}_1}\big(\|(V^{n+1}, B^{n+1},V_b^{n+1})\|_{H^{s-\f12}}+\|\widetilde{\om^{n+1}}\|_{H^{s-1}}\big)\nonumber\\
\le& \cP(\cE^{n+1}_2(t)).\label{eq:v-n+1}
\end{align}
Similarly, we have
\begin{align}
\|\na_{x,z}\widetilde{v^{n+1}_1}\|_{H^{s-1}(\overline{\cS})}\le \cP(\cE^{n+1}_2(t)).\label{eq:v1-n+1}
\end{align}

A direct calculation gives
\beno
\left\{
\begin{array}{ll}
\Delta_{x, y}\widetilde{\cD}_tv^{n+1}_1=h^{n+1}_\om\quad \textrm{in} \quad \Omega_t^{n+1},\\
\widetilde{\cD}_tv^{n+1}_1|_{y=\eta^{n+1}_1}=(\widetilde{D}_tV^{n+1}_1, \widetilde{D}_t B^{n+1}_1),\quad \widetilde{\cD}_tv^{n+1}_1|_{y=-1}=(\widetilde{D}_t^{b}V_{b,1}^{n+1},0),
\end{array}\right.
\eeno
where
\begin{align*}
-h_\om^{n+1}=&\na_{x,y}\times\widetilde{\cD}_t\om^{n+1}-\na_{x,y}v_2^{n+1}\cdot\na_{x,y}w^{n+1}+\Delta_{x,y}v_2^{n+1}\cdot\na_{x,y}v_1^{n+1}\\
&+2\pa_iv_2^{n+1}\cdot\na_{x,y}\pa_iv_1^{n+1}.
\end{align*}
We can deduce from Lemma \ref{lem:product-full}, (\ref{eq:v-n+1}) and (\ref{eq:v1-n+1}) that
\begin{align*}
\|\widetilde{h_\om^{n+1}}\|_{H^{s-2}(\overline{\cS})}\le& P_{s-\f12,\eta_1^{n+1}}\big(\|\widetilde{\widetilde{\cD}_t\om^{n+1}}\|_{H^{s-1}(\overline{\cS})}
+\|\widetilde{\om^{n+1}}\|_{H^{s-1}(\overline{\cS})}^2+\|\na_{x,z}(\widetilde{v^{n+1}},\widetilde{v_1^{n+1}})\|^2_{H^{s-1}(\overline{\cS})}\big)\\
\le & P_{s-\f12,\eta_1^{n+1}}\big(\|\widetilde{\om^{n+1}}\|_{H^{s-1}(\overline{\cS})}^2+\|\na_{x,z}(\widetilde{v^{n+1}},\widetilde{v_1^{n+1}})\|^2_{H^{s-1}(\overline{\cS})}\big)\\
\le &\cP(\cE^{n+1}_2(t)).
\end{align*}
While, we know from Lemma \ref{lem:VB-n+1-Dt} that
\begin{align*}
\|(\widetilde{D}_tV^{n+1}_1, \widetilde{D}_tB^{n+1}_1, \widetilde{D}_t^{b}V_{b,1}^{n+1})\|_{H^{s-\f12}}
\le \cA_3\big(1+\cE^{n+1}(t)\big).
\end{align*}
Then we apply Proposition \ref{prop:elliptic-boun} again to obtain
\begin{align}
\|\na_{x,z}\widetilde{\widetilde{\cD}_tv^{n+1}_1}\|_{H^{s-1}(\overline{\cS})}
\le \cP(\cA_3,\cE^{n+1}(t)).\label{eq:v1-n+1-Dt}
\end{align}

Similarly, we have
\begin{align}
\|\na_{x,z}\widetilde{\widetilde{\cD}_tv^{n+1}}\|_{H^{s-2}(\overline{\cS})}
\le \cP(\cA_3,\cE^{n+1}(t))\label{eq:v-n+1-Dt}
\end{align}
by noting that
\begin{align*}
\|(\widetilde{D}_tV^{n+1}, \widetilde{D}_tB^{n+1}, \widetilde{D}_t^{b}V_{b}^{n+1})\|_{H^{s-\f32}}
\le \cP(\cA_3,\cE^{n+1}(t)).
\end{align*}

Using (\ref{eq:v-n+1})-(\ref{eq:v1-n+1-Dt}), we can deduce that

\begin{lemma}\label{lem:vorticity-n+1-Dt}
It holds that
\beno
&&\|\widetilde{\widetilde{\cD}_t\om^{n+1}}\|_{H^{s-1}(\overline{\cS})}\le \cP(\cE^{n+1}_2(t)),\\
&&\|\widetilde{\widetilde{\cD}_t^2\om^{n+1}}\|_{H^{s-1}(\overline{\cS})}\le \cP(\cA_3,\cE^{n+1}(t)).
\eeno
\end{lemma}

For $v_\om^{n+1}$ defined by
\beno
\left\{
\begin{array}{ll}
-\Delta_{x, y} v_{\om}^{n+1}=\na_{x,y}\times\om^{n+1}\quad \textrm{in} \quad \Om_t^{n+1},\\
v_{\om}^{n+1}|_{y=\eta^{n+1}_1}=0,\quad v_\om^{n+1}|_{y=-1}=(V_{b,1}^{n+1}, 0),
\end{array}\right.
\eeno
we can deduce in a similar way that
\begin{align}
&\|\widetilde{v^{n+1}_\om}\|_{H^{s}(\overline{\cS})}\le \cP(\cE^{n+1}_2(t)),\label{eq:vom-n+1}\\
&\|\widetilde{\widetilde{\cD}_tv^{n+1}_\om}\|_{H^{s}(\overline{\cS})}\le \cP(\cA_3,\cE^{n+1}(t)),\label{eq:vom-n+1-t1}\\
&\|\widetilde{\widetilde{\cD}_tv^{n+1}_\om}\|_{X^{s-\f12}([-\f12,0])}\le \cP(\cA_2,\cE^{n+1}_2(t)),\label{eq:vom-n+1-t1-up}
\end{align}
by using the fact that
\begin{align*}
\|\widetilde{D}_t^{b}V_{b,1}^{n+1}\|_{H^{s-\f12}}
\le \cA_3,\quad \|\widetilde{D}_t^{b}V_{b,1}^{n+1}\|_{H^{s-\f32}}
\le \cA_2.
\end{align*}

On the other hand, we know that
\beno
\left\{
\begin{array}{ll}
\Delta_{x, y}\widetilde{\cD}_t^2v^{n+1}_\om=\widetilde{\cD}_th^{n+1}_\om+[\Delta_{x, y},\widetilde{\cD}_t]\widetilde{\cD}_tv^{n+1}_\om\quad \textrm{in} \quad \Omega_t^{n+1},\\
\widetilde{\cD}_t^2v^{n+1}_\om|_{y=\eta^{n+1}_1}=0,\quad \widetilde{\cD}_t^2v^{n+1}_\om|_{y=-1}=((\widetilde{D}_t^{b})^2V_{b,1}^{n+1},0),
\end{array}\right.
\eeno
where
\begin{align*}
h_\om^{n+1}=&\na_{x,y}\times\widetilde{\cD}_t\om^{n+1}-\na_{x,y}v_2^{n+1}\cdot\na_{x,y}\om^{n+1}-\na_{x,y}\om^{n+1}\cdot\na_{x,y}v_\om^{n+1}\\
&+2\pa_iv_2^{n+1}\cdot\na_{x,y}\pa_iv_\om^{n+1}.
\end{align*}
Using (\ref{eq:v1-n+1-Dt})-(\ref{eq:vom-n+1-t1}) and Lemma \ref{lem:vorticity-n+1-Dt}, we can deduce from Lemma \ref{lem:product-full} that
\beno
\|[\Delta_{x, y},\widetilde{\cD}_t]\widetilde{\cD}_tv^{n+1}_\om\|_{H^{s-2}(\overline{\cS})}+\|\widetilde{\widetilde{\cD}_th_\om^{n+1}}\|_{H^{s-2}(\overline{\cS})}\le  \cP(\cA_3,\cE^{n+1}(t)).
\eeno
And by Lemma \ref{lem:VB-n+1-Dt}, we have
\begin{align*}
\|(\widetilde{D}_t^{b})^2V_{b,1}^{n+1}\|_{H^{s-\f32}}\le& \|\pa_t\widetilde{D}_t^{b}V_{b,1}^{n+1}\|_{H^{s-\f32}}+\|V_b^{n+1}\|_{H^{s-\f32}}
\|\na^2P^n(\cdot,-1)\|_{H^{s-\f32}}\\
\le& \cP(\cA_5,\cE^{n+1}(t)).
\end{align*}
Then we can prove that
\beno
\|\na_{x,y}\widetilde{\cD}_t^2v^{n+1}_\om\|_{L^2(\Om^{n+1}_t)}\le \cP(\cA_5,\cE^{n+1}(t)).
\eeno
Thus, Proposition \ref{prop:elliptic-tan-nontame} ensures that
\begin{align}
\|\na_{x,z}\widetilde{\widetilde{\cD}_t^2v^{n+1}_\om}\|_{X^{s-\f32}([-\f12,0])}\le \cP(\cA_5,\cE^{n+1}(t)).\label{eq:vom-n+1-Dt2}
\end{align}

\subsection{Estimate of the pressure}
Recall that the pressure $P^{n+1}$ satisfies
\beno
\left\{
\begin{aligned}
&-\Delta_{x,y}P^{n+1}=\big(\pa_i( {v}^j_1)^{n+1}\pa_j ( {v}^i_1)^{n+1}\big)\circ \Phi^{n+1}_1\circ(\Phi^{n+1})^{-1}\triangleq F,\\
&P^{n+1}\mid_{y=\eta^{n+1}}=0,\quad (\pa_yP^{n+1})\mid_{y=-1}=-1.
\end{aligned}
\right.
\eeno

Let $P_1^{n+1}=P^{n+1}+y$, $I_1=[a,0]$ for $a\in (-1,0)$. First of all, we get by a similar proof of Lemma \ref{lem:pressure-L2} that
\beno
\|\na_{x,y}P_1\|_{L^2(\widetilde{\Om}^{n+1}_t)}\le P_{s-\f12,\eta^{n+1}}\|\na_{x,z}\widetilde{v^{n+1}_1}\|_{H^{s-1}(\overline{\cS})}^2,
\eeno
from which and Proposition \ref{prop:elliptic-tan-nontame}, it follows that
\begin{align}
\|\na_{x,z}\widetilde{P^{n+1}_1}\|_{X^{s-\f12}(I_1)}\le& P_{s+\f12,\eta^{n+1}}\big(\|\na_{x,y}P^{n+1}_1\|_{L^2(\widetilde{\Om}^{n+1}_t)}
+P_{s-\f12,\eta^{n+1}_1}\|\na_{x,z}\widetilde{v^{n+1}_1}\|_{H^{s-1}(\overline{\cS})}^2\big)\nonumber\\
\le& \cP(\cA_2, \cE^{n+1}(t)).\label{eq:Pn+1-up}
\end{align}
Similarly, we can deduce that
\begin{align}
\|\na_{x,z}\widetilde{P^{n+1}_1}\|_{X^{s-\f32}(I_1)}\le \cP(\cE^{n+1}_2(t)).\label{eq:Pn+1-low}
\end{align}
These ensure that there exists $y_0\in (-1,-1+h_0]$(in fact, one can take $y_0>-1+a$ with $a$ depending on $\|\eta^{n+1}\|_{H^{s-\f12}}$) so that
\beno
&&\|P_1(\cdot,y_0)\|_{H^{s+\f12}}\le \cP(\cA_2, \cE^{n+1}(t)),\\
&&\|P_1(\cdot,y_0)\|_{H^{s-\f12}}\le \cP(\cE^{n+1}_2(t)).
\eeno
Then by the elliptic estimate in the flat strip, we obtain
\begin{align}
&\|\na_{x,y}P_1^{n+1}\|_{H^{s}(\R^d\times [-1,y_0])}\le \cP(\cA_2, \cE^{n+1}_2(t)),\label{eq:Pn+1-down}\\
&\|\na_{x,y}P_1^{n+1}\|_{H^{s-1}(\R^d\times [-1,y_0])}\le \cP(\cE^{n+1}_2(t)).\label{eq:Pn+1-down-low}
\end{align}

Using Lemma \ref{lem:eta-n+1-t2}, (\ref{eq:v1-n+1}) and (\ref{eq:v1-n+1-Dt}), we can deduce that
\beno
&&\|\widetilde{\cD_tF}\|_{H^{s-2}(\overline{\cS})}\le\cP(\cA_3,\cE^{n+1}(t)).
\eeno
It is easy to show that
\beno
\|\na_{x,y}\cD_t P^{n+1}\|_{L^2(\widetilde{\Om}^{n+1}_t)}\le \cP(\cA_3,\cE^{n+1}(t)).
\eeno
Thus, we can deduce that
\begin{align}
&\|\na_{x,z}\widetilde{\cD_tP^{n+1}}\|_{X^{s-\f32}(I_1)}+\|\na_{x,y}\cD_tP^{n+1}(\cdot,-1)\|_{H^{s-\f32}}\le \cP(\cA_3,\cE^{n+1}(t)).\label{eq:Pn+1-t1}
\end{align}

It follows from (\ref{eq:Pn+1-up}), (\ref{eq:Pn+1-low}) and (\ref{eq:Pn+1-t1}) that

\begin{lemma}\label{lem:a-n+1-t2}
It holds that
\begin{align*}
&\|a^{n+1}\|_{H^{s-\f32}}\le \cP(\cE^{n+1}_2(t)),\\
&\|a^{n+1}\|_{H^{s-\f12}}\le \cP(\cA_2, \cE^{n+1}(t)),\\
&\|\pa_ta^{n+1}\|_{H^{s-\f32}}\le \cP(\cA_3,\cE^{n+1}(t)).
\end{align*}
\end{lemma}

\subsection{Estimates of the remainder of DN operator}

In this subsection, we establish some estimates for the material derivatives of the remainder of DN operator, which will be used to
estimate $R^{n}_\om$. For this purpose,
we assume that $\eta$ is a solution of the following equation
\beno
\pa_t\eta+V\cdot\na \eta=B,
\eeno
where $(V,B)=v|_{y=\eta}$.
Let $\overline{D}_t\triangleq\pa_t+V\cdot\na$ and $\cD_t\triangleq\pa_t+v\cdot\na_{x,y}$.

\medskip

Now, we state the main result in this subsection:
\begin{proposition}\label{prop:DN-remainder-Dt}
Assume that $\eta, \overline{D}_t\eta, \overline{D}_t^2\eta\in H^{s-\f12}(\R^d)$ for $s>\f d 2+5$. Then it holds that
\beno
&&\|\overline{D}_tR(\eta)f\|_{H^{s-\f32}}\le P_{s-\f12,\eta,\overline{D}_t\eta}\mathcal{V}_1(t)\big(\|f\|_{H^{s-\f12}}+\|\overline{D}_tf\|_{H^{s-\f32}}\big),\\
&&\|\overline{D}_t^2R(\eta)f\|_{H^{s-\f32}}\le P_{s-\f12,\eta,\overline{D}_t\eta,\overline{D}_t^2\eta}\mathcal{V}_2(t)\big(\|f\|_{H^{s-\f12}}+\|\overline{D}_tf\|_{H^{s-\f12}}
+\|\overline{D}_t^2f\|_{H^{s-\f32}}\big),
\eeno
where $\mathcal{V}_1(t)$ and $\mathcal{V}_2(t)$ are given by
\beno
&&\mathcal{V}_1(t)\triangleq \cP\big(\|V(t)\|_{H^{s-\f12}}, \|\na_{x,z}\widetilde{v}(t)\|_{H^{s-1}(\overline{\cS})}\big),\\
&&\mathcal{V}_2(t)\triangleq \cP\big(\|\na_{x,z}(\widetilde{v},\widetilde{\cD_t v})(t)\|_{H^{s-1}(\overline{\cS})}, \|(V, \overline{D}_tV)\|_{H^{s-\f12}}\big)
\eeno
\end{proposition}

Let $\phi(t,x,y)$ be a solution of the elliptic equation
\beno
\left\{
\begin{array}{l}
\Delta_{x,y} \phi=0\qquad\textrm{ in}\quad \Om_t=\big\{(x,y):-1<y<\eta(t,x)\big\},\\
\phi|_{y=\eta(t,x)}=f,\quad \phi|_{y=-1}=0.
\end{array}
\right.
\eeno
Proposition \ref{prop:elliptic-boun} ensures that
\ben\label{eq:DN-phi-Hs}
\|\widetilde{\phi}\|_{H^{s}(\overline{\cS})}\le P_{s-\f12,\eta}\|f\|_{H^{s-\f12}}.
\een
We next establish some estimates for $\cD_t^k\phi$ for $k=1,2$.

\begin{lemma}\label{lem:phi-Dt}
 It holds that
\begin{align*}
&\|\widetilde{\cD_t\phi}\|_{H^{s}(\overline{\cS})}\le P_{s-\f12,\eta}\big(\|\overline{D}_t f\|_{H^{s-\f12}}+\|\na_{x,z}\widetilde{v}\|_{H^{s-1}(\overline{\cS})}\|f\|_{H^{s-\f12}}\big),\\
&\|\widetilde{\cD_t^2\phi}\|_{H^{s-1}(\overline{\cS})}\le P_{s-\f12,\eta}\mathcal{V}(t)\big(\|\overline{D}_t^2f\|_{H^{s-\f32}}+\|\overline{D}_t f\|_{H^{s-\f32}}+\|f\|_{H^{s-\f32}}\big),\\
&\|(\pa_z-T_A)\widetilde{\cD_t\phi}\|_{X^{s-\f32}([-\f12,0])}\le P_{s-\f12,\eta}\big(\|\overline{D}_tf\|_{H^{s-\f32}}+\|\na_{x,z}\widetilde{v}(t)\|_{H^{s-1}(\overline{\cS})}\|f\|_{H^{s-\f12}}\big),\\
&\|(\pa_z-T_A)\widetilde{\cD_t^2\phi}\|_{X^{s-\f32}([-\f12,0])}\le P_{s-\f12,\eta}\mathcal{V}(t)\big(\|\overline{D}_t^2f\|_{H^{s-\f32}}+\|\overline{D}_t f\|_{H^{s-\f12}}+\|f\|_{H^{s-\f12}}\big).
\end{align*}
Here  $\mathcal{V}(t)\triangleq \cP\big(\|\na_{x,z}\widetilde{v}(t)\|_{H^{s-1}(\overline{\cS})}, \|\na_{x,z}\widetilde{\cD_t v}(t)\|_{H^{s-1}(\overline{\cS})}\big)$.
\end{lemma}

\no{\bf Proof.}\,A direct calculations gives
\beno
\left\{
\begin{array}{l}
\Delta_{x,y}\cD_t\phi=\Delta_{x,y}v\cdot\na_{x,y}\phi+2\na_{x,y}v\cdot\na_{x,y}^2\phi\triangleq F,\\
\cD_t\phi|_{y=\eta}=\overline{D}_t f,\quad \cD_t\phi|_{y=-1}=0.
\end{array}
\right.
\eeno
By Lemma \ref{lem:product-full}, it is easy to show that
\beno
\|\widetilde{F}\|_{H^{s-2}(\overline{\cS})}\le P_{s-\f12,\eta}\|\na_{x,z}\widetilde{v}\|_{H^{s-1}(\overline{\cS})}\|\widetilde{\phi}\|_{H^{s}(\overline{\cS})}.
\eeno
Then we deduce from Proposition \ref{prop:elliptic-boun} and (\ref{eq:DN-phi-Hs}) that
\ben\label{eq:DN-phi-Dt}
\|\widetilde{\cD_t\phi}\|_{H^{s}(\overline{\cS})}\le P_{s-\f12,\eta}\big(\|\overline{D}_t f\|_{H^{s-\f12}}+\|\na_{x,z}\widetilde{v}\|_{H^{s-1}(\overline{\cS})}\|f\|_{H^{s-\f12}}\big).
\een
Similar, we have
\beno
\|\widetilde{\cD_t\phi}\|_{H^{s-1}(\overline{\cS})}\le P_{s-\f12,\eta}\big(\|\overline{D}_t f\|_{H^{s-\f32}}
+\|\na_{x,z}\widetilde{v}\|_{H^{s-2}(\overline{\cS})}\|f\|_{H^{s-\f32}}\big).
\eeno
Then a similar argument of Proposition \ref{prop:DN-remainder} implies the third inequality.

By Lemma \ref{lem:product-full} again, we have
\begin{align*}
\|\widetilde{\cD_t F}\|_{H^{s-3}(\overline{\cS})}\le& P_{s-\f12,\eta}\|\na_{x,z}\widetilde{v}\|_{H^{s-2}(\overline{\cS})}^2\|\widetilde{\cD_t\phi}\|_{H^{s-1}(\overline{\cS})}\\
&+P_{s-\f12,\eta}\big(\|\na_{x,z}\widetilde{\cD_tv}\|_{H^{s-2}(\overline{\cS})}+\|\na_{x,z}\widetilde{v}\|_{H^{s-2}(\overline{\cS})}^2\big)
\|\widetilde{\phi}\|_{H^{s-1}(\overline{\cS})}.
\end{align*}
This implies the second inequality. We also have
\begin{align*}
\|\widetilde{\cD_t F}\|_{H^{s-2}(\overline{\cS})}\le& P_{s-\f12,\eta}\|\na_{x,z}\widetilde{v}\|_{H^{s-1}(\overline{\cS})}^2\|\widetilde{\cD_t\phi}\|_{H^{s}(\overline{\cS})}\\
&+P_{s-\f12,\eta}\big(\|\na_{x,z}\widetilde{\cD_tv}\|_{H^{s-1}(\overline{\cS})}+\|\na_{x,z}\widetilde{v}\|_{H^{s-1}(\overline{\cS})}^2\big)
\|\widetilde{\phi}\|_{H^{s}(\overline{\cS})}\\
\le &P_{s-\f12,\eta}\|\na_{x,z}\widetilde{v}\|_{H^{s-1}(\overline{\cS})}^2\big(\|\overline{D}_t f\|_{H^{s-\f12}}+\|\na_{x,z}\widetilde{v}\|_{H^{s-1}(\overline{\cS})}\|f\|_{H^{s-\f12}}\big)\\
&+P_{s-\f12,\eta}\big(\|\na_{x,z}\widetilde{\cD_tv}\|_{H^{s-1}(\overline{\cS})}+\|\na_{x,z}\widetilde{v}\|_{H^{s-1}(\overline{\cS})}^2\big)
\|f\|_{H^{s-\f12}}.
\end{align*}
Then a similar argument of Proposition \ref{prop:DN-remainder} implies the last inequality.\ef

\medskip

\begin{lemma}\label{lem:phi-Dt2}
It holds that
\beno
&&\|\overline{D}_t\widetilde{\phi}|_{z=0}\|_{H^{s-\f32}}\le P_{s-\f12,\eta}\big(\|\na_{x,z}\widetilde{v}\|_{H^{s-1}(\overline{\cS})}\|f\|_{H^{s-\f32}}+\|\overline{D}_tf\|_{H^{s-\f32}}\big),\\
&&\|\overline{D}_t^2\widetilde{\phi}|_{z=0}\|_{H^{s-\f32}}\le P_{s-\f12,\eta}\mathcal{V}(t)\big(\|f\|_{H^{s-\f32}}+\|{D}_tf\|_{H^{s-\f32}}+\|\overline{D}_t^2f\|_{H^{s-\f32}}\big).
\eeno
In addition, we have
\beno
&&\big\|\overline{D}_t\na_{x,z}\widetilde{\phi}|_{z=0}\big\|_{H^{s-\f52}}
\le P_{s-\f12,\eta,\overline{D}_t\eta}\big(\|\na_{x,z}\widetilde{v}\|_{H^{s-1}(\overline{\cS})}\|f\|_{H^{s-\f32}}+\|\overline{D}_tf\|_{H^{s-\f32}}\big),\\
&&\big\|\overline{D}_t^2\na_{x,z}\widetilde{\phi}|_{z=0}\big\|_{H^{s-\f52}}
\le P_{s-\f12,\eta,\overline{D}_t\eta, \overline{D}_t^2\eta}\mathcal{V}(t)\big(\|f\|_{H^{s-\f32}}+\|\overline{D}_tf\|_{H^{s-\f32}}+\|\overline{D}_t^2f\|_{H^{s-\f32}}\big).
\eeno
\end{lemma}

\no{\bf Proof.}\,The first two inequalities follows from Lemma \ref{lem:phi-Dt} and the fact that
\beno
\overline{D}_t^k\widetilde{\phi}|_{z=0}=\widetilde{\cD_t^k\phi}|_{z=0}.
\eeno
The last two inequalities can be deduced from Lemma \ref{lem:product}, Lemma \ref{lem:phi-Dt} and the formulas
\beno
&&\overline{D}_t\na\widetilde{\phi}|_{z=0}=\widetilde{\cD_t\na\phi}|_{z=0}+\widetilde{\cD_t\pa_y\phi}|_{z=0}\na\eta+\widetilde{\pa_y\phi}|_{z=0}\overline{D}_t\na\eta,\\
&&\overline{D}_t\pa_z\widetilde{\phi}|_{z=0}
=\widetilde{\cD_t\pa_y\phi}|_{z=0}(1+\eta+\delta|D|\eta)+\widetilde{\pa_y\phi}|_{z=0}\overline{D}_t(\eta+\delta|D|\eta).
\eeno
We omit the details.\ef

\medskip

\no{\bf Proof of Proposition \ref{prop:DN-remainder-Dt}}.\,We write
\begin{align*}
\overline{D}_tR_1(\eta)f=&\big[\overline{D}_t, T_{\zeta_1}T_A-T_{\zeta_1 A}\big]\widetilde{\phi}|_{z=0}+\big(T_{\zeta_1}T_A-T_{\zeta_1 A}\big)\widetilde{\cD_t\phi}|_{z=0}\\
=&\big(T_{\pa_t\zeta_1}T_A+T_{\zeta_1}T_{\pa_t A}-T_{\pa_t(\zeta_1 A)}\big)\widetilde{\phi}|_{z=0}+\big[V\cdot\na, T_{\zeta_1}T_A-T_{\zeta_1 A}\big]\widetilde{\phi}|_{z=0}\\
&+\big(T_{\zeta_1}T_A-T_{\zeta_1 A}\big)\widetilde{\cD_t\phi}|_{z=0}.
\end{align*}
Using Proposition \ref{prop:symbolic calculus} and Lemma \ref{lem:phi-Dt2}, and the following trick
\beno
\|\pa_tg\|_{H^{s-\f52}}\le \|\overline{D}_tg\|_{H^{s-\f52}}+\|V\|_{H^{s-\f52}}\|g\|_{H^{s-\f32}},
\eeno
we can deduce that
\beno
\|\overline{D}_tR_1(\eta)f\|_{H^{s-\f32}}\le P_{s-\f12,\eta,\overline{D}_t\eta}\mathcal{V}(t)\big(\|f\|_{H^{s-\f32}}+\|\overline{D}_tf\|_{H^{s-\f32}}\big).
\eeno

Similarly, we have
\beno
\|\overline{D}_tR_2(\eta)f\|_{H^{s-\f32}}\le P_{s-\f12,\eta,\overline{D}_t\eta}\mathcal{V}(t)\big(\|f\|_{H^{s-\f12}}+\|\overline{D}_tf\|_{H^{s-\f32}}\big).
\eeno
The same estimates holds for $\overline{D}_tR_3(\eta)$ by using the following type estimate
\begin{align*}
&\|\overline{D}_tR(f,g)\|_{H^{s-\f32}}\le \|R(\pa_t f,g)+R(f,\pa_tg)+V\cdot \na R(f,g)\|_{H^{s-\f32}}\\
&\le C\|\pa_t f\|_{H^{s-\f52}}\|g\|_{H^{s-\f32}}+C\|\pa_t g\|_{H^{s-\f52}}\|f\|_{H^{s-\f32}}+C\|V\|_{H^{s-\f32}}\|f\|_{H^{s-\f32}}\|g\|_{H^{s-\f32}}.
\end{align*}
The estimates for $\overline{D}_t^2R(\eta)f$ can be similarly proved by a very tedious computation
with the help of Lemma \ref{lem:phi-Dt} and Lemma \ref{lem:phi-Dt2}. We omit the details.\ef
\medskip

Applying Proposition \ref{prop:DN-remainder-Dt} with $v=v_1^n$ and $(V,B)=(V^n_1, B^n_1)$, we deduce that

\begin{lemma}\label{lem:DN-remainder-Dt2}
It holds that
\begin{align*}
&\|D_t^uR(\eta_1^{n})(V_1^n,B_1^n)\|_{H^{s-\f32}}\le \cP(E_1,E_2,E_4),\\
&\|(D_t^u)^2R(\eta_1^{n})(V_1^n,B_1^n)\|_{H^{s-\f32}}\le \cA_6.
\end{align*}

\end{lemma}

\subsection{Energy estimates}

\begin{proposition}\label{prop:vorticity-n+1}
It holds that
\begin{align*}
\|\widetilde{\om}^{n+1}(t)\|_{H^{s-1}(\overline{\cS})}\le \|\widetilde{\om}_0\|_{H^{s-1}(\overline{\cS})}+\int_0^t\cP(\cA_1, \cE^{n+1}(t'))dt'.
\end{align*}
\end{proposition}

\no{\bf Proof.}\,Recall that the vorticity $\om^{n+1}$ satisfies
\beno
\left\{
\begin{aligned}
&\pa_t {\om}^{n+1}+ \big((v^{n+1})^h\cdot\na+(v_1^{n+1})^{d+1}\pa_y\big){\om}^{n+1}= {\om}^{n+1}\cdot\nabla_{x,y}{v}^{n+1}_1,\\
& {\om}^{n+1}\mid_{t=0}=\om_0.
\end{aligned}
\right.
\eeno
In terms of the $(x,z)$ variable, we have
\ben
\pa_t \widetilde{\om}^{n+1}+\overline{v}^{n+1}\cdot\nabla_{x,z}\widetilde{\om}^{n+1}
=F^{n+1},
\een
where
\beno
&&\overline{v}^{n+1}=\big((\widetilde{v}^{n+1})^h,\f{1}{\pa_z \rho_{\delta,\eta^{n+1}_1}}(\widetilde{v_1}^{d+1}
-\pa_t\rho_{\delta,\eta^{n+1}_1}-(\widetilde{v}^{n+1})^h\cdot\na\rho_{\delta,\eta^{n+1}_1})\big),\\
&&F^{n+1}=(\widetilde{\om}^{n+1})^h\cdot\big(\na\widetilde{v_1}^{n+1}-\f {\na \rho_{\delta,\eta^{n+1}_1}} {\pa_z\rho_{\delta,\eta^{n+1}_1}}\pa_z\widetilde{v_1}\big)
+(\widetilde{\om}^{n+1})^{d+1}\f {\pa_z \widetilde{v_1}} {\pa_z\rho_{\delta,\eta^{n+1}_1}}.
\eeno
By Lemma \ref{lem:product-full}, it is easy to see that
\beno
\|F^{n+1}\|_{H^{s-1}(\overline{\cS})}\le \cP_{s-\f12,\eta^{n+1}_1}\|\na_{x,z}(\widetilde{v^{n+1}},\widetilde{v^{n+1}_1})\|_{H^{s-1}(\overline{\cS})}\|\widetilde{\om}^{n+1}\|_{H^{s-1}(\overline{\cS})}.
\eeno
By Lemma \ref{lem:v-trans}, we have
\beno
\|\na_{x,z}\overline{v}^{n+1}\|_{H^{s-1}(\overline{\cS})}\le \cP_{s-\f12,\eta^{n+1}_1}\|\na_{x,z}(\widetilde{v^{n+1}},\widetilde{v^{n+1}_1})\|_{H^{s-1}(\overline{\cS})}.
\eeno
Then by a similar proof of Proposition \ref{prop:vorticity}, we deduce
\begin{align}\label{eq:vorticityn-est}
\f d {dt}\|\widetilde{\om}^{n+1}\|_{H^{s-1}(\overline{\cS})}\le \cP_{s-\f12,\eta^{n+1}_1}\|\na_{x,z}(\widetilde{v^{n+1}},\widetilde{v^{n+1}_1})\|_{H^{s-1}(\overline{\cS})}.
\end{align}
This together with (\ref{eq:v-n+1}) and (\ref{eq:v1-n+1}) gives the proposition.\ef

\begin{proposition}\label{prop:VB-n+1}
It holds that
\begin{align*}
&\|(V^{n+1}, B^{n+1}, V_b^{n+1})(t)\|_{H^{s-\f12}}+\|\zeta^{n+1}(t)\|_{H^{s-\f32}}\\
&\le \Big(\|(V_0,B_0,V_{b,0})\|_{H^{s-\f12}}+\cA_3\int_0^t(1+\cE^{n+1}(t'))dt'\Big)e^{t\cA_1}.
\end{align*}

\end{proposition}

\no{\bf Proof.}\,It follows from Lemma \ref{lem:VB-n+1-Dt} that
\begin{align*}
\f12\f d {dt}\|V^{n+1}\|_{H^{s-\f12}}^2\le& \cA_3(1+\cE^{n+1}(t))\|V^{n+1}\|_{H^{s-\f12}}\\
&+\big\langle \langle D\rangle^{s-\f12}T_{V^n}\cdot\na V^{n+1}, \langle D\rangle^{s-\f12}V^{n+1}\big\rangle.
\end{align*}
We write
\begin{align*}
&\big\langle \langle D\rangle^{s-\f12}T_{V^n}\cdot\na V^{n+1}, \langle D\rangle^{s-\f12}V^{n+1}\big\rangle\\
&=\big\langle [\langle D\rangle^{s-\f12}, T_{V^n}\cdot\na]V^{n+1}, \langle D\rangle^{s-\f12}V^{n+1}\big\rangle-\f12\big\langle T_{\na\cdot V^n}\langle D\rangle^{s-\f12}V^{n+1}, \langle D\rangle^{s-\f12}V^{n+1}\big\rangle,
\end{align*}
by Lemma \ref{lem:commu-Ds}, which is bounded by
\beno
C\|V^n\|_{W^{1,\infty}}\|V^{n+1}\|_{H^{s-\f12}}^2\le C\|V^n\|_{H^{s-\f12}}\|V^{n+1}\|_{H^{s-\f12}}^2.
\eeno
This shows that
\begin{align*}
\f d {dt}\|V^{n+1}\|_{H^{s-\f12}}\le& \cA_3(1+\cE^{n+1}(t))+C\|V^n\|_{H^{s-\f12}}\|V^{n+1}\|_{H^{s-\f12}}.
\end{align*}
Then Gronwall's inequality ensures that
\beno
\|V^{n+1}(t)\|_{H^{s-\f12}}\le \Big(\|V_0\|_{H^{s-\f12}}+\cA_3\int_0^t(1+\cE^{n+1}(t'))dt'\Big)e^{t\cA_1}.
\eeno
The estimates for the other terms are similar. \ef
\medskip

Similarly, we can prove by Lemma \ref{lem:VB-n+1-Dt} and  Lemma \ref{lem:eta-n+1-t2} that

\begin{proposition}\label{prop:VB1-n+1}
It holds that
\begin{align*}
&\|(V^{n+1}_1, B^{n+1}_1, V_{b,1}^{n+1})(t)\|_{H^{s-\f12}}+\|\eta_1^{n+1}(t)\|_{H^{s-\f12}}\\
&\le \Big(\|(V_0,B_0,V_{b,0},\eta_0)\|_{H^{s-\f12}}+\cA_3\int_0^t(1+\cE^{n+1}(t'))dt'\Big)e^{C\int_0^t\cE^{n+1}_2(t')dt'}.
\end{align*}
\end{proposition}

Following the proof of Proposition \ref{prop:energy-boundary}, we can deduce that

\begin{proposition}\label{prop:U-n+1-Hs}
It holds that
\begin{align*}
\f12\f d {dt}\big(\|D_t \zeta^{n+1}\|_{H^{s-1}}^2+&\|T_{\sqrt{a^n\lambda^n}}\zeta^{n+1}\|_{H^{s-1}}^2\big)
\le \cA_4\big(\|D_t \zeta^{n+1}_1\|_{H^{s-1}}^2+\|\zeta\|_{H^{s-\f12}}^2\big)\\
&+\big\langle (f_1^n+f_2^n+D_tf_3^n+D_tf_4^n+D_tf_\om^n)_{s-1}, (D_t \zeta^{n+1})_{s-1}\big\rangle.
\end{align*}
\end{proposition}

\subsection{Nonlinear estimates}

Recall that
\beno
f^n_1=D_t\big((T_{V^n_1}-V^n_1)\cdot\na \zeta^n+[T_{\zeta^n}, T_{\lambda^n}]B_1^n\big)+(T_{a^n\lambda^n}-T_{\lambda^n}T_{a^n})\zeta^{n+1}.
\eeno
\begin{lemma}\label{lem:f1-n}
It holds that
\beno
\|f_1^n\|_{H^{s-1}}\le \cA_5.
\eeno
\end{lemma}

\no{\bf Proof.}\,We write
\begin{align*}
D_t(T_{V^n_1}-V^n_1)\cdot\na \zeta^n=-[D_t, T_{\na\zeta^n}]\cdot V^n_1+T_{\na \zeta^n}\cdot D_tV^n_1-D_tR(\na\zeta^n,V^n_1).
\end{align*}
It follows from Proposition \ref{prop:commutator-tame} that
\beno
\|[D_t, T_{\na\zeta^n}]\cdot V^n_1\|_{H^{s-1}}\le C(E_1+E_4)\|V^n_1\|_{H^{s-1}}\le \cA_4.
\eeno
And by Lemma \ref{lem:remaider}, we get
\beno
\|T_{\na \zeta^n}\cdot D_tV^n_1\|_{H^{s-1}}\le CE_1\|D_tV^n_1\|_{H^{s-1}}\le \cA_5.
\eeno
By Lemma \ref{lem:remaider} again, we get
\begin{align*}
\|D_tR(\na\zeta^n,V^n_1)\|_{H^{s-1}}\le& \|R(\na\pa_t\zeta^n, V^n_1)\|_{H^{s-1}}+\|R(\na\zeta^n, \pa_tV^n_1)\|_{H^{s-1}}
\\&+\|T_{V^n}\cdot\na R(\na\zeta^n,V^n_1)\|_{H^{s}}\le \cA_4.
\end{align*}
We write
\begin{align*}
D_t[T_{\zeta^n}, T_{\lambda^n}]B_1^n=&[T_{\pa_t\zeta^n}, T_{\lambda^n}]B_1^n+[T_{\zeta^n}, T_{\pa_t\lambda^n}]B_1^n+[T_{\zeta^n}, T_{\lambda^n}]D_tB_1^n\\
&+\big[T_{V^n}\cdot\na, [T_{\zeta^n}, T_{\lambda^n}]\big]B_1^n,
\end{align*}
which along with Proposition \ref{prop:symbolic calculus} and Corollary \ref{cor:symb} gives
\beno
\|D_t[T_{\zeta^n}, T_{\lambda^n}]B_1^n\|_{H^{s-1}}\le \cA_4.
\eeno
The proof is finished.\ef

\medskip

Recall that
\begin{align*}
f^n_2=&[D_t, T_{\lambda^n}]U^{n+1}+T_{\lambda^n}\big(h_1^n+[D_t, T_{\zeta^n}]B^{n+1}\big).
\end{align*}

\begin{lemma}\label{lem:f2-n}
It holds that
\beno
&&\|f^n_2\|_{H^{s-\f32}}\le \cP(E_1,E_2,E_4)\big(1+\cE^{n+1}_2(t)\big),\\
&&\|D_tf^n_2\|_{H^{s-\f32}}\le \cA_5\big(1+\cE^{n+1}(t)\big).
\eeno
\end{lemma}

\no{\bf Proof.}\,The first inequality is obvious. We turn to the second inequality.
It follows from Proposition \ref{prop:commutator-Dt} and Lemma \ref{lem:U-n+1-Dt} that
\begin{align*}
\|D_t[D_t, T_{\lambda^n}]U^{n+1}\|_{H^{s-\f32}}\le& \cA_5\big(\|U^{n+1}\|_{H^{s-\f12}}+\|D_tU^{n+1}\|_{H^{s-\f12}}\big)\\
\le& \cA_5\big(1+\cE^{n+1}(t)\big).
\end{align*}
We write
\begin{align*}
D_tT_{\lambda^n}\big(h_1^n+[D_t, T_{\zeta^n}]B^{n+1}\big)=&[D_t,T_{\lambda^n}]\big(h_1^n+[D_t, T_{\zeta^n}]B^{n+1}\big)\\
&+T_{\lambda^n}\big(D_th_1^n+D_t[D_t, T_{\zeta^n}]B^{n+1}\big),
\end{align*}
which along with Proposition \ref{prop:commutator-tame} and Proposition \ref{prop:symbolic calculus} gives
\begin{align*}
\|D_tT_{\lambda^n}\big(h_1^n+[D_t, T_{\zeta^n}]B^{n+1}\big)\|_{H^{s-\f32}}&\le \cA_5\big(\|h_1^n\|_{H^{s-\f12}}+\|B^{n+1}\|_{H^{s-\f12}}\\
&+\|D_th_1^n\|_{H^{s-\f12}}+\|D_tB^{n+1}\|_{H^{s-\f12}}\big).
\end{align*}
Similar to Lemma \ref{lem:f1-n}, we can prove that
\beno
\|D_th_1^n\|_{H^{s-\f32}}\le \cA_4.
\eeno
This together with Lemma \ref{lem:VB-n+1-Dt} gives
\begin{align*}
\|D_tT_{\lambda^n}\big(h_1^n+[D_t, T_{\zeta^n}]B^{n+1}\big)\|_{H^{s-\f23}}\le \cA_5\big(1+\cE^{n+1}(t)\big).
\end{align*}
Putting the above estimates together gives the second inequality.\ef
\medskip

Recall that
\beno
&&f_3^n=(\zeta^n-T_{\zeta^n})T_{\lambda^n} B^n+(\zeta^n-T_{\zeta^n})R(\eta^n_1)B^n_1,\\
&&f_4^n=R(\eta^n_1)V^n_1+T_{\zeta^n}R(\eta^n_1)B^n_1.
\eeno
Then we can deduce from Proposition \ref{prop:DN-remainder}  and Lemma \ref{lem:DN-remainder-Dt2} that
\begin{lemma}\label{lem:f3-n}
It holds that
\beno
&&\|f_3^n\|_{H^{s-1}}\le \cA_1\|\eta^n\|_{H^s},\quad \|f_3^n\|_{H^{s-\f12}}\le \cA_3,\\
&&\|D_t^uf_4^n\|_{H^{s-\f32}}\le \cP(E_1,E_2,E_4),\quad\|(D_t^u)^2f_4^n\|_{H^{s-\f32}}\le \cA_6.
\eeno
\end{lemma}

Finally, we have

\begin{lemma}\label{lem:fom-n}
It holds that
\beno
&&\|D_t^uf_\om^n\|_{H^{s-\f32}}\le \cP(E_1, E_2, E_4),\quad\|(D_t^u)^2f_\om^n\|_{H^{s-\f32}}\le \cA_6.
\eeno
\end{lemma}

\subsection{Completion of the uniform estimate}
Let
\begin{align*}
\cE_0\triangleq & \|\widetilde{\om_0}\|_{H^{s-1}(\overline{\cS})}+\|\eta_0\|_{H^{s+\f12}}+\|\widetilde{v_0}\|_{H^{s+\f12}(\overline{\cS})}
+\|D_t\zeta(0)\|_{H^{s-1}}+\|T_{\sqrt{a_0\lambda_0}}\zeta_0\|_{H^{s-1}},
\end{align*}
which is bounded by $\cP(\|\eta_0\|_{H^{s+\f12}})\big(\|v_0\|_{H^{s+\f12}(\Om_0)}+\|\eta_0\|_{H^{s+\f12}}\big)$.

Let us first assume that there exits a maximal time $T_n\in (0,T]$ so that the solution satisfies
\ben\label{ass:E-n+1}
\cE^{k}(t)\le \cP_0,\quad \cE^{k}_2(t)\le \cP_1.
\een
for $k=1,\cdots, n+1$ and $t\in [0,T_n)$ and some $\cP_0, \cP_1$ determined later.
Under this assumption, we will verify $(H1)-(H8)$ for the solution in the $(n+1)$-th iteration.
Moreover, we will show that
\ben\label{ass:E-n+1-improve}
\cE^{n+1}(t)\le \f {\cP_0} 2,\quad \cE^{n+1}_2(t)\le \f {\cP_1} 2.
\een
which implies that $T_n=T$ by a continuous argument.
\medskip

Obviously, we can take $E_1=\cP_1$. If $(\eta^{n+1},\eta_1^{n+1})$ satisfies (H8), then we deduce from (\ref{eq:v1-n+1}), (\ref{eq:vom-n+1}), Lemma \ref{lem:a-n+1-t2}
and (\ref{eq:Pn+1-down}) that
\begin{align}
&\|\na_{x,z}\widetilde{v_1^n}(t)\|_{H^{s-1}(\overline{\cS})}+\|\widetilde{v}_\om^{n+1}(t)\|_{H^s(\overline{\cS})}+\|a^{n+1}(t)\|_{H^{s-\f32}}\nonumber\\
&\quad+\|\na P^{n+1}(t,\cdot,-1)\|_{H^{s-\f32}}\le \cP(\cE^{n+1}_2(t))\le C(\cP_1)\triangleq E_2.\nonumber
\end{align}
This in turn implies that we can take $T_1$ small enough depending on  $h_0, \cP_1, E_2$ so that for $t\in \big[0,\min(T_n, T_1)\big)$,
\beno
\eta^{n+1}(t,x)+1\ge \f {h_0} 2,\quad \eta_1^{n+1}(t,x)+1\ge \f {h_0} 2.
\eeno

By Lemma \ref{lem:a-n+1-t2}, Lemma \ref{lem:energy-relation} and (\ref{eq:Pn+1-down}), we have
\begin{align}
&\|a^{n+1}(t)\|_{H^{s-\f12}}+\|\na P^{n+1}(t,\cdot,-1)\|_{H^{s-\f32}}\le \cP(\cA_2,\cE^{n+1}(t))\le C(\cP_0,\cP_1)\triangleq E_3^1,\nonumber\\
&\|\eta^{n+1}(t)\|_{H^{s+\f12}}\le \cA_2\cE^{n+1}(t)\le \cA_2\cP_0\triangleq E_3^2.\nonumber
\end{align}
Then we deduce from Proposition \ref{prop:vorticity-n+1}, Proposition \ref{prop:VB-n+1} and Proposition \ref{prop:VB1-n+1} that
\begin{align*}
\cE^{n+1}_2(t)\le& \Big(\cE^{n+1}(0)+\cA_3\int_0^t\cP(\cE^{n+1}(t'))dt'\Big)e^{\cA_1t+C\int_0^t\cE^{n+1}(t')dt'}\\
\le& C_1\big(\cE_0+tC(\cP_0,\cP_1)\big)e^{tC(\cP_0,\cP_1)},
\end{align*}
which ensures that we can take $\cP_1=2C_1\cE_0$ and $T_2$ small enough depending on $\cE_0, \cP_0$ so that
for $t\in \big[0,\min(T_n, T_1, T_2)\big)$,
\beno
\cE^{n+1}_2(t)\le \f {\cP_1} 2.
\eeno

By Lemma \ref{lem:VB-n+1-Dt}, Lemma \ref{lem:eta-n+1-t2}, Lemma \ref{lem:vorticity-n+1-Dt} and (\ref{eq:vom-n+1-t1-up}), we have
\begin{align*}
&\|D_t(V^{n+1}_1,B^{n+1}_1)(t)\|_{H^{s-\f32}}+\|\pa_t\eta^{n+1}(t)\|_{H^{s-\f12}}+\|\widetilde{D}_t^u\eta_1^{n+1}(t)\|_{H^{s-\f12}}\\
&\quad+\|\widetilde{{\cD}_t{\om}^{n+1}}(t)\|_{H^{s-1}(\overline{\cS})}
+\|\widetilde{{\cD}_tv_\om^{n+1}}(t)\|_{X^{s-\f12}([-\f12,0])}\\
&\le \cP(\cA_2,\cE^{n+1}_2(t))\le C(\cP_1)\triangleq E_4.
\end{align*}

By Lemma \ref{lem:VB-n+1-Dt}, Lemma \ref{lem:eta-n+1-t2}, Lemma \ref{lem:a-n+1-t2}, (\ref{eq:v1-n+1-Dt}) and (\ref{eq:Pn+1-t1}), we obtain
\begin{align*}
&\|\big(D_t(V^{n+1}_1, B^{n+1}_1), D_t^bV^{n+1}_{b,1}\big)(t)\|_{H^{s-\f12}}+\|\pa_t^2\eta^{n+1}(t)\|_{H^{s-\f52}}+\|\pa_ta^{n+1}(t)\|_{H^{s-\f32}}\\
&+\|\widetilde{\cD_t^2\om^n}(t)\|_{H^{s-1}(\overline{\cS})}+\|\na_{x,z}\widetilde{{\cD}_tv^{n+1}_1}(t)\|_{H^{s-1}(\overline{\cS})}+\|\na\pa_tP^{n+1}(t,\cdot,-1)\|_{H^{s-\f32}}\\
&\quad+\|(\widetilde{D}_t^u)^2\eta_1^{n+1}(t)\|_{H^{s-\f12}}\le \cP(\cA_4,\cE^{n+1}(t))\le C(\cP_0,\cP_1)\triangleq E_5.
\end{align*}
While, (H5) and the initial conditions (\ref{ass:initial-P}) ensure that there exits $T_3$ small enough depending on $E_5$ so that (H7) holds
for $t\in \big[0,\min(T_n, T_1, T_2, T_3)\big)$.

By Lemma \ref{lem:VB-n+1-Dt}, Lemma \ref{lem:vorticity-n+1-Dt} and (\ref{eq:vom-n+1-Dt2}), we have
\begin{align*}
&\|\pa_t\big(D_t(V^{n+1}_1, B^{n+1}_1),D_t^bV^{n+1}_{b,1}\big)(t)\|_{H^{s-\f32}}
+\|\widetilde{\cD_t^2v^{n+1}_\om}(t)\|_{H^{s-1}(\overline{\cS})}\\
&\le \cP(\cA_5, \cE^{n+1}(t))\le \cP(\cA_5, \cP_0)\triangleq E_6.
\end{align*}

To complete the proof, it remains to prove the first inequality of (\ref{ass:E-n+1-improve}). We know from Proposition \ref{prop:U-n+1-Hs} that
\begin{align*}
\f12\f d {dt}\big(\|D_t \zeta^{n+1}\|_{H^{s-1}}^2+&\|T_{\sqrt{a^n\lambda^n}}\zeta^{n+1}\|_{H^{s-1}}^2\big)\le \cA_4\big(\|D_t \zeta^{n+1}\|_{H^{s-1}}^2+\|\zeta\|_{H^{s-\f12}}^2\big)\\
&+\big\langle (f_1^n+f_2^n+D_tf_3^n+D_tf_4^n+f_\om^n)_{s-1}, (D_t \zeta^{n+1})_{s-1}\big\rangle.
\end{align*}
By Lemma \ref{lem:f1-n}, we have
\beno
\big\langle (f_1^n)_{s-1}, (D_t \zeta^{n+1})_{s-1}\big\rangle \le \cA_5\|D_t\zeta^{n+1}\|_{H^{s-1}}.
\eeno
We write
\begin{align*}
\big\langle (f_2^n)_{s-1}, (D_t \zeta^{n+1})_{s-1}\big\rangle=&\big\langle (f_2^n)_{s-1}, [\langle D\rangle^s, D_t]\zeta^{n+1})\big\rangle
+\big\langle (D_t)^*(f_2^n)_{s-1}, (\zeta^{n+1})_{s-1}\big\rangle\\
&+\f d {dt}\big\langle (f_2^n)_{s-1}, (\zeta^{n+1})_{s-1}\big\rangle,
\end{align*}
from which, Lemma \ref{lem:commu-Ds}, Lemma \ref{lem:f2-n} and Lemma \ref{lem:energy-relation}, we deduce that
\begin{align*}
\big\langle (f_2^n)_{s-1}, (D_t \zeta^{n+1})_{s-1}\big\rangle=&\f d {dt}\big\langle (f_2^n)_{s-1}, (\zeta^{n+1})_{s-1}\big\rangle\\
&+\cA_2\big(\|f_2^n\|_{H^{s-\f32}}+\|D_tf_2^n\|_{H^{s-\f32}}\big)\|\zeta^{n+1}\|_{H^{s-\f12}}\\
=&\f d {dt}\big\langle (f_2^n)_{s-1}, (\zeta^{n+1})_{s-1}\big\rangle+\cA_5\cE^{n+1}(t).
\end{align*}
In a similar way, we deduce from Lemma \ref{lem:fom-n} and Lemma \ref{lem:f3-n} that
\begin{align*}
&\big\langle (D_tf_\om^n)_{s-1}, (D_t \zeta^{n+1})_{s-1}\big\rangle
=\f d {dt}\big\langle(D_t^uf_\om^n)_{s-1}, (  \zeta^{n+1})_{s-1}\big\rangle+\cA_6\cE^{n+1}(t),\\
&\big\langle (D_tf_4^n)_{s-1}, (D_t \zeta^{n+1})_{s-1}\big\rangle=\f d {dt}\big\langle (D_t^uf_4^n)_{s-1}, (\zeta^{n+1})_{s-1}\big\rangle+\cA_6\cE^{n+1}(t).
\end{align*}

We write
\begin{align*}
\big\langle (D_tf_3^n)_{s-1}, (D_t\zeta^{n+1})_{s-1}&\big\rangle=\big\langle[\langle D\rangle^{s-1}, D_t]f_3^n, (D_t\zeta^{n+1})_{s-1})\big\rangle
\\&+\big\langle (f_3^n)_{s-1}, (D_t)^*(D_t\zeta^{n+1})_{s-1}\big\rangle+\f d {dt}\big\langle (f_3^n)_{s-1}, (D_t\zeta^{n+1})_{s-1}\big\rangle,
\end{align*}
which along with Lemma \ref{lem:commu-Ds} and Lemma \ref{lem:f3-n} gives
\begin{align*}
\big\langle (D_tf_3^n)_{s-1}, (D_t \zeta^{n+1})_{s-1}\big\rangle=&\f d {dt}\big\langle (f_3^n)_{s-1}, (D_t\zeta^{n+1})_{s-1}\big\rangle\\
&+\cA_2\|f_3^n\|_{H^{s-\f12}}\big(\|D_t\zeta^{n+1}\|_{H^{s-1}}+\|D_t^2\zeta^{n+1}\|_{H^{s-\f32}}\big)\\
=&\f d {dt}\big\langle (f_3^n)_{s-1}, (D_t \zeta^{n+1})_{s-1}\big\rangle+\cA_5\cE^{n+1}(t)+\cA_5.
\end{align*}

Putting the above estimates together, we obtain
\begin{align*}
&\f12\f d {dt}\big(\|D_t \zeta^{n+1}\|_{H^{s-1}}^2+\|T_{\sqrt{a^n\lambda^n}}\zeta^{n+1}\|_{H^{s-1}}^2\big)\le \cA_5\cE^{n+1}(t)^2+\cA_6\\
&\qquad+\f d {dt}\big\langle(f_2^n)_{s-1}+(f_3^n)_{s-1}+(D_t^uf_4^n)_{s-1}+(D_t^uf_\om^n)_{s-1}, (\zeta^{n+1})_{s-1}\big\rangle.
\end{align*}
Integrating in $t$ and using Lemma \ref{lem:f2-n}-Lemma \ref{lem:fom-n}, we deduce that
\begin{align*}
&\|D_t \zeta^{n+1}\|_{H^{s-1}}^2+\|T_{\sqrt{a^n\lambda^n}}\zeta^{n+1}\|_{H^{s-1}}^2
\le \cP(\cE_0)+\cA_6\int_0^t(1+\cE^{n+1}(t')^2)dt'\\
&\qquad+\cP(E_1,E_2,E_4)+\cA_1\|\eta^n\|_{H^s}^2+\epsilon\|\zeta^{n+1}\|^2_{H^{s-\f12}}\\
&\le \cP(\cE_0)+\cP(E_1,E_2,E_4)+\cA_6(1+\cP_0)^2t+\f1 {64}\cP_0^2+\f14\|T_{\sqrt{a^n\lambda^n}}\zeta^{n+1}\|_{H^{s-1}}^2.
\end{align*}
This implies that there exists $T_4>0$ depending on $\cP_0,\cP_1$ so that for $t\in [0,\min(T, T_n)]$ with $T=\min(T_1,T_2,T_3,T_4)$,
\beno
\cE^{n+1}(t)\le \cP(\cE_0)+\cP(E_1,E_2,E_4)+\cP_1+\f1 {4}\cP_0.
\eeno
Since $E_1,E_2,E_3$ are independent of $\cP_0$, we can take $\cP_0=4\big(\cP(\cE_0)+\cP(E_1,E_2,E_4)+\cP_1\big)$.
Thus, we deduce that for $t\in [0,\min(T, T_n)]$,
\beno
\cE^{n+1}(t)\le \f {\cP_0} 2.
\eeno

This completes the proof of uniform estimates.\ef

\section{Cauchy sequence and the limit system}

This section is devoted to showing that the approximate sequence constructed in last section is a Cauchy sequence.

\subsection{Set-up}
According to the uniform estimates of the approximate sequence, we may assume that
\begin{align*}
\sum_{k=0}^2\Big(&\|\pa_t^k(V^n,B^n,V^n_b,V^n_1,B^n_1,V^n_{b,1})\|_{H^{s-\f12-k}}+\|\pa_t^k\eta^n\|_{H^{s+\f12-k}}+\|\pa_t^k\eta^n_1\|_{H^{s-\f12-k}}\\
&+\|\pa_t^k\widetilde{\om^n}\|_{H^{s-1-k}(\overline{\cS})}\Big)\le \cP_0.
\end{align*}
Here and in what follows we denote by $\cP_0$ a constant depending only on $\|v_0\|_{H^s(\Om_0)}$, $\|\eta_0\|_{H^{s+\f12}}$ and $h_0, c_0$,
which may change from line to line. For a function $f(x,y)$ defined on $\big\{(x,y): -1<y<\eta(x)\big\}$,
we denote $\widetilde{f}(x,z)\triangleq f(x,\rho_\delta(x,z))$, where $\rho_{\delta}(x,z)=z+(1+z)e^{\delta z|D|}\eta$.

We introduce
\begin{align*}
\delta_{\cE^n}(t)\eqdef \delta_{\cE^n_1}(t)+\delta_{\cE^n_2}(t),
\end{align*}
where $\delta_{\cE^n_1}(t)$ and $\delta_{\cE^n_2}(t)$ are given by
\begin{align*}
\delta_{\cE^n_1}(t)=&\|D_t\delta_{\zeta^n}(t)\|_{H^{s-2}}+\|T_{\sqrt{a^n\lambda^n}}\delta_{\zeta^n}(t)\|_{H^{s-2}},\\
\delta_{\cE^n_2}(t)=&\|\big(\delta_{V^n},\delta_{B^n},\delta_{V^n_b},\delta_{V^n_1},\delta_{B^n_1},\delta_{V^n_{b,1}},\delta_{\eta^n_1},\delta_{\eta^n}\big)(t)\|_{H^{s-\f32}}
+\|\delta_{\widetilde{\om^n}}(t)\|_{H^{s-2}(\overline{\cS})}.
\end{align*}
Here $D_t=\pa_t+T_{V^{n+1}}\cdot\na$ and we denote $\delta_{f^n}\triangleq f^{n+1}-f^n$.

Throughout this section, we denote by $L^n_i(i=1,2)$  some nonlinear terms, which satisfy
\beno
\|L^n_i(t)\|_{H^{s-\f12-i}}\le \cP_0\big(\delta_{\cE^n}(t)+\delta_{\cE^{n-1}}(t)\big).
\eeno
Using Lemma \ref{lem:remaider} and Proposition \ref{prop:symbolic calculus}, it is easy to find that
\ben\label{eq:Cau-VB}
\left\{
\begin{aligned}
&D_t\big(\delta_{V^{n}},\delta_{B^{n}}\big)=L^n_1+L(\delta_{a^{n-1}}),\\
&(\pa_t+V_b^n\cdot\na)\delta_{V_b^{n}}=-\delta_{\na P^{n-1}|_{y=-1}},\\
&D_t\delta_{\zeta^{n}}=L^n_2+L^n_3+\delta_{F_2},\\
&(\pa_t+T_{V^{n+1}}\cdot\na)(\delta_{V_1^n},\delta_{B_1^n})=L^n_1+L(\delta_{a^{n-1}}),\\
&(\pa_t+V_b^{n+1}\cdot\na)\delta_{V_{b,1}^{n}}=-\delta_{\na P^{n-1}|_{y=-1}},\\
&(\pa_t+{V^{n+1}}\cdot\na)\delta_{\eta_1^n}=L^n_1,
\end{aligned}\right.
\een
where $L(\delta_{a^{n-1}})$ is a nonlinear term satisfying
\beno
\|L(\delta_{a^{n-1}})\|_{H^{s-\f32}}\le C\|\delta_{a^{n-1}}\|_{H^{s-\f32}},
\eeno
and $L^n_3$ is given by
\beno
L^n_3=(R(\eta^n_1)-R(\eta^{n-1}_1))V^{n}_1+\zeta^n(R(\eta^n_1)-R(\eta^{n-1}_1))B^{n}_1+(R^n_\om-R^{n-1}_\om).
\eeno

\subsection{Elliptic estimates with a parameter}
In order to compare the solution of the elliptic equations in two different domains and with different boundary values, we consider the following elliptic equation
in a domain $\Om_\tau(t)=\big\{(x,y): -1<y<\eta(\tau,t,x), x\in \R^d\big\}$ with a parameter $\tau\in [0,1]$:
\ben\label{eq:elliptic-para}
\left\{
\begin{array}{l}
\Delta_{x,y}v(\tau,t,x,y)=F(\tau,t,x)\quad \textrm{in} \quad \Om_\tau,\\
v|_{y=\eta(\tau,t,x)}=V(\tau,t,x),\quad v|_{y=-1}=V^b(\tau,t,x).
\end{array}
\right.
\een

In the sequel, we denote $f_\tau\triangleq\pa_\tau f(\tau,t,x,y)$. Then $v_\tau$ satisfies
\beno
\left\{
\begin{array}{l}
\Delta_{x,y}v_\tau=F_\tau\quad \textrm{in} \quad \Om_\tau,\\
v_\tau|_{y=\eta(\tau,t,x)}=V_\tau-\pa_yv\eta_\tau,\quad v_\tau|_{y=-1}=V^b_\tau.
\end{array}
\right.
\eeno
If $(\eta, V, V^b)$ and $F$ satisfy
\beno
&&\|\eta\|_{H^{s-\f12}}+\|V\|_{H^{s-\f12}}+\|V^b_\tau\|_{H^{s-\f12}}+\|\widetilde{F_\tau}\|_{H^{s-\f12}(\overline{\cS})}\le C,\\
&&\|\eta_\tau\|_{H^{s-\f32}}+\|V_\tau\|_{H^{s-\f32}}+\|V^b_\tau\|_{H^{s-\f32}}+\|\widetilde{F_\tau}\|_{H^{s-\f32}(\overline{\cS})}\le C,\\
&&\eta(\tau,t,x)+1\ge \f {h_0} 2\quad \textrm{for} \quad (\tau,t,x)\in [0,1]\times [0,T]\times \R^d,
\eeno
then we infer from Proposition \ref{prop:elliptic-tan-nontame} and Proposition \ref{prop:elliptic-boun} that for any $\sigma\in [-\f12,s-\f52]$,
\begin{align}\label{eq:Cau-ellip1}
&\|\na_{x,z}\widetilde{v_\tau}\|_{X^{\sigma}([-\f12,0])}\le C\big(\|\widetilde{F_\tau}\|_{L^2(H^{\sigma-\f12}(\overline{\cS}))}+\|(V_\tau, \eta_\tau)\|_{H^{\sigma+1}}+\|V_\tau^b\|_{H^\f12}\big),\\
&\|\na_{x,z}\widetilde{v_\tau}\|_{H^{\sigma+\f12}(\overline{\cS})}\le C\big(\|\widetilde{F_\tau}\|_{H^{\sigma-\f12}(\overline{\cS})}+\|(V_\tau, V^b_\tau, \eta_\tau)\|_{H^{\sigma+1}}\big).\label{eq:Cau-ellip2}
\end{align}

Now we take in the elliptic equation (\ref{eq:elliptic-para}):
\beno
&&\eta=\tau\eta^{n+1}_1+(1-\tau)\eta^{n}_1,\\
&&V=\tau(V^{n+1},B^{n+1})+(1-\tau)(V^{n},B^{n}),\\
&&V^b=\tau(V^{n+1}_b,0)+(1-\tau)(V^{n}_b,0),\\
&&F=\na_{x,y}\cdot\big(\tau\widetilde{\om^{n+1}}\circ\Phi^{-1}+(1-\tau)\widetilde{\om^n}\circ\Phi^{-1}\big),
\eeno
where $\Phi(\tau,x,z)=(x,\rho_\delta(\tau,x,z))$ with $\rho_\delta(\tau,x,z)=z+(1+z)e^{\delta z|D|}\eta(\tau,t,\cdot)$. Notice that
\beno
&&\delta_{\widetilde{v^n}}=\widetilde{v^{n+1}}-\widetilde{v^{n}}=\int_0^1\widetilde{v_\tau}d\tau,\quad \eta_\tau=\delta_{\eta^n},\quad
V_\tau=\delta_{(V^n,B^n)},\quad V_\tau^b=\delta_{(V^n_b,0)},\\
&&\|\widetilde{F_\tau}\|_{H^{s-3}(\overline{\cS})}\le \cP_0\delta_{\cE_2^n}(t).
\eeno
Then we can deduce from (\ref{eq:Cau-ellip2}) that
\ben\label{eq:Cau-v}
\|\na_{x,z}\delta_{\widetilde{v^n}}\|_{H^{s-2}(\overline{\cS})}\le \cP_0\delta_{\cE_2^n}(t).
\een
In a similar way, we can deduce that
\ben\label{eq:Cau-v1}
&&\|\na_{x,z}\delta_{\widetilde{v^n_1}}\|_{H^{s-2}(\overline{\cS})}\le \cP_0\delta_{\cE^n_2}(t).
\een
Let $\om(\tau,t,x)=\tau\widetilde{\om^{n+1}}\circ\Phi^{-1}+(1-\tau)\widetilde{\om^n}\circ\Phi^{-1}$ and
\beno
\widetilde{\cD}_t=\pa_t+(\widetilde{v_2^{n+1}}\circ\Phi^{-1}-\pa_t\Phi^{-1})\cdot\cA\na_{x,y}
\eeno
with $\cA=\big(\pa_i\Phi^{-1}_j\big)^{-1}$. Then $\om_\tau$ satisfies
\ben\label{eq:Cau-om-t1}
\|(\widetilde{\cD}_t\om_\tau)\circ\Phi\|_{H^{s-2}(\overline{\cS})}\le \cP_0\delta_{\cE^n_2}(t).
\een
This in turn implies that
\ben\label{eq:Cau-v1-t1}
\|\delta_{\widetilde{\widetilde{\cD}_tv^n_1}}\|_{H^{s-1}(\overline{\cS})}+\|(\widetilde{\cD}_t^2\om_\tau)\circ\Phi\|_{H^{s-2}(\overline{\cS})}
\le \cP_0\delta_{\cE^n}(t).
\een
Similarly, there holds for $\delta_{v^{n-1}_\om}$,
\ben\label{eq:Cau-vom}
&&\|\delta_{\widetilde{v^{n-1}_\om}}\|_{H^{s-1}(\overline{\cS})}+\|\delta_{\widetilde{\widetilde{\cD}_tv^{n-1}_\om}}\|_{X^{s-\f32}([-\f12,0])}\le \cP_0\delta_{\cE^{n-1}_2}(t),\\
&&\|\delta_{\widetilde{\widetilde{\cD}_t^2v^{n-1}_\om}}\|_{X^{s-\f32}([-\f12,0])}\le \cP_0\delta_{\cE^{n-1}}(t).\label{eq:Cau-vom-t2}
\een

To compare the pressure, we consider
\beno
\left\{
\begin{array}{l}
-\Delta_{x,y}P(\tau,t,x,y)=F(\tau,t,x)\quad \textrm{in} \quad \Om_\tau,\\
P|_{y=\eta(\tau,t,x)}=0,\quad \pa_y P|_{y=-1}=-1.
\end{array}
\right.
\eeno
where $\eta(\tau,t,x)=\tau\eta^{n}+(1-\tau)\eta^{n-1}$ and
\beno
F=\big(\big(\pa_i( {v_1}^j)^{n}\pa_j ( {v_1}^i)^{n}\big)\circ \Phi^{n}_1-\big(\pa_i( {v_1}^j)^{n-1}\pa_j ( {v_1}^i)^{n-1}\big)\circ \Phi^{n-1}_1\big)\circ \Phi^{-1}.
\eeno
Then $P_\tau$ satisfies
\beno
\left\{
\begin{array}{l}
-\Delta_{x,y}P_\tau=F_\tau\quad \textrm{in} \quad \Om_\tau,\\
P_\tau|_{y=\eta(\tau,t,x)}=-\pa_y P\eta_\tau,\quad \pa_y P_\tau|_{y=-1}=0.
\end{array}
\right.
\eeno
Thus, we can deduce by a similar proof of Lemma \ref{lem:a-n+1-t2} that
\begin{align}
&\|\delta_{a^{n-1}}\|_{H^{s-\f52}}+\|\delta_{\na P^{n-1}|_{y=-1}}\|_{H^{s-\f52}}\le \cP_0\delta_{\cE^{n-1}_2}(t),\label{eq:Cau-a-low}\\
&\|\delta_{a^{n-1}}\|_{H^{s-\f32}}+\|\delta_{\na P^{n-1}|_{y=-1}}\|_{H^{s-\f32}}\le \cP_0\delta_{\cE^{n-1}}(t),\label{eq:Cau-a-high}\\
&\|\delta_{\na\pa_tP^{n-1}|_{y=-1}}\|_{H^{s-\f52}}+\|\delta_{\pa_ta^{n-1}}\|_{H^{s-\f52}}\le \cP_0\delta_{\cE^{n-1}}(t).\label{eq:Cau-a-t1}
\end{align}

To compare the remainder of (DN) operators, we take in (\ref{eq:elliptic-para}):
\beno
\eta=\tau\eta^{n}_1+(1-\tau)\eta^{n-1}_1,\quad V=(V^{n}_1,B^{n}_1),\quad V^b=0,\quad F=0.
\eeno
Then $v_\tau$ satisfies
\beno
\left\{
\begin{array}{l}
\Delta_{x,y}v_\tau=0\quad \textrm{in} \quad \Om_\tau,\\
v_\tau|_{y=\eta(\tau,t,x)}=-\pa_yv\eta_\tau,\quad v_\tau|_{y=-1}=0.
\end{array}
\right.
\eeno
Notice that
\beno
\|D_t^u\eta_\tau\|_{H^{s-\f32}}\le \cP_0\delta_{\cE^{n-1}_2}(t),\quad \|(D_t^u)^2\eta_\tau\|_{H^{s-\f52}}\le \cP_0\delta_{\cE^{n-1}}(t).
\eeno
Here $D_t^u=\pa_t+V^n\cdot\na$.
Thus, we can deduce from a similar proof of Proposition \ref{prop:DN-remainder-Dt} that
\begin{align}
&\|(R(\eta_1^{n})-R(\eta_1^{n-1}))(V_1^{n},B_1^{n})\|_{H^{s-\f52}}\le \cP_0\delta_{\cE^{n-1}_2}(t),\label{eq:Cau-rem}\\
&\|D_t^u(R(\eta_1^{n})-R(\eta_1^{n-1}))(V_1^{n},B_1^{n})\|_{H^{s-\f52}}\le \cP_0\delta_{\cE^{n-1}_2}(t),\label{eq:Cau-rem-t1}\\
&\|(D_t^u)^2(R(\eta_1^{n})-R(\eta_1^{n-1}))(V_1^{n},B_1^{n})\|_{H^{s-\f52}}\le \cP_0\delta_{\cE^{n-1}}(t).\label{eq:Cau-rem-t2}
\end{align}

\subsection{Energy estimates}
Following the proofs in section 11.7 and using the estimates in subsection 12.4, we can deduce that
\begin{align*}
&\|\delta_{f_1^{n-1}}\|_{H^{s-2}}+\|\delta_{f_3^{n-1}}\|_{H^{s-\f32}}\le \cP_0\big(\delta_{\cE^{n}}(t)+\delta_{\cE^{n-1}}(t)\big),\\
&\|\delta_{f_3^{n-1}}\|_{H^{s-2}}\le \cP_0\big(\|\delta_{\eta^{n-1}}\|_{H^{s-1}}+\delta_{\cE^{n-1}_2}(t)\big),\\
&\|\delta_{f_2^{n-1}}\|_{H^{s-\f52}}+\|D_t^u\delta_{f_4^{n-1}}\|_{H^{s-\f52}}+\|D_t^u\delta_{f_\om^{n-1}}\|_{H^{s-\f52}}\le \cP_0\big(\delta_{\cE^{n-1}_2}(t)+\delta_{\cE^{n}_2}(t)\big),\\
&\|D_t\delta_{f_2^{n-1}}\|_{H^{s-\f52}}+\|(D_t^u)^2\delta_{f_4^{n-1}}\|_{H^{s-\f52}}+\|(D_t^u)^2\delta_{f_\om^{n-1}}\|_{H^{s-\f52}}\le \cP_0\big(\delta_{\cE^{n}}(t)+\delta_{\cE^{n-1}}(t)\big).
\end{align*}
By (\ref{eq:Cau-rem}) and (\ref{eq:Cau-vom}), we have
\beno
\|L_3^n\|_{H^{s-\f52}}\le \cP_0\delta_{\cE^{n-1}_2}(t).
\eeno
We also have
\ben\label{eq:Cau-eta}
D_t^2\delta_{\zeta^{n}}+ T_{a^n\lambda^n}\delta_{\zeta^{n}}=L^n_2+L^{n}_4,
\een
where $L^n_4$ is given by
\beno
L^n_4=(f_2^n+D_tf_3^n+D_tf_4^n+D_tf_\om^n)-(f_2^{n-1}+D_tf_3^{n-1}+D_tf_4^{n-1}+D_tf_\om^{n-1}).
\eeno
For $\delta_{\widetilde{\om^n}}$, we have
\beno
\|(\pa_t+\overline{v}^{n+1}\cdot\na_{x,z})\delta_{\widetilde{\om^n}}\|_{H^{s-2}(\overline{\cS})}\le \cP_0\delta_{\cE^{n-1}_2}(t).
\eeno
Making the energy estimates for the system (\ref{eq:Cau-VB}), we obtain
\beno
\delta_{\cE^n_2}(t)^2\le \cP_0\int_0^t\big(\delta_{\cE^{n}}(t')+\delta_{\cE^{n-1}}(t')\big)^2dt'.
\eeno
While, making the energy estimates for the system (\ref{eq:Cau-eta}), we obtain
\begin{align*}
\delta_{\cE^n_1}(t)^2\le& \cP_0\int_0^t\big(\delta_{\cE^{n}}(t')+\delta_{\cE^{n-1}}(t')\big)^2dt'+\f12\delta_{\cE^{n-1}}(t)^2\\
&+\cP_0\big(\delta_{\cE^{n}_2}(t)+\delta_{\cE^{n-1}_2}(t)\big)^2.
\end{align*}
Then an induction argument ensures that there exists $T>0$ depending only on $\cP_0$ so that for $t\in [0,T]$,
\beno
\delta_{\cE^n}(t)\longrightarrow 0\quad \textrm{as}\quad n\longrightarrow+\infty.
\eeno
This shows that the approximate sequence is a Cauchy sequence.

\subsection{The limit system}
Let $\big(V, B, V_b, V_1, B_1, V_{b,1}, \zeta, \eta_1, \eta, \om, v, v_1, P\big)$ be the limit of the Cauchy sequence
\beno
\big(V^n, B^n, V_b^n, V_1^n, B_1^n, V_{b,1}^n, \zeta^n, \eta^n, \eta_1^n, \om^n, v^n, v_1^n, P^n\big).
\eeno
Taking the limit for the approximate system (\ref{eq:iteration-velocity})--(\ref{eq:iteration-pressure}), we obtain the following limit system:
The boundary velocity $(V,B, V_b)$ and $\zeta$ satisfy
\begin{equation}\label{eq:iteration-velocity-limit}
\left\{
\begin{aligned}
&D_tV=T_{\zeta}a+T_{a}\zeta+R(\zeta, a)+(T_{V_1}-V_1)\cdot\nabla V,\\
&D_tB= a-1+(T_{V_1}-V_1)\cdot\nabla B,\\
&(\pa_t+V_b\cdot\na)V_b=-\na P|_{y=-1},\\
&D_t\zeta=T_{\lambda}(V+T_{\zeta}B)+(T_{V_1}-V_1)\cdot\zeta+[T_{\zeta},T_{\lam}]B_1\\
&\qquad+(\zeta-T_{\zeta})T_{\lam}B+R(\eta_1)V_1+\zeta R(\eta_1)B_1+R_\om,
\end{aligned}\right.
\end{equation}
where $D_t=\pa_t+T_V\cdot\na, a=-\pa_yP|_{y=\eta},\lambda=\lambda(\eta_1)$, and
\begin{align*}
(R_\om)^i=&\big(\pa_y (v_{\om})^i-\pa_{x_j}(v_{\om})^i\cdot \pa_{x_j}\eta_1\big)+\pa_{x_i}\eta_1\big(\pa_y (v_{\om})^{d+1}
-\pa_{x_j}\eta_1\pa_{x_j}(v_{\om})^{d+1}\big)\\
&+\big(\om_{i,d+1}-\pa_{x_j}\eta_1\om_{ij}+\pa_{x_i}\eta_1\pa_{x_j}\eta_1\om_{j,d+1}\big)\big|_{y=\eta_1}
\end{align*}
with $v_\om$ given by
\beno
\left\{
\begin{array}{ll}
-\Delta_{x, y} v_{\om}=\na_{x,y}\times\om\quad \textrm{in} \quad \Om_t=\big\{(x,y):-1<y<\eta(t,x)\big\},\\
v_{\om}|_{y=\eta_1}=0,\quad v_\om|_{y=-1}=(V_b, 0),
\end{array}\right.
\eeno
and $(V_1, B_1, V_{b,1})$ satisfies
\begin{align}\label{eq:iteration-velocity-new-limit}
\left\{
\begin{array}{l}
D_t\big(V_1, B_1\big)=D_t\big(V, B\big),\\
\big(\pa_t+V_{b}\cdot\na\big)V_{b,1}=-\na P|_{y=-1}.
\end{array}\right.
\end{align}

The free surface $(\eta,\eta_1)$ satisfies
\ben
&&-\Delta\eta+\eta=-\dive\zeta+ {\eta}_1,\\
&&(\pa_t+V\cdot\nabla) {\eta}_1=B_1.
\een

The vorticity $\om$ satisfies
\begin{equation}\label{eq:iterative-vorticity-limit}
\begin{aligned}
\pa_t {\om}+ \big(v^h\cdot\na+v_1^{d+1}\pa_y\big){\om}= {\om}\cdot\nabla_{x,y}{v}_1,
\end{aligned}
\end{equation}
where the velocity $(v, v_1)$  is given by
\begin{equation}\label{eq:iterative-velocity-limit}
\left\{
\begin{aligned}
&-\Delta_{x,y}{v}=\nabla_{x,y}\times{\om}\quad\textrm{in}\quad\Om_t,\\
& {v}\mid_{y= {\eta}_1}=(V,B),\quad v|_{y=-1}= (V_{b},0),
\end{aligned}
\right.
\end{equation}
and
\begin{equation}\label{eq:iterative-velocity-new-limit}
\left\{
\begin{aligned}
&-\Delta_{x,y}{v}_1=\nabla_{x,y}\times{\om}\quad\textrm{in}\quad\Om_t,\\
& {v}_1\mid_{y= {\eta}_1}=(V_1,B_1),\quad v_1|_{\Gamma_b}=(V_{b,1},0).
\end{aligned}
\right.
\end{equation}

Let $\widetilde{\Om}_t=\big\{(x,y):y=\eta(t,x)\big\}$. The pressure $P$ satisfies
\begin{equation}\label{eq:iteration-pressure-limit}
\left\{
\begin{aligned}
&-\Delta_{x,y}P=\big(\pa_i( {v_1}^j)\pa_j ( {v_1}^i)\big)\circ \Phi_1\circ(\Phi)^{-1}\quad \textrm{in} \quad \widetilde{\Om}_t,\\
&P\mid_{y=\eta}=0,\quad \pa_yP\mid_{y=-1}=0,
\end{aligned}
\right.
\end{equation}
where $\Phi(x,z)=(x,\rho_{\delta,\eta}(x,z))$ and $\Phi_1(x,z)=(x,\rho_{\delta,\eta_1}(x,z))$
with $\rho_{\delta,\eta}(x,z)=z+(1+z)e^{\delta z|D|}\eta$.

\section{From the limit system to the Euler equations}

The goal of this section is to  show that the limit system (\ref{eq:iteration-velocity-limit})--(\ref{eq:iteration-pressure-limit}) is equivalent to
the Euler equations (\ref{eq:euler})--(\ref{eq:euler-p}). 

First of all, it follows from (\ref{eq:iteration-velocity-new-limit}) and the third equation of (\ref{eq:iteration-velocity}) that
\beno
(V,B,V_b)=(V_1, B_1, V_{b,1}).
\eeno
Hence, $v=v_1$ since they satisfy the same elliptic equation with the same boundary conditions.
Thus, we deduce from (\ref{eq:iteration-velocity-limit}) that
\begin{equation}\label{eq:VB-limit}
\left\{
\begin{aligned}
&(\pa_t +V\cdot\nabla) V=a\zeta,\quad a=\na P|_{y=\eta}, \\
&(\pa_t +V\cdot\nabla) B =a-1,\\
&(\pa_t +V_b\cdot\nabla) V_b=-\na P|_{y=-1},\\
&(\pa_t+V\cdot\nabla)\zeta =G(\eta_1)V+ {\zeta}G(\eta_1)B+R_\om.
\end{aligned}\right.
\end{equation}
For the free surface, we have
\ben\label{eq:eta1-limit}
(\pa_t+V\cdot\nabla) {\eta}_1=B.
\een
The vorticity $\om$ satisfies
\begin{equation}\label{eq:vorticity-limit}
\left\{
\begin{aligned}
&\pa_t {\om}+ {v}\cdot\nabla_{x,y} {\om}= {\om}\cdot\nabla_{x,y}{v}\quad\textrm{in}\quad\Om_t,\\
&-\Delta_{x,y}{v}=\nabla_{x,y}\times{\om}\quad\textrm{in}\quad\Om_t,\\
&{v}\mid_{y= {\eta}_1}=(V,B),\quad{v}\mid_{y=-1}=(V_b,0).
\end{aligned}
\right.
\end{equation}

It remains to show that
\ben\label{eq:eta-om-div}
\eta=\eta_1, \quad \om=\na_{x,y}\times v,\quad \dv v=0.
\een
For this end, we introduce
\beno
&&G=v_t+v\cdot\na_{x,y} v+(\na_{x,y} (P+gy))\circ \overline{\Phi},\\
&&\delta_\om=\om-\overline{\om},\quad \delta_d=\dv v,\quad \delta_\zeta=\zeta-\na\eta_1.
\eeno
where $\overline{\om}=\na_{x,y}\times v$ and $\overline{\Phi}=\Phi\circ(\Phi_1)^{-1}$.

In what follows, we denote by $L(\cdot)$ a linear function, which may be different from line to line.

\begin{lemma}\label{lem:Phi}
It holds that
\beno
\overline{\Phi}(x,y)=(x, h(x,y)),
\eeno
where $h(x,y)$ satisfies
\beno
\|\na_{x,y}(h-y)\|_{H^1({\Om}_t)}\le C\|\eta-\eta_1\|_{H^\f 32}.
\eeno

\end{lemma}

\no{\bf Proof.}\,Let $\Phi_1^{-1}(x,y)=(x, z(x,y))$, i.e.,
\beno
y=z(x,y)+(1+z(x,y))e^{\delta z(x,y)|D|}\eta_1.
\eeno
Hence, we have
\beno
\overline{\Phi}(x,y)=\big(x, z(x,y)+(1+z(x,y))e^{\delta z(x,y)|D|}\eta\big)\triangleq (x, h(x,y)).
\eeno
Using the fact that
\beno
&&\na z(x,y)(1+e^{\delta z(x,y)|D|}\eta_1)+(1+z(x,y))(\na z(x,y)e^{\delta z(x,y)|D|}|D|\eta_1+e^{\delta z(x,y)|D|}\na\eta_1)=0,\\
&&\pa_y z(x,y)(1+e^{\delta z(x,y)|D|}\eta_1)+(1+z(x,y))\pa_y z(x,y)e^{\delta z(x,y)|D|}|D|\eta_1=1,
\eeno
we deduce that
\begin{align*}
&\na_{x,y}(h(x,y)-y)=L\big(e^{\delta z(x,y)|D|}(\eta-\eta_1), e^{\delta z(x,y)|D|}\na(\eta-\eta_1)\big),
\end{align*}
which implies that
\begin{align*}
\|\na_{x,y}(h-y)\|_{H^1(\Om_t)}\le& C\|e^{\delta z(x,y)|D|}(\eta-\eta_1)\|_{H^1(\Om_t)}+C\|e^{\delta z(x,y)|D|}\na(\eta-\eta_1)\|_{H^1(\Om_t)}\\
\le& C\|e^{\delta z|D|}(\eta-\eta_1)\|_{L^2_z(-1,0;H^2)}\le C\|\eta-\eta_1\|_{H^\f32}.
\end{align*}
The proof is finished.\ef\medskip

Let us first derive the equation of $G$.

\begin{lemma}\label{lem:G}
It holds that
\begin{align*}
\left\{
\begin{array}{l}
\Delta_{x,y}G=L\big(\delta_d, \na_{x,y}\delta_d, \na_{x,y}^2\delta_d, \delta_\om, \na_{x,y}\delta_\om, \na_{x,y}^2\delta_\om, \na_{x,y}(h-y), \na_{x,y}^2h\big),\\
G|_{y=\eta_1}=0, \quad G|_{y=-1}=0.
\end{array}\right.
\end{align*}
In particular, we have
\beno
\|G\|_{H^2(\Om_t)}\le C\big(\|(\delta_d,\delta_\om)\|_{L^2(\Om_t)}+\|\na_{x,y}(h-y)\|_{H^1({\Om}_t)}\big).
\eeno
\end{lemma}

\no{\bf Proof.}\,Thanks to (\ref{eq:vorticity-limit}), we have
 \begin{align*}
 \Delta_{x,y}v_t=-\nabla_{x,y}\times \om_t=\nabla_{x,y}\times\big(v\cdot\nabla_{x,y}\om-\om\cdot\nabla_{x,y}v\big).
 \end{align*}
 A direct calculation gives
 \begin{align*}
 \Delta_{x,y}\big(v\cdot\na_{x,y}v\big)=&-\na_{x,y}\times \big(\na_{x,y}\times(v\cdot\na v)\big)+\na_{x,y}\dv(v\cdot\na_{x,y}v)\\
 =&-\na_{x,y}\times\big(v\cdot\na_{x,y}\overline{\om}-\overline{\om}\cdot\na_{x,y}v+\overline{\om}\delta_d\big)\\
 &+\na_{x,y}(\pa_jv^i\pa_iv^j)+\na_{x,y}(v\cdot\na_{x,y}\delta_d),
 \end{align*}
 Using (\ref{eq:iteration-pressure-limit}) and Lemma \ref{lem:Phi}, we obtain
 \begin{align*}
&\Delta_{x,y}\big(\na_{x,y}P\circ \overline{\Phi}\big)=-\na_{x,y}(\pa_jv^i\pa_iv^j)+L\big(\na_{x,y}(h-y), \na^2_{x,y}h\big).
 \end{align*}
The first equation of the lemma follows by summing up the above equations. The boundary condition of $G$ follows from
(\ref{eq:VB-limit}).\ef

Now we have
\ben\label{eq:v-P}
v_t+v\cdot\na_{x,y} v+(\na_{x,y} (P+gy))\circ \overline{\Phi}=G.
\een
Take the divergence to get
\beno
\pa_t\delta_d+v\cdot\na_{x,y}\delta_d=\dv_{x,y}G+L(\na_{x,y}(h-y)),
\eeno
which along with Lemma \ref{lem:Phi} and Lemma \ref{lem:G} implies that
\begin{align}
\f d {dt}\|\delta_d\|_{H^1(\Om_t)}^2\le& C\big(\|\delta_d\|_{H^1(\Om_t)}^2+\|G\|_{H^2(\Om_t)}+\|\na_{x,y}(h-y)\|_{H^1(\Om_t)}\big)\nonumber\\
\le& C\big(\|(\delta_d,\delta_\om)\|_{H^1(\Om_t)}^2+\|\eta-\eta_1\|_{H^\f 32}\big).\label{eq:delta-d}
\end{align}
Taking the curl to (\ref{eq:v-P}), we obtain
\beno
\pa_t\overline{\om}+v\cdot\na_{x,y}\overline{\om}=\overline{\om}\cdot\na_{x,y}v+\na_{x,y}\times G+L(\na_{x,y}(h-y)),
\eeno
from which and (\ref{eq:vorticity-limit}), we infer
\beno
\pa_t\delta_\om+v\cdot\na_{x,y}\delta_\om=\delta_\om\cdot\na_{x,y}v+\na_{x,y}\times G+L(\na_{x,y}(h-y)).
\eeno
Then similar to (\ref{eq:delta-d}), we have
\begin{align}
\f d {dt}\|\delta_\om\|_{H^1(\Om_t)}^2
\le C\big(\|(\delta_d,\delta_\om)\|_{H^1(\Om_t)}^2+\|\eta-\eta_1\|_{H^\f 32}\big).\label{eq:delta-om}
\end{align}

Let $\zeta_1=\na\eta_1$. By a similar derivation of (\ref{eq:zeta-new}), we get
\beno
\pa_t\zeta_1+V\cdot\na \zeta_1=G(\eta_1)V+\zeta_1G(\eta_1)B+R_{\overline{\om}},
\eeno
where $R_{\overline{\om}}$ is given by
\begin{align*}
(R_{\overline{\om}})^i=&\big(\pa_y (v_{\om})^i-\pa_{x_j}(v_{\om})^i\cdot \pa_{x_j}\eta_1\big)+\pa_{x_i}\eta_1\big(\pa_y (v_{\om})^{d+1}
-\pa_{x_j}\eta_1\pa_{x_j}(v_{\om})^{d+1}\big)\\
&+\big(\overline{\om}_{i,d+1}-\pa_{x_j}\eta_1\overline{\om}_{ij}+\pa_{x_i}\eta_1\pa_{x_j}\eta_1\overline{\om}_{j,d+1}\big)\big|_{y=\eta_1}.
\end{align*}
Hence, $\delta_\zeta$ satisfies
\beno
\pa_t\delta_\zeta+V\cdot\na \delta_\zeta=\delta_\zeta G(\eta_1)B+R_\om-R_{\overline{\om}}.
\eeno
This ensures that
\ben\label{eq:delta-eta}
\f d {dt}\|\delta_\zeta\|_{H^\f12}^2\le C\big(\|\delta_\zeta\|_{H^\f12}^2+\|\delta_\om\|_{H^1(\Om_t)}^2\big).
\een

On the other hand, we know
\beno
-\Delta (\eta-\eta_1)+(\eta-\eta_1)=\dv(\zeta-\zeta_1)=\dv\delta_\zeta,
\eeno
which implies
\beno
\|\eta-\eta_1\|_{H^\f32}\le C\|\delta_\zeta\|_{H^\f12}.
\eeno
Then we deduce from (\ref{eq:delta-d})--(\ref{eq:delta-eta}) that
\beno
\f d {dt}\big(\|\delta_d\|_{H^1(\Om_t)}^2+\|\delta_\om\|_{H^1(\Om_t)}^2+\|\delta_\zeta\|_{H^\f12}^2\big)
\le C\big(\|\delta_d\|_{H^1(\Om_t)}^2+\|\delta_\om\|_{H^1(\Om_t)}^2+\|\delta_\zeta\|_{H^\f12}^2\big)
\eeno
together with the initial condition
\beno
\delta_d|_{t=0}=0,\quad \delta_\om|_{t=0}=0,\quad \delta_\zeta|_{t=0}=0.
\eeno
Gronwall's inequality implies that $\delta_d=0, \delta_\om=0, \delta_\zeta=0$. This proves (\ref{eq:eta-om-div}).\ef

\section{Proof of Theorem \ref{thm:local}}

This section is devoted to proving the local well-posedness of the system (\ref{eq:euler})--(\ref{eq:euler-p}) for the low regularity initial data.

\subsection{Construction of approximate smooth solution}

First of all, we can construct a sequence of smooth functions $\eta_0^n$ so that
\begin{align*}
&\|\eta_0^n-\eta_0\|_{H^{s+\f12}}\longrightarrow 0\quad \textrm{as} \quad n\longrightarrow+\infty,\\
&\eta^n_0(x)+1\ge \f 34 h_0\quad \textrm{for  } x\in \R^d.
\end{align*}
Let $\widetilde{v_0}(x,z)=v_0\circ \Phi(x,z)$, where $\Phi(x,z)=(x,\rho_\delta(x,z))$ with $\rho_\delta(x,z)=z+(1+z)e^{\delta z|D|}\eta_0$.
Then we take a sequence of smooth functions $\widetilde{v_0^n}$ so that
\beno
\|\widetilde{v_0^n}-\widetilde{v_0}\|_{H^{s+\f12}(\overline{\cS})}\longrightarrow 0\quad \textrm{as} \quad n\longrightarrow+\infty.
\eeno
Let $v_0^n(x,y)=\widetilde{v_0^n}\circ (\Phi^n)^{-1}(x,y)$, where $\Phi^n(x,z)=(x,\rho_\delta^n(x,z))$ with $\rho_\delta^n(x,z)=z+(1+z)e^{\delta z|D|}\eta_0^n$.

The pressure $P_0^n$ associated with the initial data $(v_0^n,\eta_0^n)$ is defined by
\beno
\left\{
\begin{array}{l}
-\Delta_{x,y}P_0^n=\pa_i(v_0^n)^j\pa_j(v_0^n)^i\quad \textrm{in} \quad \Om_0^n,\\
P_0^n|_{y=\eta_0(x)}=0,\quad \pa_y P_0^n|_{y=-1}=-1,
\end{array}
\right.
\eeno
where $\Om_0^n=\big\{(x,y):y=\eta_0^n(x), x\in \R^d\big\}$. Then we can show that
\beno
\|\na_{x,z}(P_0^n\circ\Phi^n-P_0\circ\Phi)\|_{X^{s-\f12}([-\f12,0])}\longrightarrow 0\quad \textrm{as} \quad n\longrightarrow+\infty.
\eeno
This ensures that for $n$ big enough, there holds
\beno
-\pa_yP_0^n\big|_{y=\eta^n_0}\ge \f 34c_0.
\eeno
Thus, we apply Theorem \ref{thm:local-smooth} to obtain a sequence of smooth solutions $(v^n,\eta^n,P^n)$ associated with the initial data $(v_0^n, \eta_0^n)$
on a maximal existence time interval $[0,T_n)$.

\subsection{Uniform estimates and existence}
We denote
\beno
(V^n,B^n)=v^n|_{y=\eta^n(t,x)},\quad V^n_b=(v^n)^h|_{y=-1}.
\eeno
We define
\beno
&&E_{s}^n(t)\eqdef\|(V^n,B^n,V_b^n)(t)\|_{H^{s}}+\|\eta^n(t)\|_{H^{s+\f12}}+\|\widetilde{\om^n}(t)\|_{H^{s-\f12}(\overline{\cS})},\\
&&E_{s,l}^n(t)\eqdef \|(V^n,B^n)(t)\|_{H^{s-\f12}}^2+\|V_b^n(t)\|_{H^s}^2+\|\eta^n(t)\|_{H^s}^2+\|\widetilde{\om^n}(t)\|_{H^{s-\f12}(\overline{\cS})}.
\eeno

The goal of this subsection is to show that there exits $T>0$ depending only on  $E_s(0)$ and $c_0, h_0, s$  such that that for any $t\in [0,\min(T,T_n))$,
there holds
\ben
E_{s}^n(t)\le \cP(E_s(0)).
\een
Here and in what follows we denote by $\cP$ an increasing function depending only on $c_0, h_0, s$, which may change from line to line.
\medskip

First of all, it follows from Proposition \ref{prop:energy-lower} that
\begin{align}
&\f d {dt}\big(\|(V^n,B^n)\|_{H^{s-\f12}}^2+\|V_b^n\|_{H^s}^2+\|\eta^n\|_{H^s}^2\big)\le \cP(E_s^n(t)).\label{eq:thm-VB-low}
\end{align}

Next, we estimate the vorticity $\om^n$, which satisfies
\beno
\pa_t \widetilde{\om^n}+\overline{v^n}\cdot \nabla_{x,z}\widetilde{\om^n}=\widetilde{\om^n}^h\cdot\big(\na\widetilde{v^n}-\f {\na \rho_\delta^n} {\pa_z\rho_\delta^n}\pa_z\widetilde{v^n}\big)
+\widetilde{\om^n}^{d+1}\f {\pa_z \widetilde{v^n}} {\pa_z\rho_\delta^n}\triangleq F^n,
\eeno
where $\overline{v^n}=\big(\widetilde{v^n},\f{1}{\pa_z \rho_\delta^n}(\widetilde{v^n}_{d+1}-\pa_t\rho_\delta^n-\widetilde{v^n}^h\cdot\na\rho_\delta^n)\big)$.
By Lemma \ref{lem:v-trans} and Lemma \ref{lem:product-full}, we have
\beno
\|\overline{v^n}\|_{H^{s+\f12}(\overline{\cS})}+\|{F^n}\|_{H^{s-\f12}(\overline{\cS})}\le \cP(E_s^n(t)).
\eeno
Let $\overline{v^n_e}$ and $F^n_e$ be the extension of $\overline{v^n}$ and $F^n$ to $\R^{d+1}$ so that
\beno
&&\|\overline{v^n_e}\|_{H^{s+\f12}(\R^{d+1})}\le C\|\overline{v^n}\|_{H^{s+\f12}(\overline{\cS})},\\
&&\|F^n_e\|_{H^{s-\f12}(\R^{d+1})}\le C\|{F^n}\|_{H^{s-\f12}(\overline{\cS})}.
\eeno
We define $\widetilde{\om_e^n}$ to be a solution of the following transport equation in $\R^{d+1}$
\beno
\pa_t\widetilde{\om_e^n}+\overline{v_e^n}\cdot\na_{x,y}\widetilde{\om_e^n}=F^n_e,\quad \widetilde{\om_e^n}(0)=\widetilde{\om_{0,e}^n}(x).
\eeno
It is obvious that $\widetilde{\om_e^n}=\widetilde{\om^n}$ in $\overline{\cS}$ by the uniqueness of the solution. By a standard $H^{s-\f12}$ energy estimate, we deduce that
\begin{align*}
\f d {dt}\|\widetilde{\om_e^n}(t)\|_{H^{s-\f12}(\R^{d+1})}\le C\|\overline{v^n_e}\|_{H^{s+\f12}(\R^{d+1})}\|\widetilde{\om_e^n}(t)\|_{H^{s-\f12}(\R^{d+1})}
+\|F^n_e\|_{H^{s-\f12}(\R^{d+1})},
\end{align*}
from which and Gronwall's inequality, it follows that
\begin{align*}
\|\widetilde{\om_e^n}(t)\|_{H^{s-\f12}(\R^{d+1})}\le \big(\|\widetilde{\om_{0,e}^n}\|_{H^{s-\f12}(\R^{d+1})}+\|F^n_e\|_{H^{s-\f12}(\R^{d+1})}\big)e^{\int_0^t\|\overline{v^n_e}\|_{H^{s+\f12}(\R^{d+1})}dt'}.
\end{align*}
This implies that
\begin{align*}
\|\widetilde{\om^n}(t)\|_{H^{s-\f12}(\overline{\cS})}\le C\big(\|\widetilde{\om_{0}^n}\|_{H^{s-\f12}(\overline{\cS})}+\cP(E^n_s(t))\big)e^{\int_0^t\cP(E^n_s(t'))dt'},
\end{align*}
which together with (\ref{eq:thm-VB-low}) leads to
\ben\label{eq:thm-VBV}
E^n_{s,l}(t)\le E^n_{s,l}(0)+\int_0^t\cP(E^n_s(t'))dt'.
\een

Now let us turn to the higher order energy estimate. For this, we need the following refined elliptic estimate from \cite{ABZ-IM}.

\begin{lemma}\label{lem:elliptic-refined}
Let $v\in H^{s+\f12}(\overline{\cS})$ be a solution of the elliptic equation (\ref{eq:ellitic-flat}) with $v(0)=f$. Then it holds that
for any $\sigma\in [-\f12,s-1]$,
\beno
\|\na_{x,z}v\|_{X^\sigma([-\f12,0])}\le \cP(\|\eta\|_{H^s})\big(\|f\|_{H^{\sigma+1}}+\|F_0\|_{Y^\sigma(I)}\big),
\eeno
if it holds for $\sigma=-\f12$.
\end{lemma}

 \begin{lemma}\label{lem:UVB-E}
It holds that
\beno
&&\|U^n\|_{H^s}+\|(V^n,B^n)\|_{H^s}+\|\eta^n\|_{H^{s+\f12}}\\
&&\le \cP(E_{s,l}^n)\big(\|D_tU^n\|_{H^{s-\f12}}+\|T_{\sqrt{a^n\lambda^n}}U^n\|_{H^{s-\f12}}\big)+\cP(E_{s,l}^n).
\eeno
\end{lemma}

\no{\bf Proof.}\,
{\bf Step 1.} Estimate for $U^n$.

Let $\epsilon\in (0,\min(s_d,1))$ with $s_d=s-\f d 2-1$. Applying Lemma \ref{lem:elliptic-refined} to $P^n$, we obtain
\begin{align}\label{eq:thm-p}
\|\na_{x,z}\widetilde{P_1^n}\|_{X^{s-1}([-\f12,0])}\le \cP(E^n_{s,l}(t)),\quad P_1^n=P^n+y.
\end{align}
Then we deduce from Proposition \ref{prop:symbolic calculus} and (\ref{eq:thm-p}) that
\begin{align*}
\|U^n\|_{H^s}\le& \|T_{(\sqrt{a^n\lambda^n})^{-1}}T_{\sqrt{a^n\lambda^n}}U^n\|_{H^s}+\|(T_{(\sqrt{a^n\lambda^n})^{-1}}T_{(\sqrt{a^n\lambda^n})}-1)U^n\|_{H^s}\\
\le& \cP(E_{s,l}^n)\|T_{\sqrt{a^n\lambda^n}}U^n\|_{H^{s-\f12}}+\cP(E_{s,l}^n)\|U^n\|_{H^{s-\epsilon}},
\end{align*}
which implies that
\ben
\|U^n\|_{H^s}\le \cP(E_{s,l}^n)\|T_{\sqrt{a^n\lambda^n}}U^n\|_{H^s}+\cP(E_{s,l}^n).\label{eq:thm-U4}
\een

{\bf Step 2. Estimate for $\eta^n$}

By the proof of Lemma \ref{lem:eta-E}, we know that
\begin{align*}
\zeta^n=T_{(a^n)^{-1}}\big(-D_tU^n+h_1^n+[D_t,T_{\zeta^n}]B^n\big)+(T_{(a^n)^{-1}}T_{a^n}-1)\zeta^n,
\end{align*}
which along with Proposition \ref{prop:symbolic calculus} and  Proposition \ref{prop:commutator-tame} gives
\begin{align*}
\|\zeta^n\|_{H^{s-\f12}}\le \cP(E_{s,l}^n)\big(&\|D_tU^n\|_{H^{s-\f12}}+\|\zeta^n\|_{H^{s-\f12-\epsilon}}+\|h_1^n\|_{H^{s-\f12}}\\
&+\|V^n\|_{B^1_{\infty,1}}\|B^n\|_{H^{s-\f12}}+\|D_t\zeta^n\|_{L^\infty}\|B^n\|_{H^{s-\f12}}\big).
\end{align*}
Recall that $h_1^n=(T_{V^n}-V^n)\cdot\nabla V^n-R(a^n,\zeta^n)+T_{\zeta^n}(T_{V^n}-V^n)\cdot\nabla B^n$. Then Lemma \ref{lem:remaider}
and (\ref{eq:thm-p}) ensure that
\beno
\|h^n_1\|_{H^{s-\f12}}\le \cP(E_{s,l}^n)\|V^n\|_{W^{1,\infty}}\|(V^n,B^n)\|_{H^{s-\f12}}+\cP(E_{s,l}^n)\|\zeta^n\|_{H^{s-\f12-\epsilon}}.
\eeno
Recall that $(\pa_t+V^n\cdot\na)\zeta^n=\na B^n-\na V^n\cdot\na \eta^n$. Then we have
\ben\label{eq:thm-zeta-Dt}
\|D_t\zeta^n(t)\|_{L^\infty}\le C\|(V^n,B^n)\|_{W^{1,\infty}}\|\eta^n\|_{W^{1,\infty}}.
\een
Thanks to $s>1+\f d 2$, we get by  the interpolation that
\ben\label{eq:thm-VB-inter}
\|(V^n,B^n)\|_{B^1_{\infty,1}}\le \|(V^n,B^n)\|_{H^{s-\f12}}^{2\epsilon}\|(V^n,B^n)\|_{H^{s}}^{1-2\epsilon}.
\een
Thus, we conclude that
\begin{align}\label{eq:thm-eta}
\|\eta^n\|_{H^{s+\f12}}\le \cP(E_{s,l}^n(t))+\cP(E_{s,l}^n(t))\|D_tU^n\|_{H^{s-\f12}}+\cP(E_{s,l}^n(t))\|V^n,B^n\|_{H^s}^{1-2\epsilon}.
\end{align}

{\bf Step 3}. Estimates for $B^n$.

By the proof of Lemma \ref{lem:VB-E}, we know that
\begin{align*}
\dv U^n=-T_{q^n}B^n-R(\eta^n)B^n+V_\om^n+T_{\dv \zeta^n}B^n,
\end{align*}
where the symbol $q^n=\lambda^n-i\zeta^n\cdot \xi$ and
\beno
V_\om^n=-\pa_y(v^n)^{d+1}_\om+\na \eta^n\cdot\na (v^n)^{d+1}_\om\big|_{y=\eta^n}+\pa_i\eta^n\om_{d+1,i}^n\big|_{y=\eta^n}.
\eeno

By using Lemma \ref{lem:elliptic-refined} and following
the proof of Proposition \ref{prop:DN-Hs}, we can deduce that
\ben
&&\|R(\eta^n)(V^n,B^n)\|_{H^{s-1}}\le \cP(E^n_{s,\ell}(t))\|\eta^n\|_{H^{s+\f12}},\nonumber\\
&&\|\na_{x,z}\widetilde{v^n_\om}\|_{X^{s-1}([-\f12,0])}\le \cP(E^n_{s,l}(t)).\label{eq:thm-vom}
\een
Thus, we infer from Lemma \ref{lem:remaider} that
\begin{align}
\|T_{q^n}B^n\|_{H^{s-1}}\le& \|U^n\|_{H^s}+\|R(\eta^n)B^n\|_{H^{s-1}}+\|V_\om^n\|_{H^{s-1}}+\|\zeta^n\|_{C^\epsilon}\|B^n\|_{H^{s-\epsilon}}\nonumber\\
\le& \|U^n\|_{H^s}+\cP(E_{s,l}^n)\|\eta^n\|_{H^{s+\f12}}+\cP(E_{s,l}^n)\|B^n\|_{H^{s-\epsilon}}.\label{eq:thm-B1}
\end{align}
Then Proposition \ref{prop:symbolic calculus} ensures that
\begin{align}
\|B^n\|_{H^s}\le& \|(T_{(q^n)^{-1}}T_{q^n}-1)B^n\|_{H^s}+\|T_{(q^n)^{-1}}T_{q^n}B^n\|_{H^s}\nonumber\\
\le& \cP(E_{s,l}^n)\|B^n\|_{H^{s-\epsilon}}+\cP(E_{s,l}^n)\|T_{q^n}B^n\|_{H^{s-1}}\nonumber\\
\le&   \cP(E_{s,l}^n)(\|B^n\|_{H^{s-\epsilon}}+\|U^n\|_{H^s}+\cP(E_{s,l}^n)\|\eta^n\|_{H^{s+\f12}}).\label{eq:thm-B2}
\end{align}

{\bf Step 4.}  Completion of the estimate

Combining the estimates (\ref{eq:thm-B2}), (\ref{eq:thm-B1}) and (\ref{eq:thm-U4}), we have that
\begin{align}
&\|U^n\|_{H^s}+\|B^n\|_{H^s}+\|\eta^n\|_{H^{s+\f12}}\nonumber\\
&\le \cP(E_{s,l}^n)\big(\|D_tU^n\|_{H^{s-\f12}}+\|T_{\sqrt{a^n\lambda^n}}U^n\|_{H^{s-\f12}}\big)+\cP(E_{s,l}^n).\label{eq:thm-B2}
\end{align}
This also gives
\begin{align*}
\|V^n\|_{H^s}\le& \|U^n\|_{H^s}+\|T_{\zeta^n}B^n\|_{H^s}\nonumber\\
\le& \cP(E_{s,l}^n)\big(\|D_tU^n\|_{H^{s-\f12}}+\|T_{\sqrt{a^n\lambda^n}}U^n\|_{H^{s-\f12}}\big)+\cP(E_{s,l}^n),
\end{align*}
which together with (\ref{eq:thm-U4}) and  (\ref{eq:thm-B2}) gives the lemma.\ef

\medskip

Proceeding the same way as in section 9.6, we can deduce that
\begin{align}
&\f d {dt}\Big(\|D_tU^n\|_{H^{s-\f12}}^2+\|T_{\sqrt{a^n\lambda^n}}U^n\|_{H^{s-\f12}}^2+\big\langle (f_\om^n)_{s-1}, U_{s}^n\big\rangle\nonumber\\
&\qquad\quad+\big\langle g_{s-1/2}^n, h_{s-1/2}^n\big\rangle+\f12\big\langle g_{s-1/2}^n, g_{s-1/2}^n\big\rangle\Big)\le\cP(E^n_s).\label{eq:thm-U1}
\end{align}
Here $g^n=[D_t,T_{\zeta^n}]B^n$ and $h^n=D_tU^n-[D_t,T_{\zeta^n}]B^n$.

It follows from Proposition \ref{prop:commutator-tame} that
\begin{align*}
&\|g^n\|_{H^{s-\f12}}\le C\big(\|V^n\|_{B^1_{\infty,1}}+\|D_t\zeta^n\|_{L^\infty}\big)\|B^n\|_{H^{s-\f12}},\\
&\|h^n\|_{H^{s-\f12}}\le \|D_tU^n\|_{H^{s-\f12}}+C\big(\|V^n\|_{B^1_{\infty,1}}+\|D_t\zeta^n\|_{L^\infty}\big)\|B^n\|_{H^{s-\f12}},
\end{align*}
and we have by (\ref{eq:thm-vom}) that
\beno
\|f_\om^n\|_{H^{s-1}}\le \cP(E^n_{s,l}(t)).
\eeno
Plugging the above estimates into (\ref{eq:thm-U1}), we obtain
\begin{align}
&\|D_tU^n(t)\|_{H^{s-\f12}}^2+\|T_{\sqrt{a^n\lambda^n}}U^n(t)\|_{H^{s-\f12}}^2\nonumber\\
&\le \cP(E_s(0))+\cP(E^n_{s,l}(t))\|U^n(t)\|_{H^s}+C\big(\|V^n(t)\|_{B^1_{\infty,1}}+\|D_t\zeta^n(t)\|_{L^\infty}\big)^2\|B^n(t)\|_{H^{s-\f12}}^2\nonumber\\
&\quad+C\big(\|V^n(t)\|_{B^1_{\infty,1}}+\|D_t\zeta^n(t)\|_{L^\infty}\big)\|B^n(t)\|_{H^{s-\f12}}\|D_tU^n\|_{H^{s-\f12}}+\int_0^t\cP(E^n_s(t'))dt',\nonumber
\end{align}
which together with (\ref{eq:thm-zeta-Dt}), (\ref{eq:thm-VB-inter}) and Lemma \ref{lem:UVB-E} implies that
\begin{align*}
E_s^n(t)\le \cP(E_s(0))+\cP(E^n_{s,l}(t))+\int_0^t\cP(E^n_s(t'))dt'.
\end{align*}
This together with (\ref{eq:thm-VBV}) ensures that there exits $T>0$ depending on $E_s(0)$ and $c_0, h_0, s$ such that for any $t\in [0,\min(T,T_n))$,
there holds
\beno
E_{s}^n(t)\le \cP(E_s(0)).
\eeno

With the uniform estimates, the existence of the solution can be deduced by a standard compact argument. Here we omit the details.

\subsection{Uniqueness of the solution}

Let $(\eta^1, v^1, P^1)$ and $(\eta^2, v^2, P^2)$ be two solutions of the system (\ref{eq:euler})--(\ref{eq:euler-p}) with the same initial data.
Assume that the solutions satisfy
\begin{align*}
\sup_{t\in [0,T]}\big(\|\eta^1(t)\|_{H^{s+\f12}}+\|\widetilde{v^1}(t)\|_{H^{s+\f12}(\overline{\cS})}+\|\eta^2(t)\|_{H^{s+\f12}}+\|\widetilde{v^2}(t)\|_{H^{s+\f12}(\overline{\cS})}\big)\le \cP_0,
\end{align*}
Here and in what follows, we denote by $\cP_0$ a constant depending only on $E_s(0), c_0, h_0, s$, which may change from line to line.

We denote $\delta_f=f^2-f^1$ and  introduce
\begin{align*}
\delta_{\cE}(t)\eqdef \delta_{\cE_1}(t)+\delta_{\cE_2}(t),
\end{align*}
where $\delta_{\cE_1}(t)$ and $\delta_{\cE_2}(t)$ are given by
\begin{align*}
\delta_{\cE_1}(t)=&\|\big(\delta_{V},\delta_{B}\big)(t)\|_{H^{s-1}}+\|\delta_{\eta}(t)\|_{H^{s-\f12}},\\
\delta_{\cE_2}(t)=&\|\big(\delta_{V},\delta_{B}\big)(t)\|_{H^{s-\f32}}+\|(\delta_{V_b},\delta_{\eta})(t)\|_{H^{s-1}}
+\|\delta_{\widetilde{\om}}(t)\|_{H^{s-\f32}(\overline{\cS})}.
\end{align*}
In what follows, we denote by $L_i(i=0,1,2)$  some nonlinear term, which satisfies
\beno
\|L_i(t)\|_{H^{s-1-\f i2}}\le \cP_0\delta_{\cE}(t).
\eeno
We denote $D_t=\pa_t+V^2\cdot\na$.

Following the proof of section 12.2, we can show that
\ben
&&\|D_t\delta_{R_\om}\|_{H^{s-2}}+\|(\delta_{R(\eta)(V,B)},\delta_{a})\|_{H^{s-\f32}}
+\|\delta_{\na P|_{y=-1}}\|_{H^{s-1}}\le \cP_0\delta_{\cE}(t),\label{eq:Un-1}\\
&&\|\delta_{R_\om}\|_{H^{s-2}}\le \cP_0\delta_{\cE_2}(t).\label{eq:Un-2}
\een

With (\ref{eq:Un-1}) and (\ref{eq:Un-2}), we can deduce that
\beno
\left\{
\begin{aligned}
&D_t\big(\delta_{V},\delta_{B}\big)=L_1,\\
&(\pa_t+V_b^2\cdot\na)\delta_{V_b}=L_0,\\
&D_t\delta_{\zeta}=L_2,\\
&D_t\delta_U+T_{{a^2\lambda^2}}\delta_U=L_1+T_{a^2}\delta_{R_\om}.
\end{aligned}\right.
\eeno
The energy estimate ensures that
\begin{align*}
&\|D_t\delta_U(t)\|_{H^{s-\f32}}+\|T_{\sqrt{a^2\lambda^2}}\delta_U(t)\|_{H^{s-\f32}}+\|\big(\delta_{V},\delta_{B}\big)(t)\|_{H^{s-\f32}}\\
&\quad+\|(\delta_{V_b},\delta_{\eta})(t)\|_{H^{s-1}}\le \cP_0\int_0^t\delta_{\cE}(t')dt'.
\end{align*}
For the vorticity, we have
\beno
\|\delta_{\widetilde{\om}}(t)\|_{H^{s-\f32}(\overline{\cS})}\le \cP_0\int_0^t\delta_{\cE}(t')dt'.
\eeno
Using the fact that
\beno
\delta_{\cE}(t)\le \cP_0\big(\delta_{\cE_2}(t)+\|D_t\delta_U(t)\|_{H^{s-\f32}}+\|T_{\sqrt{a^2\lambda^2}}\delta_U(t)\|_{H^{s-\f32}}\big).
\eeno
We conclude that
\beno
\delta_{\cE}(t)\le \cP_0\int_0^t\delta_{\cE}(t')dt',
\eeno
which implies that $\delta_{\cE}(t)=0$ for any $t\in [0,T]$.
The uniqueness is proved.\ef

\section*{Acknowledgement}

Zhifei Zhang is partially supported by NSF of China under Grant
11371037 and 11425103.

\end{document}